\documentclass[12pt]{article}

\usepackage{amssymb,amsthm,amsmath}
\makeatletter

\let\Mathrm\operator@font
\def\c@l#1{\ifmmode{\mathcal#1}\else\relax\fi}
\def\Cal{\protect\c@l}

\def\ptruebf{\fontseries{b}\fontshape{n}\selectfont}
\def\truebf{\protect\ptruebf}

\def\fr@k#1{{\@mathfrak #1}}
\def\frak{\protect\fr@k}
\def\BBB#1{{\mathbb #1}}
\def\Bbb{\protect\BBB}

\renewcommand\){{\rm )}}

\swapnumbers
\newtheorem{Theorem}{Theorem}[section]
\newtheorem{Proposition}[Theorem]{Proposition}
\newtheorem{Lemma}[Theorem]{Lemma}
\newtheorem{Corollary}[Theorem]{Corollary}

\theoremstyle{definition}

\newtheorem{example}[Theorem]{Example}
\newtheorem{Def}[Theorem]{Definition}

\theoremstyle{remark}
\newtheorem{Rem}[Theorem]{Remark}

\let\c@equation\c@Theorem

\let\c@paragraph\c@Theorem

\def\@@@sect#1#2#3#4#5#6[#7]#8{%
   \ifnum #2>\c@secnumdepth 
      \def \@svsec {}\else \refstepcounter {#1}%
      \def\@svsec{}
   \fi 
   \@tempskipa #5\relax 
   \ifdim \@tempskipa >\z@ 
     \begingroup #6\relax \@hangfrom {\hskip #3\relax 
     \@svsec}{\interlinepenalty \@M #8\par }\endgroup 
     \csname #1mark\endcsname {#7}
   \else 
   \def \@svsechd {#6\hskip #3\@svsec #8\csname #1mark\endcsname {#7}}
   \fi \@xsect {#5}}

\def\@@@startsection#1#2#3#4#5#6{%
 \if@noskipsec \leavevmode \fi \par \@tempskipa #4\relax \@afterindenttrue 
 \ifdim \@tempskipa <\z@ \@tempskipa -\@tempskipa \@afterindentfalse 
 \fi \if@nobreak \everypar {}\else \addpenalty {\@secpenalty }\addvspace 
  {\@tempskipa }\fi \@ifstar {\@ssect {#3}{#4}{#5}{#6}}{\@dblarg 
  {\@@@sect {#1}{#2}{#3}{#4}{#5}{#6}}}}

\def\theparagraph{\thesection.\arabic{paragraph}}
\def\paragraph{\@@@startsection{paragraph}{2}{\z@ }%
              {1.75ex plus .2ex minus .15ex}{-1em}{}{\truebf(\theparagraph)}}

\def\uop#1{\mathop{#1}\nolimits}
\def\standop#1{\mathop{\Mathrm #1}\nolimits}
\def\difstop#1#2{\expandafter\def\csname #1\endcsname{\standop{#2}}}
\def\defstop#1{\difstop{#1}{#1}}

\def\ev{\standop{ev}}
\defstop{Rqco}
\defstop{Dat}
\defstop{Nerve}
\defstop{LQco}
\defstop{Qco}
\defstop{EM}
\defstop{cosk}
\defstop{EqPA}
\defstop{EqAB}
\defstop{Mor}

\def\Set{\underline{\Mathrm S \Mathrm e\Mathrm t}}
\def\Ab{\underline{\Mathrm A\Mathrm b}}
\defstop{AB}
\defstop{PA}
\defstop{Mod}
\defstop{PM}
\defstop{Zar}
\def\Rng{\underline{\Mathrm R\Mathrm n\Mathrm g}}
\def\delme#1{\mathinner{\mathop{\check{\fivecdots}}\limits^{#1}}}
\def\fivecdots{\cdotp\cdotp\cdotp\cdotp\cdotp}
\def\ttimes{\mathbin{{}_{d_1}\!\times_{d_0}}}
\def\dtimes{\mathbin{{}_{d_1}\!\times_{d_1}}}

\let\setbar\mid
\newcommand{\id}{{\Mathrm i\Mathrm d}}
\def\Ext{\standop{Ext}}
\def\Sym{\standop{Sym}}
\def\ext{\uop{\textstyle\bigwedge}}
\newcommand{\Image}{\standop{Im}}
\newcommand{\Ker}{\standop{Ker}}
\def\specialarrow#1{\setbox\z@=\hbox{$\m@th\mathop{\vphantom
   {\rightarrow}}\limits^{\hspace{.5ex}{#1}\hspace{.8ex}}$}\mathord
   {\ifdim\wd\z@<1.2em\dimen\tw@1.2em\else\dimen\tw@\wd\z@\fi\copy
   \z@\kern-\wd0\hbox to\dimen\tw@{\rightarrowfill}}}
\def\specialleftarrow#1{\setbox\z@=\hbox{$\m@th\mathop{\vphantom
   {\leftarrow}}\limits^{\hspace{.7ex}{#1}\hspace{.6ex}}$}\mathord
   {\ifdim\wd\z@<1.2em\dimen\tw@1.2em\else\dimen\tw@\wd\z@\fi\copy
   \z@\kern-\wd0\hbox to\dimen\tw@{\leftarrowfill}}}
\def\suarrow#1{\uparrow\hbox to 0pt{\scriptsize$#1$\hss}}
\def\sdarrow#1{\downarrow\hbox to 0pt{\scriptsize$#1$\hss}}
\def\lrsplit#1\\#2\endlrsplit{$$\displaylines{#1\hfill\cr\hfill #2\cr}$$}

\newcommand{\Hom}{\standop{Hom}}
\newcommand{\rank}{\standop{rank}}
\def\Spec{\standop{Spec}}
\def\uSpec{\mathop{\mbox{\rm \underline{Spec}}}\nolimits}
\def\height{\standop{ht}}
\def\Coh{\standop{Coh}}
\def\SLN#1{{\em Lect. Notes Math. }{\bf #1}, Springer Verlag}

\def\depth{\standop{depth}}
\def\Id{\standop{Id}}
\def\uHom{\mathop{\mbox{\underline{$\Mathrm Hom$}}}\nolimits}
\def\suHom{\mathop{\mbox{\scriptsize\underline{$\Mathrm Hom$}}}\nolimits}
\def\uExt{\mathop{\mbox{\underline{$\Mathrm Ext$}}}\nolimits}
\def\pdim{\mathop{\mbox{$\Mathrm proj.dim$}}\nolimits}
\def\ob{\standop{ob}}
\def\triv{{\standop{triv}}}
\def\op{^{\standop{op}}}

\let\Func\func
\def\trace{\standop{tr}}
\def\red{\standop{red}}
\defstop{Sch}
\defstop{Cos}
\defstop{shift}
\defstop{qco}
\defstop{holim}
\defstop{hocolim}
\def\via{{\Mathrm via~}}
\def\mon{^{\Mathrm mon}}
\def\cc{\mathop{\check H}\nolimits}
\def\ccs{\mathop{\check{\underline H}}\nolimits}
\def\SL{{\sl{SL}}}
\def\GL{{\sl{GL}}}

\def\subcolim#1{\setbox0\hbox{$#1$}\mathop
  {\hbox{\vtop{\offinterlineskip\halign{##\cr
  \copy0\cr\hbox to\wd0{\rightarrowfill}\cr}}}}\limits}
\def\sublim#1{\setbox0\hbox{$#1$}\mathop
  {\hbox{\vtop{\offinterlineskip\halign{##\cr
  \copy0\cr\hbox to\wd0{\leftarrowfill}\cr}}}}\limits}

\let\indlim\varinjlim
\let\projlim\varprojlim
\def\uTor{\mathop{\mbox{\rm \underline{Tor}}}\nolimits}

\def\verylongrightarrow{\relbar\joinrel\longrightarrow}

\title{EQUIVARIANT TWISTED INVERSE WITHOUT
                      EQUIVARIANT COMPACTIFICATION}
\author{M{\sc itsuyasu} H{\sc ashimoto}}

\date{\normalsize
Graduate School of Mathematics, Nagoya University,
Chikusa-ku,  Nagoya 464--8602 JAPAN
}

\begin{document}

\maketitle

\begin{abstract}
An equivariant version of the twisted inverse pseudofunctor is defined,
and equivariant versions of some important properties, including the 
Grothendieck duality of proper morphisms and flat base change are proved.
As an application, a generalized version of Watanabe's theorem on
the Gorenstein property of the ring of invariants is proved.
\end{abstract}

\maketitle

\tableofcontents

\section{Introduction}

Let $S$ be a separated noetherian scheme, $G$ a flat separated 
$S$-group scheme of finite type, $X$ and $Y$ $S$-schemes separated
of finite type with $G$-actions, and $f:X\rightarrow Y$ a $G$-morphism.

The purpose of these notes is to construct an equivariant version of
the twisted inverse functor $f^!$.

One of the main motivation of the work is an application to the
invariant theory.
As an example, we give a very short proof of a generalized version
of Watanabe's theorem on the Gorenstein property of invariant subrings
\cite{Watanabe}.
Also, there might be some meaning in formulating the equivariant
duality theorem, which has the Serre duality for the representations
of reductive groups (see \cite{Jantzen}) as a special case,
in a reasonably general form for records.

In the case where $G$ is trivial, $f^!$ is defined as follows.
For a scheme $Z$, we denote the category of $\Cal O_Z$-modules by
$\Mod(Z)$.
By definition, a plump subcategory of an abelian category is a
non-empty full subcategory which is 
closed under kernels, cokernels, and extensions \cite[(1.9.1)]{Lipman}.
We denote the plump subcategory of $\Mod(Z)$ consisting of 
quasi-coherent $\Cal O_Z$-modules by $\Qco(Z)$.

By Nagata's compactification theorem, there exists some factorization
\[
X\specialarrow{i}\bar X\specialarrow{p} Y
\]
such that $p$ is proper and $i$ an open immersion.
We call such a factorization a {\em compactification}.
We define $f^!:D^+_{\Qco(Y)}(\Mod(Y))\rightarrow D^+_{\Qco(X)}(\Mod(X))$ 
to be the composite $i^*p^\times$, where 
$p^\times:D^+_{\Qco(Y)}(\Mod(Y))\rightarrow D^+_{\Qco(\bar X)}(\Mod(\bar X))$ 
is the right adjoint of $Rp_*$, and $i^*$ is the restriction.
The definition of $f^!$ is independent of choice of compactification.

In order to consider a non-trivial $G$, we need to replace
$\Qco(X)$ and $\Mod(X)$ by some appropriate categories which respect 
$G$-actions.
The category $\Qco(G,X)$ 
which corresponds to $\Qco(X)$ is fairly well-known.
It is the category of $G$-linearlized quasi-coherent $\Cal O_X$-modules
defined by Mumford \cite{GIT}.
The category $\Qco(G,X)$ is equivalent to the category of 
quasi-coherent sheaves over the diagram of schemes
\[
B_G^M(X):=\left(G\times_S G\times_S X
\begin{array}{c}
\mathop{\verylongrightarrow}\limits^{1_G\times a}\\
\mathop{\verylongrightarrow}\limits^{\mu\times 1_X}\\
\mathop{\verylongrightarrow}\limits^{p_{23}}
\end{array}
G\times_S X
\begin{array}{c}
\mathop{\verylongrightarrow}\limits^{a}\\
\mathop{\verylongrightarrow}\limits^{p_2}
\end{array}
X\right)
,
\]
where $a:G\times X\rightarrow X$ is the
action, $\mu:G\times G\rightarrow G$ the product,
and $p_{23}$ and $p_2$ are appropriate projections.
Thus it is natural to embed the category $\Qco(G,X)$ into the
category of all $\Cal O_{B_G^M(X)}$-modules $\Mod(B_G^M(X))$,
and $\Mod(B_G^M(X))$ is a good substitute of $\Mod(X)$.
As $G$ is flat, $\Qco(G,X)$ is a plump subcategory of $\Mod(B_G^M(X))$,
and we may consider $D_{\Qco(G,X)}(\Mod(B_G^M(X)))$.
However, our construction utilize an intermediate category $\LQco(G,X)$
(the category of locally quasi-coherent sheaves),
and is not an obvious interpretation of the non-equivariant case.

Note that there is a natural restriction functor $\Mod(B_G^M(X))\rightarrow
\Mod(X)$, which sends $\Qco(G,X)$ to $\Qco(X)$.
This functor is regarded as the forgetful functor, forgetting the
$G$-action.
The equivariant duality theorem which we are going to establish must be
compatible with this restriction functor, otherwise the theory 
would be something
different from the usual scheme theory and probably useless.

Most of the discussion in these notes treats more general diagrams of
schemes.
This makes the discussion easier, 
as some of the important properties are proved
by induction on the number of objects in the diagram.
Our main construction and theorems are only for the class of 
finite diagrams of schemes of certain type, which contains the diagrams
of the form $B_G^M(X)$.

In sections 2--4, we review general facts on homological algebra.
The construction of $f^!$ is divided into five steps.
The first is to analize the functoriality of sheaves over diagrams of schemes.
Sections 5--8 are devoted to this step.
The second is the derived version of the first step.
This will be done in sections 9, 10, and 13.
Note that not only the categories of all module sheaves $\Mod(X_\bullet)$ and
the category of quasi-coherent sheaves $\Qco(X_\bullet)$, 
but also the category of locally quasi-coherent sheaves $\LQco(X_\bullet)$
also plays an important role in our construction.

The third is to prove the existence of the right adjoint $p^\times_\bullet 
$ of $R(p_\bullet)_*$
for (componentwise) proper morphism $p_\bullet$ of diagrams of schemes.
This is not so difficult, and is done in section~12.
We use Neeman's existence theorem on the right adjoint of unbounded right
derived functors.
Not only to utilize Neeman's theorem, but to calculate composite of
various left and right derived functors,
it is convenient to utilize unbounded derived functors.
A short survey on unbounded derived functors is given in section~4.

The fourth step is to prove various commutativity related to the
well-definedness of the
twisted inverse pseudofunctors,
sections 11, 14, and 15.
Among them, the compatibility with restrictions 
(Proposition~\ref{twisted-restrict.thm}) is the key to our construction.
Given a $G$-morphism $f:X\rightarrow Y$ between $G$-schemes separated
of finite type over $S$, the associated morphism $B_G^M(f):B_G^M(X)
\rightarrow B_G^M(Y)$ is {\em of fiber type} (or cartesian),
see for the definition, (\ref{of-fiber-type.par}).
If we could find a compactification
\[
B_G^M(X)\specialarrow{i_\bullet} Z_\bullet\specialarrow{p_\bullet}
B_G^M(Y)
\]
such that $p_\bullet$ is proper and of fiber type, and $i_\bullet$ an
image-dense open immersion, then the construction of $f^!$ and the 
proof of various commutative diagrams would be fairly easy.
However, it seems that this is almost the same as the problem of equivariant 
compactifications.
See \cite{Sumihiro} for equivariant compactifications.
We avoid this problem, and prove the commutativity of various diagrams
without assuming that $p_\bullet$ is of fiber type.

The fifth part is the existence of factorization $f_\bullet=p_\bullet
i_\bullet$, where $p_\bullet$ is proper and $i_\bullet$ an image-dense
open immersion.
This is easily done utilizing Nagata's compactification theorem,
and is done in section 16.
This completes the basic construction of the equivariant twisted
inverse pseudofunctor.

In sections 17--24, we prove equivariant versions of most of the
known results on twisted inverse including the equivariant Grothendieck
duality and flat base change, except that the equivariant dualizing
complexes are treated later.
We also prove that the twisted inverse functor preserves quasi-coherent
cohomology groups.
As we already know the corresponding results 
on single schemes and the commutativity with restrictions, this 
consists in straightforward (but not easy) check of various commutativity
of diagrams of functors.

So far, almost all results are valid for any diagram of 
separated noetherian schemes with flat arrows over a finite ordered
category.
In sections~25--31, we prove some results which are special on the 
$G$-actions.

In section~32, we give a definition of the equivariant dualizing complexes.
As an application, we give a very short proof of a generalized
version of Watanabe's theorem 
on the Gorenstein property of invariant subrings in section~33.

\noindent
{\em Acknowledgement:}
The author is grateful to Professor Ryoshi Hotta and
Professor Joseph Lipman for valuable advice.

\section[Commutativity of diagrams]{Commutativity of diagrams 
constructed from a monoidal pair of pseudofunctors}
\label{monoidal.sec}

\paragraph Let $(?)_{*}$ be a covariant symmetric monoidal pseudofunctor
\cite[(3.6.5)]{Lipman} of closed categories on a category $\Cal S$.
Thus, for each $X\in \Cal S$, 
\[
X_{*}=(M_{X},\otimes,\Cal O_X,\alpha,\lambda,
\gamma,[?,*],\pi)
\]
is a (symmetric monoidal) closed category (see e.g., 
\cite[(3.5.1)]{Lipman}), where $M_X$ is the underlying category, 
$\otimes:M_X\times M_X\rightarrow M_X$ 
the product structure, $\Cal O_X\in M_X$ the unit object,
$\alpha:(a\otimes b)\otimes c\cong a\otimes(b\otimes c)$ 
the associativity isomorphism, $\lambda:\Cal O_X\otimes a\cong a$ 
the left unit isomorphism, $\gamma:a\otimes b\cong b\otimes a$ the
twisting (symmetry) isomorphism, $[?,-]:M_X\op\times M_X\rightarrow M_X$
the internal hom, and
\begin{equation}\label{associative-adjunction.eq}
\pi:M_X(a\otimes b,c)\cong M_X(a,[b,c])
\end{equation}
the associative adjunction isomorphism of $X$, respectively.

In this paper, various (different) adjoint pairs appears almost everywhere.
By abuse of notation, the unit (resp.\ the counit) of adjunction is usually
simply denoted by the same symbol $u$ (resp.\ $\varepsilon$).
However, the unit map and the counit map arising from the adjunction
(\ref{associative-adjunction.eq}) are denoted (by less worse abuse of notation)
by (the same symbol)
\[
\trace: a\rightarrow [b,a\otimes b]\qquad(\mbox{the trace map})
\]
and
\[
\ev: [b,c]\otimes b\rightarrow c\qquad(\mbox{the evaluation map}),
\]
respectively.

\paragraph 
Let $f:X\rightarrow Y$ be a morphism.
Then $f_*$ is a symmetric monoidal functor.
The natural map
\[
f_*a\otimes f_*b\rightarrow f_*(a\otimes b)
\]
is denoted by $m=m(f)$, and 
the map
\[
\Cal O_Y\rightarrow f_*\Cal O_X
\]
is denoted by $\eta=\eta(f)$.

\paragraph
For composable two morphisms $f$ and $g$ in $\Cal S$, the given isomorphism
$(gf)_*\cong g_*f_*$ is denoted by $c=c_{f,g}$ as in \cite[(3.6.5)]{Lipman}.
Let $gf=f'g'$ be a commutative diagram in $\Cal S$.
The composite isomorphism
\[
g_*f_*\specialarrow{c^{-1}}(gf)_*=(f'g')_*\specialarrow{c}
f'_*g'_*
\]
is also denoted by $c=c(gf=f'g')$, by abuse of notation.

\paragraph \label{H.par}
Let $f:X\rightarrow Y$ be a morphism in $\Cal S$.
The composite natural map
\[
f_*[a,b]\specialarrow{\trace}[f_*a,f_*[a,b]\otimes f_*a]
\specialarrow{\via m}
[f_*a,f_*([a,b]\otimes a)]
\specialarrow{\via \ev}[f_*a,f_*b]
\]
is denoted by $H$.
In other words, $H$ is right conjugate \cite[(3.3.5)]{Lipman} to
the map
\[
f_*[a,b]\otimes f_*a\specialarrow{m}f_*([a,b]\otimes a)
\specialarrow{f_*\ev}f_*b.
\]

\paragraph\label{Lewis.par}
G.~Lewis proved a theorem which guarantee that some diagrams
involving two symmetric monoidal closed categories and one symmetric
monoidal functor commute \cite{Lewis}.

By Lewis' result, we have that the following diagrams are commutative
for any morphism $f:X\rightarrow Y$ (also checked by direct computation).
\begin{equation}\label{Lewis-1.eq}
\begin{array}{ccc}
f_*[a,b]\otimes f_* a & \specialarrow{H\otimes 1} & [f_*a,f_*b]\otimes f_* a\\
\sdarrow{m}  & & \sdarrow{\ev}\\
f_*([a,b]\otimes a) & \specialarrow{f_*\ev} & f_*b
\end{array}
\end{equation}
\begin{equation}\label{Lewis-2.eq}
\begin{array}{ccc}
f_*a & \specialarrow{f_*\trace} & f_*[b,a\otimes b]\\
\sdarrow{\trace} &   & \sdarrow{ H}  \\
{}[f_*b,f_*a\otimes f_*b] & \specialarrow{m} & [f_*b,f_*(a\otimes b)]
\end{array}
\end{equation}
\begin{equation}\label{Lewis-3.eq}
\begin{array}{ccccccc}
f_* a & \specialarrow{\lambda^{-1}}& \Cal O_Y\otimes f_*a &
\specialarrow{\eta} & f_*\Cal O_X\otimes f_*a &
\specialarrow{\via\trace} & f_*[a,\Cal O_X\otimes a]\otimes f_* a\\
\sdarrow{1} & & & & & & \sdarrow{\via\lambda}\\
f_*a &
\multicolumn{3}{c}{
\mkern -2mu
\mathord\leftarrow
\mkern -10mu
\cleaders \hbox {$\mkern -2mu\mathord -\mkern -2mu$}
\hfill
\hbox to 0pt{\hss$
  \mathord{\mathop{\mkern -1mu\mathord-\mkern -1mu}\limits^{\hbox to
  0pt{\hss\scriptsize $\ev$\hss}}}$
   \hss}
\mkern -1mu
\cleaders \hbox {$\mkern -2mu\mathord -\mkern -2mu$}
\hfill
}
& [f_*a,f_*a]\otimes f_*a & \specialleftarrow{H\otimes 1} &
f_* [a,a]\otimes f_*a
\end{array}
\end{equation}

\begin{Lemma}\label{composition-H.thm}
Let $f:X\rightarrow Y$ and $g:Y\rightarrow Z$ be morphisms
in $\Cal S$.
Then the diagram
\[
\begin{array}{ccccc}
   g_*f_*[a,b] & \specialarrow{H} & g_*[f_*a,f_*b] & \specialarrow{H} &
   [g_*f_*a,g_*f_*b]\\
  \suarrow{c} & & & & \suarrow{[c^{-1},c]}\\
  (gf)_*[a,b] &
\multicolumn{3}{c}{
\mathord -\mkern -6mu
\cleaders \hbox {$\mkern -2mu\mathord -\mkern -2mu$}
\hfill
\hbox to 0pt{\hss$
  \mathord{\mathop{\mkern -1mu\mathord-\mkern -1mu}\limits^{\hbox to
  0pt{\hss\scriptsize $H$\hss}}}$
   \hss}
\mkern -1mu
\cleaders \hbox {$\mkern -2mu\mathord -\mkern -2mu$}
\hfill\mkern -6mu\mathord \rightarrow
} & [(gf)_*a,(gf)_*b]
\end{array}
\]
is commutative.
\end{Lemma}

\proof Consider the diagram
{\small
\[
\begin{array}{ccccccc}
g_*f_*[a,b]\otimes g_*f_* a &\specialarrow{m} &
g_*(f_*[a,b]\otimes f_*a) & \specialarrow{m} &
g_*f_*([a,b]\otimes a) & \specialarrow{\via \ev} &
g_*f_* b\\
\suarrow{c\otimes c} & &\sigma & &
\sdarrow{c^{-1}} &\tau& \sdarrow{c^{-1}}\\
(gf)_*[a,b]\otimes(gf)_*a &
\multicolumn{3}{c}{
\mathord -\mkern -6mu
\cleaders \hbox {$\mkern -2mu\mathord -\mkern -2mu$}
\hfill
\hbox to 0pt{\hss$
  \mathord{\mathop{\mkern -1mu\mathord-\mkern -1mu}\limits^{\hbox to
  0pt{\hss\scriptsize $m$\hss}}}$
   \hss}
\mkern -1mu
\cleaders \hbox {$\mkern -2mu\mathord -\mkern -2mu$}
\hfill\mkern -6mu\mathord \rightarrow
} & (gf)_*([a,b]\otimes a) &
\specialarrow{\via\ev}& (gf)_*b.
\end{array}
\]
}
The diagram $\sigma$ is commutative by the definition of monoidal 
pseudofunctor, see \cite[(3.6.7.2)]{Lipman}.
The commutativity of $\tau$ is obvious.
Taking the conjugate, we get the commutativity of the diagram in the lemma.
\qed

\paragraph\label{mon-adj-pair.par}
Let us assume that there is a left adjoint $(?)^*$ of $(?)_*$
so that $((?)^*,(?)_*)$ is a monoidal adjoint pair \cite[pp.~107--109]{Lipman}.

For composable two morphisms $f$ and $g$ in $\Cal S$, we denote
the map $f^*g^*\rightarrow (gf)^*$ conjugate to $c:(gf)_*\rightarrow g_*f_*$
by $d=d_{f,g}$.
Being the conjugate of an isomorphism, $d$ is an isomorphism.
For a commutative diagram $gf=f'g'$ in $\Cal S$, the composite map
\[
(g')^*(f')^*\specialarrow{d}(f'g')^*=(gf)^*\specialarrow{d^{-1}}f^*g^*
\]
is also denoted by $d$, by abuse of notation.
For a morphism $f:X\rightarrow Y$, the map $f^*\Cal O_Y\rightarrow \Cal O_X$
conjugate to $\eta:\Cal O_Y\rightarrow f_*\Cal O_X$ is denoted by 
$C$.
The composite map
\begin{equation}\label{def-upper*.eq}
f^*(a\otimes b)\specialarrow{\via u}
f^*(f_*f^*a\otimes f_*f^*b)
\specialarrow{\via m}
f^*f_*(f^*a\otimes f^*b)
\specialarrow{\varepsilon}
f^*a\otimes f^*b
\end{equation}
is denoted by $\Delta$.

Almost by definition, the diagrams
\begin{equation}\label{unit-otimes.eq}
\begin{array}{ccc}
   a\otimes b & \specialarrow{u\otimes u} & f_*f^*a\otimes f_*f^*b\\
  \sdarrow{u} &                         & \sdarrow{m}\\
 f_*f^*(a\otimes b) & \specialarrow{\Delta} & f_*(f^*a\otimes f^* b)
\end{array}
\end{equation}
and
\begin{equation}\label{counit-otimes.eq}
\begin{array}{ccc}
f^*(f_*a\otimes f_*b) & \specialarrow{m} & f^*f_*(a\otimes b)\\
\sdarrow{\Delta} &   & \sdarrow{\varepsilon}\\
f^*f_*a\otimes f^*f_*b & \specialarrow{\varepsilon\otimes \varepsilon} &
a\otimes b
\end{array}
\end{equation}
are commutative.

\begin{Def} 
A monoidal adjoint pair $((?)^*,(?)_*)$ is said to be {\em Lipman}
if
$\Delta:f^*(a\otimes b)\rightarrow f^*a\otimes f^*b$
and $C: f^*\Cal O_Y\rightarrow \Cal O_X$ are isomorphisms
for any morphism $f:X\rightarrow Y$ in $\Cal S$ and any $a,b\in M_Y$.
\end{Def}

\paragraph Note that $\Delta:f^*(a\otimes b)\rightarrow f^*a\otimes f^*b$ 
is a natural isomorphism if and only if its right conjugate 
(see \cite[(3.3.5)]{Lipman}) is an isomorphism.
The right conjugate is the composite
\[
f_*[f^*b,a]\specialarrow{H}[f_*f^*b,f_*a]\specialarrow{u}[b,f_*a].
\]

Let $((?)^*,(?)_*)$ be a Lipman adjoint pair of monoidal
pseudofunctors over $\Cal S$.
Then $(?)^*$ together with $\Delta^{-1}$ and $C^{-1}$ form a 
{\em covariant} symmetric monoidal pseudofunctor on $\Cal S\op$.

\paragraph
For a morphism $f:X\rightarrow Y$ in $\Cal S$, the composite map
\[
f^*[a,b]\specialarrow{\via \trace}
[f^*a,f^*[a,b]\otimes f^*a]
\specialarrow{\via \Delta^{-1}}
[f^*a,f^*([a,b]\otimes a)]
\specialarrow{\via \ev}
[f^*a,f^*b]
\]
is denoted by $P$.
We can apply Lewis' theorem to $f^*$.
In particular, the following diagrams are commutative by
(\ref{Lewis.par}) for a morphism $f:X\rightarrow Y$.

\begin{equation}\label{Lewis-1-op.eq}
\begin{array}{ccc}
f^*[a,b]\otimes f^* a & \specialarrow{P\otimes 1} & [f^*a,f^*b]\otimes f^* a\\
\sdarrow{\Delta^{-1}}  & & \sdarrow{\ev}\\
f^*([a,b]\otimes a) & \specialarrow{f^*\ev} & f_*b
\end{array}
\end{equation}
\begin{equation}\label{Lewis-2-op.eq}
\begin{array}{ccc}
f^*a & \specialarrow{f^*\trace} & f^*[b,a\otimes b]\\
\sdarrow{\trace} &   & \sdarrow{ P}  \\
{}[f^*b,f^*a\otimes f^*b] & \specialarrow{\Delta^{-1}} & [f^*b,f^*(a\otimes b)]
\end{array}
\end{equation}
\begin{equation}\label{Lewis-3-op.eq}
\begin{array}{ccccccc}
f^* a & \specialarrow{\lambda^{-1}}& \Cal O_X\otimes f^*a &
\specialarrow{C^{-1}} & f^*\Cal O_Y\otimes f^*a &
\specialarrow{\via\trace} & f^*[a,\Cal O_Y\otimes a]\otimes f^* a\\
\sdarrow{1} & & & & & & \sdarrow{\via\lambda}\\
f^*a &
\multicolumn{3}{c}{
\mkern -2mu
\mathord\leftarrow
\mkern -10mu
\cleaders \hbox {$\mkern -2mu\mathord -\mkern -2mu$}
\hfill
\hbox to 0pt{\hss$
  \mathord{\mathop{\mkern -1mu\mathord-\mkern -1mu}\limits^{\hbox to
  0pt{\hss\scriptsize $\ev$\hss}}}$
   \hss}
\mkern -1mu
\cleaders \hbox {$\mkern -2mu\mathord -\mkern -2mu$}
\hfill
}
& [f^*a,f^*a]\otimes f^*a & \specialleftarrow{P\otimes 1} &
f^* [a,a]\otimes f^*a
\end{array}
\end{equation}

\begin{Lemma}\label{composition-H-op.thm}
Let $f:X\rightarrow Y$ and $g:Y\rightarrow Z$ be morphisms
in $\Cal S$.
Then the diagram
\[
\begin{array}{ccccc}
   f^*g^*[a,b] & \specialarrow{P} & f^*[g^*a,g^*b] & \specialarrow{P} &
   [f^*g^*a,f^*g^*b]\\
  \suarrow{d^{-1}} & & & & \suarrow{[d,d^{-1}]}\\
  (gf)^*[a,b] &
\multicolumn{3}{c}{
\mathord -\mkern -6mu
\cleaders \hbox {$\mkern -2mu\mathord -\mkern -2mu$}
\hfill
\hbox to 0pt{\hss$
  \mathord{\mathop{\mkern -1mu\mathord-\mkern -1mu}\limits^{\hbox to
  0pt{\hss\scriptsize $P$\hss}}}$
   \hss}
\mkern -1mu
\cleaders \hbox {$\mkern -2mu\mathord -\mkern -2mu$}
\hfill\mkern -6mu\mathord \rightarrow
} & [(gf)^*a,(gf)^*b]
\end{array}
\]
is commutative.
\end{Lemma}

\proof Follows instantly by Lemma~\ref{composition-H.thm}.
\qed

\begin{Lemma}\label{PH-ue.thm}
The following diagrams are commutative.
\[
\begin{array}{ccc}
[a,b] & \specialarrow{u} &[a,f_*f^*b]\\
\sdarrow{u} &  & \suarrow{u}\\
f_*f^*[a,b] & \specialarrow{HP} & [f_*f^*a,f_*f^*b]
\end{array}
\]
\[
\begin{array}{ccc}
f^*f_*[a,b] & \specialarrow{PH} & [f^*f_*a,f^*f_*b]\\
\sdarrow{\varepsilon} & & \sdarrow{\varepsilon}\\
{}[a,b] & \specialarrow{\varepsilon} & [f^*f_*a,b]
\end{array}
\]
\end{Lemma}

\proof We prove the commutativity of the first diagram.
Taking the adjoint map, it suffices to prove the commutativity of
\[
\begin{array}{ccccc}
[a,b]\otimes a & 
\multicolumn{3}{c}{
\mathord -\mkern -6mu
\cleaders \hbox {$\mkern -2mu\mathord -\mkern -2mu$}
\hfill
\hbox to 0pt{\hss$
  \mathord{\mathop{\mkern -1mu\mathord-\mkern -1mu}\limits^{\hbox to
  0pt{\hss\scriptsize $\ev$\hss}}}$
   \hss}
\mkern -1mu
\cleaders \hbox {$\mkern -2mu\mathord -\mkern -2mu$}
\hfill\mkern -6mu\mathord \rightarrow} & b\\
\sdarrow{u\otimes u} & & & & \sdarrow{u}\\
f_*f^*[a,b]\otimes f_*f^*a & \specialarrow{HP\otimes 1} &
[f_*f^*a,f_*f^*b]\otimes f_*f^*a & \specialarrow{\ev} & f_*f^*b.
\end{array}
\]
This is obvious by the naturality of $u$ and the commutativity of
(\ref{Lewis-1.eq}).

The commutativity of the second diagram is similar (use (\ref{Lewis-1-op.eq})).
\qed

\paragraph
Let $X$ be an object of $\Cal S$.
We denote the composite isomorphism
\[
M_X(a,b)\specialarrow{\via\lambda} M_X(\Cal O_X\otimes a,b)
\specialarrow{\pi} M_X(\Cal O_X,[a,b])
\]
by $h_X$.

\begin{Lemma} Let $f:X\rightarrow Y$ be a morphism in $\Cal S$.
Then the composite map
\begin{multline*}
M_X(a,b)\specialarrow{h_X}
M_X(\Cal O_X,[a,b])
\specialarrow{C} M_X(f^*\Cal O_Y,[a,b])
\cong
M_Y(\Cal O_Y,f_*[a,b])\\
\specialarrow{H}
M_Y(\Cal O_Y,[f_*a,f_*b])
\specialarrow{h_Y^{-1}}
M_Y(f_*a,f_*b)
\end{multline*}
agrees with the map given by $\varphi\mapsto f_*\varphi$.
The composite map
\begin{multline*}
M_Y(a',b')\specialarrow{h_Y}M_Y(\Cal O_Y,[a',b'])
\specialarrow{f^*}
M_X(f^*\Cal O_Y,f^*[a',b'])\\
\specialarrow{M_X(C^{-1},P)}
M_X(\Cal O_X,[f^*a',f^*b'])
\specialarrow{h_X^{-1}}
M_X(f^*a',f^*b')
\end{multline*}
agrees with $f^*$.
\end{Lemma}

\proof We prove the first assertion.
All the maps are natural on $a$.
By Yoneda's lemma, we may assume that $a=b$ and it suffices to show that
the identity map $1_b$ is mapped to $1_{f_*b}$ by the map.
It is straightforward to check that $1_b$ goes to the composite map
\begin{multline*}
f_*a\specialarrow{\lambda^{-1}}\Cal O_Y\otimes f_*a\specialarrow{u}
f_*f^*\Cal O_Y\otimes f_*a\specialarrow{C}
f_*\Cal O_X\otimes f_*a\\
\specialarrow{\via\trace}
f_*[a,\Cal O_X\otimes a]\otimes f_*a
\specialarrow{\lambda}
f_*[a,a]\otimes f_*a
\specialarrow{\ev}
f_*a.
\end{multline*}
By the commutativity of (\ref{Lewis-3.eq}), we are done.

The second assertion is proved similarly, utilizing the commutativity
of (\ref{Lewis-3-op.eq}).
\qed

Let $((?)^*,(?)_*)$ be a monoidal pair over $\Cal S$ which may not
be Lipman.
Let $\sigma=(fg'=gf')$ be a commutative diagram in $\Cal S$.

\begin{Lemma}
The following composite map agree:
\begin{description}
\item[1] $g^*f_*\specialarrow{u} g^*f_*g'_*(g')^*
\specialarrow{c} g^*g_*f'_*(g')^*\specialarrow{\varepsilon}f'_*(g')^*$;
\item[2] $g^*f_*\specialarrow{u} f'_*(f')^*g^*f_*
\specialarrow{d} f'_*(g')^*f^*f_*\specialarrow{\varepsilon}f'_*(g')^*$.
\end{description}
\end{Lemma}

For the proof and more information, see \cite[(3.7.2)]{Lipman}.
We denote the composite map in the lemma by $\theta(\sigma)$ or $\theta$.

\begin{Lemma} Let $\theta=(fg'=gf')$ be a commutative diagram in 
$\Cal S$.
Then the diagram
\[
\begin{array}{ccccc}
   g^*f_*[a,b] & 
\multicolumn{3}{c}{
\mathord -\mkern -6mu
\cleaders \hbox {$\mkern -2mu\mathord -\mkern -2mu$}
\hfill
\hbox to 0pt{\hss$
  \mathord{\mathop{\mkern -1mu\mathord-\mkern -1mu}\limits^{\hbox to
  0pt{\hss\scriptsize $\theta$\hss}}}$
   \hss}
\mkern -1mu
\cleaders \hbox {$\mkern -2mu\mathord -\mkern -2mu$}
\hfill\mkern -6mu\mathord \rightarrow
} & f'_*(g')^*[a,b]\\
\sdarrow{PH} & & & & \sdarrow{HP}\\
{}[g^*f_*a,g^*f_*b] & \specialarrow{\theta} & [g^*f_*a,f'_*(g')^*b] &
\specialleftarrow{\theta} & [f'_*(g')^*a,f'_*(g')^*b]
\end{array}
\]
is commutative.
\end{Lemma}

\proof Follows from Lemma~\ref{PH-ue.thm}.
\qed

\begin{Lemma}\label{theta-delta-m.thm}
Let $\theta=(fg'=gf')$ be a commutative diagram in 
$\Cal S$.
Then the diagram
\[
 \begin{array}{ccccc}
f^*(g_*a\otimes g_*b) & \specialarrow m & f^*g_*(a\otimes b) &
   \specialarrow\theta & g'_*(f')^*(a\otimes b)\\
\sdarrow\Delta & & & & \sdarrow{\Delta}\\
f^*g_*a\otimes f^*g_*b & \specialarrow{\theta\otimes \theta} & 
   g_*'(f')^*a\otimes g_*'(f')^*b & \specialarrow m &
   g_*'((f')^*a\otimes (f')^*b)
 \end{array}
\]
is commutative.
\end{Lemma}

\proof Utilize (\ref{unit-otimes.eq}) and (\ref{counit-otimes.eq}).
\qed

\section{Sheaves on ringed sites}

\paragraph
We fix a universe $\Cal U$, and the categories of sets, abelian groups, 
rings, etc., denoted by $\Set$, $\Ab$, $\Rng$, etc., are those consisting
of sets, abelian groups, rings, etc.\ inside $\Cal U$.
In the sequel, a $\Cal U$-small set (i.e., a set which is bijective to 
an element of  $\Cal U$)
is referred simply as a small set.

\paragraph
For categories $I$ and $\Cal C$, we denote the functor category
$\Func(I\op,\Cal C)$ by $\Cal P(I,\Cal C)$.
An object of $\Cal P(I,\Cal C)$ is sometimes referred as a presheaf over $I$
with values in $\Cal C$.
We fix a universe $\Cal V$ which contains $\Cal U$ as an element.
If $\Cal C$ is a $\Cal V$-category and $I$ is $\Cal V$-svelte 
(i.e., equivalent to a $\Cal V$-small category), 
then $\Cal P(I,\Cal C)$ is a $\Cal V$-category.
The category $\Cal P(I,\Ab)$ is denoted by $\PA(I)$.

In this article, a site 
(i.e., a category with a Grothendieck topology, in the sense of 
\cite{Verdier4})
is required to be a $\Cal V$-small $\Cal U$-site defined by a pretopology.
If $\Bbb X$ is a site, then the category of sheaves on $\Bbb X$ 
with values in $\Cal C$ is denoted by $\Cal S(\Bbb X,\Cal C)$.
If $\Cal C$ is a $\Cal U$-category, then 
$\Cal S(\Bbb X,\Cal C)$ is also a $\Cal U$-category, by our requirement.
The category $\Cal S(\Bbb X,\Ab)$ is denoted by $\AB(\Bbb X)$.

Let $\Bbb X$ be a site.
The inclusion $\AB(\Bbb X)\rightarrow\PA(\Bbb X)$ is
denoted by $q(\Bbb X,\AB)$.
We denote the sheafification functor $\PA(\Bbb X)\rightarrow \AB(\Bbb X)$ 
by $a(\Bbb X,\AB)$.

\paragraph
Let $\Bbb X=(\Bbb X,\Cal O_{\Bbb X})$ be a ringed site.
Namely, let $\Bbb X$ be a site and
$\Cal O_{\Bbb X}$ a sheaf of commutative rings on $\Bbb X$.
We denote the category of presheaves (resp.\ sheaves) of $\Cal O_{\Bbb X}
$-modules
by $\PM(\Bbb X)$ (resp.\ $\Mod(\Bbb X)$).
The inclusion $\Mod(\Bbb X)\rightarrow \PM(\Bbb X)$ is denoted
by $q(\Bbb X,\Mod)$.
The sheafification $\PM(\Bbb X)\rightarrow\Mod(\Bbb X)$ is
denoted by $a(\Bbb X,\Mod)$.
The forgetful functor $\Mod(\Bbb X)\rightarrow\AB(\Bbb X)$ is denoted by
$F(\Bbb X)$.
The categories $\AB(\Bbb X)$ and 
$\Mod(\Bbb X)$ are Grothendieck in the sense that they have small sets of
generators and satisfy the (AB5) condition in \cite{Tohoku}.
They also satisfy (AB3*).

\paragraph\label{Kan-adjoint.par}
Let $f:\Bbb Y\rightarrow \Bbb X$ be a functor.
Then the pull-back $\PA(\Bbb X)\rightarrow \PA(\Bbb Y)$
is denoted by $f^{\#}_{\PA}$.
Note that $f^{\#}_{\PA}(\Cal F):=\Cal F\circ f\op$.
In general, the pull-back $\Cal P(\Bbb X,\Cal C)\rightarrow\Cal P(
\Bbb Y,\Cal C)$ is defined in a similar way, and is denoted by $f^{\#}$.
If $f$ is continuous (i.e., $f^{\#}_{\Set}$ carries sheaves to sheaves), 
then $f^{\#}_{\AB}:\AB(\Bbb X)\rightarrow \AB(\Bbb Y)$ is
defined to be the restriction of $f^{\#}_{\PA}$.
Throughout these notes, we require that a continuous functor $f:\Bbb Y
\rightarrow \Bbb X$ satisfies the following condition.
For $y\in\Bbb Y$, a covering $(y_i\rightarrow y)_{i\in I}$, and 
any $i,j\in I$, the morphisms 
$f(y_i\times_y y_j)\rightarrow
f(y_i)$ and 
$f(y_i\times_y y_j)\rightarrow
f(y_j)$ makes $f(y_i\times_y y_j)$ the fiber product $f(y_i)\times_{f(y)}
f(y_j)$.

Thanks to the re-definition of sites and continuous functors,
we have the following.

\begin{Lemma}
Let $f:\Bbb Y\rightarrow \Bbb X$ be a functor between sites.
Then $f$ is continuous if and only if the following holds.

If $(\varphi_i: y_i\rightarrow y)_{i\in I}$ is a covering,
then $(f\varphi_i: fy_i\rightarrow fy)_{i\in I}$ is a covering.
\end{Lemma}

For the proof, see \cite[(III.1.6)]{SGA4-I}.

The left adjoint of $f^{\#}_{\PA}$, which exists by Kan's lemma
(see e.g., \cite[Theorem~I.2.1]{Artin}), is denoted by $f_{\#}^{\PA}$.

For $x\in \Bbb X$, we define the category $I_x^f$ as follows.
An object of $I_x^f$ is a pair $(y,\phi)$ with $y\in\Bbb Y$ and
$\phi\in \Bbb X(x,f(y))$.
A morphism $h:(y,\phi)\rightarrow (y',\phi')$ is a morphism $h\in
\Bbb Y(y,y')$ such that $f(h)\circ \phi=\phi'$.
Note that $\Gamma(x,f_{\#}^{{\PA}}(\Cal F))=\indlim \Gamma(y,\Cal F)$,
where the colimit is taken over $(I_x^f)\op$.

The left adjoint of $f^{\#}_{\Cal C}:\Cal P(\Bbb X,\Cal C)\rightarrow
\Cal P(\Bbb X,\Cal C)$ is constructed similarly,
provided $\Cal C$ is a $\Cal U$-category with arbitrary small colimits.
The left adjoint is denoted by $f_{\#}^{\Cal C}$ or simply by $f_{\#}$.

Similarly, $f_{\#}^{\AB}:\AB(\Bbb X)\rightarrow \AB(\Bbb Y)$
and its left adjoint $f^{\#}_{\AB}$ is defined.
Note that we have $q(\Bbb Y,\AB)\circ f^{\#}_{\AB}=f^{\#}_{\PA}\circ
q(\Bbb X,\AB)$, and $f_{\#}^{\AB}=a(\Bbb X,\AB)\circ f_{\#}^{\PA}\circ
q(\Bbb Y,\AB)$.

\begin{Lemma}\label{our-site.thm}
If $(I_x^f)\op$ is pseudofiltered \(see e.g., {\rm \cite{SGA4-I,Milne}}\) 
for each 
$x\in\Bbb X$, then 
$f_{\#}^{\PA}$ is exact.
\end{Lemma}

\proof 
This is a consequence of \cite[Proposition~2.8]{SGA4-I}.
\qed

We say that $f:\Bbb Y\rightarrow \Bbb X$ 
is {\em admissible} if $f$ is continuous and the functor
$f_{\#}^{\PA}$ is exact.

\paragraph 
If $\Bbb Y$ has finite limits and $f$ preserves finite limits, then
$f_{\#}^{\PA}$ is exact by the lemma.
It follows that a continuous map induces an admissible 
continuous functor between the corresponding sites.

\paragraph\label{right-adjoint-direct.par}
The right adjoint functor of $f^{\#}_{\PA}$, which we denote by
$f_{\flat}^{\PA}$ also exists, as $\Ab\op$ has
arbitrary small colimits (i.e., $\Ab$ has small limits).
The functor $f_{\flat}^{\PA}$ 
is the left adjoint functor 
\[
\Func(\Bbb Y\op,\Ab)
\specialarrow{{\rm op}}
\Func(\Bbb Y,\Ab\op)
\rightarrow
\Func(\Bbb X,\Ab\op)\specialarrow{{\rm op}}
\Func(\Bbb X\op,\Ab)
\]
of the functor
\[
(f\op)^{\#}:\Func(\Bbb X,\Ab\op)\rightarrow\Func(\Bbb Y,\Ab\op),
\]
where $f\op=f$ is the opposite of $f$, namely, $f$ viewed as a functor 
$\Bbb Y\op\rightarrow \Bbb X\op$.

\paragraph
For $\Cal M,\Cal N\in\PM(\Bbb X)$, the presheaf tensor product is denoted
by $\Cal M\otimes^p_{\Cal O_{\Bbb X}}\Cal N$.
It is defined by
\[
\Gamma(x,\Cal M\otimes^p_{\Cal O_{\Bbb X}}\Cal N):=
\Gamma(x,\Cal M)\otimes_{\Gamma(x,\Cal O_{\Bbb X})}\Gamma(x,\Cal N)
\]
for $x\in \Bbb X$.
The sheaf tensor product $a(q\Cal M\otimes^p_{\Cal O_{\Bbb X}}q\Cal N)$ 
of $\Cal M,\Cal N\in \Mod(\Bbb X)$ 
is denoted by $\Cal M\otimes_{\Cal O_{\Bbb X}}\Cal N$.

\paragraph
We say that 
$f:(\Bbb Y,\Cal O_{\Bbb Y})\rightarrow (\Bbb X,\Cal O_{\Bbb X})$ is
a continuous ringed functor
if $f:\Bbb Y\rightarrow \Bbb X$ is a continuous functor, 
and a morphism of sheaves of commutative rings
$\Cal O_{\Bbb Y}\rightarrow f^{\#}_{\AB}(\Cal O_{\Bbb X})$ is given.
The pull-back $\PM(\Bbb X)\rightarrow\PM(\Bbb Y)$ is denoted by
$f^{\#}_{{\PM}}$, and its left adjoint is denoted by $f_{\#}^{{\PM}}$.
The left adjoint $f_{\#}^{{\PM}}$ is defined by
\[
\Gamma(x,f_{\#}^{{\PM}}\Cal M):=\indlim \Gamma(x,\Cal O_{\Bbb X})
\otimes_{\Gamma(y,\Cal O_{\Bbb Y})}\Gamma(y,\Cal M)
\]
for $x\in\Bbb X$ and $\Cal M\in\PM(\Bbb Y)$,
where the colimit is taken over the category $(I^f_x)\op$.
Similarly, $f^{\#}_{\Mod}:\Mod(\Bbb X)\rightarrow\Mod(\Bbb Y)$ and
its left adjoint $f_{\#}^{\Mod}=af^{\#}_{\PM}q$ is defined.
If $(I^f_x)\op$ is filtered for any $x\in\Bbb X$, then $f_\#^{\PA}
\Cal O_{\Bbb Y}
$ has a structure of a presheaf of rings in a natural way, and
there is a canonical isomorphism $f_\#^{\PM}\Cal M\cong \Cal O_{\Bbb X}
\otimes^p_{f_\#^{\PA}\Cal O_{\Bbb Y}}f_\#^{\PA}\Cal M$.
The right adjoint of $f^{\#}_{\PM}$, which exists as in 
(\ref{right-adjoint-direct.par}), is denoted by $f_\flat^{\PM}$.

\paragraph
Let $\Bbb X$ be a ringed site, and $x\in\Bbb X$.
The category $\Bbb X/x$ is a site with the same topology as that of $\Bbb X$.
The canonical functor $\Phi_x:\Bbb X/x\rightarrow \Bbb X$
is a continuous ringed functor, and yields the 
pull-backs $(\Phi_x)^\#_{\AB}$ and $(\Phi_x)^\#_{\PA}$, which we denote by
$(?)|_x^{\AB}$ and $(?)|_x^{\PA}$, respectively.
Their left adjoints are denoted by $L_x^{\AB}$ and $L_x^{\PA}$, respectively.
Note that 
$\Phi_x$ is admissible, see \cite[p.78]{Milne}.

\paragraph\label{L_x-exact.par}
Note that $\Bbb X/x$ is a ringed site with the structure sheaf $
\Cal O_{\Bbb X}|_x$.
Thus, $(?)|_x^{\Mod}$ and $(?)|_x^{\PM}$ are defined in an obvious way,
and their left  adjoints $L_x^{\Mod}$ and $L_x^{\PM}$ are also defined.
Note that $L_x^{\Mod}$ and $L_x^{\PM}$ are exact.

\paragraph
For a morphism $\phi:x\rightarrow y$, we have an obvious admissible
ringed continuous functor $\Phi_\phi:\Bbb X/x\rightarrow \Bbb X/y$.
The corresponding pull-back is denoted by $\phi^\star_\natural$,
and its left adjoint is denoted by $\phi_\star^\natural$,
where $\natural$ is $\AB$, $\PA$, $\Mod$ or $\PM$.

For $\Cal M,\Cal N\in\natural(\Bbb X)$, we define $
\uHom_{\natural(\Bbb X)}(\Cal M,
\Cal N)
$ to be the object of $\natural(\Bbb X)$ given by
\[
\Gamma(x,\uHom_{\natural(\Bbb X)}(\Cal M,\Cal N)):=
\Hom_{\natural(\Bbb X/x)}(\Cal M
|_x^{\natural},\Cal N|_x^{\natural}),
\]
where $\natural={\PA},{\AB},{\PM},$ or ${\Mod}$.
For $\phi:x\rightarrow y$, the restriction map 
\[
\Hom_{\natural(\Bbb X/y)}(\Cal M
|_y^{\natural},\Cal N|_y^{\natural})
\rightarrow
\Hom_{\natural(\Bbb X/x)}(\Cal M
|_x^{\natural},\Cal N|_x^{\natural})
\]
is given by $\phi^\star_{\natural}$.

The following is a slight generalization of \cite[Example~3.5.2]{Lipman}.

\begin{Lemma}\label{lipman-sites-sheaf.thm}
Let $(\Bbb X,\Cal O_{\Bbb X})$ be a ringed site.
Then we have the following.
\begin{description}
\item[1] The category $\PM(\Bbb X)$ is a closed symmetric monoidal 
category \(see {\rm \cite[p.180]{Mac Lane}}\) 
with $\otimes_{\Cal O_{\Bbb X}}^p$ the
multiplication, $q\Cal O_{\Bbb X}$ the unit object, 
$\uHom_{\PM(\Bbb X)}(?,?)$ the internal hom, etc., etc.
\item[2] The category $\Mod(\Bbb X)$ is a closed symmetric monoidal category
with $\otimes_{\Cal O_{\Bbb X}}$ the multiplication,
$\Cal O_{\Bbb X}$ the unit object,
$\uHom_{\Mod(\Bbb X)}(?,?)$ the internal hom, etc., etc.
\item[3] The inclusion $q:\Mod(\Bbb X)\rightarrow \PM(\Bbb X)$ and
the natural transformations
\[
q\Cal M\otimes^p_{\Cal O_{\Bbb X}}q\Cal N
\specialarrow{u}
qa(q\Cal M\otimes^p_{\Cal O_{\Bbb X}}q\Cal N)
=q(\Cal M\otimes_{\Cal O_{\Bbb X}}\Cal N)
\]
and 
\[
q\Cal O_{\Bbb X}\specialarrow{\id}q\Cal O_{\Bbb X}
\]
form a  symmetric monoidal functor.
\item[4] For a ringed continuous functor
$f:(\Bbb Y,\Cal O_{\Bbb Y})\rightarrow (\Bbb X,\Cal O_{\Bbb X})$,
the functor
$f^{\#}_{{\PM}}:\PM(\Bbb X)\rightarrow \PM(\Bbb Y)$,
the natural map
\begin{equation}\label{pull-back-monoidal.eq}
f^{\#}_{{\PM}}\Cal M\otimes^p_{\Cal O_{\Bbb Y}}f^{\#}_{{\PM}}\Cal N
\rightarrow
f^{\#}_{{\PM}}(\Cal M\otimes_{\Cal O_{\Bbb X}}^p \Cal N)
\end{equation}
which induces the canonical map
\[
\Gamma(fy,\Cal M)\otimes_{\Gamma(y,\Cal O_{\Bbb Y})}\Gamma(fy,\Cal N)
\rightarrow
\Gamma(fy,\Cal M)\otimes_{\Gamma(fy,\Cal O_{\Bbb X})}\Gamma(fy,\Cal N)
\]
for each $y\in\Bbb Y$,
and the map
\[
q\Cal O_{\Bbb Y}\rightarrow q f^{\#}_{{\Mod}}\Cal O_{\Bbb X}\cong
f^{\#}_{{\PM}}q\Cal O_{\Bbb X}
\]
form a symmetric monoidal functor.
\item[5] Under the same assumption as in {\bf 4}, 
the functor
$f^{\#}_{{\Mod}}:\Mod(\Bbb X)\rightarrow \Mod(\Bbb Y)$,
the natural map
\begin{multline*}
f^{\#}_{{\Mod}}\Cal M\otimes_{\Cal O_{\Bbb Y}}f^{\#}_{{\Mod}}\Cal N
=
a(q f^{\#}_{{\Mod}}\Cal M\otimes_{\Cal O_{\Bbb Y}}^p
q f^{\#}_{{\Mod}}\Cal N)
\cong
a(f^{\#}_{{\PM}}q\Cal M\otimes_{\Cal O_{\Bbb Y}}^p
f^{\#}_{{\PM}}q\Cal N)\\
\specialarrow{\via \Mathrm(\ref{pull-back-monoidal.eq})}
af^{\#}_{{\PM}}(q\Cal M\otimes_{\Cal O_{\Bbb X}}^p q\Cal N)
\specialarrow{\via \rho}
f^{\#}_{{\Mod}}a(q\Cal M\otimes_{\Cal O_{\Bbb X}}^p q\Cal N)
=f^{\#}_{{\Mod}}(\Cal M\otimes_{\Cal O_{\Bbb X}}\Cal N),
\end{multline*}
and the map
\[
\Cal O_{\Bbb Y}\rightarrow f^{\#}_{{\Mod}}\Cal O_{\Bbb X}
\]
form a symmetric monoidal functor,
where $\rho:af^{\#}_{{\PM}}\rightarrow f^{\#}_{{\Mod}}a$ is the composite
map
\[
af^{\#}_{{\PM}}\specialarrow{\via u}af^{\#}_{{\PM}}qa
\cong aqf^{\#}_{{\Mod}}a\specialarrow{\via \varepsilon}
f^{\#}_{{\Mod}}a.
\]
\item[6] Let $\Cal S$ denote the category of ringed sites \(with svelte
underlying categories\).
Then $((?)_{\#}^{{\PM}},(?)^{\#}_{{\PM}})$ and
$((?)_{\#}^{{\Mod}},(?)^{\#}_{{\PM}})$ are adjoint pairs of
monoidal pseudo-functors on $\Cal S\op$, see {\rm \cite[(3.6.7)]{Lipman}}.
\end{description}
\end{Lemma}

\paragraph
If there is no confusion, 
the $\natural$ attached to the
functors of sheaves defined above are omitted.
For example, $f^{\#}$ stands for $f^{\#}_\natural$.
Note that $f^{\#}_\natural(\Cal F)$ viewed as a presheaf of abelian
groups is independent of $\natural$.

\paragraph Let $f:\Bbb Y\rightarrow \Bbb X$ be a ringed continuous functor
between ringed sites.
The following is a restricted version of the results on cocontinuous
functors in \cite{Verdier2}.
We give a proof for convenience of readers.

\begin{Lemma}\label{relative-topology.thm}
Assume that for any
$y\in \Bbb Y$ and any covering $(x_\lambda\rightarrow fy)_{\lambda\in\Lambda}
$ of $fy$, there is
a covering $(y_\mu\rightarrow y)_{\mu\in M}$ of $y$ such that there is
a map $\phi:M\rightarrow\Lambda$ such that $f y_\mu\rightarrow fy$ factors 
through $x_{\phi(\mu)}\rightarrow fy$ for each $\mu$.
Then 
the pull-back $f^{\#}$ is compatible with the sheafification in the
sense that the canonical natural transformation
\[
\rho:a(\Bbb Y,\Mod)f^{\#}_{\PM}\rightarrow f^{\#}_{\Mod}a(\Bbb X,\Mod)
\]
is a natural isomorphism.
If this is the 
case, $f^{\#}_{\Mod}$ has the right adjoint $f_\flat^{\Mod}$,
and in particular, it preserves arbitrary limits and arbitrary colimits.
\end{Lemma}

\proof Let $\Cal M\in\PM(\Bbb X)$ and $y\in \Bbb Y$.
Then the canonical map
\begin{multline*}
\cc^0(y,f^{\#}_{\PM}\Cal M)=\indlim \cc^0((y_\mu\rightarrow y),f^{\#}_{\PM}
\Cal M)\\
=\indlim \cc^0((f y_\mu\rightarrow fy),\Cal M)
\rightarrow \cc^0(fy,\Cal M)
\end{multline*}
is an isomorphism, because the coverings of the form $(f y_\mu\rightarrow fy)$
is cofinal in the set of all coverings of $fy$ by assumption.
Thus the canonical map
\[
\ccs^0(f^{\#}_{\PM}\Cal M)\rightarrow f^{\#}_{\PM}\ccs^0(\Cal M)
\]
is an isomorphism.
The first
assertion follows immediately from the construction of the sheafification,
see \cite[(II.1)]{Artin}.

Next we show that $f^{\#}_{\Mod}$ has a right adjoint.
To prove this, it suffices to show that $f_\flat^{\PM}(\Cal M)$ is a 
sheaf if so is $\Cal M$ for $\Cal M\in\PM(\Bbb X)$.
Let $u:\Id_{\PM(\Bbb Y)}\rightarrow f^{\#}_{\PM} f_\flat^{\PM}$ be the
unit of adjunction,
$\varepsilon: f_\flat^{\PM}f^{\#}_{\PM}\rightarrow \Id_{\PM(\Bbb X)}$
the counit of adjunction,
$v(\Bbb X):q(\Bbb X,\Mod)a(\Bbb X,\Mod)\rightarrow \Id_{\PM(\Bbb X)}$ 
the unit of adjunction, and
$v(\Bbb Y):q(\Bbb Y,\Mod)a(\Bbb Y,\Mod)\rightarrow \Id_{\PM(\Bbb Y)}$ 
the unit of adjunction.
Then the diagram of functors
\[
\!\!\!\!\!\!\!\!\!
\begin{array}{ccccccc}
f_\flat^{\PM} & \specialarrow{uf_\flat^{\PM}} 
& 
f_\flat^{\PM} f^{\#}_{\PM} f_\flat^{\PM}& \specialarrow{\id}
&
f_\flat^{\PM} f^{\#}_{\PM} f_\flat^{\PM}& \specialarrow{f_\flat^{\PM}
\varepsilon}
&
f_\flat^{\PM}\\
\downarrow\hbox to 0pt{\scriptsize $v(\Bbb X)f_\flat^{\PM}$\hss}
& &
\downarrow\hbox to 0pt{\scriptsize
$f_\flat^{\PM} f^{\#}_{\PM} v(\Bbb X)f_\flat^{\PM}$\hss}
& &
\downarrow\hbox to 0pt{\scriptsize
$f_\flat^{\PM} v(\Bbb Y)f^{\#}_{\PM} f_\flat^{\PM}$\hss}
& &
\downarrow\hbox to 0pt{\scriptsize
$f_\flat^{\PM} v(\Bbb Y)$\hss}\\
qaf_\flat^{\PM} & \specialarrow{uqaf_\flat^{\PM}} &
f_\flat^{\PM}f^{\#}_{\PM}qaf_\flat^{\PM} & \specialarrow{\cong} &
f_\flat^{\PM}qa f^{\#}_{\PM}f_\flat^{\PM} & \specialarrow{f_\flat^{\PM}qa
\varepsilon} &
f_\flat^{\PM}qa
\end{array}
\]
is commutative, where $\cong$ is the {\em inverse} of the canonical map
caused by $qa f^{\#}_{\PM}\rightarrow f^{\#}_{\PM}qa$,
which exists by the first part.
As $(f^\#_{\PM},f_\flat^{\PM})$ is an adjoint pair, the
composite of the first row of the diagram is the identity.
As $v(\Bbb Y)(\Cal M)$ is an isomorphism,
the right-most vertical arrow evaluated at $\Cal M$ is an isomorphism.
Hence, $v(\Bbb X)f_\flat^{\PM}(\Cal M)$,
which is the left-most vertical arrow
evaluated at $\Cal M$, 
is a split monomorphism.
As it is a direct summand of a sheaf,
$f_\flat^{\PM}(\Cal M)$ is a sheaf, as desired.

As it is a right adjoint of $f_{\#}^{\Mod}$, 
the functor $f^{\#}_{\Mod}$ preserves arbitrary limits.
As it is a left adjoint of $f_{\flat}^{\Mod}$, 
the functor $f^{\#}_{\Mod}$ preserves arbitrary colimits.
\qed

\paragraph\label{restriction-LR.par}
Let $\Bbb X$ be a ringed site, and $x\in\Bbb X$.
It is easy to see that $\Phi_x: \Bbb X/x \rightarrow \Bbb X$ satisfies the
condition in Lemma~\ref{relative-topology.thm}.
So $(?)|^{\Mod}_x$ preserves arbitrary limits and colimits.
In particular, $(?)|^{\Mod}_x$ is exact.
Similarly, for a morphism $\phi:x\rightarrow y$ in $\Bbb X$, 
$\phi^*_{\Mod}$ preserves arbitrary limits and colimits.

\section[Derived category and derived functors of sheaves]{Derived 
category and derived functors of sheaves on ringed sites}

We utilize the notation and terminology on triangulated categories in
\cite{Verdier}.
However, 
We usually write the suspension (translation) functor of a triangulated
category by $\Sigma$ or $(?)[1]$.

Let $\Cal T$ be a triangulated category.

\begin{Lemma}
Let 
\[
(a_\lambda\specialarrow{f_\lambda} 
b_\lambda\specialarrow{g_\lambda} 
c_\lambda\specialarrow{h_\lambda}\Sigma 
a_\lambda)
\]
be a small family of distinguished triangles in $\Cal T$.
Assume that the coproducts $\bigoplus a_\lambda$, $\bigoplus b_\lambda$,
and $\bigoplus c_\lambda$ exist.
Then the triangle
\[
\bigoplus a_\lambda\specialarrow{\bigoplus f_\lambda}
\bigoplus b_\lambda
\specialarrow{\bigoplus g_\lambda}
\bigoplus c_\lambda \specialarrow{H\circ \bigoplus h_\lambda}
\Sigma (\bigoplus a_\lambda)
\]
is distinguished, where $H:\bigoplus\Sigma a_\lambda\rightarrow
\Sigma(\bigoplus a_\lambda)$ is the canonical isomorphism.
Similarly, a product of distinguished triangles is a distinguished 
triangle.
\end{Lemma}

We refer the reader to \cite[Proposition~1.2.1]{Neeman} for the proof
of the second assertion.
The proof for the first assertion is similar \cite[Remark~1.2.2]{Neeman}.

\paragraph Let $\Cal A$ be an abelian category.
The category of bounded below (resp.\ bounded above, bounded)
complexes in $\Cal A$ is denoted by $C^+(\Cal A)$ (resp.\ $C^-(\Cal A)$, 
$C^b(\Cal A)$).
The corresponding homotopy category and the derived category are
denoted by $K^?(\Cal A)$ and $D^?(\Cal A)$, where $?$ is either $+$, $-$ or
$b$.
We denote the homotopy category of complexes in $\Cal A$ with bounded below
(resp.\ bounded above, bounded) cohomology groups
by $\bar K^?(\Cal A)$.
The corresponding derived category is denoted by $\bar D^?(\Cal A)$.

For a plump subcategory $\Cal A'$ of $\Cal A$, 
we denote by $K_{\Cal A'}^?(\Cal A)$ (resp.\ $\bar K_{\Cal A'}^?(\Cal A)$)
the full subcategory of $K^?(\Cal A)$ (resp.\ $\bar K^?(\Cal A)$)
consisting of complexes whose cohomology groups are objects of $\Cal A'$.
The localization of $K^?_{\Cal A'}(\Cal A)$ by the \'epaisse subcategory
of exact complexes is denoted by $D^?_{\Cal A'}(\Cal A)$.
The category $\bar D^?_{\Cal A'}(\Cal A)$ is defined similarly.
Note that the canonical functor $D^?_{\Cal A'}(\Cal A)\rightarrow
D(\Cal A)$ is fully faithful, and hence 
$D^?_{\Cal A'}(\Cal A)$ is identified with 
the full subcategory of $D(\Cal A)$
consisting of bounded below
(resp.\ bounded above, bounded) complexes whose 
cohomology groups are in $\Cal A'$.
Note also that the canonical functor 
$D^?_{\Cal A'}(\Cal A)\rightarrow \bar D^?_{\Cal A'}(\Cal A)$ is an 
equivalence.

\paragraph Let $\Cal A$ and $\Cal B$ be abelian categories, 
and $F:K^?(\Cal A)\rightarrow K^*(\Cal B)$ a triangulated functor.
Let $\Cal C$ be a triangulated subcategory of $K^?(\Cal A)$ such that
\begin{description}
\item[1] If $c\in \Cal C$ is exact, then $Fc$ is exact.
\item[2] For any $a\in K^?(\Cal A)$, there exists some quasi-isomorphism
$a\rightarrow c$.
\end{description}
The condition {\bf 2} implies that the canonical functor
\[
i(\Cal C):\Cal C/(\Cal E\cap \Cal C)\rightarrow K^?(\Cal A)/\Cal E=D^?(\Cal A)
\]
is an equivalence, where $\Cal E$ denotes the \'epaisse subcategory of
exact complexes in $K^?(\Cal A)$, see \cite[(2.3)]{Verdier}.
We fix a quasi-inverse $p(\Cal C):D^?(\Cal A)\rightarrow \Cal C/
(\Cal E\cap \Cal C)$.
On the other hand, the composite
\begin{equation}\label{der-1.eq}
\Cal C\hookrightarrow K^?(\Cal A)\specialarrow{F} K^*(\Cal B)\specialarrow Q
D^*(\Cal B)
\end{equation}
is factorized as
\begin{equation}\label{der-2.eq}
\Cal C\specialarrow {Q_{\Cal C}} \Cal C/(\Cal E\cap \Cal C)\specialarrow{\Phi}
D^*(\Cal B)
\end{equation}
up to a unique natural isomorphism, by the universality of localization
and the condition {\bf 1} above.
Under the setting above, we have the following \cite[(I.5.1)]{Hartshorne}.

\begin{Lemma}\label{existence-derived.thm}
The composite functor
\[
RF: D^?(\Cal A)\specialarrow{p(\Cal C)}
\Cal C/(\Cal E\cap \Cal C)\specialarrow{\Phi}D^*(\Cal B)
\]
is a right derived functor of $F$.
\end{Lemma}

We denote the map $QF\rightarrow(RF)Q$ corresponding to 
$\id: RF\rightarrow RF$ by $\Xi$ or $\Xi(F)$.

\paragraph
Here we are going to review Spaltenstein's work on unbounded derived 
categories \cite{Spaltenstein}.

A chain complex $I$ of $\Cal A$ is called $K$-injective if
for any exact sequence $E$ of $\Cal A$, the complex of abelian groups
$\Hom^\bullet_{\Cal A}(E,I)$ is also exact.

A chain map $f:C\rightarrow I$ is called a $K$-injective resolution of $C$,
if $I$ is $K$-injective and $f$ is a quasi-isomorphism.

A special case of the following was proved in 
\cite{Spaltenstein} (the same proof works for the general case without
much modification).
See also \cite{AJS}.

\begin{Lemma} If $\Cal A$ is Grothendieck, then for any chain complex
$\Bbb F\in K(\Cal A)$ there is a $K$-injective resolution
$\Bbb F\rightarrow \Bbb I$ of $\Bbb F$ such that each term of $\Bbb I$
being an injective object of $\Cal A$.
\end{Lemma}

Let us denote the homotopy category of chain complexes of $\Cal A$
by $K(\Cal A)$.
A chain complex $I$ is $K$-injective if and only if $K(\Cal A)(E,I)=0$
for any exact sequence $E$.
It is easy to see that the $K$-injective complexes form an \'epaisse
subcategory $I(\Cal A)$ of $K(\Cal A)$.

\paragraph Let $F:K(\Cal A)\rightarrow K(\Cal B)$ be a triangulated 
functor, and assume that $\Cal A$ is Grothendieck.
Let $\Cal I$ be the full subcategory of $K$-injective complexes of $K(\Cal A)$.
It is easy to see that $\Cal I$ is triangulated, and $\Cal I\cap \Cal E=0$.
By (\ref{existence-derived.thm}), the composite
\[
D(\Cal A)\specialarrow{p(\Cal I)}\Cal I\specialarrow \Phi D(\Cal B)
\]
is a right derived functor $RF$ of $F$.
Note that to fix $p(\Cal I)$ is nothing but to fix 
a functorial $K$-injective resolution $\Bbb F\rightarrow pQ \Bbb F=\Bbb 
I_{\Bbb F}$.

\paragraph Let $F:K(\Cal A)\rightarrow K(\Cal B)$ be a triangulated
functor, and assume that there is a right derived functor of $F$.
For $\Bbb F\in K(\Cal A)$, we say that $\Bbb F$ is (right) $F$-acyclic
if the canonical map $QF\Bbb F\rightarrow (RF)Q\Bbb F$ is an isomorphism.

\begin{Lemma}\label{derived-adjoint.thm}
Let $\Cal A$ and $\Cal B$ be abelian categories, 
and $F:\Cal A\rightarrow \Cal B$ be an exact functor
with the right adjoint $G$.
Assume that $\Cal B$ is Grothendieck.
Then $KG: K(\Cal B)\rightarrow K(\Cal A)$ preserves $K$-injective complexes.
Moreover, $RG:D(\Cal B)\rightarrow D(\Cal A)$ is the right adjoint of
$RF=F$.
\end{Lemma}

\proof Let $\Bbb M\in K(\Cal A)$, and $\Bbb I$ a $K$-injective
complex of $K(\Cal B)$.
Then 
\begin{multline*}
\Hom_{K(\Cal A)}(\Bbb M,(KG)\Bbb I)
\cong H^0(\Hom_{\Cal A}^\bullet(\Bbb M,G\Bbb I))\\
\cong H^0(\Hom_{\Cal B}^\bullet(F\Bbb M,\Bbb I))
=\Hom_{K(\Cal B)}(F\Bbb M,\Bbb I).
\end{multline*}
If $\Bbb M$ is exact, then the last complex is zero.
This shows $(KG)\Bbb I$ is $K$-injective.

Now let $\Bbb M\in D(\Cal A)$ and $\Bbb N\in D(\Cal B)$ be arbitrary.
Then by the first part, we have a functorial isomorphism
\begin{multline*}
\Hom_{D(\Cal A)}(\Bbb M,(RG)\Bbb N)\cong
\Hom_{K(\Cal A)}(\Bbb M,(KG)\Bbb I_{\Bbb N})\\
\cong
\Hom_{K(\Cal B)}(F\Bbb M,\Bbb I_{\Bbb N})
\cong
\Hom_{D(\Cal B)}(F\Bbb M,\Bbb N).
\end{multline*}
This proves the last assertion.
\qed

\begin{Rem}
Note that for an abelian category $\Cal A$, we have
$\ob(C(\Cal A))=\ob(K(\Cal A))=\ob(D(\Cal A))$.
Thus, an object of one of the three categories is sometimes viewed as
an object of another.
\end{Rem}

\paragraph Let $\Cal A$ be a closed symmetric monoidal abelian category 
which satisfies the (AB3) and (AB3*) conditions.
Let $\otimes$ be the multiplication and $[?,?]$ be the internal hom.
For a fixed $b\in \Cal A$, $(?\otimes b,[b,?])$ is an adjoint pair.
In particular, $?\otimes b$ preserves colimits, and $[b,?]$ preserves limits.
By symmetry, $a\otimes ?$ also preserves colimits.
As we have an isomorphism
\[
\Cal A(a,[b,c])\cong \Cal A(a\otimes b,c)\cong 
\Cal A(b\otimes a,c)\cong
\Cal A(b,[a,c])\cong
\Cal A\op([a,c],b),
\]
we have $[?,c]:\Cal A\op\rightarrow \Cal A$ is right adjoint to
$[?,c]:\Cal A\rightarrow \Cal A\op$.
This shows that $[?,c]$ changes colimits to limits.

As in \cite{Hartshorne2}, we define the tensor product 
$\Bbb F\otimes^\bullet \Bbb G$
of $\Bbb F,\Bbb G\in
C(\Cal A)$ 
by
\[
(\Bbb F\otimes^\bullet\Bbb G)^n:=\bigoplus_{p+q=n}\Bbb F^p\otimes \Bbb G^q
\]
and 
\[
d^n=d_{\Bbb F}\otimes 1+(-1)^n 1\otimes d_{\Bbb G}.
\]
We have $\Bbb F\otimes^\bullet\Bbb G\in C(\Cal A)$.
Similarly, $[\Bbb F,\Bbb G]^\bullet\in C(\Cal A)$ is defined by
\[
[\Bbb F,\Bbb G]^n:=\prod_{p\in\Bbb Z}[\Bbb F^p,\Bbb G^{n+p}]
\]
and 
\[
d^n:=[d_{\Bbb F},1]+(-1)^{n+1}[1,d_{\Bbb G}].
\]

It is straightforward to prove the following.

\begin{Lemma}
Let $\Cal A$ be as above.
Then the category of chain complexes $C(\Cal A)$ is closed symmetric monoidal
with $\otimes^\bullet$ the multiplication and $[?,?]^\bullet$ the
internal hom.
The bi-triangulated functors
\[
\otimes^\bullet: K(\Cal A)\times K(\Cal A)\rightarrow K(\Cal A)
\]
and 
\[
[?,?]^\bullet: K(\Cal A)\op\times K(\Cal A)\rightarrow K(\Cal A),
\]
are induced, and $K(\Cal A)$ is a closed symmetric monoidal triangulated
category \(see {\rm \cite[(3.5), (3.6)]{Lipman}}\).
\end{Lemma}

\paragraph Let $\Cal A$ be an abelian category, and $\frak P$ a full
subcategory of $C(\Cal A)$.
An inverse system $(\Bbb F_i)_{i\in I}$ in $C(\Cal A)$ is said to be
{\em $\frak P$-special} if the following conditions are satisfied.
\begin{description}
\item[i] $I$ is well-ordered.
\item[ii] If $i\in I$ has no predecessor, then the canonical map
$I_i\rightarrow \projlim I_{j<i}$ is an isomorphism
(in particular, $I_{i_0}=0$ if $i_0$ is the minimum element of $I$).
\item[iii] If $i\in I$ has a predecessor $i-1$, then the natural chain map
$I_{i}\rightarrow I_{i-1}$ is an epimorphism, the kernel $C_i$ is
isomorphic to some object of $\frak P$, and the exact sequence
\[
0\rightarrow C_i\rightarrow I_i\rightarrow I_{i-1}\rightarrow0
\]
is semi-split.
\end{description}

Similarly, $\frak P$-special direct systems are also defined, see
\cite{Spaltenstein}.

The full subcategory of $C(\Cal A)$ consisting of inverse (resp.\ direct) 
limits of $\frak P$-special inverse (resp.\ direct) systems 
is denoted by $\sublim{\frak P}$ (resp.\ $\subcolim{\frak P}$).

\paragraph Let $(\Bbb X,\Cal O_{\Bbb X})$ be a ringed site.
Various definitions and results on unbounded complexes 
of sheaves over a ringed space by Spaltenstein \cite{Spaltenstein}
is generalized to those for ringed sites.
However, note that we can not utilize the notion related to closed subsets,
points, or stalks of sheaves.

\paragraph
We say that a complex $\Bbb F\in C(\Mod(\Bbb X))$ is {\em $K$-flat} if
$\Bbb G\otimes^\bullet \Bbb F$ is exact whenever $\Bbb G$ is an exact
complex in $\Mod(\Bbb X)$.
We say that $\Bbb A\in  C(\Mod(\Bbb X))$ is {\em weakly $K$-injective} if
$\Bbb A$ is 
$\Hom^\bullet_{\Mod(\Bbb X)}(\Bbb F,?)$-acyclic for any $K$-flat complex
$\Bbb F$.

\paragraph Let $(\Bbb X,\Cal O_{\Bbb X})$ be a ringed site.
For $x\in \Bbb X$, we define $\Cal O_x^p$ to be $L_x^{{\PM}}(\Cal O_{
\Bbb X}|_x)$, and $\Cal O_x:=L_x^{{\Mod}}(\Cal O_{\Bbb X}|x)=a\Cal O^p_x$.
For $x_1,\ldots,x_r\in \Bbb X$, we define
\[
\Cal O_{x_1,\ldots,x_r}:=\Cal O_{x_1}\otimes\cdots\otimes\Cal O_{x_r}.
\]
We denote by $\frak P_0=\frak P_0(\Bbb X,\Cal O_{\Bbb X})$ 
(resp.\ $\tilde{\frak P}_0=\tilde{\frak P}_0(\Bbb X,\Cal O_{\Bbb X})$)
the full subcategory of $C(\Mod(\Bbb X))$
consisting of complexes of the form 
$\Cal O_x[n]$ with $x\in\Bbb X$
(resp.\ $\Cal O_{x_1,\ldots,x_r}[n]$ with $x_1,\ldots,x_r\in\Bbb X$)
and $n\in\Bbb Z$.
We define $\frak P=\frak P(\Bbb X,\Cal O_{\Bbb X})$ to be 
$\subcolim{\frak P_0}$
and $\tilde{\frak P}:=\subcolim{\tilde{\frak P}_0}$.
We also define $\frak Q$ (resp.\ $\tilde{\frak Q}$) 
to be the full subcategory of $C(\Mod(\Bbb X))$
consisting of bounded above complexes whose terms are direct sums of 
copies of $\Cal O_x$ (resp.\ $\Cal O_{x_1,\ldots,x_r}$).
Note that $\tilde{\frak P}_0$, $\tilde{\frak P}$, and $\tilde{\frak Q}$
are closed under tensor products.

A complex in $\tilde{\frak P}$ is said to be {\em strongly $K$-flat}.
We say that $\Bbb A\in C(\Mod(\Bbb X))$ is $K$-limp (resp.\ strongly
$K$-limp) if $\Bbb A$ is
$\Hom_{\Mod(\Bbb X)}^\bullet(\Bbb F,?)$-acyclic for
any $\Bbb F\in\frak P$ (resp.\ $\Bbb F\in\tilde{\frak P}$).

\begin{Lemma}\label{pull-back-strong.thm}
Let $f:(\Bbb Y,\Cal O_{\Bbb Y})\rightarrow
(\Bbb X,\Cal O_{\Bbb X})$ be a ringed continuous functor.
Then we have an isomorphism
\[
f_{\#}^{\Mod}(\Cal O_{y_1,\ldots,y_r})\cong 
\Cal O_{f y_1,\ldots,f y_r}
\]
for $y_1,\ldots,y_r\in\Bbb Y$.
In particular, if $\Bbb F\in\frak P(\Bbb Y)$,
then $f_{\#}\Bbb F\in \frak P(\Bbb X)$.
If $\Bbb F$ is strongly $K$-flat, then so is $f_{\#}\Bbb F$.
\end{Lemma}

\proof Clearly, we may assume that $r=1$.
For $y\in\Bbb Y$, we denote the canonical continuous ringed functor 
\[
(\Bbb Y/y,\Cal O_{\Bbb Y}|_y)\rightarrow 
(\Bbb X/fy,\Cal O_{\Bbb X}|_{fy})
\]
by $f/y$.
We have $\Phi_{fy}\circ(f/y)=f\circ \Phi_y$.
Hence, 
\[
f_{\#}\Cal O_{y}=f_{\#}L_y(\Cal O_{\Bbb Y}|y)
\cong
L_{fy}(f/y)_{\#}(\Cal O_{\Bbb Y}|y)\cong L_{fy}(\Cal O_{\Bbb X}|fy)
=\Cal O_{fy}.
\]
The last two assertions are obvious now.
\qed

\begin{Lemma} \label{K-flat.thm}
Let $(\Bbb X,\Cal O_{\Bbb X})$ be a ringed site,
and $\Bbb F, \Bbb G\in C(\Mod(\Bbb X))$.
Then the following hold:
\begin{description}
\item[1] $\Bbb F$ is $K$-flat if and only if $\uHom^\bullet_{\Mod(\Bbb 
X)}(\Bbb F,\Bbb I)$ is $K$-injective for any $K$-injective complex
$\Bbb I$.
\item[2] If $\Bbb F$ is $K$-flat exact, then $\Bbb G
\otimes_{\Cal O_{\Bbb X}}^\bullet \Bbb F$ is exact.
\item[3] The inductive limit of a pseudo-filtered inductive system of $K$-flat
complexes is again $K$-flat.
\item[4] The tensor product of two $K$-flat complexes is again $K$-flat.
\end{description}
\end{Lemma}

See \cite{Spaltenstein} for the proof.

\begin{Proposition}
Let $(\Bbb X,\Cal O_{\Bbb X})$ be a ringed site, and $x\in\Bbb X$.
Then $\Cal O_x$ is $K$-flat.
\end{Proposition}

\proof It suffices to show that for 
any exact complex $\Cal E$, $\Cal E\otimes\Cal O_x$ is exact.
To verify this, it suffices to show that for any $K$-injective complex
$\Cal I$, the complex $\Hom^\bullet_{\Cal O_{\Bbb X}}
(\Cal E\otimes\Cal O_x,\Cal I)$ is exact.
Indeed, then if we consider the $K$-injective resolution 
$\Cal E\otimes\Cal O_x\rightarrow \Cal I$, it must be null-homotopic
and thus $\Cal E\otimes\Cal O_x$ must be exact.

Note that we have
\begin{multline*}
  \Hom^\bullet_{\Cal O_{\Bbb X}}(\Cal E\otimes\Cal O_x,\Cal I)
    \cong
  \Hom^\bullet_{\Cal O_{\Bbb X}}(\Cal O_x,
                         \uHom^\bullet_{\Cal O_{\Bbb X}}(\Cal E,\Cal I))
   \cong\\
  \Hom^\bullet_{\Mod(\Bbb X/x)}(\Cal E|_x,\Cal I|_x)
    \cong
  \Hom^\bullet_{\Cal O_{\Bbb X}}(L_x(\Cal E|_x),\Cal I).
\end{multline*}
As $(?)|_x$ and $L_x$ are exact (\ref{restriction-LR.par}), 
(\ref{L_x-exact.par}), the last complex is exact, and we are done.
\qed

\begin{Corollary}\label{strongly-K-flat.thm}
A strongly $K$-flat complex is $K$-flat.
\end{Corollary}

\proof Follows from the proposition and Lemma~\ref{K-flat.thm}.
\qed

\paragraph
Let $\Cal A$ be an abelian category.
The category of complexes of objects of $\Cal A$ is denoted by $C(\Cal A)$.
For an object
\[
\Bbb F: \cdots\rightarrow F^n
\specialarrow{d^n}F^{n+1}\specialarrow{d^{n+1}}
\rightarrow\cdots
\]
in $C(\Cal A)$, we denote the truncated complex
\[
0\rightarrow\Image{d^{n-1}}\rightarrow F^n
\specialarrow{d^n}F^{n+1}\specialarrow{d^{n+1}}
\rightarrow\cdots
\]
by $\tau^{\geq n}\Bbb F$.
Similarly, the truncated complex
\[
\cdots\rightarrow F^{n-2}\specialarrow{d^{n-2}}F^{n-1}
\specialarrow{d^{n-1}}F^n\rightarrow \Image d^{n}\rightarrow 0
\]
is denoted by $\tau^{\leq n}\Bbb F$.

\begin{Lemma}\label{strongly-K-limp.thm}
Let $(\Bbb X,\Cal O_{\Bbb X})$ be a ringed site
and $\Bbb F\in C(\AB(\Bbb X))$.
\begin{description}
\item[1] We have $\frak Q\subset\frak P\subset \tilde{\frak P}$
and $\frak Q\subset\tilde{\frak Q}\subset \tilde{\frak P}$.
\item[2] A $K$-injective complex is weakly $K$-injective,
a weakly $K$-injective complex is strongly $K$-limp, and
a strongly $K$-limp complex is $K$-limp.
\item[3] A tensor product of two strongly $K$-flat complexes is again
strongly $K$-flat.
\item[4] 
For any $\Bbb H\in C(\Mod(\Bbb X))$, there is a 
$\frak Q$-special direct system $\Bbb F_n$ and
a direct system of chain maps $(f_n:\Bbb F_n\rightarrow \tau^{\leq n}\Bbb H)$
such that $f_n$ is a quasi-isomorphism for each $n\in\Bbb N$, and
$(\Bbb F_n)_l=0$ for $l\geq n+2$.
We have $\indlim \Bbb F_n\rightarrow \Bbb H$ is a quasi-isomorphism,
and $\indlim \Bbb F_n\in\frak P$.
Moreover, $Q\Bbb H$ is the homotopy colimit of the inductive system
$(\tau^{\leq n}Q\Bbb H)$ in the category $D(\Mod(\Bbb X))$.
\item[5] The following are equivalent.
\begin{description}
\item[i] $\Bbb F$ is $K$-limp.
\item[ii] $\Bbb F$ is $K$-limp as a complex of abelian sheaves.
\item[iii] $\Bbb F$ is $\Hom^\bullet_{\Mod(\Bbb X)}(\Cal O_x,?)$-acyclic 
for $x\in\Bbb X$.
\item[iv] $\Bbb F$ is $\Gamma(x,?)$-acyclic for $x\in\Bbb X$.
\item[v] If $\Bbb G\in\frak P$ and $\Bbb G$ is exact, 
then 
$\Hom^\bullet_{\Mod(\Bbb X)}(\Bbb G,\Bbb F)$ is exact.
\end{description}
\item[6] The following are equivalent.
\begin{description}
\item[i] $\Bbb F$ is strongly $K$-limp.
\item[ii] $\Bbb F$ is strongly $K$-limp as a complex of abelian sheaves.
\item[iii] $\Bbb F$ is $\Hom_{\Mod(\Bbb X)}^\bullet(
\Cal O_{x_1,\ldots,x_r},?)$-acyclic for $x_1,\ldots,x_r\in\Bbb X$.
\item[iv] If $\Bbb G\in\tilde{\frak P}$ and $\Bbb G$ is exact, 
then
$\Hom^\bullet_{\Mod(\Bbb X)}(\Bbb G,\Bbb F)$ is exact.
\item[v] For any $\Bbb G\in\tilde{\frak P}$, 
$\uHom^\bullet_{\Mod(\Bbb X)}(\Bbb G,\Bbb F)$ is strongly $K$-limp.
\end{description}
\item[7] The following are equivalent.
\begin{description}
\item[i] $\Bbb F$ is weakly $K$-injective.
\item[ii] If $\Bbb G$ is $K$-flat exact, then $\Hom^\bullet_{\Mod(\Bbb X)}
(\Bbb G,\Bbb F)$ is exact.
\end{description}
\end{description}
\end{Lemma}

\proof {\bf 1} and {\bf 3} are trivial.
{\bf 2} follows from the definition and Corollary~\ref{strongly-K-flat.thm}.
The proof of {\bf 4} and {\bf 5} are left to the reader, see
\cite[(5.6), (5.16), (5.17), (5.21)]{Spaltenstein}.

{\bf 6} Applying Lemma~\ref{pull-back-strong.thm} to the canonical
functor
\[
(\Bbb X,a\Bbb Z) \rightarrow (\Bbb X,\Cal O_{\Bbb X}),
\]
we immediately have that {\bf i}, {\bf ii}, and {\bf iii} are equivalent.

It is also trivial that these conditions imply {\bf iv}.
We prove {\bf iv$\Rightarrow$ iii}.
Note that by {\bf 1} and {\bf 4}, we have that $\Bbb F$ is $K$-limp.
Let $\xi:\Bbb F\rightarrow \Bbb I$ be a $K$-injective resolution,
and let $\Bbb J$ be the mapping cone of $\xi$.
Take a $\frak P$-resolution $\eta:\Bbb H\rightarrow \Cal O_{x_1,\ldots,x_r}$,
and let $\Bbb G$ be the mapping cone of $\eta$.
By assumption, $\Hom^\bullet_{\Mod(\Bbb X)}(\Bbb G,\Bbb F)$ is exact.
As $\Bbb I$ is $K$-injective, we have
$\Hom^\bullet_{\Mod(\Bbb X)}(\Bbb G,\Bbb I)$ is also exact.
It follows that $\Hom^\bullet_{\Mod(\Bbb X)}(\Bbb G,\Bbb J)$ is exact.
As $\Bbb J$ is $K$-limp exact and $\Bbb H\in\frak P$, we have
$\Hom^\bullet_{\Mod(\Bbb X)}(\Bbb H,\Bbb J)$ is exact.
It follows that 
$\Hom^\bullet_{\Mod(\Bbb X)}(\Cal O_{x_1,\ldots,x_r},\Bbb J)$ is exact.
This is what we want to prove.

We prove {\bf i$\Rightarrow$v}.
Let $\Bbb H$ be an exact complex in $\tilde{\frak P}$.
Then $\Bbb H\otimes^\bullet \Bbb G$ is exact by Lemma~\ref{K-flat.thm}
and belongs to $\tilde {\frak P}$, since $\tilde{\frak P}$ is closed
under tensor products.
Hence, 
\[
\Hom^\bullet_{\Mod(\Bbb X)}(\Bbb H,\uHom^\bullet_{\Mod(\Bbb X)}
(\Bbb G,\Bbb F))
\cong
\Hom^\bullet_{\Mod(\Bbb X)}(\Bbb H\otimes^\bullet\Bbb G,\Bbb F)
\]
is exact.
By {\bf iv$\Rightarrow$i}, $\uHom^\bullet_{\Mod(\Bbb X)}(\Bbb G,\Bbb F)$
is strongly $K$-limp, as desired.

The implication {\bf v$\Rightarrow$i} is trivial, letting $\Bbb G=
\Cal O_{\Bbb X}$.

{\bf 7} is proved similarly to {\bf 6}.
\qed

\begin{Corollary}\label{skl-wki-tsuika.thm}
Let $(\Bbb X,\Cal O_{\Bbb X})$ be a ringed site and
$\Bbb F,\Bbb G\in C(\Mod(\Bbb X))$.
\begin{description}
\item[1] If $\Bbb F$ is strongly $K$-limp and $\Bbb G\in\tilde{\frak P}$,
then $\Bbb F$ is $\uHom^\bullet_{\Mod(\Bbb X)}(\Bbb G,?)$-acyclic.
\item[2] If $\Bbb F$ is weakly $K$-injective and $\Bbb G$
is $K$-flat, then $\Bbb F$ is $\uHom^\bullet_{\Mod(\Bbb X)}(\Bbb G,?)$-acyclic.
\end{description}
\end{Corollary}

\proof {\bf 1} Let $\Bbb F\rightarrow \Bbb I$ be the $K$-injective 
resolution, and
$\Bbb J$ the mapping cone.
Let 
$\rho:\Bbb H\rightarrow \uHom^\bullet_{\Mod(\Bbb X)}(\Bbb G,\Bbb J)$ be a 
$\frak P$-resolution.
As $\Bbb H\otimes^\bullet\Bbb G$ lies in $\tilde{\frak P}$ and 
$\Bbb J$ is strongly $K$-limp exact, 
\[
\Hom^\bullet_{\Mod(\Bbb X)}(\Bbb H,\uHom^\bullet_{\Mod(\Bbb X)}
(\Bbb G,\Bbb J))
\cong
\Hom^\bullet_{\Mod(\Bbb X)}(\Bbb H\otimes^\bullet\Bbb G,\Bbb J)
\]
is exact.
So $\rho$ must be null-homotopic, and hence
$\uHom^\bullet_{\Mod(\Bbb X)}(\Bbb G,\Bbb J)$ is exact.
This is what we wanted to prove.

{\bf 2} is similar, and we omit it.
\qed

\paragraph
Let $(\Bbb X,\Cal O_{\Bbb X})$ be a ringed site.
For $\Bbb G\in C(\Mod(\Bbb X))$, it is easy to see that
$\Bbb G\otimes_{\Cal O_{\Bbb X}}^\bullet ?$ induces a functor from 
$K(\Mod(\Bbb X))$ to itself.
By the lemma and \cite[Theorem~I.5.1]{Hartshorne2}, 
the derived functor $L(\Bbb G\otimes_{\Cal O_{\Bbb X}}^{\bullet}?)$ is
induced, and it is calculated using any $K$-flat resolution of $?$.
If we fix $?$, then $L(\Bbb G\otimes_{\Cal O_{\Bbb X}}^{\bullet}?)$ 
is a functor on $\Bbb G$, and it induces a bifunctor
\[
*\otimes_{O_{\Bbb X}}^{\bullet,L}?:
D(\Mod(\Bbb X))\times D(\Mod(\Bbb X))\rightarrow D(\Mod(\Bbb X)).
\]
$\Bbb G\otimes_{O_{\Bbb X}}^{\bullet,L}\Bbb F$ is calculated using
any $K$-flat resolution of $\Bbb F$ or any $K$-flat resolution of $\Bbb G$.
Note that $\otimes_{O_{\Bbb X}}^{\bullet,L}$ is a $\triangle$-functor
as in \cite[(2.5.7)]{Lipman}.

We define the hyperTor functor as follows:
\[
\uTor^{\Cal O_{\Bbb X}}_i(\Bbb F,\Bbb G)
:=H^{-i}(\Bbb F\otimes_{\Cal O_{\Bbb X}}^{\bullet,L}\Bbb G).
\]

\paragraph Let $(\Bbb X,\Cal O_{\Bbb X})$ be a ringed site.
For $\Bbb F\in C(\Mod(\Bbb X))$, the functor $\uHom_{\Cal O_{\Bbb X}}^\bullet(
\Bbb F,?)$ induces a functor from $K(\Mod(\Bbb X))$ to itself.
As $\Mod(\Bbb X)$ is Grothendieck, we can take $K$-injective resolutions,
and hence the right derived functor $R\uHom_{\Cal O_{\Bbb X}}^\bullet(
\Bbb F,\mathord ?)$ is induced.
Thus a bifunctor
\[
R\uHom_{\Cal O_{\Bbb X}}(\mathord *,\mathord ?):
D(\Mod(\Bbb X))\op\times D(\Mod(\Bbb X))\rightarrow D(\Mod(\Bbb X))
\]
is induced.
For $\Bbb F,\Bbb G\in D(\Mod(X_\bullet))$, we define the hyperExt sheaf
of $\Bbb F$ and $\Bbb G$ by
\[
\uExt^i_{\Cal O_{\Bbb X}}(\Bbb F,\Bbb G):=H^i(R\uHom^\bullet_{
\Cal O_{\Bbb X}}(\Bbb F,\Bbb G)).
\]

Similarly, the functor $\Hom_{\Cal O_{\Bbb X}}^\bullet(\mathord *,\mathord ?)$
induces
\[
R\Hom_{\Cal O_{\Bbb X}}(\mathord *,\mathord ?):
D(\Mod(\Bbb X))\op\times D(\Mod(\Bbb X))\rightarrow D(\Ab).
\]
Almost by definition, we have
\[
H^i(R\Hom_{\Cal O_{\Bbb X}}(\Bbb F,\Bbb G))\cong \Hom_{D(\Mod(\Bbb X))}
(\Bbb F,\Bbb G[i]).
\]
Sometimes we denote these groups by $\Ext^i_{\Cal O_{\Bbb X}}(\Bbb F,
\Bbb G)$.

\begin{Lemma} Let $(\Bbb X,\Cal O_{\Bbb X})$ be a ringed site.
Then $D(\Mod(\Bbb X))$ is a closed symmetric monoidal triangulated category
with $\otimes^{\bullet,L}_{\Cal O_{\Bbb X}}$ its product and
$R\uHom^\bullet_{\Cal O_{\Bbb X}}(\mathord *,\mathord ?)$ its
internal hom.
\end{Lemma}

\proof This is straightforward.

\begin{Lemma}\label{admissible-K.thm}
Let $f:(\Bbb Y,\Cal O_{\Bbb Y})\rightarrow (\Bbb X,\Cal O_{\Bbb X})$
be an admissible ringed continuous functor.
Then the following hold:
\begin{description}
\item[1] If $\Bbb I\in C(\AB(\Bbb X))$ is a $K$-injective 
\(resp.\ strongly $K$-limp, $K$-limp\) complex of
sheaves of abelian groups, then so is $f^{\#}_{{\AB}}\Bbb I$.
\item[2] If $\Bbb F\in C(\Mod(\Bbb Y))$ is strongly $K$-flat and exact, 
then $f_{\#}^{{\Mod}}\Bbb F$ is strongly $K$-flat and exact.
\item[3] If $\Bbb I\in C(\AB(\Bbb X))$ is $K$-limp and exact, then
so is $f^{\#}_{{\AB}}\Bbb I$.
\end{description}
\end{Lemma}

\proof As $f^{\#}_{{\AB}}$ has an exact left adjoint $f_{\#}^{{\AB}}$,
the assertion for $K$-injectivity in {\bf 1} is obvious.

We prove the assertion for the $K$-limp property in {\bf 1}.
Let $\Bbb P\in\frak P(\Bbb Y)$ 
be an exact complex.
As $f_{\#}^{{\AB}}$ is exact, $f_{\#}^{{\AB}}\Bbb P$ is exact and
a complex in $\frak P(\Bbb X)$ by Lemma~\ref{pull-back-strong.thm}.
Hence, 
\[
\Hom_{\AB(\Bbb Y)}^\bullet(\Bbb P,f^{\#}_{\AB}\Bbb I)
\cong
\Hom_{\AB(\Bbb X)}^\bullet(f_{\#}^{\AB}\Bbb P,\Bbb I)
\]
is exact.
This shows $f^{\#}_{\AB}\Bbb I$ is $K$-limp.
The proof for the strongly $K$-limp property is similar.

We prove {\bf 2}.
We already know that $f_{\#}\Bbb F$ is strongly $K$-flat
by Lemma~\ref{pull-back-strong.thm}.
We prove that $\Hom_{\Mod(\Bbb X)}^\bullet(f_{\#}^{{\Mod}}\Bbb F,\Bbb I)$
is exact, where $\eta:f_{\#}^{{\Mod}}\Bbb F\rightarrow 
\Bbb I$ is a $K$-injective resolution.
Then $\eta$ must be null-homotopic, and we have $f_{\#}^{{\Mod}}\Bbb F$
is exact and the proof is complete.
Clearly, $\Bbb I$ is strongly 
$K$-limp, and hence  so is $f^{\#}_{{\Mod}}\Bbb I$
by {\bf 1} and
Lemma~\ref{strongly-K-limp.thm}, {\bf 6}.
The assertion follows immediately by adjunction.

We prove {\bf 3}.
Let $\xi:\Bbb F\rightarrow f^{\#}_{\AB}\Bbb I$ be a $\frak P(\AB)$-resolution
of $f^{\#}_{\AB}\Bbb I$.
It suffices to show $\Hom^\bullet_{\AB(\Bbb Y)}(\Bbb F,f^{\#}_{\AB}\Bbb I)$ 
is exact, which is obvious since $f_{\#}^{\AB}\Bbb F\in\frak P(\AB)$.
\qed

\paragraph
By the lemma, there is a derived functor $Lf_{\#}^{{\Mod}}:D(\Mod(\Bbb Y))
\rightarrow D(\Mod(\Bbb X))$ of $f_{\#}^{{\Mod}}$ for an
admissible ringed continuous functor $f$.
It is calculated via strongly $K$-flat resolutions.
By the same lemma, we also know that a $K$-limp complex is $f^{\#}$-acyclic,
if $f$ is admissible.

Now as in \cite[section~6]{Spaltenstein} and \cite{Lipman}, 
the following is proved.

\begin{Lemma}\label{lipman-sites.thm}
Let $\Cal S$ be the category of ringed sites and admissible ringed 
continuous functors.
Then $(L(?)^{{\Mod}}_{\#},R(?)_{{\Mod}}^{\#})$ 
is a monoidal adjoint pair of $\Delta$-pseudofunctors on $\Cal S\op$.
\end{Lemma}

The proof is basically the same as that in \cite[Chapter~1--3]{Lipman},
and left to the reader.
We only remark how the structure maps (i.e., natural transformations)
which are necessary to prove the assertion are defined.
For example, the unit of adjunction $u:\Id\rightarrow (Rf^{\#})(Lf_{\#})$ 
is the unique map which makes the diagram
\[
\begin{array}{ccc}
Q & \specialarrow{uQ} & (Rf^{\#})(Lf_{\#})Q\\
\downarrow\hbox to 0pt{\scriptsize$Qu$\hss} & &
\downarrow\hbox to 0pt{\scriptsize$(Rf^{\#})\Xi$\hss}\\
Qf^{\#}f_{\#} & \specialarrow{\Xi f_{\#}} & (Rf^{\#})Qf_{\#}
\end{array}
\]
commutative.
See \cite[Chapter~3]{Lipman} for more.

\begin{Rem} The author does not know if $f_{\#}^{{\Mod}}$ preserves
$K$-flatness of complexes or exactness of $K$-flat complexes.
This is true if
$\Bbb Y$ has finite limits and $f$ preserves finite limits.
This condition is satisfied by morphisms of ringed spaces.
The author does not know if $\frak P$ is closed under tensor products.
On the other hand, strongly $K$-flat resolutions can be used to
calculate both the derived tensor products and $Lf_{\#}$.
The class of strongly $K$-flat resolutions are closed under tensor products
and preserved by $f_{\#}$, and this is enough to prove 
Lemma~\ref{lipman-sites.thm}.
\end{Rem}

\section{Sheaves over a diagram of $S$-schemes}

\paragraph
Let $S$ be a scheme, and 
$I$ a small category.
We call an object of $\Cal P(I\op,\Sch/S)$ an $I$-diagram of $S$-schemes,
where $\Sch/S$ denotes the category of $S$-schemes.
So an object of $\Cal P(I,\Sch/S)$ is referred as an $I\op$-diagram of
$S$-schemes.
Let $X_\bullet\in\Cal P(I,\Sch/S)$.
We denote $X_\bullet(i)$ by $X_i$ for $i\in I$, and
$X_\bullet(\phi)$ by $X_\phi$ for $\phi\in\Mor(I)$.
Let $\Bbb P$ be a property of schemes (e.g., quasi-compact, 
locally noetherian, regular).
We say that $X_\bullet$ satisfies $\Bbb P$ if $X_i$ satisfies $\Bbb P$ for
any $i\in I$.
Let $\Bbb Q$ be a property of morphisms of schemes (e.g., quasi-compact, 
locally of finite type, smooth).
We say that $X_\bullet$ is $\Bbb Q$ over $S$ if the structure map
$X_i\rightarrow S$ satisfies $\Bbb Q$ for any $i\in I$.
We say that $X_\bullet$ has $\Bbb Q$ arrows if $X_\phi$ satisfies
$\Bbb Q$ for any $\phi\in\Mor(I)$.

\paragraph \label{of-fiber-type.par}
Let $f_\bullet
:X_\bullet\rightarrow Y_\bullet$ be a morphism in $\Cal P(I,\Sch/S)$.
For $i\in I$, we denote $f_\bullet(i):X_i\rightarrow Y_i$ by $f_i$.
For a property $\Bbb Q$ of morphisms of schemes, we say that $f_\bullet$
satisfies $\Bbb Q$ if so does $f_i$ for any $i\in I$.
We say that $f_\bullet $ is {\em of fiber type} if the canonical map
$(f_j,X_\phi):X_j\rightarrow Y_j\times_{Y_i}X_i$ is an isomorphism for
any morphism $\phi:i\rightarrow j$ of $I$.

\paragraph
Let $S$, $I$ and $X_\bullet$ be as above.
We define the {\em Zariski site} of $X_\bullet$, denoted by $\Zar(X_\bullet)$,
as follows.
An object of $\Zar(X_\bullet)$ is a pair $(i,U)$ such that $i\in I$ and
$U$ is an open subset of $X_i$.
A morphism $(\phi,h):(j,V)\rightarrow (i,U)$ is a pair $(\phi,h)$ such that
$\phi\in I(i,j)$ and $h:V\rightarrow U$ is the restriction of $X_\phi$.
For a given morphism $\phi:i\rightarrow j$,
$U$, and $V$, such an $h$ exists if and
only if $V\subset X^{-1}_\phi(U)$, and it is unique.
We denote this $h$ by $h(\phi;U,V)$.
The composition of morphisms is defined in an obvious way.
Thus $\Zar(X_\bullet)$ is a small category.
For $(i,U)\in \Zar(X_\bullet)$, a covering of $(i,U)$ is a family of morphisms
of the form
\[
((\id_i,h(\id_i;U,U_\lambda)):(i,U_\lambda)\rightarrow (i,U))_{\lambda\in
\Lambda}
\]
such that $\bigcup_{\lambda\in \Lambda}U_\lambda=U$.
This defines a pretopology of $\Zar(X_\bullet)$, and
$\Zar(X_\bullet)$ is a site.
As we will consider only the Zariski topology, a presheaf or sheaf on
$\Zar(X_\bullet)$ will be sometimes referred as a presheaf or sheaf on
$X_\bullet$, if there is no danger of confusion.
Thus $\Cal P(X_\bullet,\Cal C)$ and $\Cal S(X_\bullet,\Cal C)$ mean
$\Cal P(\Zar(X_\bullet),\Cal C)$ and $\Cal S(\Zar(X_\bullet),\Cal C)$, 
respectively.

\paragraph\label{restrict.par}
Let $S$ and $I$ be as above, and let $\sigma:J\hookrightarrow I$ be
a subcategory of $I$.
Then we have an obvious restriction functor
$\sigma^{\#}:\Cal P(I,\Sch/S)\rightarrow \Cal P(J,\Sch/S)$, which
we denote by $(?)|_J$.
If $\ob(J)$ is finite and $I(j,i)$ is finite for each $j\in J$ and $i\in I$, 
then $(?)|_J$ has 
a right adjoint functor $\cosk^I_J=(?)\op\sigma_{\#}\op(?)\op$, 
because $\Sch/S$ has finite limits.
See \cite[pp.~9--12]{Friedlander}.

\paragraph 
Let $X_\bullet\in\Cal P(I,\Sch/S)$.
Then we have an obvious continuous functor $Q(X_\bullet,J):
\Zar((X_\bullet)|_J)\hookrightarrow \Zar(X_\bullet)$.
Note that $Q(X_\bullet,J)$ may not be admissible.
The restriction functors $Q(X_\bullet,i)^{\#}_{\AB}$ and
$Q(X_\bullet,J)^{\#}_{\PA}$ are denoted by $(?)_J^{\AB}$ and $(?)_J^{\PA}$,
respectively.
For $i\in I$, we consider that $i$ is the subcategory of $I$ whose object
set is $\{i\}$ with $\Hom_i(i,i)=\{\id\}$.
The restriction $(?)_i^\natural$ is defined.

\paragraph\label{friedlander-identification.par}
Let $\Cal F\in\PA(X_\bullet)$ and $i\in I$.
Then $\Cal F_i\in\PA(X_i)$, and thus we have a family of sheaves
$(\Cal F_i)_{i\in I}$.
Moreover, for $(i,U)\in\Zar(X_\bullet)$ and $\phi:i\rightarrow j$,
we have the restriction map
\[
\Gamma(U,\Cal F_i)=\Gamma((i,U),\Cal F)\specialarrow{\standop{res}}
\Gamma((j,X_\phi^{-1}(U)),\Cal F)=\\
\Gamma(X_\phi^{-1}(U),\Cal F_j)=
\Gamma(U,(X_\phi)_*\Cal F_j),
\]
which induces
\begin{equation}\label{beta-def.eq}
\beta_\phi(\Cal F)\in\Hom_{\PA(X_i)}(\Cal F_i,(X_\phi)_*\Cal F
_j).
\end{equation}
The corresponding map in $\Hom_{\PA(X_j)}((X_\phi)^*_{\PA}(\Cal F_i),\Cal F_j)$
is denoted by $\alpha_\phi^{\PA}(\Cal F)$.
If $\Cal F$ is a sheaf, then (\ref{beta-def.eq}) yields
\[
\alpha_\phi^{\AB}(\Cal F)\in\Hom_{\AB(X_j)}((X_\phi)^*_{\AB}(\Cal F_i),
\Cal F_j).
\]

It is straightforward to check the following.

\def\citinfo{\cite{Friedlander}}
\begin{Lemma}[\citinfo]\label{simplicial-sheaf.thm}
Let $\natural$ be either $\AB$ or $\PA$.
The following hold:
\begin{description}
\item[\truebf
1] For any $i\in I$, we have $\alpha_{\id_i}^\natural:(X_{\id_i})^*_\natural
(\Cal F_i)\rightarrow \Cal F_i$ is the canonical identification.
\item[\truebf
2] If $\phi\in I(i,j)$ and $\psi\in I(j,k)$, then the composite map
\begin{equation}\label{alpha-composition.eq}
(X_{\psi\phi})^*_\natural(\Cal F_i)\cong
(X_\psi)^*_\natural(X_\phi)^*_\natural(\Cal F_i)
\specialarrow{(X_\psi)^*_\natural\alpha_\phi^\natural}
(X_\psi)^*_\natural(\Cal F_j)
\specialarrow{\alpha_\psi^\natural}
\Cal F_k
\end{equation}
agrees with $\alpha_{\psi\phi}^\natural$.
\item[\truebf 3]
Conversely, a family $((\Cal G_i)_{i\in I},(\alpha_\phi)_{\phi\in\Mor(I)})$
such that $\Cal G_i\in \natural(X_i)$, $\alpha_\phi\in\Hom_{\natural(X_j)}
((X_\phi)^*_\natural(\Cal G_i),\Cal G_j)$ for $\phi\in I(i,j)$, 
and that the conditions corresponding to {\truebf 1,2} are satisfied
yields $\Cal G\in\natural(X_\bullet)$, and this correspondence gives
an equivalence.
\end{description}
\end{Lemma}

Similarly, a family $((\Cal G_i)_{i\in\ob(I)},
(\beta_\phi)_{\phi\in\Mor(I)})$ with
\[
\Cal G_i\in\natural(X_i) \text{ and }
\beta_\phi\in\Hom_{\natural(X_i)}(\Cal G_i,
(X_\phi)^\natural_*\Cal G_j)
\]
satisfying the conditions 
\begin{description}
\item[1'] For 
$i\in\ob(I)$, $\beta_{
\id_i}:\Cal G_i\rightarrow (\id_i)_*\Cal G_i$ is the canonical identification;
\item[2'] For $\phi\in I(i,j)$ and $\psi\in I(j,k)$, the composite
\[
\Cal F_i\specialarrow{\beta_\phi}
(X_\phi)_*(\Cal F_j)
\specialarrow{(X_\phi)_*\beta_\psi}
(X_\phi)_*(X_\psi)_*(\Cal F_k)
\cong
(X_{\psi\phi})_*(\Cal F_k)
\]
agrees with $\beta_{\psi\phi}$
\end{description}
are in one to one correspondence with $\Cal G\in\natural(X_\bullet)$.

\paragraph Let $\Cal F\in\AB(X_\bullet)$.
We say that $\Cal F$ is an {\em equivariant} abelian sheaf if
$\alpha_\phi^{\AB}$ are isomorphisms for all $\phi\in\Mor(I)$.
For $\Cal F\in\PA(X_\bullet)$, we say that $\Cal F$ is an equivariant 
abelian presheaf if $\alpha_\phi^{\PA}$ are isomorphisms for all $\phi\in\Mor
(I)$.
An equivariant sheaf may not be an equivariant presheaf.
However, an equivariant presheaf which is a sheaf is an equivariant sheaf.
We denote the category of 
equivariant sheaves and presheaves
by $\EqAB(X_\bullet)$ and $\EqPA(X_\bullet)$, respectively.
As $(X_\phi)^*_\natural$ is exact for $\natural = \AB,\PA$ and
any $\phi$, we have that $\EqAB(X_\bullet)$ is plump in $\AB(X_\bullet)$,
and $\EqPA(X_\bullet)$ is plump in $\PA(X_\bullet)$.

\paragraph
Let $X_\bullet\in\Cal P(I,\Sch/S)$.
The datum
\[
((\Cal O_{X_i})_{i\in I},(\beta_\phi:\Cal O_{X_i}\rightarrow
(X_\phi)_*\Cal O_{X_j})_{\phi\in\Mor(I)})
\]
gives a sheaf of commutative rings on $X_\bullet$, which we denote
by $\Cal O_{X_\bullet}$, and thus
$\Zar(X_\bullet)$ is a ringed site.
The categories $\PM(\Zar(X_\bullet))$ and $\Mod(\Zar(X_\bullet))$ are denoted
by $\PM(X_\bullet)$ and $\Mod(X_\bullet)$, respectively.

\paragraph\label{restrict-mod1.par}
For $J\subset I$, we have $\Cal O_{X_\bullet|_J}=
(\Cal O_{X_\bullet})_J$ by definition.
The continuous functor
\[
Q(X_\bullet,J):(\Zar(X_\bullet|_J),\Cal O_{X_\bullet|_J})
\rightarrow
(\Zar(X_\bullet),\Cal O_{X_\bullet})
\]
is actually a ringed continuous functor.

The corresponding restriction $Q(X_\bullet,J)^{\#}_\natural$ is denoted
by $(?)_J^\natural$ for $\natural=\PM, \Mod$.
For subcategories $J_1\subset J\subset I$ of $I$, we denote
the restriction $(?)_{J_1}^\natural:\natural(X_\bullet|_J)\rightarrow
\natural(X_{\bullet}|_{J_1})$ is denoted by
$(?)_{J_1,J}^\natural$, to emphasize $J$.

\paragraph
Let $\natural$ be $\PM$ or $\Mod$.
Note that $\Cal M\in\natural(X_\bullet)$ is nothing but a family
\[
\Dat(\Cal M):=
((\Cal M_i)_{i\in I}, (\alpha^\natural_\phi)_{\phi\in\Mor(I)})
\]
such that
$\Cal M_i\in\natural(X_i)$, $\alpha_\phi^\natural:(X_\phi)^*_\natural
(\Cal M_i)\rightarrow\Cal M_j$ is a morphism of $\natural(X_\bullet)$ for any 
$\phi:i\rightarrow j$, and the conditions corresponding to {\bf 1,2} in
Lemma~\ref{simplicial-sheaf.thm} are satisfied.

We say that $\Cal M\in\natural(X_\bullet)$ is {\em equivariant} if 
$\alpha_\phi^\natural$ is an isomorphism for any $\phi\in\Mor(I)$.
Note that equivariance depends on $\natural$, and is not preserved by
the forgetful functors in general.
We denote the full subcategory of $\Mod(X_\bullet)$ consisting
of equivariant objects by $\EM(X_\bullet)$.

\section[The inductions, the direct image, and the inverse image]{The 
left and right inductions and the direct and inverse images}

Let $I$ be a small category, $S$ a scheme, and $X_\bullet\in
\Cal P(I,\Sch/S)$.

\paragraph\label{restrict-mod.par}
Let $J$ be a subcategory of $I$.
The left adjoint $Q(X_\bullet,J)_{\#}^\natural$ of $(?)_J^\natural$ 
(see (\ref{restrict-mod1.par}))
is denoted by $L_J^\natural$ for $\natural=\PA, \AB, \PM, \Mod$.
The right adjoint $Q(X_\bullet,J)^\natural_{\flat}$ of $(?)^\natural_J$,
which exists by Lemma~\ref{relative-topology.thm}, is denoted by
$R^\natural_J$ for $\natural =\PA, \AB, \PM, \Mod$.
We call $L^\natural_J$ and $R^\natural_J$ the left and right
induction functor, respectively.

Let $J_1\subset J\subset I$ be subcategories of $I$.
The left and right adjoints of $(?)^\natural_{J_1,J}$ is denoted by
$L^\natural_{J,J_1}$ and $R^\natural_{J,J_1}$, respectively.
As $(?)_J^\natural$ has both a left adjoint and a right adjoint, we have

\begin{Lemma}
The functor $(?)_J^\natural$
preserves arbitrary limits and colimits
\(hence is exact\) for $\natural=\PA,\AB,\PM,\Mod$.
\end{Lemma}

The functor $\natural(X_\bullet)\rightarrow\prod_{i\in I}\natural(X_i)$
given by $\Cal F\mapsto (\Cal F_i)_{i\in I}$ is faithful for
$\natural=\PA, \AB, \PM, \Mod$.

\paragraph
Let $f_\bullet:X_\bullet\rightarrow Y_\bullet$ be a morphism
in $\Cal P(I,\Sch/S)$.
This induces an obvious morphism
\[
f_\bullet^{-1}:(\Zar(Y_\bullet),\Cal O_{Y_\bullet})\rightarrow
(\Zar(X_\bullet),\Cal O_{X_\bullet})
\]
of ringed sites.
We have $\id^{-1}=\id$, and $(g_\bullet \circ f_\bullet)^{-1}=
f_\bullet^{-1}\circ g_\bullet^{-1}$ for $g_\bullet:Y_\bullet\rightarrow
Z_\bullet$.

We define the direct image $(f_\bullet)_*^\natural
$ to be $(f_\bullet^{-1})^{\#}_\natural$,
and the inverse image $(f_\bullet)^*_\natural$ to be $(f_\bullet^{-1})_{\#}^
\natural$ for $\natural={\Mod},{\PM},{\AB},{\PA}$.

\begin{Lemma}\label{commutativity-sites.thm}
Let $f_\bullet:X_\bullet\rightarrow Y_\bullet$ be a morphism
in $\Cal P(I,\Sch/S)$, and $K\subset J\subset I$.
Then we have
\begin{description}
\item[1] $Q(X_\bullet,J)\circ Q(X_\bullet|_{J},K)=Q(X_\bullet,K)$
\item[2] $f_\bullet^{-1}\circ Q(Y_\bullet,J)=
Q(X_\bullet,J)\circ (f_\bullet|_J)^{-1}$.
\end{description}
\end{Lemma}

\paragraph Let us fix $I$ and $S$.
By Lemma~\ref{lipman-sites-sheaf.thm},
we have various natural maps between
functors on sheaves arising from the closed structures
and the monoidal pairs, involving various $J$-diagrams of schemes,
where $J$ varies subcategories of $I$.
In the sequel, many of the natural maps 
are referred as \lq the canonical maps'
or \lq the canonical isomorphisms' without any explicit definitions.
Many of them are defined in \cite{Lipman}, and various commutativity theorems
are proved there.

\begin{example}\label{lipman.ex}
Let $I$ be a small category, $S$ a scheme, and 
$f_\bullet:X_\bullet\rightarrow Y_\bullet$ and
$g_\bullet:Y_\bullet\rightarrow Z_\bullet$
are morphisms in $\Cal P(I,\Sch/S)$.
Let $K\subset J\subset I$ be subcategories, and $\natural$ denote ${\PM}$, 
${\Mod}$, ${\PA}$, or ${\AB}$.
\begin{description}
\item[1] There is a natural isomorphism
\[
c_{I,J,K}^\natural: (?)_{K,I}^\natural\cong (?)_{K,J}^\natural
\circ (?)_{J,I}^\natural.
\]
Taking the conjugate,
\[
d_{I,J,K}^\natural: L_{I,J}^\natural\circ L_{J,K}^\natural\cong 
L_{I,K}^\natural
\]
is induced.
\item[2] There is a natural isomorphism
\[
c_{J,f_\bullet}^\natural: (?)_J^\natural\circ (f_\bullet)_*^\natural 
\cong (f_\bullet|_J)_*^\natural\circ (?)_J^\natural
\]
and its conjugate
\[
d_{J,f_\bullet}^\natural: L_J^\natural\circ (f_\bullet|_J)^*_\natural
\cong (f_\bullet)^*_\natural\circ L_J^\natural.
\]
\item[3] We have
\[
(c_{K,f_\bullet|_J}^\natural(?)_J^\natural)\circ ((?)_{K,J}^\natural 
c_{J,f}^\natural)
=
((f_\bullet|_K)_*^\natural c_{I,J,K}^\natural)\circ 
c_{K,f_\bullet}^\natural\circ ((c_{I,J,K}^\natural)^{-1}(f_\bullet)_*^\natural
).
\]
\item[4] We have
\[
((g_\bullet|_J)_*^\natural c_{J,f_\bullet}^\natural)
\circ(c_{J,g_\bullet}^\natural(f_\bullet)_*^\natural)
=
(c_{f_\bullet|_J,g_\bullet|_J}^\natural(?)_J^\natural)
\circ
c_{J,g_\bullet\circ f_\bullet}^\natural
\circ((?)_J^\natural (c_{f_\bullet,g_\bullet}^\natural)^{-1}),
\]
where $c_{f_\bullet,g_\bullet}^\natural:(g_\bullet\circ f_\bullet)_*^\natural
\cong
(g_\bullet)_*^\natural\circ (f_\bullet)_*^\natural$ 
is the canonical isomorphism,
and similarly for $c_{f_\bullet|_J,g_\bullet|_J}^\natural$.
\item[5] The canonical map
\[
H:\Cal M_J\otimes_{\Cal O_{X_\bullet|_J}}\Cal N_J
\rightarrow
(\Cal M\otimes_{\Cal O_{X_\bullet}}\Cal N)_J
\]
is an isomorphism, as can be seen easily.
The canonical map
\[
\Delta:L_J(\Cal M\otimes_{\Cal O_{X_\bullet|_J}}\Cal N)\cong
(L_J \Cal M)\otimes_{\Cal O_{X_\bullet}}(L_J\Cal N).
\]
is defined, which may not be an isomorphism.
\end{description}
\end{example}

\section{Operations on sheaves via the structure data}

Let $I$ be a small category, $S$ a scheme, and $\Cal P:=\Cal P(I,\Sch/S)$.
To study sheaves on objects of $\Cal P$, it is convenient to utilize
the structure data of them, and then utilize the usual sheaf theory on
schemes.

\paragraph
Let $X_\bullet\in\Cal P$.
Let $\natural$ be any of $\PA,\AB,\PM,\Mod$,
and $\Cal M,\Cal N\in\natural(X_\bullet)$.
An element $(\varphi_i)$ in $\prod\Hom_{\natural(X_i)}
(\Cal M_i,\Cal N_i)$ is given by some
$\varphi\in\Hom_{\natural(X_\bullet)}(\Cal M,\Cal N)$
(by the canonical faithful functor $\natural(X_\bullet)\rightarrow
\prod \natural(X_i)$), if and only if
\begin{equation}\label{mor-data.eq}
\varphi_j\circ\alpha_\phi(\Cal M)=\alpha_\phi(\Cal N)\circ (X_\phi)^*_\natural
(\varphi_i)
\end{equation}
holds (or equivalently, $\beta_\phi(\Cal N)\circ\varphi_i
=
(X_{\phi})_*\varphi_j\circ \beta_\phi(\Cal M)$ holds)
for any $(\phi:i\rightarrow j)\in\Mor(I)$.

We say that a family of morphisms $(\varphi_i)_{i\in I}$ between
structure data
\[
\varphi_i:\Cal M_i\rightarrow \Cal N_i
\]
is a morphism of structure data if $\varphi_i$ is a morphism in
$\natural(X_i)$ for each $i$, and 
(\ref{mor-data.eq}) is satisfied for any $\phi$.
Thus the category of structure data of sheaves, presheaves, modules, and
premodules  on $X_\bullet$, denoted
by $\frak D_\natural(X_\bullet)$ is defined, and the equivalence
$\Dat_\natural:\natural(X_\bullet)\cong \frak D_\natural(X_\bullet)$
is given.
This is the precise meaning of Lemma~\ref{simplicial-sheaf.thm}.

\paragraph
\label{tensor-data.par}
Let $X_\bullet\in\Cal P$ and $\Cal M,\Cal N\in\Mod(X_\bullet)$.
As in Example~\ref{lipman.ex}, {\bf 5}, we have an isomorphism
\[
m_i: \Cal M_i\otimes_{\Cal O_{X_i}}\Cal N_i
     \cong 
     (\Cal M\otimes_{\Cal O_{X_\bullet}}\Cal N)_i.
\]
This is trivial for presheaves, and utilize the fact the sheafification
is compatible with $(?)_i$ for sheaves.
Under this identification, the structure map $\alpha_\phi$ of
$\Cal M\otimes\Cal N$ can be completely described via those of $\Cal M$
and $\Cal N$.
Namely, for $\phi\in I(i,j)$, $\alpha_\phi(\Cal M\otimes\Cal N)$ 
agrees with the composite map
\[
X_\phi^*(\Cal M\otimes\Cal N)_i
\specialarrow{\via m_i^{-1}}
X_\phi^*(\Cal M_i\otimes \Cal N_i)
\cong 
X_\phi^*\Cal M_i\otimes X_\phi^*\Cal N_i
\specialarrow{\alpha_\phi\otimes\alpha_\phi}
\Cal M_j\otimes\Cal N_j
\specialarrow{m_j}
(\Cal N\otimes\Cal N)_j.
\]
This can be checked directly at the presheaf level, and then sheafify them.

\paragraph\label{L_J.par}
Let $X_\bullet\in\Cal P$, and $J$ a subcategory of $I$.
The left adjoint functor $L_J^\natural=Q(X_\bullet,J)^\natural_{\#}$ 
of $(?)_J^\natural$ is given by the structure data
as follows explicitly.
For $\Cal M\in\natural(X_\bullet|_J)$ and $i\in I$, 
we have

\begin{Lemma}\label{L_J.thm} There is an isomorphism
\[
  \lambda_{J,i}:
  (L_J^\natural(\Cal M))_i^\natural
    \cong 
  \indlim (X_\phi)^*_\natural(\Cal M_j),
\]
where the colimit is taken over the subcategory 
   $(I^{(J\op\rightarrow I\op)}_i)\op$
of $I/i$
whose objects are $(\phi:j\rightarrow i)\in I/i$ with $j\in\ob(J)$ and
morphisms are morphisms $\varphi$ of $I/i$ such that 
$\varphi\in\Mor(J)$.
\end{Lemma}

In particular, we have an isomorphism

\begin{equation}\label{L_i.eq}
\lambda_{j,i}:
  (L_j^\natural(\Cal M))_i^\natural
    \cong
  \bigoplus_{\phi\in I(j,i)}(X_\phi)^*_\natural(\Cal M_j).
\end{equation}

\paragraph\label{L_J-alpha.par}
Let $\psi:i\rightarrow i'$ be a morphism.
The structure map
\[
\alpha_\psi : (X_\psi)^*_\natural((L_J^\natural(\Cal M))_i^\natural)
\rightarrow (L_J^\natural(\Cal M))_{i'}^\natural
\]
is induced by
\[
(X_\psi)^*_\natural((X_\phi)^*_\natural(\Cal M_j))
\cong (X_{\psi\phi})^*_\natural(\Cal M_j).
\]

\paragraph\label{left-induction-counit.par}
The counit map $\varepsilon :L_J(?)_J\rightarrow\Id$ is given
as a morphism of structure data as follows.
\[
\varepsilon_i:(?)_iL_J(?)_J\rightarrow (?)_i
\]
agrees with
\[
(?)_iL_J(?)_J\specialarrow{\lambda_{J,i}}
\indlim X_\phi^*(?)_j (?)_J
\cong
\indlim X_\phi^*(?)_j
\specialarrow{\alpha}(?)_i,
\]
where $\alpha$ is induced by $\alpha_\phi:X_\phi^*(?)_j\rightarrow (?)_i$.

\paragraph\label{left-induction-unit.par}
The unit map $u:\Id\rightarrow (?)_JL_J$ is also described, as
follows.
\[
u_j:(?)_j\rightarrow (?)_j(?)_JL_J
\]
agrees with
\[
(?)_j\cong \id_j^*(?)_j\rightarrow \indlim X_\phi^*(?)_k
\specialarrow{\lambda_{J,j}^{-1}}(?)_jL_J\cong (?)_j(?)_JL_J,
\]
where the colimit is taken over $(I^{(J\op\subset I\op)}_j)\op$.

\paragraph\label{right-induction.par}
Let $X_\bullet\in\Cal P$, and $J$ a subcategory of $I$.
The right adjoint functor $R_J^\natural$ of $(?)_J^\natural$
is given as follows explicitly.
For $\Cal M\in\natural(X_\bullet|_J)$ and $i\in I$,
we have
\[
\rho^{J,i}:
(R_J^\natural(\Cal M))^\natural_i=\projlim (X_\phi)_*^\natural(\Cal M_j),
\]
where the limit is taken over $I^{(J\rightarrow I)}_i$, see
(\ref{Kan-adjoint.par}) for the notation.
The description of $\alpha$, $u$, and $\varepsilon$ for the right induction
are left to the reader.

\begin{Lemma}\label{induction-restrict.thm}
Let $X_\bullet\in\Cal P$, and $J$ a full subcategory of $I$.
Then
we have the following.
\begin{description}
\item[1] The counit of adjunction $\varepsilon:(?)_J^\natural \circ 
R^\natural_J\rightarrow \Id$ is an isomorphism.
\item[2] The unit of adjunction $u:\Id\rightarrow (?)_J^\natural\circ
L^\natural_J$ is an isomorphism.
\end{description}
\end{Lemma}

\proof {\bf 1} For $i\in J$, the restriction
\[
\varepsilon_i:(?)_i^\natural(?)_J^\natural R^\natural_J\Cal M
=\projlim (X_\phi)_*^\natural(\Cal M_j)\rightarrow \Cal M_i=(?)_i\Cal M
\]
is induced by $\beta_\phi$'s, where the limit is taken over
$(\phi:i\rightarrow j)\in I^{(J\rightarrow I)}_i$.
As $J$ is a full subcategory, we have $I^{(J\rightarrow I)}_i$ equals
$i/J$, and hence it has the initial object $\id_i$.
So the limit is equal to $\Cal M_i$, and $\varepsilon_i$ is the identity map.
Since $\varepsilon_i$ is an isomorphism for each $i\in J$, we have that
$\varepsilon$ is an isomorphism.

The proof of {\bf 2} is similar, and we omit it.
\qed

Let $\Cal C$ be a small category.
A connected component of $\Cal C$ is a full subcategory of $\Cal C$
whose object set is one of the equivalence classes of $\ob(\Cal C)$
with respect to the transitive symmetric closure of the relation
$\sim$ given by
\[
c\sim c'\iff \Cal C(c,c')\neq \emptyset.
\]

\begin{Def} We say that a subcategory $J$ of $I$ is admissible if
\begin{description}
\item[1] For $i\in I$, the category $(I^{(J\subset I)}_i)\op
$ is pseudofiltered.
\item[2] For $j\in J$, we have $\id_j$ is the initial object of
one of the connected components of $I^{(J\subset I)}_j$.
\end{description}
\end{Def}

Note that for $j\in I$, the subcategory $j=(\{j\},\{\id_j\})$ of $I$
is admissible.

In Lemma~\ref{L_J.thm}, the colimit is pseudo-filtered
and hence it preserves exactness, if {\bf 1} is satisfied.
As in the proof of Lemma~\ref{induction-restrict.thm}, 
$(?)_j$ is a direct summand of $(?)_j\circ L_J$ for $j\in J$
so that $L_J$ is faithful, if {\bf 2} is satisfied.
We have the following.

\begin{Lemma}\label{admissible-subcat.thm}
Let $X_\bullet\in\Cal P(I,\Sch/S)$, and 
$K\subset J\subset I$ be admissible subcategories of $I$.
Then $L_{J,K}^{{\PA}}$ is faithful and exact.
The morphism of sites $Q(X_\bullet|_J,K)$ is admissible.
If, moreover, $X_\phi$ is flat
for any $\phi\in I(k,j)$ with $j\in J$
and $k\in K$, then $L_{J,K}^\natural$ is faithful and exact 
for $\natural={\Mod},{\PM}$.
\end{Lemma}

\paragraph \label{direct-image.par}
As in Example~\ref{lipman.ex}, {\bf 2}, we have an isomorphism

\begin{equation}
c_{i,f_\bullet}: (?)_i\circ (f_\bullet)_* \cong (f_i)_*\circ (?)_i.
\end{equation}

The translation $\alpha_\phi$ is described, as follows.
For $\phi\in I(i,j)$,
\[
\alpha_\phi (f_\bullet)_*:(Y_\phi)^*(?)_i(f_\bullet)_*\rightarrow (?)_j
(f_\bullet)_*
\]
agrees with
\begin{multline}
(Y_\phi)^*(?)_i(f_\bullet)_*
\specialarrow{c_{i,f_\bullet}}
(Y_\phi)^*(f_i)_*(?)_i
\specialarrow{\via\theta}
(f_j)_*(X_\phi)^*(?)_i\\
\specialarrow{(f_j)_*\alpha_\phi}
(f_j)_*(?)_j
\specialarrow{c_{j,f_\bullet}^{-1}}
(?)_j(f_\bullet)_*,
\label{direct-image-translate.eq}
\end{multline}
where $\theta$ is the map defined and discussed in 
\cite[(3.7.2)]{Lipman}.
One of the definitions of $\theta$ is the composite
\begin{multline}
\theta:(Y_\phi)^*(f_i)_*\specialarrow{\via u}
(Y_\phi)^*(f_i)_*(X_\phi)_*(X_\phi)^*
\cong\\
(Y_\phi)^*(Y_\phi)_*(f_j)_*(X_\phi)^*
\specialarrow{\via \varepsilon}
(f_j)_*(X_\phi)^*.
\label{direct-image-translate2.eq}
\end{multline}

\begin{Proposition}\label{crutial.thm}
Let $f_\bullet:X_\bullet\rightarrow Y_\bullet$ be a morphism in $\Cal P$, 
$J$ a subcategory of $I$, and $i\in I$.
Then the composite map
\[
(?)_iL_J (f_\bullet|_J)_*
\specialarrow{\via\theta}
(?)_i(f_\bullet)_*L_J
\specialarrow{\via c_{i,f_\bullet}}
(f_i)_*(?)_iL_J
\]
agrees with the composite map
\begin{multline*}
(?)_iL_J (f_\bullet|_J)_*
\specialarrow{\via\lambda_{J,i}}
\indlim Y_\phi^*(?)_j(f_\bullet|_J)_*
\specialarrow{\via c_{j,f_\bullet|_J}}
\indlim Y_\phi^*(f_j)_*(?)_j\\
\specialarrow{\via\theta}
\indlim (f_i)_* X_\phi^*(?)_j
\rightarrow
(f_i)_*\indlim X_\phi^*(?)_j
\specialarrow{\via\lambda_{J,i}^{-1}}
(f_i)_*(?)_i L_J.
\end{multline*}
\end{Proposition}

\proof Note that $\theta$ in the first composite map is the composite
\[
\theta=\theta(J,f_\bullet):L_J(f_\bullet|_J)_*
\specialarrow{\via u}
L_J(f_\bullet|_J)_*(?)_JL_J
\specialarrow{\cong}
L_J(?)_J(f_\bullet)_*L_J
\specialarrow{\varepsilon}(f_\bullet)_*L_J.
\]
The $\cong$ at the middle is the identity as a morphism of structure data.
The description of $u$ and $\varepsilon$ are already given, and
the proof is reduced to the iterative use of
(\ref{L_J-alpha.par}), 
(\ref{left-induction-counit.par}), (\ref{left-induction-unit.par}),
and (\ref{direct-image.par}).
The detailed argument is left to a patient reader.
The reason why the second map involves $\theta$ is (\ref{direct-image.par}).
\qed

Similarly, we have the following.

\begin{Proposition}\label{crutial2.thm}
Let $f_\bullet: X_\bullet\rightarrow Y_\bullet$
be a morphism in $\Cal P$, $J$ a subcategory of $I$, 
and $i\in I$.
Then the composite map
\[
(f_i)^*(?)_iL_J 
\specialarrow{\via\theta(f_\bullet,i)}
(?)_i(f_\bullet)^* L_J
\specialarrow{\via d_{f_\bullet,J}}
(?)_iL_J(f_\bullet|_J)^*
\]
agrees with the composite map
\begin{multline*}
(f_i)^*(?)_iL_J
\specialarrow{\via \lambda_{J,i}}
(f_i)^* \indlim Y_\phi^*(?)_j
\cong
\indlim X_\phi^*(f_j)^*(?)_j\\
\specialarrow{\via \theta(f_\bullet|_J,j)}
\indlim X_\phi^*(?)_j (f_\bullet|_J)^*
\specialarrow{\via \lambda_{J,i}^{-1}}
(?)_i L_J (f_\bullet|_J)^*.
\end{multline*}
\end{Proposition}

The proof is left to the reader.
The proof of Proposition~\ref{crutial.thm} and Proposition~\ref{crutial2.thm}
are formal, and the propositions are valid for $\natural={\PM}, {\Mod},
{\PA}, $ and ${\AB}$.

Let $f_\bullet:X_\bullet\rightarrow Y_\bullet$ be a morphism
in $\Cal P$, and $J\subset I$ a subcategory.
The inverse image $(f_\bullet)^*_{\natural}$ is compatible
with the restriction $(?)_J$.

\begin{Lemma}\label{inverse-image-restrict.thm}
The natural map
\[
\theta_\natural=
\theta_\natural(f_\bullet,J):
((f_\bullet)|_J)^*_{\natural}\circ (?)_J
\rightarrow (?)_J\circ (f_\bullet)^*_{\natural}
\]
is an isomorphism for $\natural={\PA}, {\AB}, {\PM}, {\Mod}$.
In particular, $f_\bullet^{-1}:\Zar(Y_\bullet)\rightarrow \Zar(X_\bullet)$ is
an admissible continuous functor.
\end{Lemma}

\proof Note that $\theta$ is the composite map
\[
(f_\bullet|_J)^*(?)_J\specialarrow{\via u}
(f_\bullet|_J)^*(?)_J (f_\bullet)_*f_\bullet^*
\specialarrow{\via c}
(f_\bullet|_J)^*(f_\bullet|_J)_*(?)_J f_\bullet^*
\specialarrow{\via \varepsilon}
(?)_J f_\bullet^*,
\]
see \cite[(3.7.2)]{Lipman}.

We construct a natural map 
\[
\eta: ((f_\bullet)|_J)^*_{\natural}\circ (?)_J
\rightarrow (?)_J\circ (f_\bullet)^*_{\natural},
\]
for the case where $\natural={\PM}$.

Let $\Cal M\in\PM(Y_\bullet)$, and $(j,U)\in\Zar(X_\bullet|_J)$.
We have
\[
\Gamma((j,U),(f_\bullet|_J)^*\Cal M_J)
=
\indlim \Gamma((j,U),\Cal O_{X_\bullet})
           \otimes_{\Gamma((j',V),\Cal O_{Y_\bullet})}
        \Gamma((j',V),\Cal M),
\]
where the colimit is taken over $(j',V)\in 
(I^{(f_\bullet|_J)^{-1}}_{(j,U)})\op$.
On the other hand, we have
\[
\Gamma((j,U),(?)_J f_\bullet^*\Cal M)
=
\indlim \Gamma((j,U),\Cal O_{X_\bullet})
           \otimes_{\Gamma((i,V),\Cal O_{Y_\bullet})}
        \Gamma((i,V),\Cal M),
\]
where the colimit is taken over $(i,V)\in 
(I^{f_\bullet^{-1}}_{(j,U)})\op$.
There is an obvious map from the first to the second.

To verify that this is an isomorphism, it suffices to show that
the category $(I^{f_\bullet^{-1}|_J}_{(j,U)})\op$
is final in
the category $(I^{f_\bullet^{-1}}_{(j,U)})\op$.
In fact, any $(\phi,h):(j,U)\rightarrow (i,f_i^{-1}(V))$ 
with $(i,V)\in\Zar(Y_\bullet)$ factors through 
\[(\id_j,h):(j,U)\rightarrow (j,f_j^{-1}(Y_\phi)^{-1}(V)).
\]

Hence, $\eta_\natural$ is an isomorphism for $\natural={\PM}$.
The construction for the case where $\natural={\PA}$ is similar.
As $(?)_J$ is compatible with the sheafification
by Lemma~\ref{relative-topology.thm}, an isomorphism $\eta$
for $\natural={\Mod}, {\AB}$ is also defined.

In any case, 
it is straightforward to check that the diagram
\begin{equation}\label{f_*-f^*-unit-restrict.eq}
\begin{array}{ccccc}
(?)_J  & \multicolumn{3}{c}{
\mathord -\mkern -6mu
\cleaders \hbox {$\mkern -2mu\mathord -\mkern -2mu$}
\hfill
\hbox to 0pt{\hss$
  \mathord{\mathop{\mkern -1mu\mathord-\mkern -1mu}\limits^{\hbox to
  0pt{\hss\scriptsize $\id$\hss}}}$
   \hss}
\mkern -1mu
\cleaders \hbox {$\mkern -2mu\mathord -\mkern -2mu$}
\hfill\mkern -6mu\mathord \rightarrow
} & (?)_J\\
\downarrow\hbox to 0pt{\scriptsize$\via u$\hss} & & & &
\downarrow\hbox to 0pt{\scriptsize$\via u$\hss}\\
(f_\bullet|_J)_*(f_\bullet|_J)^*(?)_J &
\specialarrow{\via \eta}&
(f_\bullet|_J)_*(?)_J f_\bullet^*&
\specialarrow{\via c^{-1}}&
(?)_J (f_\bullet)_* f_\bullet^*
\end{array}
\end{equation}
is commutative.
By the definition of $\theta$ and the identity 
$(\varepsilon (f_\bullet|_J)^*)\circ ((f_\bullet|_J)^* u)=\id$, it is easy
to check that $\eta^{-1}\theta=\id$.
This shows $\theta=\eta$, and we are done.
\qed

\begin{Lemma}\label{lipman.thm}
Let the notation be as above, and 
$\Cal M,\Cal N\in\natural(Y_\bullet)$.
Then the diagram
{\small
\begin{equation}\label{lipman.eq}
\begin{array}{ccccc}
   (f_\bullet|_J)^*_\natural(\Cal M_J\otimes \Cal N_J) & \specialarrow{m}&
   ((f_\bullet|_J)^*_\natural((\Cal M\otimes \Cal N)_J) & \specialarrow\theta&
   ((f_\bullet)^*_\natural(\Cal M\otimes\Cal N))_J\\
   \sdarrow\Delta & & & & \sdarrow{(?)_J\Delta}\\
   (f_\bullet|_J)^*_\natural\Cal M_J \otimes (f_\bullet|_J)^*_\natural\Cal N_J&
   \specialarrow{\theta\otimes\theta}&
   ((f_\bullet)^*_\natural\Cal M)_J\otimes ((f_\bullet)^*_\natural\Cal N)_J&
   \specialarrow m &
   ((f_\bullet)^*_\natural\Cal M \otimes (f_\bullet)^*_\natural\Cal N)_J
 \end{array}
\end{equation}
}
is commutative.
\end{Lemma}

\proof This is an immediate consequence of
Lemma~\ref{theta-delta-m.thm}.
\qed

\begin{Corollary} 
The adjoint pair $((?)^*_{\Mod},(?)_*^{\Mod})$ 
over the category $\Cal P(I,\Sch/S)$ is Lipman.
\end{Corollary}

\proof Let $f_\bullet:X_\bullet \rightarrow Y_\bullet$ 
be a morphism of $\Cal P(I,\Sch/S)$.
The morphism $\eta:f^*\Cal O_{Y_\bullet}
\rightarrow \Cal O_{X_\bullet}$ is an isomorphism, like any
such morphism arising from a ringed continuous functor.
Let us consider $\Cal M,\Cal N\in\natural(Y_\bullet)$.
To verify that $\Delta$ is an isomorphism, 
it suffices to show that
\[
(?)_i\Delta: 
  (f^*_\bullet(\Cal M\otimes\Cal N))_i
  \rightarrow
  (f^*_\bullet\Cal M\otimes f^*_\bullet\Cal N)_i
\]
is an isomorphism for any $i\in\ob(I)$.
Now consider the diagram (\ref{lipman.eq}) for $J=i$.
Horizontal maps in the diagram are isomorphisms by (\ref{tensor-data.par})
and (\ref{inverse-image-restrict.thm}).
The left $\Delta$ is an isomorphism, since $f_i$ is a morphism of 
single schemes.
By Lemma~\ref{lipman.thm}, $(?)_i\Delta$ is also an isomorphism.
\qed

\paragraph \label{inv-im-translate.par}
The description of the translation map $\alpha_\phi$ is as follows.
For $\phi\in I(i,j)$, 
\[
\alpha_\phi: X_\phi^* (?)_ i f^*_\bullet\rightarrow (?)_jf^*_\bullet
\]
is the composite
\[
X_\phi^*(?)_i f^*_\bullet\specialarrow{X_\phi^*\theta^{-1}}
X_\phi^*f_i^*(?)_i\cong
f_j^* Y_\phi^*(?)_i
\specialarrow{f_j^*\alpha_\phi}
f_j^* (?)_j
\specialarrow{\theta}(?)_jf^*_\bullet
.
\]

\paragraph Let $X_\bullet\in\Cal P$, and $\Cal M,\Cal N\in\natural(X_\bullet)$.
Although there is a canonical map
\[
(?)_i:\uHom_{\natural(X_\bullet)}(\Cal M,\Cal N)_i
\rightarrow\uHom_{\natural(X_i)}(\Cal M_i,\Cal N_i)
\]
arising from the closed structure for $i\in I$,
this may not be an isomorphism. 
However, we have the following.

\begin{Lemma}\label{uhom.thm}
Let $i\in I$.
If $\Cal M$ is equivariant, then the canonical map
\[
H_i :\uHom_{\natural(X_\bullet)}(\Cal M,\Cal N)_i
\rightarrow\uHom_{\natural(X_i)}(\Cal M_i,\Cal N_i)
\]
is an isomorphism.
\end{Lemma}

\proof It suffices to prove that
\[
H_i:\Hom_{\natural(\Zar(X_\bullet)/(i,U))}(\Cal M|_{(i,U)},
\Cal N|_{(i,U)})\rightarrow
\Hom_{\natural(U)}(\Cal M_i|_U,\Cal N_i|_U)
\]
is an isomorphism for any Zariski open set $U$ in $X_i$.

To give an element of $\varphi\in
\Hom_{\natural(\Zar(X_\bullet)/(i,U))}(\Cal M|_{(i,U)},\Cal N|_{(i,U)})$
is the same as to give a family $(\varphi_\phi)_{\phi:i\rightarrow j}$ with
\[
\varphi_\phi\in\Hom_{\natural(X_\phi^{-1}(U))}
   (\Cal M_j|_{X_\phi^{-1}(U)},\Cal N_j|_{X_\phi^{-1}(U)})
\]
such that 
\[
\varphi_\phi\circ (\alpha_\phi(\Cal M))|_{X_\phi^{-1}(U)}
= (\alpha_\phi(\Cal N)
)|_{X_\phi^{-1}(U)}\circ ((X_\phi)|_{X_\phi^{-1}(U)})^*_\natural(
\varphi_{\id_i}).
\]
As $\alpha_\phi(\Cal M)$ is an isomorphism for any $\phi:i\rightarrow j$,
we have that such a $(\varphi_\phi)$ is uniquely determined by
$\varphi_{\id_i}$.
Hence $H_i$ is bijective, as desired.
\qed

\begin{Lemma}\label{uhom2.thm}
Let $J$ be a subcategory of $I$.
If $\Cal M$ is equivariant, then the canonical map
\[
H_J:\uHom_{\natural(X_\bullet)}(\Cal M,\Cal N)_J
\rightarrow\uHom_{\natural(X_J)}(\Cal M_J,\Cal N_J)
\]
is an isomorphism.
\end{Lemma}

\proof It suffices to show that 
\[
(H_J)_i:
(\uHom_{\natural(X_\bullet)}(\Cal M,\Cal N)_J)_i
\rightarrow
\uHom_{\natural(X_J)}(\Cal M_J,\Cal N_J)_i
\]
is an isomorphism for each $i\in J$.
By Lemma~\ref{composition-H.thm}, the composite map
\begin{multline*}
\uHom_{\natural(X_\bullet)}(\Cal M,\Cal N)_i
\cong
(\uHom_{\natural(X_\bullet)}(\Cal M,\Cal N)_J)_i\\
\specialarrow{(H_J)_i}
\uHom_{\natural(X_J)}(\Cal M_J,\Cal N_J)_i
\specialarrow{H_i}
\uHom_{\natural(X_i)}(\Cal M_i,\Cal N_i)
\end{multline*}
agrees with $H_i$.
As $\Cal M_J$ is also equivariant, we have that the two $H_i$ are
isomorphisms by Lemma~\ref{uhom.thm}, and hence
$(H_J)_i$ is an isomorphism for any $i\in J$.
\qed

\paragraph\label{hom-structure.par}
By the lemma, the sheaf $\uHom_{\natural(X_\bullet)}(\Cal M,\Cal N)$ is
given by the collection 
\[
(\uHom_{\natural(X_i)}(\Cal M_i,\Cal N_i))_{i\in I}
\]
provided $\Cal M$ is equivariant.
The structure map is the canonical composite map
\begin{multline*}
\alpha_\phi:(X_\phi)_\natural^*\uHom_{\natural(X_i)}(\Cal M_i,\Cal N_i)
\specialarrow{P}
\uHom_{\natural(X_j)}((X_\phi)_\natural^*\Cal M_i,(X_\phi)_\natural^*\Cal N_i)
\\
\specialarrow{\suHom_{\natural(X_j)}(\alpha_\phi^{-1},\alpha_\phi)}
\uHom_{\natural(X_j)}(\Cal M_j,\Cal N_j).
\end{multline*}

Similarly, the following is also easy to prove.

\begin{Lemma}\label{initial-equivariant.thm}
Let $i\in I$ be an initial object of $I$.
Then the following hold:
\begin{description}
\item[\truebf 1] If $\Cal M\in\natural(X_\bullet)$ is equivariant,
then
\[
(?)_i:\Hom_{\natural(X_\bullet)}(\Cal M,\Cal N)\rightarrow
\Hom_{\natural(X_i)}(\Cal M_i,\Cal N_i)
\]
is an isomorphism.
\item[\truebf 2] $(?)_i:\EM(X_\bullet)\rightarrow\Mod(X_i)$ is an equivalence,
whose quasi-inverse is $L_i$.
\end{description}
\end{Lemma}

The fact that $L_i(\Cal M)$ is equivariant for $\Cal M\in\Mod(X_i)$
is checked directly from the definition.

\section{Quasi-coherent sheaves over a diagram of schemes}

Let $I$ be a small category, $S$ a scheme, and $X_\bullet\in\Cal P(I,
\Sch/S)$.

\paragraph Let $\Cal M\in\Mod(X_\bullet)$.
We say that $\Cal M$ is locally quasi-coherent (resp.\ locally coherent)
if $\Cal M_i$ is quasi-coherent (resp.\ coherent) for any $i\in I$.
We say that $\Cal M$ is {\em quasi-coherent}
if for any $(i,U)\in\Zar(X_\bullet)$ with $U=\Spec A$ being affine, 
there exists an exact sequence
in $\Mod(\Zar(X_\bullet)/(i,U))$ of the form
\begin{equation}\label{quasi-coherent.eq}
\Cal (O_{X_\bullet}|_{(i,U)})^{(T)}\rightarrow
\Cal (O_{X_\bullet}|_{(i,U)})^{(\Sigma)}\rightarrow
\Cal M|_{(i,U)}\rightarrow 
0,
\end{equation}
where $T$ and $\Sigma$ are arbitrary small sets.

\begin{Lemma}\label{quasi-coherent.thm}
Let $\Cal M\in\Mod(X_\bullet)$.
Then the following are equivalent.
\begin{description}
\item[1] $\Cal M$ is quasi-coherent.
\item[2] $\Cal M$ is locally quasi-coherent and equivariant.
\item[3] For any morphism $(\phi,h):(j,V)\rightarrow (i,U)$ in
$\Zar(X_\bullet)$ such that $V=\Spec B$ and $U=\Spec A$ are affine,
the canonical map $B\otimes_A \Gamma((i,U),\Cal M)\rightarrow \Gamma((j,V),
\Cal M)$ is an isomorphism.
\end{description}
\end{Lemma}

\proof {\bf 1$\Rightarrow$2} Let $i\in I$ and $U$ an affine open subset
of $X_i$.
Then there is an exact sequence of the form (\ref{quasi-coherent.eq}).
Applying the restriction functor $\Mod(\Zar(X_\bullet)/(i,U))\rightarrow
\Mod(U)$,
we get an exact sequence
\[
\Cal O_U^{(T)}\rightarrow
\Cal O_U^{(\Sigma)}\rightarrow
\Cal (\Cal M_i)|_U\rightarrow 
0,
\]
which shows that $\Cal M_i$ is quasi-coherent for any $i\in I$.
We prove that $\alpha_\phi(\Cal M)$ is an isomorphism for any $\phi:i
\rightarrow j$, to show that $\Cal M$ is equivariant.
Take an affine open covering $(U_\lambda)$ of $X_i$, and we prove that
$\alpha_\phi(\Cal M)$ is an isomorphism over $X_\phi^{-1}(U_\lambda)$ for
each $\lambda$.
But this is obvious by the existence of an
exact sequence of the form (\ref{quasi-coherent.eq}) and the five lemma.

{\bf 2$\Rightarrow$3}
Set $W:=X_\phi^{-1}(U)$, and let $\iota:V\hookrightarrow W$ be the
inclusion map.
Obviously, we have $h=(X_\phi)|_W\circ\iota$.
As $\Cal M$ is equivariant, we have that the canonical map
\[
\alpha_\phi|_W(\Cal M):{(X_\phi)|_{W}}^*_{\Mod}(\Cal M_i)|_U 
\rightarrow(\Cal M_j)|_W
\]
is an isomorphism.
Applying $\iota^*_{\Mod}$ to the isomorphism, we have that
$h^*_{\Mod}((\Cal M_i)|_U)\cong (\Cal M_j)|_V$.
The assertion follows from the assumption that $\Cal M_i$ is quasi-coherent.

{\bf 3$\Rightarrow$1}
Let $(i,U)\in\Zar(X_\bullet)$ with $U=\Spec A$ affine.
There is a presentation of the form
\[
A^{(T)}\rightarrow A^{(\Sigma)}\rightarrow \Gamma((i,U),\Cal M)\rightarrow 0.
\]
It suffices to prove that the induced sequence (\ref{quasi-coherent.eq})
is exact.
To verify this, it suffices to prove that the sequence is exact
after taking the section at $((\phi,h):(j,V)\rightarrow (i,U))\in\Zar(
X_\bullet)/(i,U)$ with $V=\Spec B$ being affine.
We have a commutative diagram
\[
\begin{array}{ccccccc}
B\otimes_A A^{(T)} & \rightarrow & B\otimes_A A^{(\Sigma)}&
\rightarrow &B\otimes\Gamma((i,U),\Cal M) & \rightarrow & 0\\
\downarrow\hbox to 0pt{\small$\cong$\hss} & &
\downarrow\hbox to 0pt{\small$\cong$\hss} & &
\downarrow\hbox to 0pt{\small$\cong$\hss} & & \\
B^{(T)} & \rightarrow & B^{(\Sigma)} & \rightarrow &\Gamma((j,V),\Cal M)&
\rightarrow & 0
\end{array}
\]
whose first row is exact and vertical arrows are isomorphisms.
Hence, the second row is also exact, and (\ref{quasi-coherent.eq}) 
is exact.
\qed

\begin{Def} We say that $\Cal M\in\Mod(X_\bullet)$ is {\em coherent}
if it is equivariant and $\Cal M_i$ is coherent for any $i\in I$.
We denote the full subcategory of $\Mod(X_\bullet)$ consisting of
coherent objects by $\Coh(X_\bullet)$.
\end{Def}

\paragraph Let $J\subset I$ be a subcategory.
We say that $J$ is big in $I$ if for any $(\psi:j\rightarrow k)\in\Mor(I)$,
there exists some $(\phi:i\rightarrow j)\in\Mor(J)$ such that
$\psi\circ\phi\in\Mor(J)$.
Note that $\ob(J)=\ob(I)$ if $J$ is big in $I$.
Let $\Bbb Q$ be a property of morphisms of schemes.
We say that $X_\bullet$ has $\Bbb Q$ $J$-arrows if
$(X_\bullet)|_J$ has $\Bbb Q$-arrows.

\begin{Lemma}\label{quasi-coherent-thick.thm}
Let $J\subset I$ be a subcategory, and $\Cal M\in
\Mod(X_\bullet)$.
\begin{description}
\item[\truebf 1] 
The full subcategory $\LQco(X_\bullet)$ of $\Mod(X_\bullet)$ 
consisting of locally quasi-coherent objects is a plump subcategory.
\item[\truebf 2]
If $\Cal M$ is equivariant \(resp.\ locally quasi-coherent, quasi-coherent\),
then so is $\Cal M_J^{\Mod}$.
\item[\truebf 3]
If $J$ is big in $I$ and $\Cal M|_J^{\Mod}$ is equivariant \(resp.\ 
locally quasi-coherent, quasi-coherent\),
then so is $\Cal M$.
\item[\truebf 4]
If $J$ is big in $I$ and $X_\bullet$ has flat $J$-arrows, then
the full subcategory $\EM(X_\bullet)$ \(resp.\ $\Qco(X_\bullet)$\) 
of $\Mod(X_\bullet)$ 
consisting of equivariant (resp.\ quasi-coherent) 
objects is a plump subcategory.
\item[\truebf 5] If $J$ is big in $I$, then $(?)_J$ is faithful and exact.
\end{description}
\end{Lemma}

\proof {\bf 1} and {\bf 2} are trivial.

We prove {\bf 3}.
The assertion for the local quasi-coherence is obvious, because we have
$\ob(J)=\ob(I)$.
By Lemma~\ref{quasi-coherent.thm}, it remains to show the
assertion for the equivariance.
Let us assume that $\Cal M|_J^{\Mod}$ is equivariant and 
$\psi:j\rightarrow
k$ be a morphism in $I$, and take $\phi:i\rightarrow j$ such that
$\phi,\psi\phi\in\Mor(J)$.
Then the composite map
\[
(X_{\psi\phi})^*_{\Mod}(\Cal M_i)\cong
(X_\psi)^*_{\Mod}(X_\phi)^*_{\Mod}(\Cal M_i)
\specialarrow{(X_\psi)^*_{\Mod}\alpha_\phi^{\Mod}}
(X_\psi)^*_{\Mod}(\Cal M_j)
\specialarrow{\alpha_\psi^{\Mod}}
\Cal M_k,
\]
which agrees with $\alpha_{\psi\phi}^{\Mod}$, is an isomorphism by assumption.
As we have $\alpha_\phi^{\Mod}$ is also an isomorphism,
we have that $\alpha_\psi^{\Mod}$ is an isomorphism.
Thus $\Cal M$ is equivariant.

{\bf 4} By {\bf 1} and Lemma~\ref{quasi-coherent.thm}, it suffices to
prove the assertion only for the equivariance.
Let 
\[
\Cal M_1\rightarrow \Cal M_2\rightarrow \Cal M_3\rightarrow \Cal M_4
\rightarrow \Cal M_5
\]
be an exact sequence in $\Mod(X_\bullet)$, and assume that
$\Cal M_i$ is equivariant for $i=1,2,4,5$.
We prove that $\Cal M_3$ is equivariant.
The sequence remains exact after applying the functor $(?)_J^{\Mod}$.
By {\bf 3}, replacing $I$ by $J$ and 
$X_\bullet$ by $(X_\bullet)|_J$, we may assume that $X_\bullet$ has
flat arrows.
Now the assertion follows easily from the five lemma.

The assertion {\bf 5} is obvious, because $\ob(J)=\ob(I)$.
\qed

\begin{Lemma}\label{indlim-qco.thm}
Let $(\Cal M_\lambda)$ be a diagram in 
$\Mod(X_\bullet)$.
If each $\Cal M_\lambda$ is locally quasi-coherent \(resp.\ equivariant,
quasi-coherent\),
then so is $\indlim \Cal M_\lambda$.
\end{Lemma}

\proof As $(?)_i$ preserves colimits, the assertion for local quasi-coherence
is trivial.
Assume that each $\Cal M_\lambda$ is equivariant.
For $(\phi:i\rightarrow j)\in \Mor(I)$, $\alpha_\phi(\Cal M_\lambda)$
is an isomorphism.
As $\alpha_\phi(\indlim\Cal M_\lambda)$ is nothing but the composite
\[
(X_\phi)^*_{\Mod}((\indlim \Cal M_\lambda)_i)
\cong \indlim (X_\phi)^*_{\Mod}(\Cal M_\lambda)_i
\specialarrow{\indlim \alpha_\phi(\Cal M_\lambda)}
\indlim (\Cal M_\lambda)_j
\cong(\indlim \Cal M_\lambda)_j,
\]
it is an isomorphism.
The rest of the assertions follow.
\qed

By Lemma~\ref{L_J.par}, we have the following.

\begin{Lemma} Let $J\subset I$ be a subcategory, and $\Cal M\in
\LQco(X_\bullet|_J)$.
Then we have $L_J^{\Mod}(\Cal M)\in\LQco(X_\bullet)$.
\end{Lemma}

Similarly, we have the following.

\begin{Lemma}\label{R_i-lqco.thm}
Let $j\in I$.
Assume that $X_\bullet$ has separated quasi-compact arrows,
and that $I(i,j)$ is finite for any $i \in I$.
If $\Cal M\in \Qco(X_j)$, then we have $R_j\Cal M\in\LQco(X_\bullet)$.
\end{Lemma}

The following is also proved easily, using (\ref{tensor-data.par}).

\begin{Lemma}\label{tensor-product-qco.thm}
Let $\Cal M$ and $\Cal N$ be
locally quasi-coherent \(resp.\ 
equivariant, quasi-coherent\) $\Cal O_{X_\bullet}$-modules.
Then $\Cal M\otimes_{\Cal O_{X_\bullet}}
\Cal N$ is also locally quasi-coherent
\(resp.\
equivariant, quasi-coherent\).
\end{Lemma}

The following is a consequence of the observation in (\ref{hom-structure.par}).

\begin{Lemma}\label{hom-qco.thm}
Let $\Cal M$ be a coherent $\Cal O_{X_\bullet}$-module, and
$\Cal N$ a locally quasi-coherent $\Cal O_{X_\bullet}$-module.
Then $\uHom_{{\Mod}(X_\bullet)}(\Cal M,\Cal N)$ is locally quasi-coherent.
If, moreover, there is a big subcategory $J$ of $I$ such that
$X_\bullet$ has flat $J$-arrows and $\Cal N$ is quasi-coherent, then
$\uHom_{{\Mod}(X_\bullet)}(\Cal M,\Cal N)$ is quasi-coherent.
\end{Lemma}

\paragraph Let $f_\bullet:X_\bullet\rightarrow Y_\bullet$ be 
a morphism in $\Cal P$.

As $\theta$ in (\ref{direct-image-translate2.eq}), which appears in
(\ref{direct-image-translate.eq})
is not an isomorphism in general, $(f_\bullet)_*^\natural(\Cal M)$
need not be equivariant even if $\Cal M$ is equivariant.
However, we have

\begin{Lemma}\label{equivariant-direct-image.thm}
Let $f_\bullet:X_\bullet\rightarrow Y_\bullet$ be a morphism in $\Cal P$,
and $J$ a big subcategory of $I$.
Then we have the following:
\begin{description}
\item[\truebf 1] $f_\bullet$ is of fiber type if and only if $(f_\bullet)|_J$
is of fiber type.
\item[\truebf 2] If $f_\bullet$ is quasi-compact separated and
$\Cal M\in\LQco(X_\bullet)$, then $(f_\bullet)_*(\Cal M)\in\LQco(Y_\bullet)$.
\item[\truebf 3] If $Y_\bullet$ has flat $J$-arrows, $f_\bullet$ is
of fiber type and quasi-compact separated, and $\Cal M\in\Qco(X_\bullet)$,
then we have $(f_\bullet)_*(\Cal M)\in\Qco(Y_\bullet)$.
\end{description}
\end{Lemma}

\proof {\truebf1}
Assume that $f_\bullet|_J$ is of fiber type, and
let $\psi:j\rightarrow k$ be a morphism in $I$.
Take $\phi:i\rightarrow j$ such that $\phi,\psi\phi\in\Mor(J)$.
Consider the commutative diagram
\[
\begin{array}{ccccc}
X_k & \specialarrow{X_\psi} & X_j & \specialarrow{X_\phi} & X_i\\
\downarrow\hbox to 0pt{\small$f_k$\hss} & \mbox{(a)} &
\downarrow\hbox to 0pt{\small$f_j$\hss} & \mbox{(b)} &
\downarrow\hbox to 0pt{\small$f_i$\hss} \\
Y_k & \specialarrow{Y_\psi} & Y_j & \specialarrow{Y_\phi} & Y_i.
\end{array}
\]
By assumption, the square (b) and the whole rectangle ((a)$+$(b)) are
fiber squares.
Hence (a) is also a fiber square.
This shows that $f_\bullet$ is of fiber type.
The converse is obvious.

The assertion {\bf 2} is obvious by the isomorphism
$(f_\bullet)_*(\Cal M)_i\cong (f_i)_*(\Cal M_i)$ for $i\in I$, see
\cite[(9.2.2)]{EGA-I}.

We prove {\bf 3}.
By Lemma~\ref{quasi-coherent-thick.thm}, we may assume that $J=I$.
Then $(f_\bullet)_*(\Cal M)$ is locally quasi-coherent by {\bf 2}.
As $\Cal M$ is equivariant and $\theta$ in 
(\ref{direct-image-translate2.eq}) is
an isomorphism by \cite[(1.4.15)]{EGA-III}, we have that
$(f_\bullet)_*(\Cal M)$ is equivariant.
Hence by Lemma~\ref{quasi-coherent.thm}, $(f_\bullet)_*(\Cal M)$ is
quasi-coherent.
\qed

\paragraph
If $f_\bullet$ is quasi-compact separated of fiber type and $Y_\bullet$
has flat $J$-arrows for some big subcategory $J$ of $I$, 
then $(f_\bullet)_*^{\Qco}:\Qco(X_\bullet)\rightarrow\Qco(Y_\bullet)$
is defined as the restriction of $(f_\bullet)_*^{\Mod}$.

Similarly, we have the following.

\begin{Lemma}\label{fiber-type-1.thm}
Let $f_\bullet:X_\bullet\rightarrow Y_\bullet$ and
$g_\bullet:Y_\bullet\rightarrow Z_\bullet$ be morphisms in $\Cal P$.
Then the following hold.
\begin{description}
\item[\truebf 0] An isomorphism is a morphism of fiber type.
\item[\truebf 1] If $f_\bullet$ and $g_\bullet$ are of fiber type, then
so is $g_\bullet\circ f_\bullet$.
\item[\truebf 2] If $g_\bullet$ and $g_\bullet\circ f_\bullet$ are
of fiber type, then so is $f_\bullet$.
\item[\truebf 3] If $f_\bullet$ is faithfully flat of fiber type and
$g_\bullet\circ f_\bullet$ is of fiber type, then $g_\bullet$ is
of fiber type.
\end{description}
\end{Lemma}

\begin{Lemma}\label{fiber-type-2.thm}
Let $f_\bullet:X_\bullet \rightarrow Y_\bullet$ and
$g_\bullet : Y_\bullet'\rightarrow Y_\bullet$ be morphisms in $\Cal P$.
Let $f_\bullet':X'_\bullet\rightarrow Y_\bullet'$ be the
base change of $f_\bullet$ by $g_\bullet$.
\begin{description} 
\item[\bf 1] If $f_\bullet$ is of fiber type, then so is $f_\bullet'$.
\item[\bf 2] If $f_\bullet'$ is of fiber type and $g_\bullet$ is
faithfully flat, then $f_\bullet$ is of fiber type.
\end{description}
\end{Lemma}

\paragraph Let $f:X\rightarrow Y$ be a morphism of schemes.
If $f$ is quasi-compact separated, then $f_*^{{\Qco}}$ 
is compatible with filtered inductive limits.

\def\citinfo{\cite[p.641, Proposition~6]{Kempf}}
\begin{Lemma}[\citinfo]\label{direct-image-flim-schemes.thm}
Let $f:X\rightarrow Y$ be a quasi-compact separated morphism of schemes, and
$(\Cal M_i)$ a pseudo-filtered inductive system of $\Cal O_X$-modules.
Then the canonical map
\[
\indlim f_*\Cal M_i\rightarrow f_*\indlim \Cal M_i
\]
is an isomorphism.
\end{Lemma}

By the lemma, the following follows immediately.

\begin{Lemma}\label{direct-image-flim.thm}
Let $f_\bullet:X_\bullet\rightarrow Y_\bullet$ be a
morphism in $\Cal P(I,\Sch/S)$.
If $f_\bullet$ is quasi-compact and separated, 
then $(f_\bullet)_*^{{\Mod}}$ and
$(f_\bullet)_*^{{\LQco}}$ preserve pseudo-filtered inductive limits.
\end{Lemma}

\begin{Lemma}\label{direct-im-left-ind.thm}
Let $f_\bullet:X_\bullet\rightarrow Y_\bullet$ be a morphism in $\Cal P$.
Let $J$ be an admissible subcategory of $I$.
If $Y_\bullet$ has flat arrows and $f_\bullet$ is
of fiber type and quasi-compact separated, then the canonical map
\[
\theta(J,f_\bullet): L_J\circ (f_J)_*\rightarrow (f_\bullet)_* \circ L_J
\]
is an isomorphism of functors from $\LQco(X_\bullet|_J)$ to
$\LQco(Y_\bullet)$.
\end{Lemma}

\proof This is obvious by Proposition~\ref{crutial.thm}, 
Lemma~\ref{direct-image-flim-schemes.thm}, and
\cite[(1.4.15)]{EGA-III}.
\qed

The following is obvious by Lemma~\ref{inverse-image-restrict.thm}
and (\ref{inv-im-translate.par}).

\begin{Lemma}\label{inv-im-qco.thm}
Let $f_\bullet:X_\bullet\rightarrow Y_\bullet$ be a
morphism in $\Cal P$.
If $\Cal M\in\Mod(Y_\bullet)$ is equivariant \(resp.\ locally quasi-coherent, 
quasi-coherent\), then so is $(f_\bullet)^*_{\Mod}(\Cal M)$.
If $\Cal M\in\Mod(Y_\bullet)$, $f_\bullet$ is faithfully flat,
and $(f_\bullet)^*_{\Mod}(\Cal M)$ is equivariant, 
then we have $\Cal M$ is equivariant.
\end{Lemma}

The restriction $(f_\bullet)^*:\Qco(Y_\bullet)\rightarrow\Qco(X_\bullet)$
is sometimes denoted by $(f_\bullet)^*_{\Qco}$.

\section[Derived functors]{Derived functors 
of functors on sheaves of modules over diagrams of schemes}

Let $I$ be a small category, and $S$ a scheme.
Set $\Cal P:=\Cal P(I,\Sch/S)$, and let $X_\bullet\in\Cal P$.

\begin{Proposition}\label{K-limp-restrict.thm}
Let $X_\bullet\in \Cal P$, and $\Bbb I\in K(\Mod(X_\bullet))$.
We have $\Bbb I$ is $K$-limp if and only if so is $\Bbb I_i$ for $i\in I$.
\end{Proposition}

\proof 
The only if part is proved in Lemma~\ref{admissible-K.thm}.

We prove the if part.
Let $\Bbb I\rightarrow \Bbb J$ be a $K$-injective resolution,
and let $\Bbb C$ be the mapping cone.
Note that $\Bbb C_i$ is exact for each $i$.

Let $(U,i)\in \Zar(X_\bullet)$.
We have an isomorphism
\[
\Gamma((U,i),\Bbb C)\cong \Gamma(U,\Bbb C_i).
\]
As $\Bbb C_i$ is $K$-limp by the only if part, these are exact
for each $(U,i)$.
It follows that $\Bbb I$ is $K$-limp.
\qed

\begin{Corollary}\label{derived-direct-restrict.thm}
Let $J$ be a subcategory of $I$, and $f_\bullet:X_\bullet\rightarrow
Y_\bullet$ a morphism in $\Cal P(I,\Sch/S)$.
Then there is an isomorphism
\[
c(J,f_\bullet): (?)_J R(f_\bullet)_*\cong R(f_\bullet|_J)_* (?)_J.
\]
\end{Corollary}

\paragraph Let $X$ be a scheme, $x\in X$, and $M$ an $\Cal O_{X,x}$-module.
We define $\xi_x(M)\in\Mod(X)$ by $\Gamma(U,\xi_x(M))=M$ if $x\in U$,
and zero otherwise.
The restriction maps are defined in an obvious way.
For an exact complex $\Bbb H$ of $\Cal O_{X,x}$-modules,
$\xi_x(\Bbb H)$ is exact not only as a complex of sheaves,
but also as a complex of presheaves.
For a morphism of schemes $f:X\rightarrow Y$,
we have that $f_*\xi_x(M)\cong \xi_{f(x)}(M)$.

\begin{Lemma}\label{K-flat2.thm}
Let $\Bbb F\in C(\Mod(X_\bullet))$.
The following are equivalent.
\begin{description}
\item[1] $\Bbb F$ is $K$-flat.
\item[2] $\Bbb F_i$ is $K$-flat for $i\in\ob(I)$.
\item[3] $\Bbb F_{i,x}$ is a $K$-flat complex of 
$\Cal O_{X,x}$-modules for any $i\in \ob(I)$.
\end{description}
\end{Lemma}

\proof
{\bf 3$\Rightarrow$1} Let $\Bbb G\in C(\Mod(X_\bullet))$
be exact.
We are to prove that $\Bbb F\otimes_{\Cal O_{X_\bullet}}\Bbb G$ is exact.
For $i\in \ob(I)$ and $x\in X$, we have
\[
(\Bbb F\otimes_{\Cal O_{X_\bullet}}\Bbb G)_{i,x}
\cong
(\Bbb F_i\otimes_{\Cal O_{X_i}}\Bbb G_i)_x
\cong
\Bbb F_{i,x}\otimes_{\Cal O_{X_i,x}}\Bbb G_{i,x}.
\]
Since $\Bbb G_{i,x}$ is exact, $(\Bbb F\otimes_{\Cal O_{X_i,x}}
\Bbb G)_{i,x}$ is exact.
So $\Bbb F\otimes_{\Cal O_{X_\bullet}}\Bbb G$ is exact.

{\bf 1$\Rightarrow$3}
Let $\Bbb H\in C(\Mod(\Cal O_{X_i,x}))$ be an exact complex, and
we are to prove that $\Bbb F_{i,x}\otimes_{\Cal O_{X_i,x}}\Bbb H$ is
exact.
For each $j\in\ob(I)$, 
\[
(?)_j R_i\xi_x(\Bbb H)\cong
\prod_{\phi\in I(j,i)}(X_\phi)_*\xi_x(\Bbb H)
\cong
\prod_{\phi\in I(j,i)} \xi_{X_\phi(x)}(\Bbb H)
\]
is exact, since a direct product of exact complexes of presheaves is
exact.
So $R_i\xi_x(\Bbb H)$ is exact.
It follows that $\Bbb F\otimes_{\Cal O_{X_\bullet}}R_i\xi_x(\Bbb H)$ is
exact.
Hence
\[
(\Bbb F\otimes_{\Cal O_{X_\bullet}}R_i\xi_x(\Bbb H))_i
\cong
\Bbb F_i\otimes_{\Cal O_{X_i}}\prod_{\phi\in I(i,i)}\xi_{X_\phi(x)}\Bbb H
\]
is also exact.
So $\Bbb F_i\otimes_{\Cal O_{X_i}}\xi_{\id_X(x)}(\Bbb H)
=
\Bbb F_i\otimes_{\Cal O_{X_i}}\xi_x\Bbb H$ is exact.
So
\[
(\Bbb F_i\otimes_{\Cal O_{X_i}}\xi_x\Bbb H)_x
\cong
\Bbb F_{i,x}\otimes_{\Cal O_{X_i,x}}(\xi_x\Bbb H)_x
\cong
\Bbb F_{i,x}\otimes_{\Cal O_{X_i,x}}\Bbb H
\]
is also exact.

Applying {\bf 1$\Leftrightarrow$3}, which has already been proved, to the
complex $\Bbb F_i$ over the single scheme $X_i$, 
we get {\bf 2$\Leftrightarrow$3}.
\qed

Hence by \cite{Spaltenstein}, we have the following.

\begin{Lemma} Let $f_\bullet:X_\bullet\rightarrow Y_\bullet$ be 
a morphism in $\Cal P$.
Then we have the following.
\begin{description}
\item[1] If $\Bbb F\in C(\Mod(Y_\bullet))$ is $K$-flat, then so is
$f_\bullet^*\Bbb F$.
\item[2] If $\Bbb F \in C(\Mod(Y_\bullet))$ is $K$-flat exact, then
so is $f_\bullet^*\Bbb F$.
\item[3] If $\Bbb I\in C(\Mod(X_\bullet))$ is weakly $K$-injective, then
so is $(f_\bullet)_*\Bbb I$.
\end{description}
\end{Lemma}

\paragraph
By the lemma, the left derived functor $Lf^*_\bullet$, 
which we already know its
existence by Lemma~\ref{inverse-image-restrict.thm}, 
can also be calculated by $K$-flat resolutions.

\begin{Lemma}\label{theta-f-J.thm}
Let $J$ be a subcategory of $I$, and 
$f_\bullet:X_\bullet\rightarrow Y_\bullet$ a morphism in $\Cal P$.
Then we have the following.
\begin{description}
\item[1] 
The canonical map
\[
\theta(f_\bullet,J):L(f_\bullet|_J)^*(?)_J\rightarrow (?)_J Lf_\bullet^*
\]
is an isomorphism.
\item[2] The diagram
\[
\begin{array}{ccccc}
(?)_J  & \multicolumn{3}{c}{
\mathord -\mkern -6mu
\cleaders \hbox {$\mkern -2mu\mathord -\mkern -2mu$}
\hfill
\hbox to 0pt{\hss$
  \mathord{\mathop{\mkern -1mu\mathord-\mkern -1mu}\limits^{\hbox to
  0pt{\hss\scriptsize $\id$\hss}}}$
   \hss}
\mkern -1mu
\cleaders \hbox {$\mkern -2mu\mathord -\mkern -2mu$}
\hfill\mkern -6mu\mathord \rightarrow
} & (?)_J\\
\downarrow\hbox to 0pt{\scriptsize$\via u$\hss} & & & &
\downarrow\hbox to 0pt{\scriptsize$\via u$\hss}\\
R(f_\bullet|_J)_*L(f_\bullet|_J)^*(?)_J &
\specialarrow{\via \theta}&
R(f_\bullet|_J)_*(?)_J L f_\bullet^*&
\specialarrow{\via c^{-1}}&
(?)_J R (f_\bullet)_* L f_\bullet^*
\end{array}
\]
is commutative.
\item[3] The diagram
\[
\begin{array}{ccccc}
(?)_J  & \multicolumn{3}{c}{
\mathord -\mkern -6mu
\cleaders \hbox {$\mkern -2mu\mathord -\mkern -2mu$}
\hfill
\hbox to 0pt{\hss$
  \mathord{\mathop{\mkern -1mu\mathord-\mkern -1mu}\limits^{\hbox to
  0pt{\hss\scriptsize $\id$\hss}}}$
   \hss}
\mkern -1mu
\cleaders \hbox {$\mkern -2mu\mathord -\mkern -2mu$}
\hfill\mkern -6mu\mathord \rightarrow
} & (?)_J\\
\uparrow\hbox to 0pt{\scriptsize$\via \varepsilon$\hss} & & & &
\uparrow\hbox to 0pt{\scriptsize$\via \varepsilon$\hss}\\
L(f_\bullet|_J)^*R(f_\bullet|_J)_*(?)_J &
\specialarrow{\via c^{-1}}&
L(f_\bullet|_J)^*(?)_J R (f_\bullet)_*&
\specialarrow{\via \theta}&
(?)_J L(f_\bullet)^* R (f_\bullet)_*
\end{array}
\]
is commutative.
\end{description}
\end{Lemma}

\proof 
Since $(?)_J$ preserves $K$-flat complexes by Lemma~\ref{K-flat2.thm},
we have
\[
L(f_\bullet|_J)^*(?)_J\cong L((f_\bullet|_J)\circ (?)_J).
\]
On the other hand, it is obvious that we have $(?)_J L f_\bullet^*
\cong L((?)_J f_\bullet^*)$.
By Lemma~\ref{inverse-image-restrict.thm}, we have a composite isomorphism
\[
\theta': L(f_\bullet|_J)^*(?)_J
  \cong 
L((f_\bullet|_J)^*\circ(?)_J)
  \specialarrow{L\theta}
L((?)_J\circ f_\bullet^*)
  \cong
(?)_J L(f_\bullet)^*.
\]

It is straightforward 
to check the diagram in {\bf 2} ($\theta$ there should be
replaced by $\theta'$) is commutative
utilizing
(\ref{f_*-f^*-unit-restrict.eq}) (note that $\eta$ there equals $\theta$).
As in the proof of Lemma~\ref{inverse-image-restrict.thm}, 
we have $\theta=\theta'$, and {\bf 1} and {\bf 2} are proved.
The proof of {\bf 3} is formal, and we omit it.
\qed

\begin{Lemma} Let $X_\bullet\in\Cal P$, and $\Bbb F,\Bbb G
\in D(\Mod(X_\bullet))$.
Then we have the following.
\begin{description}
\item[1] $\Bbb F_J\otimes_{\Cal O_{X_\bullet|_J}}^{\bullet,L}\Bbb G_J
\cong 
(\Bbb F\otimes_{\Cal O_{X_\bullet}}^{\bullet,L}\Bbb G)_J$
for any subcategory $J\subset I$.
\item[2] If $\Bbb F$ and $\Bbb G$ has locally quasi-coherent cohomology
groups, then $\uTor^{\Cal O_{X_\bullet}}_i(\Bbb F,\Bbb G)$ is also
locally quasi-coherent for any $i\in\Bbb Z$.
\item[3] Assume that there exists some big subcategory $J$ of $I$
such that $X_\bullet$ has flat $J$-arrows.
If both $\Bbb F$ and $\Bbb G$ has equivariant (resp.\ quasi-coherent) 
cohomology groups, 
then $\uTor^{\Cal O_{X_\bullet}}_i(\Bbb F,\Bbb G)$ is also
equivariant (resp.\ quasi-coherent).
\end{description}
\end{Lemma}

\proof The assertion {\bf 1} is an immediate consequence of
Lemma~\ref{K-flat2.thm} and Example~\ref{lipman.ex}, {\bf 5}.

{\bf 2} In view of {\bf 1}, we may assume that $X=X_\bullet$ is
a single scheme.
As the question is local, we may assume that $X$ is even affine.

We may assume that $\Bbb F=\indlim \Bbb F_n$, where $(\Bbb F_n)$
is the $\frak P(X_\bullet)$-special direct system such that 
each $\Bbb F_n$ is bounded above and has locally quasi-coherent cohomology
groups as in Lemma~\ref{strongly-K-limp.thm}, {\bf 4}.
Similarly, we may assume that $\Bbb G=\indlim \Bbb G_n$.
As filtered inductive limits are exact and compatible with tensor products,
and the colimit of locally quasi-coherent sheaf is locally quasi-coherent,
we may assume that both $\Bbb F$ and $\Bbb G$ are bounded above, flat, and
has locally quasi-coherent cohomology groups.
By \cite[Proposition~I.7.3]{Hartshorne2}, we may assume that
both $\Bbb F$ and $\Bbb G$ are single quasi-coherent sheaves.
This case is trivial.

{\bf 3} 
In view of {\bf 1}, we may assume that $J=I$ and $X_\bullet$ has
flat arrows.
By {\bf 2}, it suffices to show the assertion for equivariance.
Assuming that $\Bbb F$ and $\Bbb G$ are $K$-flat with equivariant
cohomology groups, we prove that $\Bbb F\otimes \Bbb G$ has equivariant
cohomology groups.
This is enough.

Let $\phi:i\rightarrow j$ be a morphism of $I$.
As $X_\phi$ is flat and $\Bbb F$ and $\Bbb G$ have
equivariant cohomology groups,
$\alpha_\phi: X_\phi^*\Bbb F_i\rightarrow \Bbb F_j$ 
and $\alpha_\phi:X_\phi^*\Bbb G_i\rightarrow \Bbb G_j$ are quasi-isomorphisms.
The composite
\[
X_\phi^*(\Bbb F\otimes_{\Cal O_{X_\bullet}}^\bullet\Bbb G)_i\cong
X_\phi^*\Bbb F_i\otimes_{\Cal O_{X_j}}^\bullet X_\phi^*\Bbb G_i
\specialarrow{\alpha_\phi\otimes\alpha_\phi}
\Bbb F_j\otimes_{\Cal O_{X_j}}^\bullet \Bbb G_j
\cong
(\Bbb F\otimes_{\Cal O_{X_\bullet}}^\bullet \Bbb G)_j
\]
is a quasi-isomorphism, since $\Bbb F_j$ and $X_\phi^*\Bbb F_i$ are
$K$-flat.
By (\ref{tensor-data.par}), $\alpha_\phi(\Bbb F\otimes^\bullet_{\Cal O_{
X_\bullet}}\Bbb G)$ is a quasi-isomorphism.

As $X_\bullet$ has flat arrows, this shows that $\Bbb F\otimes^\bullet
_{\Cal O_{X_\bullet}}\Bbb G$ has equivariant cohomology groups.
\qed

\paragraph Let $X_\bullet\in \Cal P$, and $J$ an admissible subcategory
of $I$.
By Lemma~\ref{admissible-subcat.thm}, the left derived functor
\[
LL_J^{\Mod}:D(\Mod(X_\bullet|_J))\rightarrow D(\Mod(X_\bullet))
\]
of $L_J^{\Mod}$ 
is defined, since $Q(X_\bullet,J)$ is admissible.
This is also calculated using $K$-flat resolutions.
Namely, 

\begin{Lemma}
Let $X_\bullet$ and $J$ be as above.
If $\Bbb F\in K(\Mod(X_\bullet|_J))$ is $K$-flat, then
so is $L_J\Bbb F$.
If $\Bbb F$ is $K$-flat exact, then so is $L_J\Bbb F$.
\end{Lemma}

\proof This is trivial by (\ref{L_J.thm}).
\qed

\begin{Corollary}\label{w-k.thm}
Let $X_\bullet\in \Cal P(I,\Sch/S)$, $J$ an admissible subcategory of $I$,
and $\Bbb I\in K(\Mod(X_\bullet))$.
If $\Bbb I$ is weakly $K$-injective, then
$\Bbb I_J$ is weakly $K$-injective.
\end{Corollary}

\proof Let $\Bbb F$ be a $K$-flat exact complex in $K(\Mod(X_\bullet|_J))$.
Then, 
\[
\Hom^\bullet_{\Mod(X_\bullet|_J)}(\Bbb F,\Bbb I_J)\cong
\Hom^\bullet_{\Mod(X_\bullet)}(L_J\Bbb F,\Bbb I)
\]
is exact by the lemma.
By Lemma~\ref{strongly-K-limp.thm}, {\bf 7}, we are done.
\qed

\begin{Lemma}\label{L_i-lqco.thm}
Let $X_\bullet$ and $J$ be as above.
Let $\Bbb F\in D_{\LQco(X_\bullet|J)}(\Mod(X_\bullet|J))$.
Then, $L_J\Bbb F\in D_{\LQco(X_\bullet)}(\Mod(X_\bullet))$.
\end{Lemma}

\proof This is trivial by the standard spectral sequence argument, if
$\Bbb F$ is bounded above.
The general case follows immediately by Lemma~\ref{strongly-K-limp.thm}, 
{\bf 4}.
\qed

\paragraph Let $I$ be a small category, and $S$ a scheme.
Set $\Cal P:=\Cal P(I,\Sch/S)$.
As we have seen, for a morphism $f_\bullet: X_\bullet\rightarrow Y_\bullet$,
$f_\bullet^{-1}:\Zar(Y_\bullet)\rightarrow \Zar(X_\bullet)$ is an admissible
ringed continuous functor by Lemma~\ref{inverse-image-restrict.thm}.
Moreover, if $J$ and $K$ are admissible subcategories of $I$ such that
$J\subset K$, then
$Q(X_\bullet|_J,K):\Zar(X_\bullet|_K)\rightarrow \Zar(X_\bullet|_J)$ is
also admissible.
By Lemma~\ref{lipman-sites.thm} and Lemma~\ref{commutativity-sites.thm},
we have the following.

\begin{example}\label{derived-lipman.ex}
Let $I$ be a small category, $S$ a scheme, and 
$f_\bullet:X_\bullet\rightarrow Y_\bullet$ and
$g_\bullet:Y_\bullet\rightarrow Z_\bullet$
are morphisms in $\Cal P(I,\Sch/S)$.
Let $K\subset J\subset I$ be admissible subcategories.
Then we have the following.
\begin{description}
\item[1] There is a natural isomorphism
\[
c_{I,J,K}: (?)_{K,I}\cong (?)_{K,J}\circ (?)_{J,I}.
\]
Taking the conjugate,
\[
d_{I,J,K}: LL_{I,J}\circ LL_{J,K}\cong LL_{I,K}
\]
is induced.
\item[2] There are natural isomorphism
\[
c_{J,f_\bullet}: (?)_J\circ R(f_\bullet)_* \cong R(f_\bullet|_J)_*\circ (?)_J
\]
and its conjugate
\[
d_{J,f_\bullet}: LL_J\circ L(f_\bullet|_J)^*\cong L(f_\bullet)^*\circ LL_J.
\]
\item[3] We have
\[
(c_{K,f_\bullet|_J}(?)_J)\circ ((?)_{K,J}c_{J,f})
=
(R(f_\bullet|_K)_* c_{I,J,K})\circ 
c_{K,f_\bullet}\circ (c_{I,J,K}^{-1}R(f_\bullet)_*).
\]
\item[4] We have
\[
(R(g_\bullet|_J)_*c_{J,f_\bullet})\circ(c_{J,g_\bullet}R(f_\bullet)_*)
=
(c_{f_\bullet|_J,g_\bullet|_J}(?)_J)
\circ
c_{J,g_\bullet\circ f_\bullet}\circ((?)_J c_{f_\bullet,g_\bullet}^{-1}),
\]
where $c_{f_\bullet,g_\bullet}:R(g_\bullet\circ f_\bullet)_*\cong
R(g_\bullet)_*\circ R(f_\bullet)_*$ is the canonical isomorphism,
and similarly for $c_{f_\bullet|_J,g_\bullet|_J}$.
\item[5] The adjoint pair $(L(?)^*_{\Mod},R(?)_*^{\Mod})$ 
over the category $\Cal P(I,\Sch/S)$ is Lipman.
\end{description}
\end{example}

It is now quite formal to prove the following.

\begin{Lemma}\label{theta-J-f.thm}
Let the notation be as in {\rm Example~\ref{derived-lipman.ex}}.
Let $\theta=\theta(J,f_\bullet)$ denotes the natural map
\[
LL_J R(f_\bullet|_J)_*\rightarrow R(f_\bullet)_* LL_J
\]
defined in {\rm\cite[(3.7.2)]{Lipman}}.
Then we have the following.
\begin{description}
\item[1] The diagram
\[
\begin{array}{ccccc}
LL_J  & \multicolumn{3}{c}{
\mathord -\mkern -6mu
\cleaders \hbox {$\mkern -2mu\mathord -\mkern -2mu$}
\hfill
\hbox to 0pt{\hss$
  \mathord{\mathop{\mkern -1mu\mathord-\mkern -1mu}\limits^{\hbox to
  0pt{\hss\scriptsize $\id$\hss}}}$
   \hss}
\mkern -1mu
\cleaders \hbox {$\mkern -2mu\mathord -\mkern -2mu$}
\hfill\mkern -6mu\mathord \rightarrow
} & LL_J\\
\downarrow\hbox to 0pt{\scriptsize$\via u$\hss} & & & &
\downarrow\hbox to 0pt{\scriptsize$\via u$\hss}\\
LL_J R(f_\bullet|_J)_*L(f_\bullet|_J)^*&
\specialarrow{\via \theta}&
R(f_\bullet)_* LL_J L (f_\bullet|_J)^*&
\specialarrow{\via d}&
R (f_\bullet)_* L f_\bullet^*LL_J
\end{array}
\]
is commutative.
\item[2] The diagram
\[
\begin{array}{ccccc}
LL_J  & \multicolumn{3}{c}{
\mathord -\mkern -6mu
\cleaders \hbox {$\mkern -2mu\mathord -\mkern -2mu$}
\hfill
\hbox to 0pt{\hss$
  \mathord{\mathop{\mkern -1mu\mathord-\mkern -1mu}\limits^{\hbox to
  0pt{\hss\scriptsize $\id$\hss}}}$
   \hss}
\mkern -1mu
\cleaders \hbox {$\mkern -2mu\mathord -\mkern -2mu$}
\hfill\mkern -6mu\mathord \rightarrow
} & LL_J\\
\uparrow\hbox to 0pt{\scriptsize$\via \varepsilon$\hss} & & & &
\uparrow\hbox to 0pt{\scriptsize$\via \varepsilon$\hss}\\
LL_J L(f_\bullet|_J)^*R(f_\bullet|_J)_*&
\specialarrow{\via d}&
L(f_\bullet)^* LL_J R (f_\bullet|_J)_*&
\specialarrow{\via \theta}&
L(f_\bullet)^* R (f_\bullet)_*LL_J
\end{array}
\]
is commutative.
\item[3] Let $h_\bullet:Y'_\bullet\rightarrow Y_\bullet$ be a morphism
in $\Cal P$.
Set $X'_\bullet:=Y'_\bullet\times_{Y_\bullet}X_\bullet$.
Let $h_\bullet':X'_\bullet\rightarrow X_\bullet$ be the second projection,
and $f_\bullet':X'_\bullet\rightarrow Y_\bullet'$ be the first projection.
Then the diagram
\[
\begin{array}{ccccc}
LL_J L(h_\bullet|_J)^* R(f_\bullet|_J)_*
& \specialarrow{d} & L(h_\bullet)^* LL_J R(f_\bullet|_J)_* 
& \specialarrow{\theta} & L(h_\bullet)^* R(f_\bullet)_* LL_J\\
\downarrow\hbox to 0pt{\scriptsize$\theta$\hss} & & & &
\downarrow\hbox to 0pt{\scriptsize$\theta$\hss}\\
LL_J R(f_\bullet'|_J)_* L(h_\bullet'|_J)^*
& \specialarrow{\theta} & R(f_\bullet')_* LL_J L(h_\bullet'|_J)^* &
\specialarrow{d} & R(f_\bullet')_* L(h_\bullet')^* LL_J
\end{array}
\]
is commutative.
\end{description}
\end{Lemma}

\proof Left to the reader.

\section[The right adjoint of the derived direct image functor]{The 
right adjoint of the derived direct image functor
of a morphism of quasi-compact separated schemes}

Let $X$ be a quasi-compact separated scheme.

\begin{Lemma}\label{I-B-N.thm}
The following hold:
\begin{description}
\item[1] There is  a right adjoint $\qco(X):\Mod(X)\rightarrow \Qco(X)$ 
of the canonical inclusion $F_X: \Qco(X)\hookrightarrow \Mod(X)$.
\item[1'] $\qco(X)$ preserves $K$-injective complexes.
$R\qco(X):D(\Mod(X))\rightarrow D(\Qco(X))$ is right adjoint
to $F_X:D(\Qco(X))\rightarrow D(\Mod(X))$.
\item[2] The functor
$F_X:D(\Qco(X))\rightarrow D(\Mod(X))$
is full and faithful, and induces an equivalence
$D(\Qco(X))\rightarrow D_{\Qco(X)}(\Mod(X))$.
\item[3] The unit of adjunction $u:\Id\rightarrow R\qco(X)F_X$ is
an isomorphism, and
$\varepsilon: F_X R\qco(X)\rightarrow \Id$ is an isomorphism 
on $D_{\Qco(X)}(\Mod(X))$.
\end{description}
\end{Lemma}

\proof See \cite{Illusie} for {\bf 1}.
The assertion 
{\bf 1'} follows from {\bf 1} and Lemma~\ref{derived-adjoint.thm}.
For {\bf 2}, see \cite{BN}.
The assertion {\bf 3} follows from {\bf 2}.
\qed

By the lemma, the following follows immediately.

\begin{Corollary}\label{B-N-cor.thm}
Let $X$ be a quasi-compact separated scheme.
Then $\Qco(X)$ has enough injectives and has direct products.
For any $\Bbb F\in K(\Qco(X))$, there is a $K$-injective resolution
$\Bbb F\rightarrow \Bbb I$ such that each term of $\Bbb I$
being an injective object of $\Qco(X)$.
\end{Corollary}

In fact, if $F_X\Bbb F\rightarrow \Bbb J$ is a $K$-injective resolution
such that each term of $\Bbb J$ is injective, then
the corresponding $\Bbb F\rightarrow \qco\Bbb J$ is still a quasi-isomorphism
by the lemma, and this is a desired $K$-injective resolution.

\paragraph
Let $f:X \rightarrow Y$ be a quasi-compact separated morphism of schemes.
Then $f_*^{{\Qco}}:\Qco(X)\rightarrow \Qco(Y)$ is defined,
and we have $F_Y\circ f_*^{{\Qco}}\cong f_*^{{\Mod}}\circ F_X$, where
$F_Y$ and $F_X$ are the forgetful functors.
By the corollary, there is a right derived functor $Rf_*^{{\Qco}}$ of
$f_*^{{\Qco}}$.

\begin{Lemma}\label{qco-mod-direct.thm}
Let $f:X\rightarrow Y$ be a 
morphism of schemes.
Assume that both $X$ and $Y$ are quasi-compact separated.
Then, we have the following.
\begin{description}
\item[1]
The canonical maps
\[
F_Y\circ Rf_*^{{\Qco}}\cong
R(F_Y\circ f_*^{\Qco})
\cong
R(f_*^{\Mod}\circ F_X)
\rightarrow Rf_*^{{\Mod}}\circ
F_X
\]
are all isomorphisms.
\item[2] 
There are a left adjoint $F$ of $Rf_*^{\Qco}$ and an isomorphism
$F_XF\cong Lf^*_{\Mod}F_Y$.
\item[3] There is an isomorphism
\[
Rf_*^{\Qco}R\qco\cong R\qco Rf_*^{\Mod}.
\]
\end{description}
\end{Lemma}

\proof We prove {\bf 1}.
It suffices to show that, if $\Bbb I\in K(\Qco(X))$ is $K$-injective,
then $F_X\Bbb I$ is $f_*^{\Mod}$-acyclic.
By Corollary~\ref{B-N-cor.thm}, we may assume that each term of
$\Bbb I$ is injective.
By \cite[(3.9.3.5)]{Lipman}, it suffices to show that an injective
object $\Cal I$ of $\Qco(X)$ is $f_*^{\Mod}$-acyclic.
Let $g:Z\rightarrow X$ be an affine
faithfully flat morphism such that $Z$ is affine.
Such a morphism exists, as $X$ is quasi-compact separated.
Then it is easy to see that $\Cal I$ is a direct summand of $g_*\Cal J$
for some injective object $\Cal J$ in $\Qco(Z)$.
So it suffices to show that $F_Xg_*\Cal J$ is $f_*^{\Mod}$-acyclic.
Let $F_Z\Cal J\rightarrow \Bbb J$ be an injective resolution.
As $g$ is affine and $\Cal J$ is quasi-coherent, 
we have $R^if_*^{\Mod}F_Z\Cal J=0$ for $i>0$.
Hence $g_*F_Z\Cal J\rightarrow g_*\Bbb J$ is a quasi-isomorphism,
and hence is a $K$-limp resolution of $F_Xg_*\Cal J\cong g_*F_Z\Cal J$.
As $f\circ g$ is also affine, 
$f_*g_*F_Z\Cal J\rightarrow f_*g_*\Bbb J$ is still a quasi-isomorphism,
and this shows that $F_Xg_*\Cal J$ is $f_*$-acyclic.

{\bf 2}
Define $F:D(\Qco(Y))\rightarrow D(\Qco(X))$ by
$F:=\Rqco Lf^*_{\Mod}F_Y$.
As we have $Lf^*_{\Mod}F_Y(D(\Qco(Y)))\subset D_{\Qco(X)}(D(\Mod(X)))$, 
we have
\[
F_X F=F_X R\qco Lf^*_{\Mod}F_Y\specialarrow{\via\varepsilon} Lf^*_{\Mod}F_Y
\]
is an isomorphism by Lemma~\ref{I-B-N.thm}.
Hence
$$\displaylines{
\quad
  \Hom_{D(\Qco(X))}(F\Bbb F,\Bbb G)
  \cong
  \Hom_{D(\Mod(X))}(F_X F\Bbb F,F_X \Bbb G)
  \cong
\hfill\cr
\quad
  \Hom_{D(\Mod(X))}(Lf^*_{\Mod}F_Y\Bbb F,F_X \Bbb G)
  \cong
  \Hom_{D(\Mod(Y))}(F_Y\Bbb F,Rf_*^{\Mod}F_X \Bbb G)
  \cong
\hfill\cr
\quad
  \Hom_{D(\Mod(Y))}(F_Y\Bbb F,F_YRf_*^{\Qco} \Bbb G)
  \cong
  \Hom_{D(\Qco(Y))}(\Bbb F,Rf_*^{\Qco}\Bbb G).
\hfill
}
$$
This shows that $F$ is left adjoint to $Rf_*^{\Qco}$.

{\bf 3} Take the conjugate of {\bf 2}.
\qed

\begin{Rem} By the lemma, we know that there is no difference
between
$Rf_*^{\Qco}$ on $D(\Qco(X))$ and $Rf_*^{\Mod}$ on $D_{\Qco(X)}(\Mod(X))$
when we consider quasi-compact separated schemes.
We are to consider the latter.
\end{Rem}

\paragraph Let $\Cal C$ be an additive category, and $c\in\Cal C$.
We say that $c$ is a compact object, if for any small family of objects
$(t_\lambda)_{\lambda\in\Lambda}$ of $\Cal C$ such that the coproduct
(direct sum) $\bigoplus_{\lambda \in \Lambda}t_\lambda$ exists, 
the canonical map
\[
\bigoplus_\lambda\Hom_{\Cal C}(c,t_\lambda)
\rightarrow
\Hom_{\Cal C}(c,\bigoplus_\lambda t_\lambda)
\]
is an isomorphism.

A triangulated category $\Cal T$ is said to be compactly generated, if
$\Cal T$ has small coproducts, and there is a small set $C$ of 
compact objects of $\Cal T$
such that $\Hom_{\Cal T}(c,t)=0$ for all $c\in C$ implies $t=0$.
The following was proved by A.~Neeman \cite{Neeman2}.

\begin{Theorem}\label{brown-representativity.thm}
Let $\Cal S$ be a compactly generated triangulated category, 
$\Cal T$ any triangulated category, and $F:\Cal S\rightarrow \Cal T$ 
a triangulated functor.
Suppose that $F$ preserves coproducts, that is to say, 
for any small family of objects $(s_\lambda)$ of $\Cal S$,
the canonical maps $F(s_\lambda)\rightarrow F(\bigoplus_\lambda s_\lambda)$
make $F(\bigoplus_\lambda s_\lambda)$ the coproduct of $F(s_\lambda)$.
Then $F$ has a right adjoint $G:\Cal T\rightarrow \Cal S$.
\end{Theorem}

Note that a right adjoint $G$ of a triangulated functor $F$ is triangulated.
Note also that, if both $\Cal S$ and $\Cal T$ have
$t$-structures and $F$ is way-out left, then $G$ is way-out right.

The following was also proved by A.~Neeman \cite{Neeman2}.

\begin{Theorem}\label{c-generation-scheme.thm}
Let $X$ be a quasi-compact separated scheme.
Then $c\in D(\Qco(X))$ is a compact object 
if and only if $c$ is isomorphic to a
perfect complex, where we say that $\Bbb C\in C(\Qco(X))$ is perfect
if $\Bbb C$ is bounded, and each term of $\Bbb C$ is locally free of
finite rank.
Moreover, $D(\Qco(X))$ is compactly generated.
\end{Theorem}

By the theorem, $D_{\Qco(X)}(\Mod(X))$ is also compactly generated.

\begin{Lemma}\label{respects-coproducts.thm}
Let $f:X\rightarrow Y$ be a quasi-compact 
separated morphism of schemes.
Then $Rf_*^{{\Mod}}:D_{\Qco(X)}(\Mod(X))\rightarrow D(\Mod(Y))$ 
preserves coproducts.
\end{Lemma}

\proof See \cite{Neeman2} or \cite{Lipman}.
\qed

\paragraph Let $f:X\rightarrow Y$ be a morphism between quasi-compact
separated schemes.
Then by Theorem~\ref{brown-representativity.thm}, 
Theorem~\ref{c-generation-scheme.thm} and Lemma~\ref{respects-coproducts.thm},
there is a right adjoint
\[
f^\times: D(\Mod(Y))\rightarrow D_{\Qco(X)}(\Mod(X))
\]
of $Rf_*^{{\Mod}}$.

By restriction, $f^\times:D_{\Qco(Y)}(\Mod(Y))\rightarrow D_{\Qco(X)}(\Mod(X))
$ is a right adjoint of $Rf_*^{\Mod}:D_{\Qco(X)}(\Mod(X))\rightarrow
D_{\Qco(Y)}(\Qco(Y))$.
By \cite[(3.6.8.2)]{Lipman}, $(R(?)_*,(?)^\times)$ is an adjoint pair of
$\Delta$-pseudofunctors on the category of quasi-compact separated
schemes.

\section{The Composition of two pseudofunctors}
\label{composition.sec}

\begin{Def}\label{composition-data.def}
We say that
$\Cal C=(\Cal A,\Cal F,\Cal P,\Cal I,\Cal D,\Cal D^+,(?)^{\#},(?)^{\flat},
\zeta)$ is a {\em composition data of contravariant pseudofunctors}
if the following sixteen conditions are satisfied:
\begin{description}
\item[1] $\Cal A$ is a category with fiber products.
\item[2] $\Cal P$ and $\Cal I$ are classes of morphisms of $\Cal A$.
\item[3] Any isomorphism in $\Cal A$ is in $\Cal P\cap \Cal I$.
\item[4] The composite of two morphisms in $\Cal P$ (resp.\ 
$\Cal I$) is again a morphism in $\Cal P$ (resp.\ $\Cal I$).
\item[5] A base change of a morphism in $\Cal P$ 
is again a morphism in $\Cal P$.
\item[6] If $f$ is a morphism in $\Cal A$, $g\in\Cal P$,
and $g\circ f$ is defined and a morphism in $\Cal P$,
then $f\in\Cal P$.
\item[7] Any $f\in\Mor(\Cal A)$ admits a factorization
$f=pi$ such that $p\in\Cal P$ and $i\in\Cal I$.
\end{description}
Before we state the remaining conditions, we give some definitions
for convenience.
\begin{description}
\item[i] Let $\Cal C$ be a class of morphisms in $\Cal A$ containing all
identity maps and closed under compositions.
We define $\Cal A_{\Cal C}$ by $\ob(\Cal A_{\Cal C}):=
\ob(\Cal A)$ and $\Mor(\Cal A_{\Cal C}):=\Cal C$.
In particular, the subcategories $\Cal A_{\Cal P}$ and $\Cal A_{\Cal I}$ of
$\Cal A$ are defined.
\item[ii] We call a commutative diagram of the form $p\circ i=i'\circ p'$
with $p,p'\in\Cal P$ and $i,i'\in\Cal I$ a {\em pi-square}.
We denote the set of all pi-squares by $\Pi$.
\end{description}
\begin{description}
\item[8] $\Cal D=(\Cal D(X))_{X\in\ob(\Cal A)}$ is a family of categories.
\item[9] $(?)^{\#}$ is a contravariant pseudofunctor on $\Cal A_{\Cal P}$,
$(?)^{\flat}$ is a contravariant pseudofunctor on $\Cal A_{\Cal I}$,
and we have $X^{\#}=X^{\flat}=\Cal D(X)$ for each $X\in\ob(\Cal A)$.
\item[10] $\zeta=(\zeta(\sigma))_{\sigma=(pi=jq)\in\Pi}$ is a family
of natural transformations
\[
\zeta(\sigma):i^\flat p^{\#}\rightarrow q^{\#}j^\flat.
\]
\item[11] If
\[
\begin{array}{ccccc}
U_1 & \specialarrow{j_1} & V_1 & \specialarrow{i_1} & X_1 \\
\downarrow\hbox to 0pt{\scriptsize$p_U$\hss} & \sigma' &
\downarrow\hbox to 0pt{\scriptsize$p_V$\hss} & \sigma  &
\downarrow\hbox to 0pt{\scriptsize$p_X$\hss}\\
U & \specialarrow{j} & V & \specialarrow{i} & X
\end{array}
\]
is a commutative diagram in $\Cal A$ such that $\sigma,\sigma'\in\Pi$,
then the composite map
\[
(i_1j_1)^\flat p_X^{\#}\cong j_1^\flat i_1^\flat p_X^{\#}
\specialarrow{\via\zeta(\sigma)}j_1^\flat p_V^{\#}i^\flat
\specialarrow{\via\zeta(\sigma')}p_U^{\#}j^\flat i^\flat
\cong p_U^{\#}(ij)^\flat
\]
agrees with $\zeta(\sigma'\sigma)$,
where $\sigma'\sigma$ is the square $p_X(i_1j_1)=(ij)p_U$.
\item[12] If
\[
\begin{array}{ccc}
U_2 & \specialarrow{i_2} & X_2\\
\downarrow\hbox to 0pt{\scriptsize$q_U$\hss} & \sigma'  &
\downarrow\hbox to 0pt{\scriptsize$q_X$\hss}\\
U_1 & \specialarrow{i_1} & X_1\\
\downarrow\hbox to 0pt{\scriptsize$p_U$\hss} & \sigma  &
\downarrow\hbox to 0pt{\scriptsize$p_X$\hss}\\
U & \specialarrow{i} & X
\end{array}
\]
is a commutative diagram in $\Cal A$ such that $\sigma,\sigma'\in\Pi$,
then the composite
\[
i_2^\flat(p_Xq_X)^{\#}\cong
i_2^\flat q_X^{\#}p_X^{\#}
\specialarrow{\via\zeta(\sigma')}
q_U^{\#}i_1^\flat p_X^{\#}
\specialarrow{\via\zeta(\sigma)}
q_U^{\#}p_U^{\#}i^\flat
\cong 
(p_Uq_U)^{\#}i^\flat
\]
agrees with $\zeta(((p_Xq_X)i_2=i(p_Uq_U)))$.
\item[13] $\Cal F$ is a subcategory of $\Cal A$, and 
any isomorphism in $\Cal A$ 
between objects of $\Cal F$ is in $\Mor(\Cal F)$.
\item[14] $\Cal D^+=(\Cal D^+(X))_{X\in\ob(\Cal F)}$ is a family of categories
such that $\Cal D^+(X)$ is a full subcategory of $\Cal D(X)$ for each
$X\in\ob(\Cal F)$.
\item[15] If $f:X\rightarrow Y$ is a morphism in $\Cal F$, 
$f=p\circ i$, $p\in\Cal P$ and $i\in\Cal I$, then
we have $i^\flat p^!(\Cal D^+(Y))\subset \Cal D^+(X)$.
\item[16] If
\[
\begin{array}{ccccccc}
V & \specialarrow j & U_1 & \specialarrow{i_1} & X_1 &  &  \\
  &                 &
\downarrow\hbox to 0pt{\scriptsize$p_U$\hss} & \sigma  &
\downarrow\hbox to 0pt{\scriptsize$p_X$\hss} &   &  \\
  &                 & U & \specialarrow{i}&  X & \specialarrow{q} & Y
\end{array}
\]
is a diagram in $\Cal A$ such that $\sigma\in \Pi$, $V,U,Y\in\ob(\Cal F)$, 
$p_U j\in \Mor(\Cal F)$ and $qi\in\Mor(\Cal F)$, then
\[
j^\flat \zeta(\sigma)q^{\#}:j^\flat i_1^\flat p_X^{\#}q^{\#}
\rightarrow j^\flat p_U^{\#} i^\flat q^{\#}
\]
is an isomorphism between functors from $\Cal D^+(Y)$ to $\Cal D^+(V)$.
\end{description}
\end{Def}

\paragraph Let
$\Cal C=(\Cal A,\Cal F,\Cal P,\Cal I,\Cal D,\Cal D^+,(?)^{\#},(?)^{\flat},
\zeta)$ be a composition data of contravariant pseudofunctors.
We call a commutative diagram of the form $f=pi$ with $p\in\Cal P$, 
$i\in\Cal I$ and $f\in\Mor(\Cal A)$ a {\em compactification}.
We call a commutative diagram of the form $pi=qj$ with $p,q\in\Cal P$,
$i,j\in\Cal I$ and $pi=qj\in\Mor(\Cal F)$ an {\em independence diagram}.

\begin{Lemma}\label{independence-isom.thm}
Let
\[
\begin{array}{ccc}
U & \specialarrow{i_1} & X_1\\
\downarrow\hbox to 0pt{\scriptsize$i$\hss}& \tau &
\downarrow\hbox to 0pt{\scriptsize$p_1$\hss}\\
X & \specialarrow{p} & Y
\end{array}
\]
be an independence diagram.
Then the following hold:
\begin{description}
\item[1] There is a diagram of the form
\[
\begin{array}{ccccc}
U & \specialarrow j & Z & \specialarrow{q_1} & X_1\\
& & \downarrow\hbox to 0pt{\scriptsize$q$\hss}&  &
\downarrow\hbox to 0pt{\scriptsize$p_1$\hss}\\
& & X & \specialarrow{p} & Y
\end{array}
\]
such that $qj=i$, $q_1j=i_1$, $pq=p_1q_1$, $q,q_1\in\Cal P$ and $j\in\Cal I$.
\item[2] $\zeta(qj=i1_U)p^{\#}:j^\flat q^{\#}p^{\#}\rightarrow
i^\flat p^{\#}$ is an isomorphism 
between functors from $\Cal D^+(Y)$
to $\Cal D^+(U)$.
\item[3] $\zeta(q_1j_1=i_11_U)p_1^{\#}$ is also an isomorphism
between functors from $\Cal D^+(Y)$ to $\Cal D^+(U)$.
\item[4] The composite isomorphism
\[
\rho(\tau):
i^\flat p^{\#}\specialarrow{(\zeta(qj=i1_U)p^{\#})^{-1}}j^\flat q^{\#}p^{\#}
\cong
j^\flat q_1^{\#}p_1^{\#}
\specialarrow{\zeta(q_1j_1=i_11_U)p_1^{\#}}
i_1^\flat p_1^{\#}
\]
\(between functors defined over $\Cal D^+(Y)$, not over $\Cal D(Y)$\)
depends only on $\tau$.
\item[5] If $\tau'=(p_1i_1=p_2i_2)$ is an independence diagram, then we have
\[
\rho(\tau')\circ \rho(\tau)=\rho(pi=p_2i_2).
\]
\end{description}
\end{Lemma}

The proof is left to the reader.
We call $\rho(\tau)$ the independence isomorphism of $\tau$.

\paragraph Any $f\in\Mor(\Cal A)$ has a compactification by assumption.
We fix a family of compactifications $\Cal T:=(\tau(f):(f=p(f)\circ i(f))
)_{f\in\Mor(\Cal A)}$
so that $p(1_X)=i(1_X)=1_X$ for $X\in\ob(\Cal A)$.
We denote the target of $i(f)$, which agrees with the source of $p(f)$
by $M(f)$.

For $X\in\ob(\Cal F)$, we define $X^!:=\Cal D^+(X)$.
For a morphism $f:X\rightarrow Y$ in $\Cal F$, 
we define $f^!:=i(f)^\flat p(f)^{\#}$,
which is a functor from $X^!$ to $Y^!$ by assumption.

Let $f:X\rightarrow Y$ and $g:Y\rightarrow Z$ be morphisms in $\Cal F$.
Set $q:=p(i(g)p(f))$ and $j:=i(i(g)p(f))$.
Then by {\bf 16} in Definition~\ref{composition-data.def}
and Lemma~\ref{independence-isom.thm},
the composite map
\begin{multline*}
(gf)^!=i(gf)^\flat p(gf)^{\#}
\specialarrow{\rho(p(gf)i(gf)=(p(g)q)(ji(f)))}
(j\circ i(f))^\flat (p(g)\circ q)^{\#}\\
\cong
i(f)^\flat j^\flat q^{\#}p(g)^{\#}
\specialarrow{i(f)^\flat\zeta(qj=i(g)p(f))p(g)^{\#}}
i(f)^\flat p(f)^{\#}i(g)^\flat p(g)^{\#}=f^!g^!
\end{multline*}
is an isomorphism.
We define $d_{f,g}:f^!g^!\rightarrow (gf)^!$ to be the inverse of 
this composite.

\begin{Lemma} The definition of $d_{f,g}$ is independent of choice of
$q$ and $j$ above.
\end{Lemma}

The proof of the lemma is left to the reader.

\begin{Proposition} Let the notation be as above.
\begin{description}
\item[1] $(?)^!$ and $(d_{f,g})$ form a contravariant pseudofunctor
on $\Cal F$.
\item[2] For $j\in\Cal I\cap \Mor(\Cal F)$, 
define $\rho:j^!\rightarrow j^\flat$ to be
$\rho(p(j)i(j)=j\id)$.
Then $\rho:(?)^!\rightarrow (?)^\flat$ is an isomorphism of
pseudofunctors on $\Cal A_{\Cal I}\cap \Cal F$.
\item[3] For $q\in\Cal P\cap \Mor(\Cal F)$, 
define $\rho:q^!\rightarrow q^{\#}$ to be
$\rho(p(q)i(q)=\id q)$.
Then $\rho:(?)^!\rightarrow (?)^{\#}$ is an isomorphism of 
pseudofunctors on $\Cal A_{\Cal P}\cap \Cal F$.
\item[4] Let us take another family of compactifications $(f=
p_1(f)i_1(f))_{f\in\Mor(\Cal F)}$, and let $(?)^\star$ be the resulting
pseudofunctor defined by $f^\star=i_1(f)^\flat p_1(f)^{\#}$.
Then $\rho: f^!\rightarrow f^\star$ induces an isomorphism of 
pseudofunctors $(?)^!\cong (?)^\star$.
\end{description}
\end{Proposition}

The proof is left to the reader.
We call $(?)^!$ the composition of $(?)^{\#}$ and $(?)^\flat$.
The composition is uniquely defined up to isomorphisms of pseudofunctors
on $\Cal F$.
The discussion above has an obvious triangulated version.
Composition data of contravariant triangulated pseudofunctors are
defined appropriately, and the composition of two triangulated 
pseudofunctors is obtained as a triangulated pseudofunctor.

\paragraph\label{scheme-composition.par}
Let $S$ be a separated noetherian scheme.
Let $\Cal A$ be the category whose objects are separated noetherian 
$S$-schemes and morphisms are morphisms of finite type.
Set $\Cal F=\Cal A$.
Let $\Cal I$ be the class of open immersions.
Let $\Cal P$ be the class of proper morphisms.
Set $\Cal D(X)=\Cal D^+(X)=D^+_{\Qco(X)}(\Mod(X))$ for $X\in\Cal A$.
Let $(?)^\flat:=(?)^*$, the derived inverse image pseudofunctor
for morphisms in $\Cal I$,
where $X^\flat:=\Cal D(X)$.
Let $(?)^{\#}:=(?)^\times$, the twisted inverse pseudofunctor
for morphisms in $\Cal P$,
where $X^{\#}:=\Cal D(X)$ again.
Note that
the left adjoint $R(?)_*$ is way-out left for morphisms in $\Cal P$
so that $(?)^\times$ is way-out right, and $(?)^{\#} $ is well-defined.
The conditions {\bf 1--9} in Definition~\ref{composition-data.def}
hold.
Note that {\bf 7} is nothing but Nagata's compactification theorem
\cite{Nagata}.
The conditions {\bf 13--15} are trivial.

Let $\sigma_0:pi=jq$ be a pi-diagram, which is also a fiber square.
Then the canonical map
\[
\theta: j^*(Rp_*)\rightarrow (Rq_*)i^*
\]
is an isomorphism of triangulated functors, see \cite[(3.9.5)]{Lipman}.
Hence, taking the inverse of the conjugate, we have an isomorphism
\begin{equation}\label{verdier.eq}
\xi=\xi(\sigma_0):(Ri_*)q^\times\cong p^\times(Rj_*).
\end{equation}
So we have a morphism of triangulated functors
\[
\zeta_0(\sigma_0):i^*p^\times\specialarrow{\via u}
i^*p^\times(Rj_*)j^*\specialarrow{\via\xi^{-1}}
i^*(Ri_*)q^\times j^*
\specialarrow{\via\varepsilon}
q^\times j^*,
\]
which is an {\em isomorphism}, see \cite{Verdier3}.
The statements corresponding to {\bf 10--12} only for {\em 
fiber square pi-diagrams},
is readily proved.

In particular, for a closed open immersion $\eta:U\rightarrow X$, 
we have an isomorphism
\[
v(\eta):=\zeta_0(\pi 1_U=1_U\pi)^{-1}: \eta^*\rightarrow \eta^\times.
\]

Let $\sigma=(pi=qj)$ be an arbitrary pi-diagram.
Let $j_1$ be the base change of $j$ by $p$, and let
$p_1$ be the base change of $p$ by $j$.
Let $\eta$ be the unique morphism such that $q=p_1\eta$ and $i=j_1\eta$.
Note that $\eta$ is a closed open immersion.
Define $\zeta(\sigma)$ to be the composite isomorphism
\[
i^*p^\times\cong \eta^*j_1^*p^\times\specialarrow{\via\zeta_0}
\eta^*p_1^\times j^*\specialarrow{\via v(\eta)}
\eta^\times p_1^\times j^*\cong q^\times j^*.
\]

Now the proof of conditions {\bf 10--12} consists in diagram chasing
arguments, while {\bf 16} is trivial, since $\zeta(\sigma)$ is
always an isomorphism.
Thus the twisted inverse triangulated pseudofunctor $(?)^!$ on
$\Cal A$ is defined to be the composite of $(?)^\times$ and $(?)^*$.

\section[Morphism of diagrams]{The right adjoint of the derived direct image
functor of a morphism of diagrams}

Let $I$ be a small category, $S$ a scheme, and 
$X_\bullet\in \Cal P(I,\Sch/S)$.

\begin{Lemma} Assume that $X_\bullet$ is quasi-compact separated.
That is, $X_i$ is quasi-compact separated for each $i\in I$.
Let $C_i$ be a small set of compact generators of $D(\Qco(X_i))$.
Then 
\[
C:=\{LL_ic\setbar i\in I,\;c\in C_i\}
\}
\]
is a small 
set of compact generators of $D_{\LQco(X_\bullet)}(\Mod(X_\bullet))$.
In particular, the category 
$D_{\LQco(X_\bullet)}(\Mod(X_\bullet))$ is compactly generated.
\end{Lemma}

\proof Let $t\in D_{\LQco(X_\bullet)}(\Mod(X_\bullet))$ and assume that
\[
\Hom_{D(\Mod(X_\bullet))}(LL_ic,t)\cong
\Hom_{D(\Mod(X_i))}(c,t_i)=0
\]
for any $i\in I$ and any $c\in C_i$.
Then, $t_i=0$ for all $i$.
This shows $t=0$.
It is easy to see that $LL_ic$ is compact, and $C$ is small.
So $C$ is a small set of compact generators.
As $D_{\LQco(X_\bullet)}(\Mod(X_\bullet))$ has coproducts,
it is compactly generated.
\qed

\begin{Lemma} Let $f_\bullet:X_\bullet\rightarrow Y_\bullet$ be a 
quasi-compact separated morphism in $\Cal P(I,\Sch/S)$.
Then
\[
R(f_\bullet)_*:D_{\LQco(X_\bullet)}(\Mod(X_\bullet))\rightarrow
D_{\LQco(Y_\bullet)}(\Mod(Y_\bullet))
\]
preserves coproducts.
\end{Lemma}

\proof Let $(t_\lambda)$ be a small family of objects in 
$D_{\LQco(X_\bullet)}(\Mod(X_\bullet))$, and consider the canonical
map
\[
\bigoplus_\lambda R(f_\bullet)_* t_\lambda
\rightarrow R(f_\bullet)_*(\bigoplus_\lambda t_\lambda).
\]
For each $i\in I$, apply $(?)_i$ to the map.
As $(?)_i$ obviously preserves 
coproducts and we have a canonical isomorphism
\[
(?)_i R(f_\bullet)_*\cong R(f_i)_* (?)_i,
\]
the result is  the canonical map
\[
\bigoplus_\lambda R(f_i)_* (t_\lambda)_i
\rightarrow
R(f_i)_*(\bigoplus_\lambda (t_\lambda)_i),
\]
which is an isomorphism by
Lemma~\ref{respects-coproducts.thm}.
Hence, $R(f_\bullet)_*$ preserves coproducts.
\qed

By Theorem~\ref{brown-representativity.thm}, we have:

\begin{Theorem}
Let $I$ be a small category, $S$ a scheme, and $f_\bullet:X_\bullet\rightarrow
Y_\bullet$ a morphism in $\Cal P(I,\Sch/S)$.
If $X_\bullet$ and $Y_\bullet$ are quasi-compact separated, then
\[
R(f_\bullet)_*:D_{\LQco(X_\bullet)}(\Mod(X_\bullet))\rightarrow
D_{\LQco(Y_\bullet)}(\Mod(Y_\bullet))
\]
has a right adjoint $f_\bullet^\times$.
\end{Theorem}

\begin{Lemma}\label{theta-g-f-isom.thm}
Let $S$ be a scheme, $I$ a small category, and
$f_\bullet:X_\bullet\rightarrow Y_\bullet$ and
$g_\bullet:Y_\bullet'\rightarrow Y_\bullet$ be
morphisms in $\Cal P(I,\Sch/S)$.
Set $X_\bullet':=Y_\bullet'\times_{Y_\bullet} X_\bullet$.
Let $f_\bullet': X'_\bullet\rightarrow Y_\bullet'$ be the first
projection, and $g_\bullet': X'_\bullet\rightarrow X_\bullet$ be the
second projection.
Assume that $f_\bullet$ is quasi-compact separated, and $g_\bullet$ is
flat.
Then the canonical map
\[
\theta(g_\bullet,f_\bullet)
: (g_\bullet)^* R(f_\bullet)_*\rightarrow R(f_\bullet')_* (g_\bullet')^*
\]
is an isomorphism of functors from $D_{\LQco(X_\bullet)}(\Mod(X_\bullet))$ 
to $D_{\LQco(Y_\bullet')}(\Mod(Y_\bullet'))$.
\end{Lemma}

\proof It suffices to show that for each $i\in I$, 
\[
(?)_i\theta: (?)_i(g_\bullet)^* R(f_\bullet)_*
            \rightarrow 
             (?)_i R(f_\bullet')_* (g_\bullet')^*
\]
is an isomorphism.
By Lemma~\ref{theta-f-J.thm}, {\bf 2, 3}, 
it is easy to verify that the diagram
\begin{equation}\label{theta-f-J-restrict.eq}
\begin{array}{ccccc}
   g_i^* R(f_i)_* (?)_i & \specialarrow{c^{-1}} & g_i^*(?)_i R(f_\bullet)_*
  & \specialarrow{\theta}& (?)_i g_\bullet^* R(f_\bullet)_*\\
\downarrow\hbox to 0pt{\scriptsize $\theta(g_i,f_i)$\hss} & & & &
\downarrow\hbox to 0pt{\scriptsize $(?)_i\theta(g_\bullet,f_\bullet)$\hss}\\
R(f_i')_* (g_i')^* (?)_i & \specialarrow{\theta} & R(f_i')_*(?)_i(g_\bullet')^*
& \specialarrow{c^{-1}} & (?)_i R(f_\bullet')_* (g_\bullet')^*\\
\end{array}
\end{equation}
is commutative.
As the horizontal maps and $\theta(g_i,f_i)$ are isomorphisms
by \cite[(3.9.5)]{Lipman}, $(?)_i\theta(g_\bullet,f_\bullet)$ is
also an isomorphism.
\qed

\section{B\"okstedt--Neeman resolutions and hyperExt sheaves}

\paragraph Let $\Cal T$ be a triangulated category with small 
direct products.
Note that a direct product of distinguished triangles is again a
distinguished triangle.

Let
\begin{equation}\label{t_i.eq}
\cdots\rightarrow t_3\specialarrow{s_3}t_2\specialarrow{s_2}t_1
\end{equation}
be a sequence of morphisms in $\Cal T$.
We define $d:\prod_{i\geq 1}t_i\rightarrow \prod_{i\geq 1}t_i$ by
$p_i\circ d=p_i-s_{i+1}\circ p_{i+1}$, where $p_i:\prod_i t_i\rightarrow t_i$ 
is the projection.
Consider a distinguished triangle of the form
\[
M\specialarrow{m}\prod_{i\geq 1}t_i\specialarrow p \prod_{i\geq 1}t_i
\specialarrow{q}\Sigma M,
\]
where $\Sigma$ denotes the suspension.

We call $M$, which is determined uniquely up to isomorphisms, the
homotopy limit of (\ref{t_i.eq}) and denote it by $\holim t_i$.

\paragraph Dually, homotopy colimit is defined and denoted by $\hocolim$,
if $\Cal T$ has small coproducts.

\paragraph Let $\Cal A$ be an abelian category which satisfies (AB3*).
Let $(\Bbb F_\lambda)_{\lambda\in\Lambda}$ be a small family of objects
in $K(\Cal A)$.
Then for any $\Bbb G\in K(\Cal A)$, we have that
\begin{multline*}
\Hom_{K(\Cal A)}(\Bbb G,\prod_\lambda\Bbb F_\lambda)
=H^0(\Hom_{\Cal A}^\bullet(\Bbb G,\prod_\lambda \Bbb F_\lambda))
\cong
H^0(\prod_\lambda\Hom_{\Cal A}^\bullet(\Bbb G,\Bbb F_\lambda))\\
\cong
\prod_\lambda H^0(\Hom_{\Cal A}^\bullet(\Bbb G,\Bbb F_\lambda))
=\prod_\lambda\Hom_{K(\Cal A)}(\Bbb G,\Bbb F_\lambda).
\end{multline*}
That is, the direct product $\prod_\lambda\Bbb F_\lambda$ in $C(\Cal A)$ is
also a direct product in $K(\Cal A)$.

\paragraph Let $\Cal A$ be a Grothendieck abelian category which satisfies
(AB3*),
and $(t_\lambda)$ a small family of objects of $D(\Cal A)$.
Let $(\Bbb F_\lambda)$ be a family of $K$-injective objects of $K(\Cal A)$ 
such that $\Bbb F_\lambda$ represents $t_\lambda$ for each $\lambda$.
Then $Q(\prod_{\lambda}\Bbb F_\lambda)$ is a direct product of $t_\lambda$ 
in $D(\Cal A)$.
Hence $D(\Cal A)$ has small products.

\begin{Lemma}\label{bn-holim.thm}
Let $I$ be a small category, $S$ be a scheme, and 
let $X_\bullet\in\Cal P(I,\Sch/S)$.
Let $\Bbb F$ be an object of $C(\Mod(X_\bullet))$.
Assume that $\Bbb F$ has locally quasi-coherent cohomology groups.
Then the following hold.
\begin{description}
\item[1] Let $\frak I$ denote the
full subcategory of $C(\Mod(X_\bullet))$ consisting of bounded below 
complexes of injective objects of $\Mod(X_\bullet)$ 
with locally 
quasi-coherent cohomology groups.
There is an
$\frak I$-special inverse
system $(I_n)_{n\in\Bbb N}$
with the index set $\Bbb N$ and an inverse system of
chain maps $(f_n: \tau^{\geq -n}\Bbb F\rightarrow I_n)$ such that
\begin{description}
\item[i] $f_n$ is a quasi-isomorphism for any $n\in\Bbb N$.
\item[ii] $(I_n)_i=0$ for $i\leq -n-1$.
\end{description}
\item[2] If $(I_n)$ and $(f_n)$ are as in {\bf 1}, then the following hold.
\begin{description}
\item[i] For each $i\in \Bbb Z$,
the canonical map $H^i(I_n)\rightarrow H^i(\projlim I_n)$ 
is an isomorphism for $n\geq\max(1,-i)$, 
where the projective limit is taken in the
category $C(\Mod(X_\bullet))$,
and 
$H^i(?)$ denotes the $i$th cohomology sheaf of a complex of sheaves.
\item[ii] 
$\projlim f_n: \Bbb F\rightarrow \projlim I_n$ is a quasi-isomorphism.
\item[iii] The projective limit $\projlim I_n$, viewed as an object of
$K(\Mod(X))$, is the homotopy limit of $(I_n)$.
\item[iv] $\projlim I_n$ is $K$-injective.
\end{description}
\item[3] $Q\Bbb F$ is the projective limit of the projective
system $(\tau^{\geq -n}Q\Bbb F)$ in the category $D(\Mod(X_\bullet))$.

\end{description}
\end{Lemma}

\proof The assertion {\bf 1} is \cite[(3.7)]{Spaltenstein}.

We prove {\bf 2}, {\bf i}.
Let $j\in\ob(I)$ and $U$ an affine open subset of $X$.
Then for any $n\geq 1$, $I_n^i$ and $H^i(I_n)$ are 
$\Gamma((j,U),?)$-acyclic for each $i\in\Bbb Z$.
As $I_n$ is bounded below, each $Z^i(I_n)$ and $B^i(I_n)$ are also
$\Gamma((j,U),?)$-acyclic, and the sequence
\begin{equation}\label{zib.eq}
0\rightarrow \Gamma((j,U),Z^i(I_n))\rightarrow \Gamma((j,U),I_n^i)\rightarrow
\Gamma((j,U),B^{i+1}(I_n))\rightarrow 0
\end{equation}
and 
\begin{equation}\label{bzh.eq}
0\rightarrow \Gamma((j,U),B^i(I_n))\rightarrow \Gamma((j,U),Z^i(I_n))
\rightarrow \Gamma((j,U),H^i(I_n))\rightarrow 0
\end{equation}
are exact for each $i$, as can be seen easily,
where $B^i$ and $Z^i$ respectively denote the $i$th coboundary and the
cocycle sheaves.

In particular, the inverse system $(\Gamma((j,U),B^i(I_n)))$ is a
Mittag-Leffler inverse system of abelian groups by (\ref{zib.eq}),
since $(\Gamma((j,U),I^i_n))$ is.
On the other hand, as we have $H^i(I_n)\cong H^i(\Bbb F)$ for $n\geq \max(1,
-i)$, the inverse system $(\Gamma((j,U),H^i(I_n)))$ stabilizes,
and hence we have $(\Gamma((j,U),Z^i(I_n)))$ is also Mittag-Leffler.

Passing through the projective limit, 
\[
0\rightarrow\Gamma((j,U),Z^i(\projlim I_n))
\rightarrow \Gamma((j,U),\projlim I_n)
\rightarrow \Gamma((j,U),\projlim B^{i+1}(I_n))\rightarrow 0
\]
is exact.
Hence, the canonical map $B^i(\projlim I_n)\rightarrow \projlim B^i(I_n)$
is an isomorphism, since $(j,U)$ with $U$ an affine open subset of $X_j$
generates the topology of $\Zar(X_\bullet)$.

Taking the projective limit of (\ref{bzh.eq}), we have
\[
0\rightarrow \Gamma((j,U),B^i(\projlim I_n))\rightarrow
\Gamma((j,U),Z^i(\projlim I_n))
\rightarrow \Gamma((j,U),\projlim H^i(I_n))\rightarrow
0
\]
is an exact sequence for any affine open subset $U$ of $X$.

Hence, the canonical maps
\[
\Gamma((j,U),H^i(I_n))\cong \Gamma((j,U),\projlim H^i(I_n))\leftarrow 
\Gamma((j,U),H^i(\projlim I_n))
\]
are all isomorphisms for $n\geq \max(1,-i)$, and we have
$H^i(I_n)\cong H^i(\projlim I_n)$ for $n\geq \max(1,-i)$.

The assertion {\bf ii} is now trivial.

The assertion {\bf iii} is now a consequence of
\cite[(2.3.1)]{BN} (one can work at the presheaf level where we have
the (AB4*) property).
The assertion {\bf iv} is now obvious.
\qed

Let $I$ be a small category, $S$ a scheme, and
$X_\bullet\in \Cal P(I,\Sch/S)$.

\begin{Lemma}\label{H_J-isom.thm}
Assume that $X_\bullet$ has flat arrows.
Let $J$ be a subcategory of $I$, and
let $\Bbb F\in D_{\EM(X_\bullet)}(\Mod(X_\bullet))$ and
$\Bbb G\in D(\Mod(X_\bullet))$.
Then the canonical map
\[
H_J: (?)_J R\uHom_{\Mod(X_\bullet)}^\bullet(\Bbb F,\Bbb G)
        \rightarrow
     R\uHom_{\Mod(X_J)}^\bullet(\Bbb F_J,\Bbb G_J)
\]
is an isomorphism.
\end{Lemma}

\proof First, we prove the lemma for the case where $\Bbb F$ has
bounded above cohomology groups.

By Lemma~\ref{composition-H.thm}, we may assume that $J=i$ for
an object $i$ of $I$.
By the way-out lemma \cite[Proposition~I.7.1]{Hartshorne2}, we may 
assume that $\Bbb F$ is a single equivariant $\Cal O_{X_\bullet}$-module,
and $\Bbb G$ is a single injective object of $\Mod(X_\bullet)$.

As $X_\bullet$ has flat arrows, the restriction $(?)_i$ has an exact
left adjoint $L_i$.
It follows that $\Bbb G_i$ is an injective $\Cal O_{X_i}$-module.
The assertion is now clear by Lemma~\ref{uhom2.thm}.

Now consider the general case.
As the functors in question on $\Bbb F$ changes coproducts to products
and a direct sum of complexes with equivariant cohomology has
equivariant cohomology, the map in question is an isomorphism if
$\Bbb F$ is a direct sum of objects in $D^-_{\EM(X_\bullet)}(X_\bullet)$.
In particular, the map is an isomorphism if $\Bbb F$ is
a homotopy colimit of objects of
$D^-_{\EM(X_\bullet)}(X_\bullet)$.
As any object $\Bbb F$ of $D_{\EM(X_\bullet)}(X_\bullet)$ is 
the homotopy colimit of $(\tau^{\leq n}\Bbb F)$,  we are done.
\qed

\begin{Lemma}\label{H_J-isom2.thm}
Assume that $X_\bullet$ has flat arrows.
Let $J$ be a subcategory of $I$, and
let $\Bbb F\in D_{\EM(X_\bullet)}(\Mod(X_\bullet))$ and
$\Bbb G\in D_{\LQco(X_\bullet)}(\Mod(X_\bullet))$.
Then the canonical map
\[
H_J: (?)_J R\uHom_{\Mod(X_\bullet)}^\bullet(\Bbb F,\Bbb G)
        \rightarrow
     R\uHom_{\Mod(X_J)}^\bullet(\Bbb F_J,\Bbb G_J)
\]
is an isomorphism.
\end{Lemma}

\proof By Lemma~\ref{bn-holim.thm}, $\Bbb G$ is the homotopy limit of
objects of $D^+(\Mod(X_\bullet))$.
The functors in question on $\Bbb G$ preserve direct products and
triangles.
In particular, they preserve homotopy limits.
The assertion follows immediately by Lemma~\ref{H_J-isom.thm}.
\qed

\begin{Lemma}\label{ext-equivariant.thm}
Let $I$ be a small category,
$S$ a scheme, and $X_\bullet
\in\Cal P(I,\Sch/S)$.
Assume that $X_\bullet$ has flat arrows and is locally noetherian.
Let $\Bbb F\in D^-_{\Coh(X_\bullet)}(\Mod(X_\bullet))$
and $\Bbb G\in D^+_{\LQco(X_\bullet)}(\Mod(X_\bullet))$.
Then $\uExt^i_{\Cal O_{X_\bullet}}(\Bbb F,\Bbb G)$ is locally quasi-coherent
for $i\in\Bbb Z$.
If moreover $\Bbb G$ has quasi-coherent cohomology groups,
then then
$\uExt^i_{\Cal O_{X_\bullet}}(\Bbb F,\Bbb G)$ is quasi-coherent
for $i\in\Bbb Z$.
\end{Lemma}

\proof 
We prove the assertion for the local quasi-coherence.
By Lemma~\ref{H_J-isom.thm}, we may assume that $X_\bullet$ is a
single scheme.
This case is \cite[Proposition~II.3.3]{Hartshorne2}.

We prove the assertion for the quasi-coherence, assuming
that $\Bbb G$ has quasi-coherent cohomology groups.
By \cite[Proposition~I.7.3]{Hartshorne2}, we may assume that
$\Bbb F$ is a single coherent sheaf, and $\Bbb G$ is an injective resolution
of a single quasi-coherent sheaf.

As $X_\bullet$ has flat arrows and the restrictions are exact, 
it suffices to show that
\[
\alpha_\phi:X_\phi^*(?)_i\uHom^\bullet_{\Mod(X_\bullet)}(\Bbb F,\Bbb G)
\rightarrow (?)_j\uHom^\bullet_{\Mod(X_\bullet)}(\Bbb F,\Bbb G)
\]
is a quasi-isomorphism for any morphism $\phi:i\rightarrow j$ in $I$.

As $X_\phi$ is flat, $\alpha_\phi:X_\phi^*\Bbb F_i\rightarrow \Bbb F_j$ 
and $\alpha_\phi:X_\phi^*\Bbb G_i\rightarrow \Bbb G_j$ are quasi-isomorphisms.
In particular, the latter is a $K$-injective resolution.

By the derived version of (\ref{hom-structure.par}), 
it suffices to show that
\[
P: X_\phi^* R\uHom_{\Cal O_{X_i}}^\bullet(\Bbb F_i,\Bbb G_i)
\rightarrow
R\uHom_{\Cal O_{X_j}}^\bullet(X_\phi^*\Bbb F_i,X_\phi^*\Bbb G_i)
\]
is an isomorphism.
This is \cite[Proposition~II.5.8]{Hartshorne2}.
\qed

\section{Commutativity of twisted inverse with restrictions}

\paragraph\label{derived-theta.par}
Let $S$ be a scheme, $I$ a small category, and
$f_\bullet:X_\bullet\rightarrow Y_\bullet$ be a morphism in $\Cal P(I,\Sch/S)$.
Let $J$ be an admissible subcategory of $I$.
Assume that $f_\bullet$  is quasi-compact and separated.
Then there is a natural map
\begin{equation}\label{derived-theta.eq}
\theta(J,f_\bullet): 
LL_J\circ R(f_\bullet|_J)_*\rightarrow R(f_\bullet)_*\circ LL_J
\end{equation}
between functors from $D_{\LQco(X_\bullet|_J)}(\Mod(X_\bullet|_J))$
to $D_{\LQco(Y_\bullet)}(\Mod(Y_\bullet))$,
see \cite[(3.7.2)]{Lipman}.

\paragraph\label{xi-def.par}
Let $S$, $I$ and $f_\bullet$ be as in (\ref{derived-theta.par}).
We assume that $X_\bullet$ and $Y_\bullet$ are 
quasi-compact $\Bbb Z$-separated,
so that the right adjoint functor
\[
f_\bullet^\times:D_{\LQco(Y_\bullet)}(\Mod(Y_\bullet))
\rightarrow D_{\LQco(X_\bullet)}(\Mod(X_\bullet))
\]
of
\[
R(f_\bullet)_*:D_{\LQco(X_\bullet)}(\Mod(X_\bullet))
\rightarrow
D_{\LQco(Y_\bullet)}(\Mod(Y_\bullet))
\]
exists.
Let $J$ be a subcategory of $I$ which may not be admissible.

We define the natural transformation
\[
\xi(J,f_\bullet): (?)_{J}\circ f_\bullet^\times \rightarrow
f_{J}^\times\circ (?)_{J}
\]
to be the composite
\[
(?)_{J} f_\bullet^\times 
  \specialarrow{u}
(f_\bullet|_J)^\times R(f_\bullet|_J)_*(?)_J f_\bullet^\times
  \specialarrow{c}
(f_\bullet|_J)^\times (?)_J R(f_\bullet)_* f_\bullet^\times
  \specialarrow{\varepsilon}
(f_\bullet|_J)^\times (?)_J.
\]
By definition, $\xi$ is the 
conjugate map of $\theta(J,f_\bullet)$ in (\ref{derived-theta.eq}) if
$J$ is admissible.

\begin{Lemma}\label{restrict-xi-map.thm}
Let $J_2\subset J_1\subset I$ be subcategories of $I$.
Let $S$ be a scheme, $f_\bullet:X_\bullet\rightarrow Y_\bullet$ be
a morphism in $\Cal P(I,\Sch/S)$.
Assume both $X_\bullet$ and $Y_\bullet$ are quasi-compact separated.
Then the composite map
\begin{multline*}
(?)_{J_2}f_\bullet^\times\cong
(?)_{J_2}(?)_{J_1}f_\bullet^\times
\specialarrow{\via\xi(J_1)}
(?)_{J_2}(f_\bullet|_{J_1})^\times(?)_{J_1}\\
\specialarrow{\via\xi(J_2)}
(f_\bullet|_{J_2})^\times(?)_{J_2}(?)_{J_1}
=(f_\bullet|_{J_2})^\times(?)_{J_2}
\end{multline*}
is equal to $\xi(J_2)$.
\end{Lemma}

\proof Straightforward (and tedious) diagram drawing.
\qed

\begin{Lemma}\label{easy-commutativity.thm}
Let $S$ and $f_\bullet:X_\bullet\rightarrow Y_\bullet$ 
be as in Lemma~\ref{restrict-xi-map.thm}.
Let $J$ be a subcategory of $I$.
Assume that $Y_\bullet$ has flat arrows and $f_\bullet$ is of 
fiber type.
Then $\xi(J,f_\bullet): (?)_J f_\bullet^\times
\rightarrow (f_\bullet|_J)^\times (?)_J$
is an isomorphism between functors from $D_{\LQco(Y_\bullet)}(\Mod(Y_\bullet))
$ to $D_{\LQco(X_\bullet|_J)}(\Mod(X_\bullet|_J))$.
\end{Lemma}

\proof In view of Lemma~\ref{restrict-xi-map.thm}, we may assume
that $J=i$ for an object $i$ of $I$.
Then, as $i$ is an admissible subcategory of $I$ and
$\xi(i,f_\bullet)$ is an conjugate map of $\theta(i,f_\bullet)$,
it suffices to show that $(?)_j\theta(i,f_\bullet)$ is an isomorphism
for any $j\in\ob(I)$.
As $Y_\bullet$ has flat arrows, $L_i:\Mod(Y_i)\rightarrow \Mod(Y_\bullet)$ 
is exact.
As $f_\bullet$ is of fiber type, $L_i:\Mod(X_i)\rightarrow \Mod(X_\bullet)$ 
is also exact.

By Proposition~\ref{crutial.thm}, the composite
\[
(?)_jL_iR(f_i)_*\specialarrow{(?)_j\theta}(?)_jR(f_\bullet)_*L_i
\specialarrow{c}R(f_j)_*(?)_jL_i
\]
agrees with the composite
\begin{multline*}
(?)_jL_iR(f_i)_*\specialarrow{\lambda_{i,j}}
\bigoplus_{\phi\in I(i,j)}X_\phi^* R(f_i)_*
\specialarrow{\oplus\theta}
\bigoplus_{\phi}R(f_j)_*Y_\phi^*\\
\specialarrow{C}
R(f_j)_*\left(
  \bigoplus_\phi Y_\phi^*\right)
\specialarrow{\lambda_{i,j}^{-1}}
R(f_j)_*(?)_jL_i,
\end{multline*}
where $C$ is the canonical map.
By Lemma~\ref{respects-coproducts.thm}, $C$ is an isomorphism.
As $f_\bullet$ is of fiber type, $\theta:X_\phi^*R(f_i)_*\rightarrow
R(f_j)_*Y_\phi^*$ is an isomorphism for each $\phi\in I(i,j)$
by \cite[(3.9.5)]{Lipman}.
Hence the second composite is an isomorphism.
As the first composite is an isomorphism and $c$ is also an isomorphism,
we have that $(?)_j\theta(i,f_\bullet)$ is an isomorphism.
\qed

\paragraph Let $I$ be a small category, $S$ a scheme, and 
$f_\bullet:X_\bullet\rightarrow Y_\bullet$ a morphism in $\Cal P(I,\Sch/S)$.
Assume that $X_\bullet$ and $Y_\bullet$ are quasi-compact separated.

\begin{Lemma}\label{Rf-twisted-restrict.thm}
Let $J$ be a subcategory of $I$.
Then the following hold:
\begin{description}
\item[1] The composite map
\begin{multline*}
(?)_J\specialarrow{(?)_J u} (?)_J f_\bullet^\times R(f_\bullet)_*
\specialarrow{\via\xi(J,f_\bullet)} 
\\
(f_\bullet|_J)^\times (?)_J R(f_\bullet)_*
\specialarrow{\via c_{J,f_\bullet}}(f_\bullet|_J)^\times R(f_\bullet|_J)_*
(?)_J
\end{multline*}
agrees with $u(?)_J$, where $u:\Id\rightarrow
(f_\bullet|_J)^\times R(f_\bullet|_J)_*$ is the unit of adjunction.
\item[2] The composite map
\begin{multline*}
(?)_J R(f_\bullet)_* f_\bullet^\times  \specialarrow{\via c_{J,f_\bullet}}
R(f_\bullet|_J)_*(?)_J f_\bullet^\times
\specialarrow{\via\xi(J,f_\bullet)}
\\
R(f_\bullet|_J)_*(f_\bullet|_J)^\times(?)_J
\specialarrow{\varepsilon (?)_J}(?)_J
\end{multline*}
agrees with $(?)_J\varepsilon$.
\end{description}
\end{Lemma}

\proof This consists in straightforward diagram drawings.
\qed

\paragraph \label{well-behaved.par} Let $I$ be a small category.
For $i,j\in \ob(I)$, we say that $i\leq j$ if $I(i,j)\neq\emptyset$.
This definition makes $\ob(I)$ a pseudo-ordered set.
We say that $I$ is {\em ordered} if $\ob(I)$ is an ordered set
with the pseudo-order structure above, and $I(i,i)=\{\id\}$ for $i\in I$.

\begin{Lemma}\label{simplicial-restrict.thm}
Let $I$ be an ordered small category.
Let $J_0$ and $J_1$ be full subcategories of $I$, such that
$\ob(J_0)\cup \ob(J_1)=\ob(I)$, $\ob(J_0)\cap\ob(J_1)=\emptyset$, and
$I(j_1,j_0)=\emptyset$ for $j_1\in J_1$ and $j_0\in J_0$.
Let $X_\bullet\in\Cal P(I,\Sch/S)$.
Then, we have the following.
\begin{description}
\item[1] The unit of adjunction $u:\Id_{\Mod(X_\bullet|_{J_1})}\rightarrow
(?)_{J_1}\circ L_{J_1}$ is an isomorphism.
\item[2] $(?)_{J_0}\circ L_{J_1}$ is zero.
\item[3] $L_{J_1}$ is exact, and $J_1$ is an admissible subcategory of $I$.
\item[4] For $\Cal M\in\Mod(X_\bullet)$, $\Cal M_{J_0}=0$ if and only if
$\varepsilon: L_{J_1}\Cal M_{J_1}\rightarrow \Cal M$ is an isomorphism.
\item[5] The counit of adjunction $
(?)_{J_0}\circ R_{J_0} \rightarrow \Id_{\Mod(X_\bullet|_{J_0})}$
is an isomorphism.
\item[6] $(?)_{J_1}\circ R_{J_0}$ is zero.
\item[7] $R_{J_0}$ is exact and preserves local-quasi-coherence.
\item[8] For $\Cal M\in\Mod(X_\bullet)$, $\Cal M_{J_1}=0$ if and only if
$u:\Cal M\rightarrow R_{J_0}\Cal M_{J_0}$ is an isomorphism.
\item[9] The sequence
\[
0\rightarrow L_{J_1}(?)_{J_1}\specialarrow{\varepsilon}\Id\specialarrow
{u}R_{J_0}(?)_{J_0}\rightarrow0
\]
is semisplit 
exact, and induces a distinguished triangle in $D(\Mod(X_\bullet))$.
\end{description}
\end{Lemma}

\proof {\bf 1}
This is obvious by Lemma~\ref{induction-restrict.thm}.

{\bf 2} The category
$I^{(J_1\hookrightarrow I)}_{j_0}$ is empty, if
$j_0\in J_0$, since $I(j_1,j_0)=\emptyset$ if $j_1\in J_1$ and $j_0\in J_0$.
It follows that $(?)_{j_0}\circ L_{J_1}=0$ if $j_0\in J_0$.
Hence, $(?)_{J_0}\circ L_{J_1}=0$.

{\bf 3} This is trivial by {\bf 1,2} and their proof.

{\bf 4} The \lq if' part is trivial by {\bf 2}.
We prove the \lq only if' part.
By assumption and {\bf 2}, both $\Cal M_{J_0}$ and $(?)_{J_0}L_{J_1}\Cal M_
{J_1}$ are zero, and $(?)_{J_0}\varepsilon$ is an isomorphism.
On the other hand, 
$((?)_{J_1}\varepsilon)(u(?)_{J_1})=\id$, and $u$ is an isomorphism by
{\bf 1}.
Hence, 
\[
(?)_{J_1}\varepsilon:(?)_{J_1}L_{J_1}\Cal M_{J_1}\rightarrow (?)_{J_1}\Cal M
\]
is an isomorphism.
Hence, $\varepsilon$ is an isomorphism.

The assertions {\bf 5,6,7,8,9} are similar, and we omit the proof.
\qed

\begin{Lemma}\label{stupid-xi.thm}
Let $I$, $S$, 
$J_1$ and $J_0$ be as in Lemma~\ref{simplicial-restrict.thm}.
Let $f_\bullet:X_\bullet\rightarrow Y_\bullet$ be a morphism in $\Cal P(I,
\Sch/S)$.
Assume that $X_\bullet$ and $Y_\bullet$ are quasi-compact separated.
Then we have that $\theta(J_1,f_\bullet)$ and 
$\xi(J_1,f_\bullet)$ is an isomorphism.
\end{Lemma}

\proof Note that $J_1$ is admissible by Lemma~\ref{simplicial-restrict.thm},
{\bf 3}, and hence $\theta(J_1,f_\bullet)$ and $\xi(J_1,f_\bullet)$ are
defined.
Since we have $\xi(J_1,f_\bullet)$ is the conjugate of $\theta(
J_1,f_\bullet)$ by definition, it suffices to show that
$\theta(J_1,f_\bullet)$ is an isomorphism.
It suffices to show that 
\[
(?)_i LL_{J_1} R(f_\bullet|_{J_1})_*
   \specialarrow{(?)_i\theta}
(?)_i R(f_\bullet)_*  LL_{J_1}
   \cong 
R(f_i)_* (?)_iLL_{J_1}
\]
is an isomorphism for any $i\in I$.

If $i\in J_0$, then both hand sides are zero functors, and it is an 
isomorphism.
On the other hand, if $i\in J_1$, then the map in question is
equal to the composite isomorphism
\[
(?)_iLL_{J_1}R(f_\bullet|_{J_1})_*
\cong
(?)_i R(f_\bullet|_{J_1})_*
\specialarrow{c_{i,f_\bullet|_{J_1}}}
R(f_i)_*(?)_i
\cong
R(f_i)_*(?)_iLL_{J_1}
\]
by Proposition~\ref{crutial.thm}.
Hence $\theta(J_1,f_\bullet)$ is an isomorphism, as desired.
\qed

\paragraph Let $S$ be a scheme, $I$ an ordered small category,
$i\in I$, and $X_\bullet\in\Cal P(I,\Sch/S)$.
Let $J_1$ be a filter of $\ob(I)$ such that $i$ is a minimal element of $J_1$
(e.g., $[i,\infty)$), and
set $\Gamma_i:=L_{I,J_1}\circ R_{J_1,i}$.
Then we have
$(?)_j\Gamma_i=0$ if $j\neq i$ and $(?)_i\Gamma_i=\Id$.
It follows that $\Gamma_i$ preserves arbitrary limits and colimits (hence
is exact).
Assume that $X_i$ is quasi-compact and separated.
Then $D_{\Qco(X_i)}(\Mod(X_i))$ is compactly generated, and the derived 
functor
\[
\Gamma_i: D_{\Qco(X_i)}(\Mod(X_i))\rightarrow D_{\LQco(X_\bullet)}(
\Mod(X_\bullet))
\]
preserves coproducts.
It follows that there is a right adjoint
\[
\Sigma_i: D_{\LQco(X_\bullet)}(\Mod(X_\bullet))\rightarrow D_{\Qco(X_i)}
(\Mod(X_i)).
\]
As $\Gamma_i$ is obviously way-out left, we have $\Sigma_i$ is way-out right.

\paragraph \label{D^+.par}
Let $S$ be a scheme, $I$ a small category, 
and $X_\bullet\in \Cal P(I,\Sch/S)$.
We define $\Cal D^+(X_\bullet)$ (resp.\ $\Cal D^-(X_\bullet)$) 
to be the full subcategory 
of $D(\Mod(X_\bullet))$ consisting of $\Bbb F\in D(\Mod(X_\bullet))$
such that $\Bbb F_i$ is bounded below (resp.\ above) and has quasi-coherent
cohomology groups for each $i\in I$.
For a plump full subcategory $\Cal A$ of $\LQco(X_\bullet)$, we denote the
triangulated subcategory of $\Cal D^+(X_\bullet)$ 
(resp.\ $\Cal D^-(X_\bullet)$) consisting of objects all of whose 
cohomology groups belong to $\Cal A$ by
$\Cal D^+_{\Cal A}(X_\bullet)$
(resp.\ $\Cal D^-_{\Cal A}(X_\bullet)$).

\paragraph
Let $P$ be an ordered set.
We say that $P$ is upper Jordan-Dedekind (UJD for short) 
if for any $p\in P$, the subset
\[
[p,\infty):=\{q\in P\setbar q\geq p\}
\]
is finite.
We say that an ordered small category $I$ is UJD if the ordered
set $\ob(I)$ is UJD, and $I(i,j)$ is finite for $i,j\in I$.

\begin{Proposition}\label{twisted-restrict.thm}
Let $I$ be an ordered UJD small category.
Let $S$ be a scheme, and $g_\bullet:U_\bullet\rightarrow X_\bullet$
and $f_\bullet:X_\bullet\rightarrow Y_\bullet$ be morphisms in $\Cal P(I,
\Sch/S)$.
Assume that $Y_\bullet$ is separated noetherian
which has flat arrows, $f_\bullet$ is proper,
$g_\bullet$ is an open immersion such that $g_i(U_i)$ is dense in 
$X_i$ for each $i\in I$, and $f_\bullet\circ g_\bullet$ is
of fiber type.
Then $g_\bullet$ is of fiber type, and for any
$i\in I$ the composite natural map
\[
(?)_i g_\bullet^*f_\bullet^\times\specialarrow{\via \theta^{-1}}
g_i^*(?)_if_\bullet^\times
\specialarrow{\via \xi(i)}
g_i^*f_i^\times(?)_i
\]
is an isomorphism between functors $\Cal D^+(Y_\bullet)\rightarrow
D_{\Qco(U_i)}(\Mod(U_i))$,
where $\theta:g_i^*(?)_i\rightarrow (?)_i g_\bullet^*$ is the canonical
isomorphism.
\end{Proposition}

\proof Note that $U_\bullet$ has flat arrows, since $Y_\bullet$ has
flat arrows and $f_\bullet\circ g_\bullet$ is of fiber type.

We prove that $g_\bullet$ is of fiber type.
Let $\phi:i\rightarrow j$ be a morphism in $I$.
Then, the canonical map $(U_\phi,f_jg_j):
U_j\rightarrow U_i\times_{Y_i}Y_j$ is an 
isomorphism by assumption.
This map factors through $(U_\phi,g_j):U_j\rightarrow U_i\times_{X_i}X_j$,
and hence $(U_\phi,g_j)$ is a closed immersion.
On the other hand, it is an image dense open immersion, as can be seen 
easily, and hence it is an isomorphism.
So $g_\bullet$ is of fiber type.

Set $J_1:=[i,\infty)$ and $J_0:=\ob(I)\setminus J_1$.
By Lemma~\ref{stupid-xi.thm}, $\xi(J_1,f_\bullet)$ is an isomorphism.
By Lemma~\ref{restrict-xi-map.thm}, we may replace $I$ by $J_1$, and
we may assume that $I$ is an ordered finite category, and $i$ is a
minimal element of $\ob(I)$.
Now we have $\Cal D^+(Y_\bullet)$ agrees with $D^+_{\LQco(Y_\bullet)}(
\Mod(Y_\bullet))$.
Since we have $\ob(I)$ is finite, it is easy to see that $R(f_\bullet)_*$
is way-out both.
It follows that $f_\bullet^\times$ is way-out right.

It suffices to show
\[
g_i^*\xi(i): g_i^*(?)_i f_\bullet^\times\rightarrow
g_i^*f_i^\times(?)_i
\]
is an isomorphism of functors from $D^+_{\LQco(Y_\bullet)}(\Mod(Y_\bullet))$
to $D^+_{\Qco(U_i)}(\Mod(U_i))$.
As $g_i^*R(g_i)_*\cong\Id$, it suffices to show that
$R(g_i)_*g_i^*\xi(i)$
is an isomorphism.
This is equivalent to say that for any perfect complex $\Bbb P\in C(\Qco(X_i))
$,
we have
\begin{multline*}
R(g_i)_* g_i^*\xi(i):
\Hom_{D(\Mod(X_i))}(\Bbb P,R(g_i)_*g_i^*(?)_if_\bullet^\times)\\
\rightarrow
\Hom_{D(\Mod(X_i))}(\Bbb P,R(g_i)_*g_i^*f_i^\times(?)_i)
\end{multline*}
is an isomorphism.
By \cite[Lemma~2]{Verdier3}, this is equivalent to say that
the canonical map
\begin{multline*}
\indlim \Hom_{D(\Mod(Y_\bullet))}(R(f_\bullet)_* LL_i(\Bbb P
\otimes_{\Cal O_{X_i}}^{\bullet,L}\Cal J^n),?)\\
\rightarrow
\indlim \Hom_{D(\Mod(Y_\bullet))}(LL_iR(f_i)_*(\Bbb P
\otimes_{\Cal O_{X_i}}^{\bullet,L}\Cal J^n),?)
\end{multline*}
induced by the conjugate $\theta(i,f_\bullet)$ 
of $\xi(i)$ is an isomorphism,
where $\Cal J$ is a defining ideal sheaf of the closed 
subset $X_i\setminus U_i$ in $X_i$.

As $I$ is ordered and $\ob(I)$ is finite, we may label 
\[\ob(I)=\{i=i(0),i(1),i(2),\ldots\}
\]
so that $I(i(s),i(t))\neq\emptyset$ implies that $s\leq t$.
Let $J(r)$ denote the full subcategory of $I$ whose object set is
$\{i(r),i(r+1),\ldots\}$.

By descending induction on $t$, we prove that the map
\begin{multline*}
\via\theta(i,f_\bullet): 
\indlim \Hom_{D(\Mod(Y_\bullet))}(L_{J(t)}(?)_{J(t)}R(f_\bullet)_* LL_i(\Bbb P
\otimes_{\Cal O_{X_i}}^{\bullet,L}\Cal J^n),?)\\
\rightarrow
\indlim \Hom_{D(\Mod(Y_\bullet))}(L_{J(t)}(?)_{J(t)} LL_iR(f_i)_*(\Bbb P
\otimes_{\Cal O_{X_i}}^{\bullet,L}\Cal J^n),?)
\end{multline*}
is an isomorphism.
This is enough to prove the proposition, since $L_{J(1)}(?)_{J(1)}=\Id$.

Since the sequence
\[
0\rightarrow L_{J(t+1)}(?)_{J(t+1)}\specialarrow{\via\varepsilon}
L_{J(t)}(?)_{J(t)}\rightarrow \Gamma_{i(t)}(?)_{i(t)}\rightarrow 0
\]
is an exact sequence of exact functors, it suffices to prove that
the map
\begin{multline*}
\via\theta(i,f_\bullet): 
\indlim \Hom_{D(\Mod(Y_\bullet))}(\Gamma_{i(t)}(?)_{i(t)}
R(f_\bullet)_* LL_i(\Bbb P
\otimes_{\Cal O_{X_i}}^{\bullet,L}\Cal J^n),?)\\
\rightarrow
\indlim \Hom_{D(\Mod(Y_\bullet))}(\Gamma_{i(t)}(?)_{i(t)} LL_iR(f_i)_*(\Bbb P
\otimes_{\Cal O_{X_i}}^{\bullet,L}\Cal J^n),?)
\end{multline*}
is an isomorphism by induction assumption and the five lemma.
By Proposition~\ref{crutial.thm}, this is equivalent to say that the map
\begin{multline*}
\via\theta:
\indlim \Hom_{D(\Mod(Y_{i(t)}))}(R(f_{i(t)})_* \bigoplus_\phi LX_\phi^*(\Bbb P
\otimes_{\Cal O_{X_i}}^{\bullet,L}\Cal J^n),\Sigma_i(?))\\
\rightarrow
\indlim \Hom_{D(\Mod(Y_{i(t)}))}(\bigoplus_\phi Y_\phi^*R(f_i)_*(
\Bbb P
\otimes_{\Cal O_{X_i}}^{\bullet,L}\Cal J^n),\Sigma_i(?))
\end{multline*}
is an isomorphism, where the sum is taken over the {\em finite} set
$I(i,i(t))$.
It suffices to prove that the map
\begin{multline*}
\via\xi:
\indlim \Hom_{D(\Mod(X_i))}(\Bbb P
\otimes_{\Cal O_{X_i}}^{\bullet,L}\Cal J^n,
R(X_\phi)_*f_{i(t)}^\times \Sigma_i(?)
)\\
\rightarrow
\indlim \Hom_{D(\Mod(X_i))}(\Bbb P
\otimes_{\Cal O_{X_i}}^{\bullet,L}\Cal J^n,
f_i^\times R(Y_\phi)_* \Sigma_i(?)
)
\end{multline*}
induced by the map $\xi:R(X_\phi)_* f_{i(t)}^\times\rightarrow 
f_i^\times R(Y_\phi)_*$, which is conjugate to 
\[
\theta: Y_\phi^* R(f_i)_*
\rightarrow R(f_{i(t)})_* LX_\phi^*, 
\]
is an isomorphism for
$\phi\in I(i,i(t))$.
Since both $R(X_\phi)_*f_{i(t)}^\times \Sigma_i$ and
$f_i^\times R(Y_\phi)_*\Sigma_i$ are way-out right, it suffices to show
the canonical map
\[
g_i^*\xi(Y_\phi f_{i(t)}=f_iX_\phi): g_i^*R(X_\phi)_*f_{i(t)}^\times
\rightarrow
g_i^*f_i^\times R(Y_\phi)_*
\]
is an isomorphism between functors from $D^+_{\Qco(Y_{i(t)})}(\Mod(Y_{i(t)}))$
to $D^+_{\Qco(U_i)}(\Mod(U_i))$.
Let $X':=X_i\times_{Y_i}Y_{i(t)}$, $p_1:X'\rightarrow X_i$ be the
first projection, $p_2:X'\rightarrow Y_{i(t)}$ the second projection,
and $\pi: X_{i(t)}\rightarrow X'$ be the map $(X_\phi,f_{i(t)})$.
It is easy to see that $\xi(Y_\phi f_{i(t)}=f_i X_\phi)$ equals
the composite map
\[
R(X_\phi)_*f_{i(t)}^\times\cong R(p_1)_*R\pi_*\pi^\times p_2^\times
\specialarrow{\varepsilon}R(p_1)_*p_2^\times \cong f_i^\times R(Y_\phi)_*.
\]
Note that the last map is an isomorphism since $Y_\phi$ is flat.
As we have
\[
U_{i(t)}\rightarrow U_i\times_{Y_i}Y_{i(t)}\cong
U_i\times_{X_i}X'
\]
is an isomorphism and $g_i$ is an open immersion
by assumption, the canonical map
\[
g_i^* R(p_1)_*\rightarrow  R(U_{\phi})_* (\pi\circ g_{i(t)})^*
\]
is an isomorphism.
So it suffices to prove that
\[
(\pi\circ g_{i(t)})^* R\pi_*\pi^\times\specialarrow{\via\varepsilon}
(\pi\circ g_{i(t)})^* 
\]
is an isomorphism.

Consider the fiber square
\[
\begin{array}{ccc}
        U_{i(t)} & \specialarrow{g_{i(t)}} & X_{i(t)}\\
        \downarrow\hbox to 0pt{\scriptsize$\id$\hss} & \sigma &
        \downarrow\hbox to 0pt{\scriptsize$\pi$\hss}\\
        U_{i(t)} & \specialarrow{\pi \circ g_{i(t)}} & X'.
\end{array}
\]
By \cite{Verdier3}, $\zeta_0(\sigma):g_{i(t)}^*\pi^\times\rightarrow 
(\pi\circ g_{i(t)})^*$ is an isomorphism.
By definition (\ref{scheme-composition.par}), $\zeta_0$ is the composite map
\[
g_{i(t)}^*\pi^\times\specialarrow{u} \id^\times R\id_* g_{i(t)}^*\pi^\times
\specialarrow{\cong}\id_*(\pi\circ g_{i(t)})^*R\pi_*\pi^\times
\specialarrow{\varepsilon}(\pi\circ g_{i(t)})^*.
\]
Since the first and the second maps are isomorphisms, the third map
is an isomorphism.
This was what we wanted to prove.
\qed

\begin{Corollary}\label{fifteen.thm}
Under the same assumption as in the proposition,
we have $g_\bullet^*f_\bullet^\times(\Cal D^+(Y_\bullet))\subset
\Cal D^+(U_\bullet)$.
\end{Corollary}

\proof This is because $g_i^*f_i^\times(D^+_{\Qco(Y_i)}(\Mod(Y_i)))
\subset D^+_{\Qco(U_i)}(\Mod(U_i))$ for each $i\in I$.

\section{Open immersion base change}

\paragraph\label{define-zeta.par}
Let $S$ be a scheme, $I$ a small category, and
\[
\begin{array}{ccc}
X_\bullet' & \specialarrow{g_\bullet'} & X_\bullet \\
\downarrow\hbox to 0pt{\scriptsize$f_\bullet'$\hss}& \sigma &
\downarrow\hbox to 0pt{\scriptsize$f_\bullet$\hss}\\
Y_\bullet' & \specialarrow{g_\bullet} & Y_\bullet
\end{array}
\]
a fiber square in $\Cal P(I,\Sch/S)$.
Assume that the all diagrams of schemes in the diagram are quasi-compact
separated, and $g_\bullet$ is flat.
By Lemma~\ref{theta-g-f-isom.thm}, the canonical map
\[
\theta(g_\bullet,f_\bullet)
: g_\bullet^* R(f_\bullet)_*\rightarrow R(f_\bullet')_* (g_\bullet')^*
\]
is an isomorphism of functors from $D_{\LQco(X_\bullet)}(\Mod(X_\bullet))$ 
to $D_{\LQco(Y_\bullet')}(\Mod(Y_\bullet'))$.
We define $\zeta(\sigma)=\zeta(g_\bullet,f_\bullet)$ to be the composite map
\[
\zeta(\sigma): (g_\bullet')^*f_\bullet^\times
\specialarrow{u}
(f_\bullet')^\times R(f_\bullet')_*(g_\bullet')^*f_\bullet^\times
\specialarrow{\theta^{-1}}
(f_\bullet')^\times g_\bullet^* R(f_\bullet)_* f_\bullet^\times
\specialarrow{\varepsilon}
(f_\bullet')^\times g_\bullet^*.
\]

\begin{Lemma}\label{zeta-restrict.thm}
Let $\sigma$ be as above, and $J$ a subcategory
of $I$.
Then the diagram
\[
\begin{array}{ccccc}
  (?)_J(g_\bullet')^*f_\bullet^\times & \specialarrow{\theta^{-1}} &
  (g_\bullet'|_J)^*(?)_Jf_\bullet^\times & \specialarrow{\xi} &
  (g_\bullet'|_J)^*(f_\bullet|_J)^\times(?)_J\\
  \downarrow\hbox to 0pt{\scriptsize$\zeta$\hss} & & & &
  \downarrow\hbox to 0pt{\scriptsize$\zeta$\hss}\\
  (?)_J(f_\bullet')^\times g_\bullet^* & \specialarrow{\xi} &
  (f_\bullet'|_J)^\times(?)_J g_\bullet^* & \specialarrow{\theta^{-1}} &
  (f_\bullet'|_J)^\times (g_\bullet|_J)^*(?)_J
\end{array}
\]
is commutative.
\end{Lemma}

\proof Follows immediately from Lemma~\ref{Rf-twisted-restrict.thm} and
the commutativity of (\ref{theta-f-J-restrict.eq}).
\qed

\begin{Theorem}\label{sixteen.thm}
Let $S$ be a scheme, $I$ an ordered UJD small category, 
and
\[
\begin{array}{ccccccc}
V_\bullet & \specialarrow j_\bullet & U'_\bullet & \specialarrow{i'_\bullet} 
    & X'_\bullet &  &  \\
  &                 &
\downarrow\hbox to 0pt{\scriptsize$p^U_\bullet$\hss} & \sigma  &
\downarrow\hbox to 0pt{\scriptsize$p^X_\bullet$\hss} &   &  \\
  &                 & U_\bullet & \specialarrow{i_\bullet}&  X_\bullet 
                    & \specialarrow{q_\bullet} & Y_\bullet
\end{array}
\]
be a diagram in $\Cal P(I,\Sch/S)$.
Assume the following.
\begin{description}
\item[1] $Y_\bullet$ is separated noetherian.
\item[2] $j_\bullet$, $i'_\bullet$ and $i_\bullet$ are image dense
open immersions.
\item[3] $q_\bullet$, $p_\bullet^X$ and $p_\bullet^U$ are proper.
\item[4] $V_\bullet$, $U_\bullet$ and $Y_\bullet$ have flat arrows.
\item[5] $p^X_\bullet i'_\bullet= i_\bullet p^U_\bullet$.
\item[6] $q_\bullet i_\bullet$ and $p_\bullet^U j_\bullet$ are of fiber type.
\end{description}
Then $\sigma$ is a fiber square, and the composite map
\[
j_\bullet^* \zeta(\sigma) q_\bullet^\times: j_\bullet^*(i'_\bullet)^*
(p^X_\bullet)^\times q_\bullet^\times
\rightarrow j_\bullet^* (p^U_\bullet)^\times i_\bullet^* q_\bullet^\times
\]
is an isomorphism of functors from $\Cal D^+(Y_\bullet)$ to
$\Cal D^+(V_\bullet)$.
\end{Theorem}

\proof The square $\sigma$ is a fiber square, since
the canonical map $U'_\bullet\rightarrow U_\bullet\times_{X_\bullet}X'_\bullet$
is an image dense  closed open immersion, and is an isomorphism.

To prove the theorem, it suffices to show that the map in question is an
isomorphism after applying $(?)_i$ for any $i\in I$.
By Proposition~\ref{twisted-restrict.thm}, 
Lemma~\ref{zeta-restrict.thm}, and \cite[(3.7.2), (iii)]{Lipman},
the problem is reduced to the flat base change theorem (in fact open
immersion base change theorem is enough) for schemes
\cite[Theorem~2]{Verdier3}, and we are done.
\qed

\section[The existence of compactification and composition data]{The 
existence of compactification and composition data for
finite diagrams of schemes}

\paragraph\label{composition-setting.par}
Let $S$ be a noetherian separated scheme, and
$I$ an ordered finite category which is non-empty.
Let $\Cal A$ denotes the category of
$I\op$-diagrams of $S$-schemes separated of
finite type over $S$.
Let $\Cal P$ denote the class of proper morphisms in $\Mor(\Cal A)$.
Let $\Cal I$ denote the class of image dense open immersions 
in $\Mor(\Cal A)$.

Define $\Cal D(X_\bullet):=D_{\LQco(X_\bullet)}(\Mod(X_\bullet))$ 
for $X_\bullet\in\ob(\Cal A)$.
Define a pseudofunctor $(?)^{\#}$ on $\Cal A_{\Cal P}$
to be $(?)^\times$, where 
$X^{\#}_\bullet=\Cal D(X_\bullet)$ for $X_\bullet\in\ob(\Cal A_{\Cal P})$.
Define a pseudofunctor $(?)^\flat$ on $\Cal A_{\Cal I}$
to be $(?)^*$, where 
$X^{\flat}_\bullet=\Cal D(X_\bullet)$ for $X_\bullet\in\ob(\Cal A_{\Cal I})$.

For a pi-square $\sigma$, define $\zeta(\sigma)$ to be the natural map
defined in (\ref{define-zeta.par}).

\begin{Lemma}\label{eleven.thm}
Let the notation be as above.
Conditions {\bf 1--6} and {\bf 8--12} in 
Definition~\ref{composition-data.def} are satisfied.
Moreover, any pi-square is a fiber square.
\end{Lemma}

\proof This is easy.
\qed

\begin{Proposition}\label{finite-compactification.thm}
Let the notation be as in 
{\rm (\ref{composition-setting.par})}.
Then the condition {\bf 7} in {\rm Definition~\ref{composition-data.def}} is
satisfied.
That is, for any morphism $f_\bullet$ in $\Cal A$, 
there is a factorization $f_\bullet=p_\bullet j_\bullet$ 
with $p_\bullet\in\Cal P$ and $j_\bullet\in\Cal I$.
\end{Proposition}

\proof Label the object set $\ob(I)$ of $I$ as
$\{i(1),\cdots,i(n)\}$ so that $I(i(s),i(t))=\emptyset$ if $s>t$.
Set $J(r)$ to be the full subcategory of
$I$ with $\ob(J(r))=\{i(1),\ldots,i(r)\}$.
By induction on $r$, we construct morphisms $h_\bullet(r):X_\bullet|_{J(r)}
\rightarrow Z_\bullet(r)$ in $\Cal P(J(r),\Sch/S)$ such that
\begin{description}
\item[1] $Z_\bullet(r)$ is proper over $S$;
\item[2] $h_\bullet(r)_i$ is an open immersion, whose scheme theoretic
image is $Z_\bullet(r)_i$ for $i\in J(r)$;
\item[3] $Z_\bullet(r)|_{J(j)}=Z_\bullet(j)$ and
$h_\bullet(r)|_{J(j)}=h_\bullet(j)$ for $j<r$.
\end{description}
We may assume that the construction is done for $j<r$.

First, consider the case where $i(r)$ is a minimal element in $\ob(I)$.
By Nagata's compactification theorem \cite{Nagata}, there is 
an $S$-morphism
\[
h: X_{i(r)}\rightarrow Z
\]
such that $Z$ is proper over $S$ and $h$ is an open immersion whose
scheme theoretic image is $Z$.
Now define $Z_\bullet(r)_{i(r)}:=Z$ and $Z_\bullet(r)_{\id_{i(r)}}=\id_Z$.
Defining the other structures after {\bf 3}, we get $Z_\bullet(r)$,
since $I(i(r),i(s))=\emptyset=I(i(s),i(r))$ for $s<r$ 
and $I(i(r),i(r))=\{\id\}$
by assumption.
Define $h_\bullet(r)_{i(r)}:=h$ and by {\bf 3}, we get a morphism
$h_\bullet(r): X_\bullet|_{J(r)}\rightarrow Z_\bullet(r)$ by the same reason,
and {\bf 1,2,3} are satisfied by induction assumption.
So this case is OK.

Now assume that $i(r)$ is not minimal so that $\bigcup_{i<r}I(i(j),i(r))\neq
\emptyset$.
Consider the $S$-morphism
\[
\varphi:
X_{i(r)}\specialarrow{(X_\phi)}\prod_{j<r}\prod_{\phi\in I(i(j),i(r))}X_{i(j)}
\specialarrow{(h_j)}
\prod_{j<r}\prod_{\phi\in I(i(j),i(r))}Z_j,
\]
where $Z_j:=(Z_\bullet(j))_{i(j)}$, and $h_j:=(h_\bullet(j))_{i(j)}$.
By Nagata's compactification theorem \cite{Nagata}, 
There is a factorization
\[
X_{i(r)}\specialarrow{h_r}Z_r\specialarrow{p_r}\prod_{j<r}\prod_{\phi\in
I(i(j),i(r))}Z_j
\]
such that $h_r$ is an open immersion $S$-morphism 
such that the scheme theoretic image of $h_r$ agrees with
$Z_r$, and $p_r$ is a proper $S$-morphism.
Set $Z_\bullet(r)_{i(j)}:=Z_j$, and $Z_\bullet(r)_{i(r)}:=Z_r$.
For $j\leq j'<r$ and $\phi\in I(i(j),i(j'))$, define $(Z_\bullet(r))_\phi:=
Z_\bullet(r-1)_\phi$.
For each $j<r$ and $\phi\in I(i(j),i(r))$, 
define $(Z_\bullet(r))_\phi:Z_r\rightarrow Z_j$ to be the composite
\[
Z_r\specialarrow{p_r}\prod_{j<r}\prod_{\phi}Z_j\specialarrow{\rm 
projection} Z_j.
\]
Finally, define $(Z_\bullet(r))_{\id_{i(r)}}:=\id_{Z_r}$.
Almost by definition, we have $h_j\circ X_\phi=(Z_\bullet(r))_\phi
\circ h_r$ for $j<r$ and
$\phi\in I(i(j),i(r))$.
As $h_r$ is scheme theoretically dense, it is easy to see that $Z_\bullet(r)$
is a $J(r)\op$ diagram of $S$-schemes, and the conditions {\bf 1,2,3} are
checked easily.

Hence there is a morphism $h_\bullet:X_\bullet\rightarrow Z_\bullet$ 
of $\Cal P(I,\Sch/S)$ such that $Z_\bullet$ is proper over $S$, and
$h_\bullet$ is an image-dense open immersion.

Now let $f_\bullet:X_\bullet\rightarrow Y_\bullet$ be a morphism in $\Cal A$.
Consider the morphism
\[
\psi_\bullet:
X_\bullet\specialarrow{\Gamma(f_\bullet)}X_\bullet\times Y_\bullet
\specialarrow{h_\bullet\times 1}Z_\bullet\times Y_\bullet.
\]
We define $\bar X_\bullet$ as follows.
For each $i\in I$, define $\bar X_i$ to be the scheme theoretic image of
$\psi_i$, and $\bar X_\phi$ is the restriction of $Z_\phi\times Y_\phi$.
As $\psi_\bullet$ is an immersion, we have a factorization
\[
X_\bullet\specialarrow{i_\bullet}\bar X_\bullet\specialarrow{q_\bullet}
Z_\bullet\times Y_\bullet
\]
of $\psi_\bullet$ such that $i_\bullet\in \Cal I$ and $q_\bullet$ is
a closed immersion.

Let $p_2:Z_\bullet\times Y_\bullet\rightarrow Y_\bullet$ be the
second projection.
Set $p_\bullet:=p_2\circ q_\bullet$.
Then, as we have $p_2\circ \psi_\bullet=f_\bullet$, we have
$f_\bullet=p_\bullet i_\bullet$, and $p_\bullet\in\Cal P$.
This proves the lemma.
\qed

\begin{Theorem}\label{finite-data.thm}
Let the notation be as 
in {\rm (\ref{composition-setting.par})}.
Set $\Cal F$ to be the subcategory of $\Cal A$ whose 
objects are objects of $\Cal A$ with flat arrows, and
whose morphisms are morphisms in $\Cal A$ of fiber type.
Define $\Cal D^+$ by
$\Cal D^+(X_\bullet):=D^+_{\LQco(X_\bullet)}(\Mod(X_\bullet))$.
Then $(\Cal A,\Cal F,\Cal P,\Cal I,\Cal D,\Cal D^+,(?)^{\#},(?)^{\flat},
\zeta)$ is a composition data of contravariant pseudofunctors.
\end{Theorem}

\proof Conditions {\bf 1--12} in Definition~\ref{composition-data.def}
have already been checked.
{\bf 13} follows from Lemma~\ref{fiber-type-1.thm}.
{\bf 14} is trivial.
Since $I$ is finite, 
the definition of $\Cal D^+(X_\bullet)$ is consistent with
that in (\ref{D^+.par}).
Hence {\bf 15} is Corollary~\ref{fifteen.thm}.
{\bf 16} is Theorem~\ref{sixteen.thm}.
\qed

We call the composite of $(?)^{\#}$ and $(?)^\flat$ defined by the
composition data in the theorem the {\em equivariant twisted inverse
pseudofunctor}, and denote it by $(?)^!$.

\section{Flat base change}

Let the notation be as in Theorem~\ref{finite-data.thm}.
Let $f_\bullet:X_\bullet\rightarrow Y_\bullet$ be a morphism
in $\Cal F$, and $J$ a subcategory of $I$.
Let $f_\bullet=p_\bullet i_\bullet$ be a compactification.

\begin{Lemma} The composite map
\[
(?)_J f_\bullet^!
  \specialarrow{\via\rho}
(?)_J i_\bullet^* p_\bullet^\times
\specialarrow{\theta^{-1}}
i_J^*(?)_J p_\bullet^\times
\specialarrow{\xi}
i_J^*p_J^\times (?)_J
\specialarrow{\via\rho}
f_J^!(?)_J
\]
is independent of choice of compactification 
$f_\bullet=p_\bullet i_\bullet$,
where $\rho$'s are independence isomorphisms.
\end{Lemma}

The proof utilizes 
Lemma~\ref{independence-isom.thm}, and left to the reader.
We denote 
by $\bar\xi=\bar\xi(i,f_\bullet)$ the composite map in the lemma.

\begin{Lemma}\label{restrict-bar-xi.thm}
Let $f_\bullet:X_\bullet\rightarrow Y_\bullet$ be a morphism
in $\Cal F$, and $K\subset J\subset I$ be subcategories.
Then the composite map
\[
(?)_Kf_\bullet^!\cong (?)_K(?)_J f_\bullet^!
  \specialarrow{\bar\xi}
(?)_K f_J^! (?)_J
  \specialarrow{\bar\xi}
f_K^! (?)_K (?)_J\cong f_K^!(?)_K
\]
agrees with $\bar\xi(K,f_\bullet)$.
\end{Lemma}

\proof Follows easily from Lemma~\ref{restrict-xi-map.thm}.
\qed

\begin{Lemma}\label{twisted-restrict2.thm}
Let $f_\bullet:X_\bullet\rightarrow Y_\bullet$ be a
morphism in $\Cal F$, and $J$ a subcategory of $I$.
Then, $\bar\xi(J,f_\bullet)$ is an isomorphism.
\end{Lemma}

\proof 
It suffices to show that $(?)_i\bar\xi(J,f_\bullet)$ is an isomorphism for
any $i\in\ob(J)$.
By Lemma~\ref{restrict-bar-xi.thm}, 
we have
\[
\bar\xi(i,f_J)\circ ((?)_i\bar\xi(J,f_\bullet))=\bar\xi(i,f_\bullet).
\]
By Proposition~\ref{twisted-restrict.thm},
we have $\bar\xi(i,f_J)$ and $\bar\xi(i,f_\bullet)$ are isomorphisms.
Hence the natural map $(?)_i\bar\xi(J,f_\bullet)$ is also an isomorphism.
\qed

\begin{Lemma}\label{flat-base-change-dense.thm}
Let $f:X\rightarrow Y$ be a flat morphism of locally noetherian
schemes, and $U$ a dense open subset of $Y$.
Then $f^{-1}(U)$ is a dense open subset of $X$.
\end{Lemma}

\proof The question is local both on $Y$ and $X$, and hence we may assume
that both $Y=\Spec A$ and $X=\Spec B$ are affine.
Let $I$ be the radical ideal of $A$ defining the closed subset $Y\setminus U$.
By assumption, $I$ is not contained in any minimal prime of $A$.
Assume that $f^{-1}(U)$ is not dense in $X$.
Then, there is a minimal prime $P$ of $B$ which contains $IB$.
As we have $I\subset IB\cap A\subset P\cap A$ and $P\cap A$ is minimal
by the going-down theorem (see \cite[Theorem~9.5]{CRT}), this is a
contradiction.
\qed

\paragraph\label{fiber.par}
Let $b: S'\rightarrow S$ be a morphism of noetherian separated
schemes.
Let $\Cal A'$ be the full subcategory of $\Cal P(I,\Sch/S')$ consisting
of objects separated of finite type over $S'$.
Let $\Cal F'$ be the category of objects of $\Cal A'$ with flat arrows
and morphisms in $\Cal A'$ of fiber type.
Any $S'$-scheme is viewed as an $S$-scheme via $b$.

Let
\[
\begin{array}{ccc}
X_\bullet' & \specialarrow{f_\bullet'} & Y_\bullet'\\
\sdarrow{g_\bullet^X} & \sigma & \sdarrow{g_\bullet}\\
X_\bullet & \specialarrow{f_\bullet} & Y_\bullet
\end{array}
\]
be a diagram in $\Cal P(I,\Sch/S)$ such that
\begin{description}
\item[1] $X_\bullet'$ and $Y_\bullet'$ are in $\Cal F'$,
and $f_\bullet'$ is a morphism in $\Cal F'$.
\item[2] $f_\bullet$ is a morphism in $\Cal F$;
\item[3] $\sigma$ is a fiber square;
\item[4] $g_\bullet$ is flat.
\end{description}

By assumption, there is a diagram
\begin{equation}\label{compactification-flat-base-change.eq}
\begin{array}{ccccc}
X_\bullet' & \specialarrow{i_\bullet'} & Z_\bullet' & 
    \specialarrow{p_\bullet'}& Y_\bullet'\\
\sdarrow{g_\bullet^X} & \sigma_1 & \sdarrow{g_\bullet^Z} & \sigma_2 & 
     \sdarrow{g_\bullet}\\
X_\bullet & \specialarrow{i_\bullet} & Z_\bullet & 
    \specialarrow{p_\bullet}& Y_\bullet\\
\end{array}
\end{equation}
such that $f_\bullet=p_\bullet i_\bullet$ is a compactification, 
$\sigma_1$ and $\sigma_2$ are fiber squares, and the whole rectangle
$\sigma_1\sigma_2$ equals $\sigma$.
By Lemma~\ref{flat-base-change-dense.thm}, we have that
$f_\bullet'=p_\bullet' i'_\bullet$ is a compactification.

\begin{Lemma} The composite map
\[
(g_\bullet^X)^* f_\bullet^!
  \specialarrow{\rho}
(g_\bullet^X)^* i_\bullet^* p_\bullet^\times
  \specialarrow{d}
(i_\bullet')^*(g_\bullet^Z)^*p_\bullet^\times
  \specialarrow{\zeta}
(i_\bullet')^*(p_\bullet')^\times g_\bullet^*
  \specialarrow{\rho}
(f_\bullet')^! g_\bullet^*
\]
is independent of choice of the diagram
{\rm (\ref{compactification-flat-base-change.eq})}, 
and depends only on $\sigma$,
where $\rho$'s are independence isomorphisms.
\end{Lemma}

\proof Obvious by Lemma~\ref{independence-isom.thm}.
\qed

We denote the composite map in the lemma by $\bar\zeta=\bar\zeta(\sigma)$.

\begin{Theorem}\label{flat-base-change.thm}
Let the notation be as above.
Then we have:
\begin{description}
\item[1] Let $J$ be an admissible subcategory of $I$. 
Then the diagram
\begin{equation}\label{zeta-restrict2.eq}
\begin{array}{ccccc}
(?)_J(g_\bullet^X)^*f_\bullet^! & \specialarrow{\theta^{-1}} &
 (g_J^X)^*(?)_Jf_\bullet^! & \specialarrow{\bar\xi} &
 (g_J^X)^*f_J^!(?)_J\\
\sdarrow{(?)_J\bar\zeta(\sigma)} & & & & \sdarrow{\bar\zeta(\sigma_J)(?)_J}\\
(?)_J(f_\bullet')^!(g_\bullet)^* & \specialarrow{\bar\xi} &
 (f_J')^!(?)_J(g_\bullet)^* & \specialarrow{\theta^{-1}} &
 (f_J')^!(g_J)^*(?)_J
\end{array}
\end{equation}
is commutative.
\item[2] $\bar\zeta(\sigma)$ is an isomorphism.
\end{description}
\end{Theorem}

\proof {\bf 1} is an immediate consequence of
Lemma~\ref{zeta-restrict.thm} and \cite[(3.7.2)]{Lipman}.

{\bf 2} Let $i$ be an object of $I$.
By Lemma~\ref{twisted-restrict2.thm}, the horizontal arrows
in the diagram (\ref{zeta-restrict2.eq}) are isomorphisms.
By Verdier's flat base change theorem \cite[Theorem~2]{Verdier3}, 
we have that $\bar\zeta(\sigma_i)$ is an isomorphism.
Hence, we have that $(?)_i\bar\zeta(\sigma)$ is an isomorphism for
any $i\in I$ by {\bf 1} applied to $J=i$, and the assertion follows.
\qed

\section{Preservation of Quasi-coherent cohomology}

\paragraph Let the notation be as in (\ref{composition-setting.par}).
Let $\Cal F$ be as in Theorem~\ref{finite-data.thm}.

\begin{Lemma}\label{coherence-schemes.thm}
Let $f:X\rightarrow Y$ be a morphism of finite type
between separated noetherian schemes.
If $\Bbb F\in \Cal D^+_{\Coh(Y)}(\Mod(Y))$, then $f^!\Bbb F\in
\Cal D^+_{\Coh(X)}(\Mod(X))$.
\end{Lemma}

\proof We may assume that both $Y$ and $X$ are affine.
So we may assume that $f$ is either smooth or a closed immersion.
The case where $f$ is smooth is obvious by \cite[Theorem~3]{Verdier3}.
The case where $f$ is a closed immersion is also obvious by
Proposition~III.6.1 and Theorem~III.6.7 in \cite{Hartshorne2}.
\qed

\begin{Proposition} Let $f_\bullet:X_\bullet\rightarrow Y_\bullet$ be 
a morphism in $\Cal F$, and $\phi:i\rightarrow j$ a morphism in $I$.
Then the composite map
\[
X_\phi^*(?)_i f_\bullet^!
  \specialarrow{\alpha_\phi}
(?)_j f_\bullet^!
  \specialarrow{\bar\xi}
f_j^! (?)_j
\]
agrees with the composite map
\[
X_\phi^*(?)_i f_\bullet^!
  \specialarrow{\bar\xi}
X_\phi^* f_i^! (?)_i
  \specialarrow{\bar\zeta}
f_j^! Y_\phi^*(?)_i
  \specialarrow{\alpha_\phi}
f_j^! (?)_j.
\]
\end{Proposition}

\proof By Lemma~\ref{restrict-bar-xi.thm},
we may assume that $I$ is the ordered category given by 
$\ob(I)=\{i,j\}$ and $I(i,j)=\{\phi\}$.

Then it is easy to see that there is a compactification
\[
X_\bullet\specialarrow{i_\bullet}Z_\bullet\specialarrow{p_\bullet}Y_\bullet
\]
of $f_\bullet$ such that $p_\bullet$ is of fiber type.
Note that $i_\bullet$ is of fiber type, and $Z_\bullet$ has flat arrows.

By the definition of $\bar\xi$ and $\bar\zeta$, it suffices to prove that
the composite map
\[
X_\phi^*(?)_i i_\bullet^* p_\bullet^\times
  \specialarrow{\theta^{-1}}
X_\phi^* i_i^*(?)_i p_\bullet^\times
  \specialarrow{d}
i_j^* Z_\phi^*(?)_i p_\bullet^\times
  \specialarrow{\xi}
i_j^* Z_\phi^* p_i^\times (?)_i
  \specialarrow{\zeta}
i_j^*p_j^\times Y_\phi^*(?)_i
  \specialarrow{\alpha_\phi}
i_j^*p_j^\times (?)_j
\]
agrees with
\[
X_\phi^*(?)_i i_\bullet^* p_\bullet^\times
  \specialarrow{\alpha_\phi}
(?)_j i_\bullet^* p_\bullet^\times
  \specialarrow{\theta^{-1}}
i_j^* (?)_j p_\bullet^\times
  \specialarrow{\xi}
i_j^* p_j^\times(?)_j.
\]

By the ``derived version'' of (\ref{inv-im-translate.par}),
the composite map
\[
X_\phi^*(?)_i i_\bullet^* p_\bullet^\times
  \specialarrow{\theta^{-1}}
X_\phi^* i_i^*(?)_i p_\bullet^\times
  \specialarrow{d}
i_j^* Y_\phi^*(?)_i p_\bullet^\times
  \specialarrow{\alpha_\phi}
i_j^*(?)_j p_\bullet^\times
\]
agrees with
\[
X_\phi^*(?)_i i_\bullet^* p_\bullet^\times
  \specialarrow{\alpha_\phi}
(?)_j i_\bullet^* p_\bullet^\times
  \specialarrow{\theta^{-1}}
i_j^* (?)_j p_\bullet^\times.
\]
Hence it suffices to prove the map
\[
Z_\phi^*(?)_i p_\bullet^\times
  \specialarrow{\xi}
Z_\phi^*p_i^\times(?)_i
  \specialarrow{\zeta}
p_j^\times Y_\phi^*(?)_i
  \specialarrow{\alpha_\phi}
p_j^\times (?)_j
\]
agrees with
\[
Z_\phi^*(?)_i p_\bullet^\times
  \specialarrow{\alpha_\phi}
(?)_j p_\bullet^\times
  \specialarrow{\xi}
p_j^\times (?)_j.
\]

Now the proof consists in a straightforward diagram drawing utilizing 
the derived version of (\ref{direct-image.par}).
\qed

\begin{Corollary}\label{preservation.thm}
Let $f_\bullet:X_\bullet\rightarrow Y_\bullet$ be 
a morphism in $\Cal F$.
Then we have $f_\bullet^!(\Cal D^+_{\Qco}(Y_\bullet))\subset
\Cal D^+_{\Qco}(X_\bullet)$ and
$f_\bullet^!(\Cal D^+_{\Coh}(Y_\bullet))\subset
\Cal D^+_{\Coh}(X_\bullet)$.
\end{Corollary}

\proof Let $\phi:i\rightarrow j$ be a morphism in $I$.
By Theorem~\ref{flat-base-change.thm}, {\bf 2}, 
Lemma~\ref{twisted-restrict2.thm}
and the proposition, we have $\alpha_\phi:X_\phi^*(?)_if_\bullet^!
\rightarrow (?)_j f_\bullet^!$ is an isomorphism
if $\alpha_\phi:Y_\phi^*(?)_i\rightarrow (?)_j$ is an isomorphism.
So $f_\bullet^!$ preserves equivariance of cohomology groups,
and the first assertion follows.

On the other hand, by Lemma~\ref{coherence-schemes.thm}
and Proposition~\ref{twisted-restrict.thm}, $f^!$ preserves local
coherence of cohomology groups.
Hence it also preserves the coherence of cohomology groups, by the
first paragraph.
\qed

\section{Compatibility with derived direct image}

\paragraph Let the notation be as in (\ref{fiber.par}).
Consider the diagram (\ref{compactification-flat-base-change.eq}).
Lipman's theta $\theta(\sigma_2): g_\bullet^*R(p_\bullet)_*
\rightarrow
R(p_\bullet')_* (g_\bullet^Z)^*$ 
induces the conjugate map
\[
\xi(\sigma_2): R(g_\bullet^Z)_* (p_\bullet')^\times
\rightarrow
p_\bullet^\times R(g_\bullet)_*.
\]
As $\sigma_2$ is a fiber square, $\theta(\sigma_2)$ is an isomorphism.
Hence $\xi(\sigma_2)$ is also an isomorphism.
Note that
\[
\theta: i_\bullet^*R(g_\bullet^Z)_*\rightarrow R(g_\bullet^X)_*(i_\bullet')^*
\]
is an isomorphism, since $\sigma_1$ is a fiber square.
We define $\bar \xi: R(g_\bullet^X)_* (f_\bullet')^!\rightarrow f_\bullet^!
R(g_\bullet)_*$ to be the composite
\[
R(g_\bullet^X)_* (f_\bullet')^!
\specialarrow{\rho}
R(g_\bullet^X)_* (i_\bullet')^*(p_\bullet')^\times
\specialarrow{\theta^{-1}}
i_\bullet^*R(g_\bullet^Z)_*(p_\bullet')^\times
\specialarrow{\xi}
i_\bullet^*p_\bullet^\times R(g_\bullet)_*.
\specialarrow{\rho}
f_\bullet^! R(g_\bullet)_*.
\]

As all the maps in the composition are isomorphisms, we have

\begin{Lemma} $\bar\xi$ is an isomorphism.
\end{Lemma}

\begin{Lemma} For any subcategory $J$, the composite
\[
(?)_J R(g_\bullet^Z)_* (p_\bullet')^\times
\specialarrow{\xi}
(?)_J p_\bullet^\times R(g_\bullet)_*
\specialarrow{\xi}
(p_\bullet|_J)^\times (?)_J R(g_\bullet)_*
\specialarrow{c}
(p_\bullet|_J)^\times R(g_\bullet|_J)_* (?)_J
\]
agrees with the composite
\begin{multline*}
(?)_J R(g_\bullet^Z)_* (p_\bullet')^\times
\specialarrow{c}
R(g_\bullet^Z|_J)_*(?)_J (p_\bullet')^\times
\specialarrow{\xi}\\
R(g_\bullet^Z|_J)_*(p_\bullet'|_J)^\times(?)_J
\specialarrow{\xi}
(p_\bullet|_J)^\times R(g_\bullet|_J)_* (?)_J.
\end{multline*}
\end{Lemma}

\proof Follows immediately from Lemma~\ref{Rf-twisted-restrict.thm}.
\qed

\section{Compatibility with derived right inductions}

\paragraph Let the notation be as in Theorem~\ref{finite-data.thm}.
Let $X_\bullet\in\Cal F$, and $J$ a subcategory of $I$.
As $I$ is assumed to be a finite category, $I_i^{(J\rightarrow I)}$ is
a finite category for any $i\in\ob(I)$.
Hence for any $\Cal M\in \LQco(X_\bullet|_J)$, we have that
$R_J\Cal M\in\LQco(X_\bullet)$ by (\ref{right-induction.par}).
It is easy to see that $RR_J: D_{\LQco(X_\bullet|_J)}^+(\Mod(X_\bullet|_J))
\rightarrow D_{\LQco(X_\bullet)}^+(\Mod(X_\bullet))$ is the right adjoint
of $(?)_J$.

Let $f_\bullet:X_\bullet\rightarrow Y_\bullet$ be a morphism in $\Cal A$.
As
\[
c: (?)_J R(f_\bullet)_*\rightarrow R(f_\bullet|_J)_*(?)_J
\]
is an isomorphism, its conjugate map
\[
c': RR_J (f_\bullet|_J)^\times\rightarrow f_\bullet^\times RR_J
\]
is also an isomorphism.

Let $g_\bullet:U_\bullet\rightarrow X_\bullet$ be an image dense open
immersion of fiber type in $\Cal A$.
Let $\sigma$ be the canonical map
\[
g_\bullet^*RR_J
\specialarrow{u}
g_\bullet^*RR_J R(g_\bullet|_J)_* (g_\bullet|_J)^*
\specialarrow{\xi^{-1}}
g_\bullet^* R (g_\bullet)_*RR_J(g_\bullet|_J)^*
\specialarrow{\varepsilon}
RR_J(g_\bullet|_J)^*,
\]
where $\xi:R (g_\bullet)_*RR_J \rightarrow RR_J R(g_\bullet|_J)_*$
is the conjugate of the isomorphism
$\theta: (g_\bullet)|_J^*(?)_J\rightarrow (?)_J g_\bullet^*$.

\begin{Lemma}\label{rind-xi.thm}
Let the notation be as above.
Then $\sigma:g_\bullet^*RR_J\rightarrow RR_J(g_\bullet|_J)^*$ 
is an isomorphism of functors from $D(\Mod(X_\bullet|_J))$ to
$D(\Mod(U_\bullet))$.
\end{Lemma}

\proof As $g_\bullet^*$, $(g_\bullet|_J)^*$, $R_J$, $(g_\bullet)_*$, 
and $(g_\bullet|_J)_*$ have exact left adjoints,
it suffices to show that
\[
(?)_i\sigma: (?)_i g_\bullet^* R_J\Bbb I\rightarrow (?)_i 
R_J (g_\bullet|_J)^*\Bbb I
\]
is an isomorphism for any $K$-injective complex in $C(\Mod(X_\bullet|_J))$
and $i\in\ob(I)$.
This map agrees with
\begin{multline*}
(?)_i g_\bullet^* R_J\Bbb I\specialarrow{\theta^{-1}}
g_i^*(?)_i R_J\Bbb I
\cong
g_i^*\projlim (X_\phi)_*\Bbb I_j
\cong
\projlim g_i^*(X_\phi)_*\Bbb I_j\\
\specialarrow{\theta}
\projlim (U_\phi)_* g_j^*\Bbb I_j
\specialarrow{\theta}
\projlim(U_\phi)_* (?)_j (g_\bullet|_J)^*\Bbb I
\cong
R_J (g_\bullet|_J)^*\Bbb I,
\end{multline*}
which is obviously an isomorphism.
\qed

\paragraph Let $f_\bullet:X_\bullet\rightarrow Y_\bullet$ be a 
morphism in $\Cal F$, and $f_\bullet=p_\bullet i_\bullet$ a compactification.
We define $\bar c: f_\bullet^! RR_J\rightarrow RR_J(f_\bullet|_J)^!$ 
to be the composite
\[
f_\bullet^!RR_J\specialarrow{\rho}
i_\bullet^*p_\bullet^\times RR_J
\specialarrow{c'}
i_\bullet^*RR_J (p_\bullet|_J)^\times
\specialarrow{\sigma}
RR_J (i_\bullet|_J)^*(p_\bullet|_J)^\times
\specialarrow{\rho}
RR_J (f_\bullet|_J)^!.
\]

By Lemma~\ref{rind-xi.thm}, we have

\begin{Lemma}\label{bar-xi-rind.thm}
$\bar c: f_\bullet^! RR_J\rightarrow RR_J(f_\bullet|_J)^!$
is an isomorphism of functors from $D_{\LQco(Y_\bullet|_J)}^+
(\Mod(X_\bullet))$ 
to $D_{\LQco(X_\bullet)}^+(\Mod(X_\bullet))$.
\end{Lemma}

\section{Equivariant Grothendieck's duality}

\begin{Theorem}[Equivariant Grothendieck's duality]\label{egd.thm}
Let $I$ be a small category, 
and $f_\bullet:X_\bullet\rightarrow Y_\bullet$ a morphism in $\Cal P
(I,\Sch/\Bbb Z)$.
If $Y_\bullet$ is separated noetherian with flat arrows
and $f_\bullet$ is proper of
fiber type, then the composite
\begin{multline*}
\Theta(f_\bullet):
R(f_\bullet)_*R\uHom_{\Mod(X_\bullet)}^\bullet(\Bbb F,f_\bullet^\times
  \Bbb G)
\specialarrow{H}
R\uHom_{\Mod(Y_\bullet)}^\bullet(R(f_\bullet)_*\Bbb F,R(f_\bullet)_*
  f_\bullet^\times \Bbb G)\\
\specialarrow{\varepsilon}
R\uHom_{\Mod(Y_\bullet)}^\bullet(R(f_\bullet)_*\Bbb F,\Bbb G)
\end{multline*}
is an isomorphism for $\Bbb F\in \Cal D_{\Coh}^-(X_\bullet)$ and
$\Bbb G\in \Cal D^+(Y_\bullet)$.
\end{Theorem}

\proof It suffices to show that $(?)_i \Theta(f_\bullet)$ is 
an isomorphism for $i\in\ob(I)$.
By Lemma~\ref{composition-H.thm} and 
Lemma~\ref{Rf-twisted-restrict.thm}, {\bf 2}, it is easy to see that 
$(?)_i\Theta(f_\bullet)$  agrees with the composite
\begin{multline*}
(?)_iR(f_\bullet)_*
R\uHom_{\Mod(X_\bullet)}^\bullet(\Bbb F,f_\bullet^\times\Bbb G)
\specialarrow{c}
R(f_i)_*(?)_i
R\uHom_{\Mod(X_\bullet)}^\bullet(\Bbb F,f_\bullet^\times\Bbb G)\\
\specialarrow{H_i}
R(f_i)_*R\uHom_{\Mod(X_i)}^\bullet(\Bbb F_i,(?)_if_\bullet^\times \Bbb G)
\specialarrow{\xi}
R(f_i)_*R\uHom_{\Mod(X_i)}^\bullet(\Bbb F_i,f_i^\times \Bbb G_i)\\
\specialarrow{\Theta(f_i)}
R\uHom_{\Mod(Y_i)}^\bullet(R(f_i)_*\Bbb F_i,\Bbb G_i).
\end{multline*}
By Lemma~\ref{H_J-isom.thm}, $H_i$ is an isomorphism.
By Lemma~\ref{easy-commutativity.thm}, $\xi$ is an isomorphism.
By Grothendieck's duality for usual schemes \cite[section~6]{Neeman2},
$\Theta(f_i)$ is an isomorphism.
Hence, $(?)_i\Theta(f_\bullet)$ is an isomorphism.
\qed

\section{Morphisms of finite flat dimension}

\paragraph Let $((?)^*,(?)_*)$ be a monoidal adjoint pair over
a category $\Cal S$.
For a morphism $f:X\rightarrow Y$ in $\Cal S$, we define the projection
morphism $\Pi=\Pi(f)$ to be the composite
\[
f_*a\otimes b\specialarrow{u}f_*f^*(f_*a\otimes b)
\specialarrow{\Delta}f_*(f^*f_*a\otimes f^*b)
\specialarrow{\varepsilon}f_*(a\otimes f^*b),
\]
where $a\in M_X$ and $b\in M_Y$ (see section~\ref{monoidal.sec} for 
notation).

\begin{Lemma} Let the notation be as above, and $f:X\rightarrow Y$ and
$g:Y\rightarrow Z$ be morphisms in $\Cal S$.
For $x\in M_X$ and $z\in M_Z$, the composite
\[
(gf)_*x\otimes z \specialarrow{c}g_*(f_*x)\otimes z
\specialarrow{\Pi(g)}g_*(f_*x\otimes g^*z)
\specialarrow{\Pi(f)}g_*f_*(x\otimes f^*g^*z)
\specialarrow{\cong}(gf)_*(x\otimes(gf)^*z)
\]
agrees with $\Pi(gf)$.
\end{Lemma}

\proof Left to the reader.
\qed

\paragraph Let $I$ be a small category, $S$ a scheme, and $f_\bullet:
X_\bullet\rightarrow Y_\bullet$ a morphism in $\Cal P(I,\Sch/S)$.

\begin{Lemma}[Projection Formula] \label{projection-formula.thm}
Assume that $X_\bullet$ and $Y_\bullet$ are quasi-compact
and separated.
Then the natural map
\[
\Pi=\Pi(f_\bullet):
(Rf_\bullet)_*\Bbb F\otimes^{\bullet,L}_{\Cal O_{Y_\bullet}}\Bbb G
\rightarrow
(Rf_\bullet)_*(\Bbb F\otimes^{\bullet,L}_{\Cal O_{X_\bullet}}Lf^*_\bullet
\Bbb G)
\]
is an isomorphism for $\Bbb F\in D_{\LQco(X_\bullet)}(\Mod(X_\bullet))$
and $\Bbb G\in D_{\LQco(Y_\bullet)}(\Mod(Y_\bullet))$.
\end{Lemma}

\proof For each $i\in \ob(I)$, the composite
\lrsplit
R(f_i)_*\Bbb F_i\otimes^L_{\Cal O_{Y_i}}\Bbb G_i
\specialarrow{\Pi(f_i)}
R(f_i)_*(\Bbb F_i\otimes^L_{\Cal O_{X_i}}Lf^*_i\Bbb G_i)
\specialarrow{\theta(f_\bullet,i)}
R(f_i)_*(\Bbb F_i\otimes^L_{\Cal O_{X_i}}(?)_iLf_\bullet^*\Bbb G)\\
\specialarrow{m}
R(f_i)_*(?)_i(\Bbb F\otimes^{\bullet,L}_{\Cal O_{X_\bullet}}Lf_\bullet^*\Bbb G)
\specialarrow{c}
(?)_iR(f_\bullet)_*(\Bbb F\otimes^{\bullet,L}_{\Cal O_{X_\bullet}}
Lf_\bullet^*\Bbb G)
\endlrsplit
is an isomorphism by \cite[Proposition~5.3]{Neeman2}
and Lemma~\ref{theta-f-J.thm}, {\bf 1}.
On the other hand, it is straightforward to check that 
this composite isomorphism agrees with the composite
\begin{multline*}
R(f_i)_*\Bbb F_i\otimes^L_{\Cal O_{Y_i}}\Bbb G_i
\specialarrow{c}
(?)_i R(f_\bullet)_*\Bbb F\otimes^L_{\Cal O_{Y_i}}\Bbb G_i
\specialarrow{m}
(?)_i(R(f_\bullet)_*\Bbb F\otimes^{\bullet,L}_{\Cal O_{Y_\bullet}}\Bbb G)\\
\specialarrow{(?)_i\Pi(f_\bullet)}
(?)_iR(f_\bullet)_*(\Bbb F\otimes^{\bullet,L}_{\Cal O_{X_\bullet}}
Lf_\bullet^*\Bbb G).
\end{multline*}
It follows that $(?)_i\Pi(f_\bullet)$ is an isomorphism for any $i\in\ob(I)$.
Hence, $\Pi(f_\bullet)$ is an isomorphism.
\qed

\paragraph\label{Phi-def.par}
Let $f_\bullet: X_\bullet\rightarrow Y_\bullet$ be a 
morphism in $\Cal P(I,\Sch/S)$, and assume that both $X_\bullet$
and $Y_\bullet$ are quasi-compact separated.
Define $\Phi=\Phi(f_\bullet)$ to be the composite
\begin{multline*}
f_\bullet^\times\Bbb F\otimes^{\bullet,L}_{\Cal O_{X_\bullet}}
Lf_\bullet^*\Bbb G
\specialarrow u
f^\times_\bullet R(f_\bullet)_*(f_\bullet^\times\Bbb F
\otimes^{\bullet,L}_{\Cal O_{X_\bullet}}Lf_\bullet^*\Bbb G)\\
\specialarrow{\Pi(f_\bullet)^{-1}}
f^\times_\bullet(R(f_\bullet)_*f_\bullet^\times\Bbb F
\otimes^{\bullet,L}_{\Cal O_{Y_\bullet}}\Bbb G)
\specialarrow{\varepsilon}
f^\times_\bullet(\Bbb F
\otimes^{\bullet,L}_{\Cal O_{Y_\bullet}}\Bbb G),
\end{multline*}
where $\Bbb F,\Bbb G\in D_{\LQco(Y_\bullet)}(\Mod(Y_\bullet))$.

Utilizing the commutativity 
as in the proof of Lemma~\ref{projection-formula.thm}
and Lemma~\ref{Rf-twisted-restrict.thm}, it is not so difficult to
show the following.

\begin{Lemma}\label{Phi-restrict.thm}
Let $f_\bullet:X_\bullet\rightarrow Y_\bullet$
be as in {\rm (\ref{Phi-def.par})}.
For a subcategory $J$ of $I$, the composite
\begin{multline*}
(?)_J f_\bullet^\times \Bbb F\otimes_{\Cal O_{X_\bullet|_J}}
^{\bullet,L}(?)_J Lf_\bullet^*\Bbb G
\specialarrow{\xi\otimes \theta^{-1}}
(f_\bullet|_J)^\times\Bbb F_J\otimes_{\Cal O_{X_\bullet|_J}}
^{\bullet,L}L(f_\bullet|_J)^*\Bbb G_J\\
\specialarrow{\Phi(f_\bullet|_J)}
(f_\bullet|_J)^\times(\Bbb F_J\otimes_{\Cal O_{Y_\bullet|_J}} \Bbb G_J)
\specialarrow{m}
(f_\bullet|_J)^\times(?)_J(\Bbb F\otimes_{\Cal O_{Y_\bullet}} \Bbb G)
\end{multline*}
agrees with
\begin{multline*}
(?)_J f_\bullet^\times \Bbb F\otimes_{\Cal O_{X_\bullet|_J}}
^{\bullet,L}(?)_J Lf_\bullet^*\Bbb G
\specialarrow m
(?)_J(f_\bullet^\times \Bbb F\otimes_{\Cal O_{X_\bullet}}
^{\bullet,L}Lf_\bullet^*\Bbb G)\\
\specialarrow{\Phi(f_\bullet)}
(?)_J f_\bullet^\times 
(\Bbb F\otimes^{\bullet,L}_{\Cal O_{Y_\bullet}}
\Bbb G)
\specialarrow{\xi}
(f_\bullet|_J)^\times (?)_J
(\Bbb F\otimes^{\bullet,L}_{\Cal O_{Y_\bullet}}
\Bbb G).
\end{multline*}
\end{Lemma}

\begin{Lemma}\label{fgh.thm}
Let $f_\bullet:X_\bullet\rightarrow Y_\bullet$
be as in {\rm (\ref{Phi-def.par})}.
The composite
\begin{multline*}
f_\bullet^\times \Bbb F\otimes^{\bullet,L}_{\Cal O_{X_\bullet}}
L  f_\bullet^*(\Bbb G\otimes^{\bullet,L}_{\Cal O_{Y_\bullet}}\Bbb H)
\specialarrow{\Delta}
f_\bullet^\times \Bbb F\otimes^{\bullet,L}_{\Cal O_{X_\bullet}}
L  f_\bullet^*\Bbb G\otimes^{\bullet,L}_{\Cal O_{X_\bullet}}
L  f_\bullet^*\Bbb H\\
\specialarrow{\Phi}
f_\bullet^\times(\Bbb F\otimes^{\bullet,L}_{\Cal O_{Y_\bullet}}\Bbb G)
  \otimes^{\bullet,L}_{\Cal O_{X_\bullet}}L f_\bullet^*\Bbb H
\specialarrow{\Phi}
f_\bullet^\times(\Bbb F\otimes^{\bullet,L}_{\Cal O_{Y_\bullet}}\Bbb G
  \otimes^{\bullet,L}_{\Cal O_{Y_\bullet}}\Bbb H)
\end{multline*}
agrees with $\Phi$.
\end{Lemma}

\begin{Lemma}\label{Phi-base-change.thm}
Let $S$, $I$ and $\sigma$ be as in {\rm (\ref{define-zeta.par})}.
For $\Bbb F\in D_{\LQco(Y_\bullet)}(\Mod(Y_\bullet))$, the composite
\begin{multline*}
(g'_\bullet)^*f_\bullet^\times \Bbb F\otimes_{\Cal O_{X_\bullet'}}
^{\bullet,L}(g'_\bullet)^*L f_\bullet^*\Bbb G
\specialarrow{\Delta^{-1}}
(g'_\bullet)^*(f_\bullet^\times \Bbb F\otimes_{\Cal O_{X_\bullet}}
^{\bullet,L}L f_\bullet^*\Bbb G)\\
\specialarrow{\Phi}
(g'_\bullet)^* f_\bullet^\times 
(\Bbb F\otimes^{\bullet,L}_{\Cal O_{Y_\bullet}}\Bbb G)
\specialarrow{\zeta(\sigma)}
(f'_\bullet)^\times g_\bullet^*
(\Bbb F\otimes^{\bullet,L}_{\Cal O_{Y_\bullet}}\Bbb G)
\end{multline*}
agrees with
\begin{multline*}
(g'_\bullet)^*f_\bullet^\times \Bbb F\otimes_{\Cal O_{X_\bullet'}}
^{\bullet,L}(g'_\bullet)^*L f_\bullet^*\Bbb G
\specialarrow{\zeta(\sigma)\otimes d}
(f_\bullet')^\times g_\bullet^*\Bbb F\otimes_{\Cal O_{X_\bullet'}}
^{\bullet,L}L (f'_\bullet)^*g_\bullet^*\Bbb G\\
\specialarrow{\Phi}
(f_\bullet')^\times(g_\bullet^*\Bbb F
  \otimes_{\Cal O_{Y_\bullet'}}^{\bullet,L}g_\bullet^*\Bbb G)
\specialarrow{\Delta^{-1}}
(f_\bullet')^\times g_\bullet^*(\Bbb F\otimes_{\Cal O_{Y_\bullet}}^{\bullet,L}
\Bbb G).
\end{multline*}
\end{Lemma}

\begin{Lemma}\label{Phi-composition.thm}
Let $f_\bullet:X_\bullet\rightarrow Y_\bullet$
and $g_\bullet:Y_\bullet\rightarrow Z_\bullet$ be
morphisms in $\Cal P(I,\Sch/S)$.
Assume that $X_\bullet$, $Y_\bullet$ and $Z_\bullet$
are quasi-compact separated.
Then the composite
\begin{multline*}
(g_\bullet f_\bullet)^\times\Bbb F\otimes^{\bullet,L}_{\Cal O_{
X_\bullet}}(g_\bullet f_\bullet)^*\Bbb G
\cong
f_\bullet^\times g_\bullet^\times \Bbb F\otimes^{\bullet,L}_{\Cal O_{
X_\bullet}}f_\bullet^*g_\bullet^*\Bbb G
\specialarrow{\Phi(f_\bullet)}
f_\bullet^\times(g_\bullet^\times\Bbb F\otimes^{\bullet,L}_{\Cal O_{
Y_\bullet}}g_\bullet^*\Bbb G)\\
\specialarrow{\Phi(g_\bullet)}
f_\bullet^\times g_\bullet^\times(\Bbb F\otimes^{\bullet,L}_{\Cal O_{
Z_\bullet}}\Bbb G)
\cong
(g_\bullet f_\bullet)^\times(\Bbb F\otimes^{\bullet,L}_{\Cal O_{
Z_\bullet}}\Bbb G)
\end{multline*}
agrees with $\Phi(g_\bullet f_\bullet)$.
\end{Lemma}

The proof of the lemmas above are left to the reader.

\paragraph Let the notation be as in Theorem~\ref{finite-data.thm}.
Let $i_\bullet:U_\bullet\rightarrow X_\bullet$ be a morphism in $\Cal I$, 
and $p_\bullet:X_\bullet\rightarrow Y_\bullet$ a morphism in $\Cal P$.

We define $\bar\Phi=\bar\Phi(p_\bullet,i_\bullet)$ to be the composite
\begin{multline*}
\bar\Phi:
i_\bullet^*p_\bullet^\times
\Bbb F
\otimes^{\bullet,L}_{\Cal O_{U_\bullet}}
L (p_\bullet i_\bullet)^*\Bbb G
\specialarrow{d^{-1}}
i_\bullet^*p_\bullet^\times
\Bbb F
\otimes^{\bullet,L}_{\Cal O_{U_\bullet}}
i_\bullet^*L p_\bullet^*\Bbb G\\
\specialarrow{\Delta^{-1}}
i_\bullet^*(p_\bullet^\times
\Bbb F\otimes^{\bullet,L}_{\Cal O_{X_\bullet}}
L p_\bullet^*\Bbb G)
\specialarrow{i_\bullet^*\Phi}
i_\bullet^*p_\bullet^\times(\Bbb F\otimes^{\bullet,L}_{\Cal O_{Y_\bullet}}
\Bbb G).
\end{multline*}

\begin{Lemma}\label{Phi-independent.thm}
Let $f_\bullet$ be a morphism in $\Cal F$, and
$f_\bullet=p_\bullet i_\bullet=q_\bullet j_\bullet$ an independence square.
Then the composite
\[
i_\bullet^*p_\bullet^\times
\Bbb F
\otimes^{\bullet,L}_{\Cal O_{U_\bullet}}
L (p_\bullet i_\bullet)^*
\Bbb G
\specialarrow{\bar\Phi}
i_\bullet^*p_\bullet^\times(\Bbb F\otimes^{\bullet,L}_{\Cal O_{Y_\bullet}}
\Bbb G)
\specialarrow{\rho(p_\bullet i_\bullet=q_\bullet j_\bullet)}
j_\bullet^*q_\bullet^\times(\Bbb F\otimes^{\bullet,L}_{\Cal O_{Y_\bullet}}
\Bbb G)
\]
agrees with
\[
i_\bullet^*p_\bullet^\times
\Bbb F
\otimes^{\bullet,L}_{\Cal O_{U_\bullet}}
L (p_\bullet i_\bullet)^*
\Bbb G
\specialarrow{\rho\otimes 1}
j_\bullet^*q_\bullet^\times
\Bbb F
\otimes^{\bullet,L}_{\Cal O_{U_\bullet}}
L (q_\bullet j_\bullet)^*
\Bbb G
\specialarrow{\bar\Phi}
j_\bullet^*q_\bullet^\times(\Bbb F\otimes^{\bullet,L}_{\Cal O_{Y_\bullet}}
\Bbb G).
\]
\end{Lemma}

\proof As $\rho$ is constructed from $\zeta$ and $d$ by definition,
the assertion follows easily from
Lemma~\ref{Phi-base-change.thm} and Lemma~\ref{Phi-composition.thm}.
\qed

\paragraph
Let $f_\bullet:X_\bullet\rightarrow Y_\bullet$ be a morphism in $\Cal F$.
We define $\bar \Phi(f_\bullet)$ to be
$\bar\Phi(p_\bullet,i_\bullet)$, where $f_\bullet=p_\bullet i_\bullet$ 
is the (fixed) compactification of $f_\bullet$.

By Lemma~\ref{Phi-independent.thm}, 
$\bar\Phi(f_\bullet)$ is an isomorphism if and only if 
there exists some compactification $f_\bullet=q_\bullet j_\bullet$ such
that $\bar \Phi(q_\bullet,j_\bullet)$ is an isomorphism.

\begin{Lemma}\label{bar-Phi-restrict.thm}
Let the notation be as in {\rm Theorem~\ref{finite-data.thm}},
and $f_\bullet:X_\bullet\rightarrow Y_\bullet$ and
$g_\bullet: Y_\bullet\rightarrow Z_\bullet$ morphisms in $\Cal F$.
Then the composite
\begin{multline*}
(g_\bullet f_\bullet)^!\Bbb F\otimes^{\bullet,L}_{\Cal O_{
X_\bullet}} L (g_\bullet f_\bullet)^*\Bbb G
\cong
f_\bullet^! g_\bullet^! \Bbb F\otimes^{\bullet,L}_{\Cal O_{
X_\bullet}}L f_\bullet^*L g_\bullet^*\Bbb G
\specialarrow{\bar\Phi(f_\bullet)}
f_\bullet^!(g_\bullet^!\Bbb F\otimes^{\bullet,L}_{\Cal O_{
Y_\bullet}}L g_\bullet^*\Bbb G)\\
\specialarrow{\bar\Phi(g_\bullet)}
f_\bullet^! g_\bullet^!(\Bbb F\otimes^{\bullet,L}_{\Cal O_{
Z_\bullet}}\Bbb G)
\cong
(g_\bullet f_\bullet)^!(\Bbb F\otimes^{\bullet,L}_{\Cal O_{
Z_\bullet}}\Bbb G)
\end{multline*}
agrees with $\bar\Phi(g_\bullet f_\bullet)$.
\end{Lemma}

\begin{Theorem}\label{ffd.thm}
Let the notation be as in {\rm Theorem~\ref{finite-data.thm}}, and
$f_\bullet:X_\bullet\rightarrow Y_\bullet$ a morphism in $\Cal F$.
If $f_\bullet$ is of finite flat dimension, then
\[
\bar\Phi(f_\bullet):
f_\bullet^!\Bbb F\otimes^{\bullet,L}_{\Cal O_{Y_\bullet}}L f_\bullet^*\Bbb G
\rightarrow
f_\bullet^!(\Bbb F\otimes^{\bullet,L}_{\Cal O_{X_\bullet}}\Bbb G)
\]
is an isomorphism for $\Bbb F,\Bbb G\in D^+_{\LQco(Y_\bullet)}(\Mod(Y_\bullet))
$.
\end{Theorem}

\proof Let $f_\bullet=p_\bullet i_\bullet$ be a compactification of $f_\bullet
$.
It suffices to show that $\bar\Phi(p_\bullet,i_\bullet)$ is an isomorphism.
In view of Proposition~\ref{twisted-restrict.thm} and 
Lemma~\ref{Phi-restrict.thm}, it suffices to show that
\[
\bar\Phi(f_j): i_j^* p_j^\times\Bbb F_j\otimes^{\bullet,L}_{\Cal O_{X_j}}
L f_j^*\Bbb G_i\rightarrow
i_j^*p_j^\times(\Bbb F_j\otimes^{\bullet,L}_{\Cal O_{Y_j}}\Bbb G_j)
\]
is an isomorphism for any $j\in\ob(I)$.
So we may assume that $I=j$.

By the flat base change theorem and Lemma~\ref{Phi-base-change.thm}, 
the question is local on $Y_j$.
Clearly, the question is local on $U_j$.
Hence we may assume that $Y_j$ and $U_j$ are affine.
Set $f=f_j$, $Y=Y_j$, and $U=U_j$.
Note that $f$ is a closed immersion defined by an ideal of finite
projective dimension, followed by an affine $n$-space.

By Lemma~\ref{bar-Phi-restrict.thm}, it suffices to prove that
$\bar\Phi(f)$ is an isomorphism if $f$ is a closed immersion defined
by an ideal of finite projective dimension or an affine $n$-space.

By Lemma~\ref{fgh.thm}, we may assume that $\Bbb F=\Cal O_{Y}$.
Both cases are proved easily, using \cite[Theorem~5.4]{Neeman2}.
\qed

\section{Finite morphisms of fiber type}

\paragraph\label{finite-setting.par}
Let $I$ be a small category, $S$ a scheme, and
$f_\bullet:X_\bullet\rightarrow Y_\bullet$ a morphism in $\Cal P(I,\Sch/S)$.
Let $Z$ denote the ringed site $(\Zar(Y_\bullet),(f_\bullet)_*
(\Cal O_{X_\bullet}))$.
Assume that $Y_\bullet$ is locally noetherian.
There are obvious admissible ringed continuous functors 
$i:\Zar(Y_\bullet)\rightarrow Z$ and
$g:Z\rightarrow \Zar(X_\bullet)$ such that 
$gi=f_\bullet^{-1}$.
It is easy to see that $g_{\#}:\Mod(Z)\rightarrow\Mod(X_\bullet)$ is
an exact functor.

\begin{Lemma} If $f_\bullet$ is affine, then the counit 
\[
\varepsilon: g_{\#} R g^{\#}\Bbb F\rightarrow \Bbb F
\]
is an isomorphism for $\Bbb F\in D^+_{\LQco(X_\bullet)}(\Mod(X_\bullet))$.
\end{Lemma}

\proof The construction of $\varepsilon$ is compatible with restrictions.
So we may assume that $f_\bullet=f:X\rightarrow Y$ is an affine morphism of 
single schemes.
Further, the question is local on $Y$, and hence we may assume that
$Y=\Spec A$ is affine.
As $f$ is affine, $X=\Spec B$ is affine.
By the way-out lemma, we may assume that $\Bbb F$ is a single quasi-coherent
sheaf, say $\tilde M$, where $M$ is a $B$-module.
As $R^if_*\tilde M=0$ for $i>0$, it suffices to show that
$\varepsilon:g_{\#}g^{\#}\tilde M\rightarrow \tilde M$ is an isomorphism.
As $g_{\#}g^{\#}$ is exact on $\tilde M$ and preserves direct sums, 
we may assume that $M=B$, which
case is trivial.
\qed

\paragraph 
Let $I$, $S$, $f_\bullet:X_\bullet\rightarrow Y_\bullet$,
$Z$, $g$, and $i$ be
as in (\ref{finite-setting.par}).
Assume that $f_\bullet$ is finite of fiber type.

We say that an $\Cal O_Z$-module $\Cal M$ is locally quasi-coherent
(resp.\ quasi-coherent, coherent) if $i^{\#}\Cal M$ is.
The corresponding full subcategory of $\Mod(Z)$ is denoted by
$\LQco(Z)$ (resp.\ $\Qco(Z)$, $\Coh(Z)$).

\begin{Lemma} Let the notation be as above.
Then an $\Cal O_Z$-module $\Cal M$ is locally quasi-coherent
if and only if for any $i\in\ob(I)$ and any affine open subscheme
$U$ of $Y_i$, there exists an exact sequence of $((\Cal O_Z)_i)|_U$-modules
\[
\Cal (((\Cal O_Z)_i)|_U)^{(T)}\rightarrow
\Cal (((\Cal O_Z)_i)|_U)^{(\Sigma)}\rightarrow
\Cal M_i|_{U}\rightarrow 
0.
\]
\end{Lemma}

\proof As we assume that $f_\bullet$ is finite of fiber type, 
$\Cal O_Z$ is coherent.
Hence the existence of such exact sequences implies that $\Cal M$ is
locally quasi-coherent.

We prove the converse.
Let $i\in\ob(I)$ and $U$ an affine open subset of $Y_i$.
Set $C:=\Gamma(U,(\Cal O_Z)_i)=\Gamma(f_i^{-1}(U),\Cal O_{X_i})$ 
and $M:=\Gamma(U,\Cal M_i)$.
There is a canonical map 
$(g_i|_{f_i^{-1}(U)})^{\#}(\tilde M)\rightarrow
\Cal M_i|_U$,
where $\tilde M$ is the quasi-coherent sheaf over $\Spec C\subset X_i$
associated with the $C$-module $M$.
When we apply $i^{\#}$ to this map, we get $\tilde M_0
\rightarrow ((i^{\#}\Cal M)_i)_U$, where $M_0$ is $M$ viewed as a 
$\Gamma(U,\Cal O_{Y_i})$-module.
This is an isomorphism, since $((i^{\#}\Cal M)_i)_U$ is quasi-coherent and
$U$ is affine.
As $i^{\#}$ is faithful and exact, 
we have
$(g_i|_{f_i^{-1}(U)})^{\#}(\tilde M)\cong
\Cal M_i|_U$.

Take an exact sequence of the form
\[
C^{(T)}\rightarrow C^{(\Sigma)}\rightarrow M\rightarrow 0.
\]
Applying 
the exact functor $(g_i|_{f_i^{-1}(U)})^{\#}\circ \tilde{\mathord ?}$,
we get an exact sequence of the desired type.
\qed

\begin{Corollary}\label{g-qco.thm}
Under the same assumption as in the lemma,
the functor $g_{\#}$ preserves local quasi-coherence.
\end{Corollary}

\proof As $g_{\#}$ is compatible with restrictions, we may assume that
$I$ consists of one object and one morphism.
Further, as the question is local, we may assume that $Y_\bullet=Y$ is
an affine scheme.
By the lemma, it suffices to show that $g_{\#}\Cal O_Z$ is quasi-coherent,
since $g_{\#}$ is exact and preserves direct sums.
As $g_{\#}\Cal O_Z=g_{\#}g^{\#}\Cal O_X\cong \Cal O_X$, we are done.
\qed

\begin{Lemma}
Let the notation be as above.
The unit of adjunction $u:\Bbb F\rightarrow Rg^{\#}g_{\#}\Bbb F$ is
an isomorphism for $\Bbb F\in D^+_{\LQco(Z)}(\Mod(Z))$.
\end{Lemma}

\proof We may assume that $I$ consists of one object and one morphism,
and $Y_\bullet=Y$ is affine.
By the way-out lemma, we may assume that $\Bbb F$ is a single 
quasi-coherent sheaf.
By Corollary~\ref{g-qco.thm}, we may assume that $\Bbb F=\Cal O_Z
=g^{\#}\Cal O_X$ and it suffices to prove that
$ug^{\#}:g^{\#}\Cal O_X\rightarrow g^{\#}g_{\#}g^{\#}\Cal O_X$ is
an isomorphism.
As $\id=(g^{\#}\varepsilon)(ug^{\#})$ and $\varepsilon$ is an isomorphism,
we are done.
\qed

For $\Cal N\in\Mod(Y_\bullet)$ and $\Cal M\in\Mod(Z)$, 
the sheaf 
$\uHom_{\Cal O_{Y_\bullet}}(\Cal M,\Cal N)$ on $Y_\bullet$ 
has a structure of
$\Cal O_Z$-module, and it belongs to $\Mod(Z)$.
There is an obvious isomorphism of functors
\[
\kappa: i^{\#}\uHom_{\Cal O_{Y_\bullet}}(\Cal M,\Cal N)
\cong \uHom_{\Cal O_{Y_\bullet}}(i^{\#}\Cal M,\Cal N).
\]
For $\Cal M,\Cal M'\in\Mod(Z)$, there is a natural map
\[
\upsilon: \uHom_{\Mod(Z)}(\Cal M,\Cal M')\rightarrow
\uHom_{\Cal O_{Y_\bullet}}(\Cal M,i^{\#}\Cal M').
\]
Note that the composite
\[
i^{\#}\uHom_{\Mod(Z)}(\Cal M,\Cal M')\specialarrow{i^{\#}\upsilon}
i^{\#}\uHom_{\Cal O_{Y_\bullet}}(\Cal M,i^{\#}\Cal M')
\specialarrow{\kappa}
\uHom_{\Cal O_{Y_\bullet}}(i^{\#}\Cal M,i^{\#}\Cal M')
\]
agrees with $H$.

\paragraph
Let $I$, $S$, $f_\bullet:X_\bullet\rightarrow Y_\bullet$,
$Z$, $g$, and $i$ be
as in (\ref{finite-setting.par}).
Assume that $f_\bullet$ is finite of fiber type,
and $Y_\bullet$ has flat arrows.
Define $f_\bullet^\natural:\Cal D^+(Y_\bullet)
\rightarrow D(\Mod(X_\bullet))$ by
\[
f_\bullet^\natural(\Bbb F):=g_{\#}R\uHom^\bullet_{\Cal O_{Y_\bullet}}
(\Cal O_Z,\Bbb F).
\]
As $f_\bullet$ is finite of fiber type, $\Cal O_Z$ is coherent.
By Lemma~\ref{ext-equivariant.thm},
$i^{\#}R\uHom^\bullet_{\Cal O_{Y_\bullet}}(\Cal O_Z,\Bbb F)
\in \Cal D^+(Y_\bullet)$.
It follows that $f_\bullet^\natural(\Bbb F)\in \Cal D^+(X_\bullet)$,
and $f_\bullet^\natural$ is a functor from $\Cal D^+(Y_\bullet)$ to
$\Cal D^+(X_\bullet)$.

Define $\varepsilon: R(f_\bullet)_* f_\bullet^\natural\rightarrow
\Id_{\Cal D^+(Y_\bullet)}$ by
\lrsplit
R(f_\bullet)_*f_\bullet^\natural\Bbb F=
i^{\#} Rg^{\#} g_{\#}R\uHom^\bullet_{\Cal O_{Y_\bullet}}(\Cal O_Z,\Bbb F)
\specialarrow{u^{-1}}
i^{\#}R\uHom^\bullet_{\Cal O_{Y_\bullet}}(\Cal O_Z,\Bbb F)\\
\specialarrow{\kappa}
R\uHom^\bullet_{\Cal O_{Y_\bullet}}(i^{\#}\Cal O_Z,\Bbb F)
=
R\uHom^\bullet_{\Cal O_{Y_\bullet}}((f_\bullet)_*\Cal O_{X_\bullet},\Bbb F)
\specialarrow{\eta}
R\uHom^\bullet_{\Cal O_{Y_\bullet}}(\Cal O_{Y_\bullet},\Bbb F)
\cong
\Bbb F.
\endlrsplit
Define $u: \Id_{\Cal D^+(X_\bullet)}\rightarrow f_\bullet^\natural 
R(f_\bullet)_*$ by
\lrsplit
\Bbb F\cong R\uHom^\bullet_{\Cal O_{X_\bullet}}(\Cal O_{X_\bullet},\Bbb F)
\specialarrow{\varepsilon^{-1}}g_{\#} R g^{\#}
R\uHom^\bullet_{\Cal O_{X_\bullet}}(\Cal O_{X_\bullet},\Bbb F)\\
\specialarrow{H}
g_{\#}R\uHom^\bullet_{\Mod(Z)}(g^{\#}\Cal O_{X_\bullet},R g^{\#}\Bbb F)
\specialarrow{\upsilon}
g_{\#}R\uHom^\bullet_{\Cal O_{Y_\bullet}}
(\Cal O_Z,R(f_\bullet)_*\Bbb F)
=
f_\bullet^\natural R(f_\bullet)_*\Bbb F.
\endlrsplit

\begin{Theorem}\label{finite-duality.thm}
Let the notation be as above.
Then $f_\bullet^\natural$ is right adjoint to $R(f_\bullet)_*$, and
$\varepsilon$ and $u$ defined above are the counit and unit of adjunction,
respectively.
In particular, if moreover $Y_\bullet$ is quasi-compact separated,
then $f_\bullet^\natural$ is isomorphic to $f_\bullet^\times$.
\end{Theorem}

\proof It is easy to see that the composite
$$\displaylines{
R\uHom_{\Cal O_{Y_\bullet}}^\bullet(\Cal O_{Z},\Bbb F)
\cong
R\uHom_{\Mod(Z)}^\bullet(\Cal O_Z,
  R\uHom_{\Cal O_{Y_\bullet}}^\bullet(\Cal O_{Z},\Bbb F))\hfill\cr
\specialarrow{\upsilon}
R\uHom_{\Cal O_{Y_\bullet}}^\bullet(\Cal O_Z,
i^{\#}  R\uHom_{\Cal O_{Y_\bullet}}^\bullet(\Cal O_{Z},\Bbb F))
\specialarrow{\kappa}
R\uHom_{\Cal O_{Y_\bullet}}^\bullet(\Cal O_Z,
R\uHom_{\Cal O_{Y_\bullet}}^\bullet(i^{\#}\Cal O_{Z},\Bbb F))\hfill\cr
\hfill
\specialarrow{\eta}
R\uHom_{\Cal O_{Y_\bullet}}^\bullet(\Cal O_Z,
R\uHom_{\Cal O_{Y_\bullet}}^\bullet(\Cal O_{Y_\bullet},\Bbb F))
\cong
R\uHom_{\Cal O_{Y_\bullet}}^\bullet(\Cal O_{Z},\Bbb F)
\crcr}
$$
is the identity.
Utilizing this, $(f_\bullet^\natural \varepsilon)\circ(u f_\bullet^\natural)
=\id$ and $ (\varepsilon R(f_\bullet)_*)\circ(R(f_\bullet)_*u)=\id$ are
checked directly.

The last assertion is obvious, as the right adjoint functor is unique.
\qed

\section{Regular embeddings and Smooth morphisms of fiber type}

\paragraph 
Let $I$ be a small category, $S$ a scheme, and 
$X_\bullet\in\Cal P(I,\Sch/S)$.
An $\Cal O_{X_\bullet}$-module sheaf $\Cal M\in\Mod(X_\bullet)$ is
said to be locally free (resp.\ invertible) if $\Cal M$ is coherent 
and $\Cal M_i$ is locally free (resp.\ invertible) for any $i\in\ob(I)$.
A perfect complex of $X_\bullet$ is a bounded complex in $C^b(\Mod(X_\bullet))
$ each of whose terms is locally free.

A point of $X_\bullet$ is a pair $(i,x)$ such that $i\in\ob(I)$ and
$x\in X_i$.
A stalk of a sheaf $\Cal M\in\AB(X_\bullet)$ at the point $(i,x)$ is
defined to be $(\Cal M_i)_x$, and we denote it by $\Cal M_{i,x}$.

A connected component of $X_\bullet$ is an equivalence class
with respect to the equivalence relation of the set of points of $X_\bullet$
generated by the following relations.
\begin{description}
\item[1] $(i,x)$ and $(i',x')$ are equivalent if $i=i'$ and $x$ and $x'$
belong to the same connected component of $X_i$.
\item[2] $(i,x)$ and $(i',x')$ are equivalent if there exists some
$\phi:i\rightarrow i'$ such that $X_\phi(x')=x$.
\end{description}
We say that $X_\bullet$ is $d$-connected if $X_\bullet$ consists of
one connected component (note that the word \lq connected' is reserved
for componentwise connectedness).
If $X_\bullet$ is locally noetherian, then a connected component of 
$X_\bullet$ is a closed open subdiagram of schemes in a natural way.
If this is the case, the rank function $(i,x)\mapsto \rank_{\Cal O_{X_i,x}}
\Cal F_{i,x}$ 
of a locally free sheaf $\Cal F$ is
constant on a connected component of $X_\bullet$.

\begin{Lemma} \label{perfect-dual.thm}
Let $I$ be a small category, $S$ a scheme, and 
$X_\bullet\in\Cal P(I,\Sch/S)$.
Let $\Bbb F$ be
 a perfect complex of $X_\bullet$.
Then we have
\begin{description}
\item[1] The canonical map
\[
H_J: (?)_J R\uHom_{\Cal O_{X_\bullet}}^\bullet(\Bbb F,\Bbb G)
\rightarrow R\uHom_{\Cal O_{X_\bullet|_J}}^\bullet(\Bbb F_J,\Bbb G_J)
\]
is an isomorphism for $\Bbb G\in D(\Mod(X_\bullet))$.
\item[2] The canonical map
\[
R\uHom_{\Cal O_{X_\bullet}}^\bullet
(\Bbb F,\Bbb G)\otimes^{\bullet,L}_{\Cal O_{X_\bullet}}\Bbb H
\rightarrow
R\uHom_{\Cal O_{X_\bullet}}^\bullet
(\Bbb F,\Bbb G\otimes^{\bullet,L}_{\Cal O_{X_\bullet}}\Bbb H)
\]
is an isomorphism for $\Bbb G,\Bbb H\in D(\Mod(X_\bullet))$.
\end{description}
\end{Lemma}

\proof {\bf 1} 
It suffices to show that 
\[
H_J: (?)_J\uHom_{\Cal O_{X_\bullet}}^\bullet(\Bbb F,\Bbb G)
\rightarrow \uHom_{\Cal O_{X_\bullet|_J}}^\bullet(\Bbb F_J,\Bbb G_J)
\]
is an isomorphism of complexes if $\Bbb G$ is a $K$-injective
complex in $C(\Mod(X_\bullet))$, since $\Bbb F_J$ is $K$-flat and
$\Bbb G_J$ is weakly $K$-injective.
The assertion follows immediately by Lemma~\ref{uhom2.thm}.

{\bf 2} We may assume that $\Bbb F$ is a single locally free sheaf.
By {\bf 1}, we may assume that $X=X_\bullet$ is a single scheme.
We may assume that $X$ is affine and $\Bbb F=\Cal O_X^n$ for some $n$.
This case is trivial.
\qed

\paragraph 
Let $I$ be a small category, $S$ a scheme, and 
$X_\bullet\in\Cal P(I,\Sch/S)$.
An $\Cal O_{X_\bullet}$-module $\Cal M$ is said to be locally of finite
projective dimension if $\Cal M_{i,x}$ is of finite projective dimension
as an $\Cal O_{X_i,x}$-module for any point $(i,x)$ of $X_\bullet$.
We say that $\Cal M$ has finite projective dimension if there exists 
some non-negative integer $d$ such that $\pdim_{\Cal O_{X_i,x}}\Cal M_{i,x}
\leq d$ for any point $(i,x)$ of $X_\bullet$.

\begin{Lemma} \label{fpd-dual.thm}
Let $I$ be a small category, $S$ a scheme, and 
$X_\bullet\in\Cal P(I,\Sch/S)$.
Assume that $X_\bullet$ has flat arrows and is locally noetherian.
If $\Bbb F$ is a complex in $C(\Mod(X_\bullet))$
with bounded coherent cohomology groups which 
have finite projective dimension, then
the canonical map
\[
R\uHom_{\Cal O_{X_\bullet}}^\bullet
(\Bbb F,\Bbb G)\otimes^{\bullet,L}_{\Cal O_{X_\bullet}}\Bbb H
\rightarrow
R\uHom_{\Cal O_{X_\bullet}}^\bullet
(\Bbb F,\Bbb G\otimes^{\bullet,L}_{\Cal O_{X_\bullet}}\Bbb H)
\]
is an isomorphism for $\Bbb G,\Bbb H\in D(\Mod(X_\bullet))$.
\end{Lemma}

\proof We may assume that $\Bbb G=\Cal O_{X_\bullet}$.
By the way-out lemma, we may assume that $\Bbb F$ is a single coherent sheaf
which has finite projective dimension, say $d$.
By Lemma~\ref{H_J-isom2.thm},
it is easy to see that $\uExt_{\Cal O_{X_\bullet}}^i(\Bbb F,G)=0$
($i>d$) for $G\in\Mod(X_\bullet)$.
In particular, $R\uHom_{\Cal O_{X_\bullet}}^\bullet(\Bbb F,?)$ is
way-out in both directions.
On the other hand, as $R\uHom_{\Cal O_{X_\bullet}}(\Bbb F,\Cal O_{X_\bullet})$
has finite flat dimension, and hence
$R\uHom_{\Cal O_{X_\bullet}}(\Bbb F,\Cal O_{X_\bullet})\otimes^{\bullet,L}
_{\Cal O_{X_\bullet}}?$ is also way-out in both directions.
By the way-out lemma, we may assume that $\Bbb H$ is a single $\Cal 
O_{X_\bullet}$-module.
By Lemma~\ref{H_J-isom.thm}, 
we may assume that 
$X=X_\bullet$ is a single scheme.
The question is local, and we may assume
that $X=\Spec A$ is affine.
Moreover, we may assume
that $\Bbb F$ is a complex of sheaves associated with a finite projective 
resolution of a single finitely generated module.
As $\Bbb F$ is perfect, the result follows from
Lemma~\ref{perfect-dual.thm}.
\qed

\paragraph Let $S$, $I$, and $X_\bullet$ be as above.
For a locally free sheaf $\Cal F$ over $X_\bullet$, we denote
$\uHom_{\Cal O_{X_\bullet}}(\Cal F,\Cal O_{X_\bullet})$ by $\Cal F^\vee$.
It is easy to see that $\Cal F^\vee$ is again locally free.
If $\Cal L$ is an invertible sheaf, then
\[
\Cal O_{X_\bullet}\specialarrow{\trace}
\uHom_{\Cal O_{X_\bullet}}(\Cal L,\Cal L)
\cong \Cal L^\vee\otimes_{\Cal O_{X_\bullet}}\Cal L
\]
are isomorphisms.

\paragraph\label{closed-immersion-ideal.par}
Let $I$ be a small category, $S$ a scheme, and 
$i_\bullet:Y_\bullet\rightarrow X_\bullet$ 
a closed immersion in $\Cal P(I,\Sch/S)$.
Then the canonical map $\eta: \Cal O_{X_\bullet}\rightarrow (i_\bullet)_*
\Cal O_{Y_\bullet}$ is an epimorphism in $\LQco(X_\bullet)$.
Set $\Cal I:=\Ker\eta$.
Then $\Cal I$ is a locally quasi-coherent ideal of $\Cal O_{X_\bullet}$.
Conversely, if $\Cal I$ is a given locally 
quasi-coherent ideal of $\Cal O_{X_\bullet}$, then
\[
Y_\bullet:=\uSpec_\bullet \Cal O_{X_\bullet}/\Cal I\specialarrow{i_\bullet}
X_\bullet
\]
is defined appropriately, and $i_\bullet$ is a closed immersion.
Thus the isomorphism classes of closed immersions to $X_\bullet$ 
in the category $\Cal P(I,\Sch/S)/X_\bullet$ and locally quasi-coherent
ideals of $\Cal O_{X_\bullet}$ are in one-to-one correspondence.
We call $\Cal I$ the defining ideal sheaf of $Y_\bullet$.

Obviously, $i_\bullet$ is of fiber type if and only if $\Cal I$ is
equivariant.

\paragraph Let $X_\bullet$ be locally noetherian.
A morphism $i_\bullet : Y_\bullet\rightarrow X_\bullet$ is said to be
a regular embedding, if $i_\bullet$ is a closed immersion of fiber type
such that $i_j:Y_j\rightarrow X_j$ is a regular embedding for each
$j\in\ob(I)$, or equivalently, $\Cal I$ is coherent and $\Cal I_{j,x}$ is
a complete intersection ideal of $\Cal O_{X_j,x}$ for any $j\in\ob(I)$ and
$x\in X_j$.
If this is the case, we say that $\Cal I$ is a local complete intersection
ideal sheaf.

A closed immersion of fiber type 
$i_\bullet:Y_\bullet\rightarrow X_\bullet$ with $X_\bullet$
locally noetherian is a regular embedding
if and only if $i_\bullet^*\Cal I$ is locally free and
$\Cal I_{i,x}$ is of finite projective dimension as an 
$\Cal O_{X_i,x}$-module.

Note that $i_\bullet^*\Cal I\cong \Cal I/\Cal I^2$, and we have
\[
\height_{\Cal O_{X_i,x}}\Cal I_x=\rank_{\Cal O_{Y_i,y}}(i_\bullet^*\Cal I)_{
i,y}
\]
for any point $(i,y)$ of $Y_\bullet$, where $x=i_i(y)$.
We call these numbers the codimension of $\Cal I$ at $(i,y)$.

\begin{Proposition}\label{regular-embedding.thm}
Let $I$ be a small category, $S$ a scheme, and 
$i_\bullet:Y_\bullet\rightarrow X_\bullet$ a morphism in $\Cal P(I,\Sch/S)$.
Assume that $X_\bullet$ is locally noetherian with flat arrows
and $i_\bullet$ is
a regular embedding.
Let $\Cal I$ be the defining ideal of $Y_\bullet$, and assume that
$Y_\bullet$ has a constant codimension $d$.
Then we have the following.
\begin{description}
\item[1] 
$\uExt^i_{\Cal O_{X_\bullet}}((i_\bullet)_*\Cal O_{Y_\bullet},
\Cal O_{X_\bullet})=0$ for $i\neq d$.
\item[2] The canonical map
\[
\uExt^d_{\Cal O_{X_\bullet}}((i_\bullet)_*\Cal O_{Y_\bullet},
\Cal O_{X_\bullet})
\rightarrow
\uExt^d_{\Cal O_{X_\bullet}}(
(i_\bullet)_*\Cal O_{Y_\bullet},
(i_\bullet)_*\Cal O_{Y_\bullet})
\]
is an isomorphism.
\item[3] The Yoneda algebra
\[
\uExt^\bullet_{\Cal O_{Y_\bullet}}(
(i_\bullet)_*\Cal O_{Y_\bullet},
(i_\bullet)_*\Cal O_{Y_\bullet})
:=\bigoplus_{j\geq 0}
\uExt^j_{\Cal O_{Y_\bullet}}(
(i_\bullet)_*\Cal O_{Y_\bullet},
(i_\bullet)_*\Cal O_{Y_\bullet})
\]
is isomorphic to the exterior algebra $(i_\bullet)_*
\ext^\bullet (i_\bullet^*\Cal I)^\vee$ 
as graded $\Cal O_{X_\bullet}$-algebras.
\item[4] There is an isomorphism
\[
i_\bullet^\natural\Cal O_{X_\bullet}\cong \ext^d (i_\bullet^*\Cal I)^\vee
[-d].
\]
\item[5] For $\Bbb F\in\Cal D^+(X_\bullet)$, there is a functorial 
isomorphism
\[
i_\bullet^\natural\Bbb F\cong
\ext^d (i_\bullet^*\Cal I)^\vee\otimes^{\bullet,L}_{\Cal O_{Y_\bullet}}
Li_\bullet^*\Bbb F[-d].
\]
\end{description}
\end{Proposition}

\proof {\bf 1} is trivial, since $\Cal I_{i,x}$ is a complete intersection
ideal of the local ring $\Cal O_{X_i,x}$ of codimension $d$ 
for any point $(i,x)$ of $X_\bullet$.

{\bf 2} Note that $(i_\bullet)_*\Cal O_{Y_\bullet}\cong \Cal O_{X_\bullet}/
\Cal I$.
From the short exact sequence
\[
0\rightarrow \Cal I\rightarrow\Cal O_{X_\bullet}\rightarrow
\Cal O_{X_\bullet}/\Cal I\rightarrow 0,
\]
we get an isomorphism
\[
\uExt^1_{\Cal O_{X_\bullet}}((i_\bullet)_*\Cal O_{Y_\bullet},
(i_\bullet)_*\Cal O_{Y_\bullet})\cong
\uHom_{\Cal O_{X_\bullet}}(\Cal I,\Cal O_{X_\bullet}/\Cal I)
\cong (i_\bullet)_*(i_\bullet^*\Cal I)^\vee.
\]
The canonical map
\[
(i_\bullet)_*(i_\bullet^*\Cal I)^\vee\cong
\uExt^1_{\Cal O_{X_\bullet}}((i_\bullet)_*\Cal O_{Y_\bullet},
(i_\bullet)_*\Cal O_{Y_\bullet})\hookrightarrow
\uExt^\bullet_{\Cal O_{X_\bullet}}((i_\bullet)_*\Cal O_{Y_\bullet},
(i_\bullet)_*\Cal O_{Y_\bullet})
\]
is uniquely extended to an $\Cal O_{X_\bullet}$-algebra map
\[
T_\bullet((i_\bullet)_*(i_\bullet^*\Cal I)^\vee)\rightarrow
\uExt^\bullet_{\Cal O_{X_\bullet}}((i_\bullet)_*\Cal O_{Y_\bullet},
(i_\bullet)_*\Cal O_{Y_\bullet}),
\]
where $T_\bullet$ denotes the tensor algebra.
It suffices to prove that this map is an epimorphism, which induces
an isomorphism
\[
\ext^\bullet((i_\bullet)_*(i_\bullet^*\Cal I)^\vee)\rightarrow
\uExt^\bullet_{\Cal O_{X_\bullet}}((i_\bullet)_*\Cal O_{Y_\bullet},
(i_\bullet)_*\Cal O_{Y_\bullet}).
\]
In fact, the exterior algebra is compatible with base change, and
\begin{multline*}
\ext^\bullet((i_\bullet)_*(i_\bullet^*\Cal I)^\vee)\cong
(i_\bullet)_*i_\bullet^*
\ext^\bullet((i_\bullet)_*(i_\bullet^*\Cal I)^\vee)\\
\cong
(i_\bullet)_*\ext^\bullet ((i_\bullet^*(i_\bullet)_*)(i_\bullet^*\Cal I)^\vee)
\cong
(i_\bullet)_*\ext^\bullet (i_\bullet^*\Cal I)^\vee.
\end{multline*}
To verify this, we may assume that $i_\bullet:Y_\bullet\rightarrow 
X_\bullet$ is a morphism of single schemes, $X_\bullet=\Spec A$ affine,
and $\Cal I=\tilde I$ generated by an $A$-sequence.
The proof for this case is essentially the same as 
\cite[Lemma~IV.1.1.8]{Hashimoto}, and we omit it.

{\bf 4} Let $Z$ denote the ringed site $(\Zar(X_\bullet),(i_\bullet)_*
\Cal O_{Y_\bullet})$, and $g:Z\rightarrow \Zar(Y_\bullet)$ the associated
admissible ringed continuous functor.
By {\bf 2--3}, there is a sequence of isomorphisms of $\Coh(X_\bullet)$
\[
(i_\bullet)_*\ext^d(i_\bullet^*\Cal I)^\vee
\cong
\uExt^d_{\Cal O_{Y_\bullet}}(
(i_\bullet)_*\Cal O_{Y_\bullet},
(i_\bullet)_*\Cal O_{Y_\bullet})
\cong
\uExt^d_{\Cal O_{X_\bullet}}((i_\bullet)_*\Cal O_{Y_\bullet},
\Cal O_{X_\bullet}).
\]
In view of {\bf 1}, there is an isomorphism
\[
R(g_\bullet)^{\#}\ext^d(i_\bullet^*\Cal I)^\vee
\cong R\uHom_{\Cal O_{X_\bullet}}^\bullet(\Cal O_Z,\Cal O_{X_\bullet})[d]
\]
in $D^b_{\Coh}(\Mod(Z))$.
Applying $(g_\bullet)_{\#}$ to both sides, we get
\[
\ext^d(i_\bullet^*\Cal I)^\vee\cong i_\bullet^\natural\Cal O_{X_\bullet}[d].
\]

{\bf 5} is an immediate consequence of {\bf 4} and Lemma~\ref{fpd-dual.thm}.
\qed

\paragraph \label{Kahler.par}
Let $I$ and $S$ be as in (\ref{closed-immersion-ideal.par}).
Let $f_\bullet:X_\bullet\rightarrow Y_\bullet$ be a morphism in 
$\Cal P(I,\Sch/S)$.
Assume that $f_\bullet$ is separated so that the diagonal
$\Delta_{X_\bullet/Y_\bullet}: X_\bullet\rightarrow X_\bullet\times_{Y_\bullet}
X_\bullet$ is a closed immersion.
Define $\Omega_{X_\bullet/Y_\bullet}:=i_\bullet^*\Cal I$, where
$\Cal I:=\Ker(\eta:\Cal O_{X_\bullet\times_{Y_\bullet}X_\bullet}\rightarrow
(\Delta_{X_\bullet/Y_\bullet})_*\Cal O_{X_\bullet})$.
Note that $(\Delta_{X_\bullet/Y_\bullet})_*\Omega_{X_\bullet/
Y_\bullet}\cong \Cal I/\Cal I^2$.

\begin{Lemma} Let the notation be as above.
Then we have
\begin{description}
\item[1] $\Omega_{X_\bullet/Y_\bullet}$ is locally quasi-coherent.
\item[2] If $f_\bullet$ is of fiber type, then $\Omega_{X_\bullet/Y_\bullet}$ 
is quasi-coherent.
\item[3] For $i\in\ob(I)$, 
there is a canonical isomorphism $\Omega_{X_i/Y_i}\cong (\Omega_{X_\bullet/
Y_\bullet})_i$.
\end{description}
\end{Lemma}

\proof Easy.
\qed

\begin{Theorem}\label{smooth-twisted.thm}
Let $I$ be a finite category, 
$S$ a separated noetherian scheme,
and $f_\bullet:X_\bullet\rightarrow Y_\bullet$ a 
morphism in $\Cal P(I,\Sch/S)$.
Assume that $Y_\bullet$ is separated 
of finite type over $S$, and $f_\bullet$ is separated smooth of
finite type and fiber type.
Assume that $f_\bullet$ has a constant relative dimension $d$.
Then for any $\Bbb F\in D^+_{\LQco(Y_\bullet)}(\Mod(X_\bullet))$,
there is a functorial isomorphism
\[
\ext^d\Omega_{X_\bullet/Y_\bullet}[d]\otimes^\bullet_{\Cal O_{X_\bullet}}
f^*_\bullet\Bbb F\cong f_\bullet^!\Bbb F,
\]
where $[d]$ denotes the shift of degree.
\end{Theorem}

\proof In view of Theorem~\ref{ffd.thm}, it suffices to show that
there is an isomorphism $f_\bullet^!\Cal O_{Y_\bullet}
\cong \ext^d\Omega_{X_\bullet/Y_\bullet}[d]$.
Consider the commutative diagram
\[
\begin{picture}(360,65)(-38,35)
\put(50,80){\makebox(0,0)[b]{$X_\bullet$}}
\put(60,83){\vector(1,0){25}}
\put(114,80){\makebox(0,0)[b]{$X_\bullet\times_{Y_\bullet}X_\bullet$}}
\put(114,75){\vector(0,-1){25}}
\put(114,47){\makebox(0,0)[t]{$X_\bullet$}}
\put(53,76){\vector(3,-2){50}}
\put(143,83){\vector(1,0){25}}
\put(179,80){\makebox(0,0)[b]{$X_\bullet$}}
\put(179,75){\vector(0,-1){25}}
\put(124,41){\vector(1,0){44}}
\put(179,47){\makebox(0,0)[t]{$Y_\bullet$.}}
\put(72,86){\makebox(0,0)[b]{\small$\Delta$}}
\put(155,86){\makebox(0,0)[b]{\small$p_2$}}
\put(182,63){\makebox(0,0)[l]{\small$f_\bullet$}}
\put(144,44){\makebox(0,0)[b]{\small$f_\bullet$}}
\put(117,63){\makebox(0,0)[l]{\small$p_1$}}
\put(74,61){\makebox(0,0)[tr]{\small$\id$}}
\end{picture}
\]
By Lemma~\ref{fiber-type-2.thm} and Lemma~\ref{fiber-type-1.thm},
all the morphisms in the diagrams are of fiber type.
As $p_1$ is smooth of finite type of relative dimension $d$, 
$\Delta$ is a regular embedding of the constant codimension $d$.

By Theorem~\ref{flat-base-change.thm}, 
Theorem~\ref{finite-duality.thm}, and
Proposition~\ref{regular-embedding.thm},
we have
\begin{multline*}
\Cal O_{X_\bullet}\cong
f_\bullet^*\Cal O_{Y_\bullet}
\cong
\Delta^!p_1^!f_\bullet^*\Cal O_{Y_\bullet}
\cong
\Delta^\natural p_2^*f_\bullet^!\Cal O_{Y_\bullet}
\cong\\
\ext^d \Omega_{X_\bullet/Y_\bullet}^\vee[-d]
\otimes^{\bullet,L}_{\Cal O_{X_\bullet}}
L\Delta^*p_2^*f_\bullet^!\Cal O_{Y_\bullet}
\cong
\ext^d\Omega_{X_\bullet/Y_\bullet}^\vee[-d]
\otimes^{\bullet,L}_{\Cal O_{X_\bullet}}
f_\bullet^!\Cal O_{Y_\bullet}.
\end{multline*}
As $\ext^d\Omega_{X_\bullet/Y_\bullet}$ is an invertible sheaf,
we are done.
\qed

\section{Simplicial objects}

\paragraph For $n\in\Bbb Z$ with $n\geq -1$, we define 
$[n]$ to be the totally ordered finite set
$\{0<1<\ldots<n\}$.
Thus, $[-1]=\emptyset$, $[0]=\{0\}$, $[1]=\{0<1\}$, and so on.
We define $(\Delta^+)$ to be the small category given by
$\ob(\Delta^+):=\{[n]\setbar n\in\Bbb Z,\;n\geq -1\}$ and
\[
\Mor(\Delta^+):=\{\mbox{monotone maps}\}.
\]
For a subset $S$ of $\{-1,0,1,\ldots\}$, 
we define $(\Delta^+)_S$ to be the full subcategory of $(\Delta^+)$
such that $\ob((\Delta^+)_S)=\{[n]\setbar n\in S\}$.
We define $(\Delta):=(\Delta^+)_{[0,\infty)}$.
If $-1\notin S$, then $(\Delta^+)_S$ is also denoted by $(\Delta)_S$.

We define $(\Delta^+)\mon$ to be the subcategory of $(\Delta^+)$
by $\ob((\Delta^+)\mon):=\ob(\Delta^+)$ and
\[
\Mor ((\Delta^+)\mon):=\{\mbox{injective monotone maps}\}.
\]
For $S\subset \{-1,0,1,\ldots\}$, 
the full subcategories 
$(\Delta^+)\mon_S$ and $(\Delta)\mon_S$ of $(\Delta^+)\mon$ 
are defined similarly.

We denote $(\Delta)\mon_{\{0,1,2\}}$ and $(\Delta^+)\mon_{\{-1,0,1,2\}}$
by $\Delta_M$ and $\Delta^+_M$, respectively.

Let $\Cal C$ be a category.
We call an object of
$\Cal P((\Delta^+),\Cal C)$
(resp.\ $\Cal P((\Delta),\Cal C)$,
an augmented simplicial object 
(resp.\ simplicial object)
of $\Cal C$.

For a subcategory $\Cal D$ of $(\Delta^+)$ 
and an object $X_\bullet\in\Cal P(\Cal D,\Cal C)$,
we denote $X_{[n]}$ by $X_n$.

As $[-1]$ is the initial object of $\Delta^+$,
an augmented simplicial object
$X_\bullet$ of $\Cal C$
with $X_{-1}=c$ is identified with a simplicial object 
of $\Cal C/c$.

We define some particular morphisms in $\tilde \Delta$.
The unique map $[-1]\rightarrow[n]$ is denoted by $\varepsilon(n)$.
The unique injective monotone map $[n-1]\rightarrow[n]$ such that
$i$ is not in the image is denoted by $\delta_i(n)$ for $i\in[n]$.
The unique surjective monotone map $[n+1]\rightarrow[n]$ such that
$i$ has two inverse images is denoted by $\sigma_i(n)$ for $i\in[n]$.
The unique map $[0]\rightarrow[n]$ such that $i$ is in the image
is denoted by $\rho_i(n)$.
The unique map $[n]\rightarrow [0]$ is denoted by $\lambda_n$.

Let $\Cal D$ be a subcategory of $(\Delta^+)$.
For $X_\bullet\in\Cal P(\Cal D,\Cal C)$, we denote $X_\bullet(\varepsilon(n))$
(resp.\ $X_\bullet(\delta_i(n))$, $X_\bullet(\sigma_i(n))$, $X_\bullet
(\rho_i(n))$, and $X_\bullet(\lambda_n)$) 
by $e(n,X_\bullet)$ (resp.\ $d_i(n,X_\bullet)$, $s_i(n,X_\bullet)$,
$r_i(n,X_\bullet)$, and $l_n(X_\bullet)$),
or simply by $e(n)$ (resp.\ $d_i(n)$, $s_i(n)$, $r_i(n)$, $l_n$),
if there is no confusion.

Note that $\Delta$ is generated by $\delta_i(n)$, $\sigma_i(n)$ 
for various $i$ and $n$.

\paragraph Note that $\Delta^+([m],[n])$ is  a finite set for any
$m,n$.
Assume that $\Cal C$ has finite limits and
let $f:X\rightarrow Y$ be a morphism in $\Cal C$.
Then the {\em \v Cech nerve} is defined to be
$\Nerve(f):=\cosk_{(\Delta^+)_{\{0,1\}}}^{\Delta^+}(f)$,
where $\cosk_{(\Delta^+)_{\{0,1\}}}^{\Delta^+}$ is the right adjoint of
the restriction, which exists.
It is described as follows.
$\Nerve(f)_n=X\times_Y\times\cdots\times_Y X$ ($(n+1)$-fold fiber product)
for $n\geq 0$, and $\Nerve(f)_{-1}=Y$.
Note that $d_i(n)$ is given by
\[
d_i(n)(x_n,\ldots,x_1,x_0)=(x_n,\delme{i},x_1,x_0),
\]
and $s_i(n)$ is given by
\[
s_i(n)(x_n,\ldots,x_1,x_0)=(x_n,\ldots,x_{i+1},x_{i},x_{i},x_{i-1},
\ldots,x_1,x_0)
\]
if $\Cal C=\Set$.

\paragraph Let $S$ be a scheme.
A simplicial object (resp.\ augmented simplicial object)
in $\Sch/S$, in other words, 
an object of $\Cal P((\Delta),\Sch/S)$ (resp.\ $\Cal P((\Delta^+),\Sch/S)$),
is called a simplicial (resp.\ augmented simplicial) $S$-scheme.
The following is well-known.

\begin{Lemma}\label{simplicial-M.thm}
Let $X_{\bullet}\in\Cal P((\Delta),\Sch/S)$.
Then the restriction $(?)_{\Delta_M}:
\EM(X_{\bullet})\rightarrow
\EM(X_{\bullet}|_{\Delta_M})$ is an equivalence.
With the equivalence, quasi-coherent
sheaves correspond to quasi-coherent sheaves.
\end{Lemma}

\proof 
We define a third category $\Cal A$ as follows.
An object of $\Cal A$ is a pair $(\Cal M_0,\varphi)$ such that, 
$\Cal M_0\in\Mod(X_{0})$, $\varphi\in\Hom_{\Mod(X_{1})}
(d_0^*(\Cal M_0),d_1^*(\Cal M_0))$, $\varphi$ an isomorphism, and that
$d_1^*(\varphi)=d_2^*(\varphi)\circ d_0^*(\varphi)$.
Note that applying $l_2^*$ to the last equality, we get
$l_1^*(\varphi)=l_1^*(\varphi)\circ l_1^*(\varphi)$.
As $\varphi$ is an isomorphism, we get $l_1^*(\varphi)=\id$.

A morphism $\gamma_0:(\Cal M_0,\varphi)\rightarrow(\Cal N_0,\psi)$ is an
element of
\[
\gamma_0\in
\Hom_{\Mod(X_{0})}(\Cal M_0,\Cal N_0)
\]
such that
\[
\psi\circ d_0^*(\gamma_0)=d_1^*(\gamma_0)\circ \varphi.
\]

We define a functor $\Phi:\EM(X_{\bullet}|_{\Delta_M})
\rightarrow \Cal A$ by
\[
\Phi(\Cal M):=(\Cal M_0,
\alpha_{d_1(1)}^{-1}\circ\alpha_{d_0(1)}).
\]
It is easy to verify that this gives a well-defined functor.

Now we define a functor $\Psi:\Cal A\rightarrow \EM(X_{\bullet})$.
Note that an object $\Cal M$ of $\EM(X_{\bullet})$ is identified with
a family $(\Cal M_n,\alpha_w)_{[n]\in(\Delta),\;w\in\Mor((\Delta))}$
such that $\Cal M_n\in\Mod(X_{n})$,
\[
\alpha_w\in\Hom_{\Mod(X_{n})}((X_{w})^*_{\Mod}(\Cal M_m),\Cal M_n)
\]
for $w\in\Delta(m,n)$,
$\alpha_w$ is an isomorphism,
and
\begin{equation}\label{cocycle.eq}
\alpha_{ww'}=\alpha_{w'}\circ X_{w'}^*\alpha_w
\end{equation}
whenever $ww'$ is defined, see (\ref{friedlander-identification.par}).

For $(\Cal M_0,\varphi)\in\Cal A$, we define
$\Cal M_{n,i}:=(r_i(n))^*(\Cal M_0)$, and $\Cal M_n:=\Cal M_{n,0}$
(as we have $r_0(0)=\id$, this will not cause a double-meaning of 
$\Cal M_0$) for $n\geq 0$ and $0\leq i\leq n$.
We define $\psi_i(n):\Cal M_{n,i+1}\rightarrow \Cal M_{n,i}$
to be $(X_{q(i,n)})^*(\varphi)$
for $n\geq 1$ and $0\leq i<n$,
where $q(i,n):[1]\rightarrow[n]$
is the unique injective monotone map with $\{i,i+1\}=\Image q(i,n)$.
We define $\varphi_i(n):\Cal M_{n,i}\cong\Cal M_n$ to be the composite map
\[
\varphi_i(n):=\psi_0(n)\circ \psi_1(n)\circ\cdots\circ\psi_{i-1}(n)
\]
for $n\geq 0$ and $0\leq i\leq n$.

Now we define
\[
\alpha_w\in\Hom_{\Mod(X_{n})}(X_{w}^*(\Cal M_m),\Cal M_n)
\]
to be the map
\[
X_{w}^*\Cal M_m=X_{w}^* r_0(m)^*
\Cal M_0=r_{w(0)}(n)^*\Cal M_0
=\Cal M_{n,w(0)}\specialarrow{\varphi_{w(0)}(n)}
\Cal M_n
\]
for $w\in\Delta(m,n)$.

Thus $(\Cal M_0,\varphi)$ yields a family $(\Cal M_n,\alpha_w)$,
and this gives the definition of $\Psi:\Cal A\rightarrow
\EM(X_{\bullet})$.
The details of the proof of the well-definedness is left to the reader.

It is also straightforward to check that $(?)_{\Delta_M}$, $\Phi$,
and $\Psi$ give the equivalence of these three categories.
The proof is also left to the reader.

The last assertion is obvious from the construction.
\qed

\section{Descent theory}

Let $S$ be a scheme.

\paragraph
Consider the functor $\shift:(\Delta^+)\rightarrow(\Delta)$
given by $\shift[n]:=[n+1]$, $\shift(\delta_i(n)):=\delta_{i+1}(n+1)$,
$\shift(\sigma_i(n)):=\sigma_{i+1}(n+1)$, and $\shift(\varepsilon(0)):
=\delta_1(1)$.
We have a natural transformation $(\delta_0^+):\Id_{(\Delta^+)}\rightarrow
\iota\circ\shift$ given by $(\delta_0^+)_n:=\delta_0(n+1)$ for $n\geq 0$ and
$(\delta_0^+)_{-1}:=\varepsilon(0)$,
where $\iota:(\Delta)\hookrightarrow(\Delta^+)$ is the inclusion.
We denote $(\delta_0^+)\iota$ by $(\delta_0)$.
Note that $(\delta_0)$ can be viewed as a natural map
$(\delta_0):\Id_{(\Delta)}\rightarrow\shift\iota$.

Let $X_\bullet\in\Cal P((\Delta),\Sch/S)$.
We define $X_\bullet'$ to be the augmented simplicial scheme
$\shift^{\#}(X_\bullet)=X_\bullet\shift$.
The natural map
\[
X_\bullet(\delta_0):X_\bullet'|_{(\Delta)}=X_\bullet\shift\iota
\rightarrow X_\bullet
\]
is denoted by $(d_0)(X_\bullet)$ or $(d_0)$.
Similarly, if $Y_\bullet\in\Cal P((\Delta^+),\Sch/S)$, then
\[
(d_0^+)(Y_\bullet):(Y_\bullet|_{(\Delta)})'=
Y_\bullet\iota\shift\specialarrow{Y_\bullet(\delta_0^+)}
Y_\bullet
\]
is defined as well.

\paragraph We say that $X_\bullet\in\Cal P((\Delta),\Sch/S)$ is a simplicial
groupoid of $S$-schemes if there is a faithfully flat morphism of $S$-schemes
$g:Z\rightarrow Y$ such that there is a faithfully 
flat morphism
$f_\bullet:Z_\bullet\rightarrow X_\bullet$ of fiber type of $\Cal P((\Delta),
\Sch/S)$, where $Z_\bullet=\Nerve(g)|_{(\Delta)}$.

\begin{Lemma}\label{groupoid.thm}
Let $X_\bullet\in\Cal P((\Delta),\Sch/S)$.
\begin{description}
\item[\bf 1] If $X_\bullet\cong\Nerve(g)|_{(\Delta)}$ 
for some faithfully flat morphism $g$ of $S$-schemes, then
$X_\bullet$ is a simplicial groupoid.
\item[\bf 2] If $f_\bullet:Z_\bullet\rightarrow X_\bullet$ is a faithfully
flat morphism of fiber type of simplicial $S$-schemes
and $Z_\bullet$ is a simplicial groupoid, then we have $X_\bullet$ is also
a simplicial groupoid.
\item[\bf 3] $X_\bullet$ is a simplicial groupoid if and only if 
$(d_0):X_\bullet'|_{(\Delta)}\rightarrow X_\bullet$ is of fiber type, 
the canonical unit map
\[
X_\bullet'\rightarrow \Nerve(d_1(1))
\]
is an isomorphism, and $d_0(1)$ and $d_1(1)$ are flat.
\item[\bf 4] If $f_\bullet:Z_\bullet\rightarrow X_\bullet$ is a morphism of
fiber type of simplicial $S$-schemes and $X_\bullet$ is a simplicial groupoid,
then $Z_\bullet$ is a simplicial groupoid.
\item[\bf 5] A simplicial groupoid has flat $(\Delta)\mon$-arrows.
\item[\bf 6] If $X_\bullet$ is a simplicial groupoid of $S$-schemes
such that $d_0(1)$ and $d_1(1)$ are separated \(resp.\ quasi-compact, 
of finite type, smooth, \'etale\), then $X_\bullet$ has
separated \(resp.\ quasi-compact, of finite type, smooth,\'etale\) 
$(\Delta)\mon$-arrows, and $(d_0):X_\bullet'|_{(\Delta)}\rightarrow X_\bullet$
is separated \(resp.\ quasi-compact, of finite type, smooth, \'etale\).
\end{description}
\end{Lemma}

\proof {\bf 1} and {\bf 2} are obvious by definition.

We prove {\bf 3}.
We prove the `if' part.
As $d_0(1)s_0(0)=\id=d_1(1)s_0(0)$, we have that $d_0(1)$ and $d_1(1)$ are
faithfully flat by assumption.
As $d_0(1)=(d_0)_0$ is faithfully flat and $(d_0)$ is of fiber type, 
it is easy to see that $(d_0)$ is also faithfully flat.
So this direction is obvious.

We prove the `only if' part.
As $X_\bullet$ is a simplicial groupoid,
there is a faithfully flat $S$-morphism $g:Z\rightarrow Y$ and
a faithfully flat morphism $f_\bullet:Z_\bullet\rightarrow X_\bullet$
of simplicial $S$-schemes of fiber type, where
$Z_\bullet=\Nerve(g)|_{(\Delta)}$.
It is easy to see that $(d_0):Z_\bullet'|_{(\Delta)}\rightarrow Z_\bullet$
is nothing but the base change by $g$, and it is faithfully flat
of fiber type.
It is also obvious that $Z_\bullet'\cong\Nerve(d_1(1)(Z_\bullet))$ and
$d_0(1)(Z_\bullet)$ and $d_1(1)(Z_\bullet)$ are flat.
It is obvious that $f_\bullet':Z_\bullet'\rightarrow X_\bullet'$ is
faithfully flat of fiber type.
Now by Lemma~\ref{fiber-type-1.thm}, $(d_0)(X_\bullet)$ is of fiber type.
As $f_\bullet$ is faithfully flat of fiber type and $d_0(1)(Z_\bullet)$
and $d_1(1)(Z_\bullet)$ are flat, we have that $d_0(1)(X_\bullet)$ and
$d_1(1)(X_\bullet)$ are flat.
When we base change $X_\bullet'\rightarrow\Nerve(d_1(1)(X_\bullet))$
by $f_0:Z_0\rightarrow X_0$, then we have the isomorphism
$Z_\bullet'\cong\Nerve(d_1(1)(Z_\bullet))$.
As $f_0$ is faithfully flat, we have that $X_\bullet'\rightarrow
\Nerve(d_1(1))$ is also an isomorphism.

The assertions {\bf 4}, {\bf 5} and {\bf 6} are proved easily.
\qed

\paragraph Let $X_\bullet\in\Cal P((\Delta),\Sch/S)$.
Then we define $F:\Zar(X_\bullet')\rightarrow \Zar(X_\bullet)$ by
$F(i,U)=(\shift i,U)$.
The corresponding pull-back $F^{\#}_{\Mod}$ is denoted by $(?)'$.
It is easy to see that $(?)'$ have a left and a right adjoint.
It also preserves equivariant and locally quasi-coherent sheaves.

Let $\Cal M\in\Mod(X_\bullet)$.
Then we define $(\alpha):(d_0)^*\Cal M\rightarrow \Cal M_{(\Delta)}'$
by 
\[
(\alpha)_n:((d_0)^*\Cal M)_n=d_0(n)^*\Cal M_n\specialarrow
{\alpha_{\delta_0(n)}}\Cal M_{n+1}=\Cal M'_n.
\]
It is easy to see that $(\alpha):(d_0)^*\rightarrow(?)_{(\Delta)}\circ(?)'$
is a natural map.
Similarly, for $Y_\bullet\in\Cal P((\Delta^+),\Sch/S)$, 
\[
(\alpha^+):(d_0^+)^*\rightarrow (?)'\circ(?)_{(\Delta)}
\]
is defined.

\paragraph
Let $X_{\bullet}\in\Cal P((\Delta^+),\Sch/S)$, 
and $\Cal M\in\Mod(X_{\bullet}|_{(\Delta)})$.
Then, we have a cosimplicial object $\Cos(\Cal M)$
of $\Mod(X_{-1})$ (i.e., a simplicial object of $\Mod(X_{-1})\op$).
We have $\Cos(\Cal M)_n:=e(n)_*(\Cal M_{n})$,
and\[
w(\Cos(\Cal M)):= e(m)_*(\beta_{w}):
 e(m)_*(\Cal M_{m})\rightarrow
 e(m)_*(w_*(\Cal M_{n}))\cong
 e(n)_*(\Cal M_{n})
\]
for a morphism $w:m\rightarrow n$ in $\Delta$.
Similarly, the augmented cosimplicial object $\Cos^+(\Cal N)$
of $\Mod(X_{-1})$ is defined for $\Cal N\in\Mod(X_\bullet)$.

By (\ref{right-induction.par}), it is easy to see that
$\Cal M^+:=R_{(\Delta)}^{\Mod}(\Cal M)$ is
$\Cal M$ on $X_\bullet|_{(\Delta)}$, $\Cal M^+_{-1}$
is $\projlim \Cos(\Cal M)$, and
$
\beta_{\varepsilon(n)}(\Cal M^+)
$
is nothing but the canonical map
\[
\projlim \Cos(\Cal M)\rightarrow
\Cos(\Cal M)_n= e(n)_*(\Cal M_{n})
= e(n)_*(\Cal M^+_{n}).
\]

Note also that we have
\begin{equation}\label{plim-simp.eq}
\projlim\Cos(\Cal M)=
\Ker( e(0)_*\beta_{\delta_0(1)}- e(0)_*\beta_{\delta_1(1)}),
\end{equation}
which is determined only by $\Cal M_{(\Delta)_{\{0,1\}}}$.

\begin{Lemma}\label{inv-im-right-ind.thm}
Let $f_{\bullet}:X_{\bullet}\rightarrow
Y_{\bullet}$ be a morphism of $\Cal P((\Delta^+),\Sch/S)$.
If $f_{\bullet}|_{(\Delta^+)_{\{-1,0,1\}}}$ is of fiber type, 
flat quasi-compact separated and $Y_{\bullet}$ has flat 
$(\Delta^+)_{\{-1,0,1\}}$-arrows,
then the canonical map
\[
f_{\bullet}^*\circ R_{(\Delta)}\rightarrow
R_{(\Delta)}\circ (f_{\bullet}|_{(\Delta)})^*
\]
is an isomorphism of functors from $\LQco(Y_{\bullet}|_{(\Delta)})$ to
$\LQco(X_{\bullet})$.
\end{Lemma}

\proof 
To prove that the map in question is an isomorphism,
it suffices to show that the map is an isomorphism after applying the
functor $(?)_n$ for $n\geq -1$.
This is trivial if $n\geq 0$.
On the other hand, if $n=-1$, the map restricted at $-1$ and
evaluated by $\Cal M\in\LQco(Y_\bullet)$ is nothing but
\[
f_{-1}^*(H^0(\Cos(\Cal M)))
\cong
H^0(f_{-1}^*(\Cos(\Cal M)))
\rightarrow
H^0(\Cos((f_\bullet)^*(\Cal M))).
\]
The first map is an isomorphism as $f_{-1}$ is flat.
Although the map $f_{-1}^*(\Cos(\Cal M))\rightarrow \Cos((f_\bullet
)^*(\Cal M))$
may not be a chain isomorphism, it is an isomorphism at the degrees
$-1,0,1$, and it induces the isomorphism of $H^0$.
\qed

\begin{Lemma}\label{unit-isom.thm}
Let $X_\bullet\in\Cal P((\Delta),\Sch/S)$, and $\Cal M\in\Mod(X_\bullet)$.
Then the \(associated chain complex of the\) augmented cosimplicial object
$\Cos^+(\Cal M')$ of $\Cal M'\in\Mod(X_\bullet')$ is split exact.
In particular, the unit map $u:\Cal M'\rightarrow R_{(\Delta)}
\Cal M'_{(\Delta)}$ is an isomorphism.
\end{Lemma}

\proof 
Define
$s_n:\Cos^+(\Cal M')_{n}\rightarrow\Cos^+(\Cal M')_{n-1}$
to be
\begin{multline*}
\Cos^+(\Cal M')_{n}=(r_0)(n+1)_*(\Cal M_{n+1})
\specialarrow{(r_0)(n+1)_*\beta{\sigma_0(n)}}\\
(r_0)(n+1)_*s_0(n)_*(\Cal M_{n})
=(r_0)(n)_*(\Cal M_{n})
\Cos^+(\Cal M')_{n-1}
\end{multline*}
for $n\geq 0$, and $s_{-1}:\Cos^+(\Cal M')_{-1}\rightarrow 0$ to be $0$.
It is easy to verify that $s$ is a chain deformation of $\Cos^+(\Cal M')$.
\qed

\begin{Corollary}\label{equiv-descent.thm}
Let the notation be as in the lemma.
Then there is a functorial isomorphism
\begin{equation}
\label{descent-eq.eq}
R_{(\Delta)}(d_0)^*(\Cal M)\rightarrow \Cal M'
\end{equation}
for $\Cal M\in\EM(X_\bullet)$.
In particular, there is a functorial isomorphism
\begin{equation}\label{descent-eq2.eq}
(R_{(\Delta)}(d_0)^*(\Cal M))_{-1}\rightarrow \Cal M_0.
\end{equation}
\end{Corollary}

\proof The first map (\ref{descent-eq.eq}) is defined to be the composite
\[
R_{(\Delta)}(d_0)^*
\specialarrow{R_{(\Delta)}(\alpha)}
R_{(\Delta)}(?)'_{(\Delta)}
\specialarrow{u^{-1}}
(?)'.
\]
As $(\alpha)(\Cal M)$ is an isomorphism if $\Cal M$ is equivariant,
this is an isomorphism.
The second map (\ref{descent-eq2.eq}) is obtained from (\ref{descent-eq.eq}),
applying $(?)_{-1}$.
\qed

The following well-known theorem in descent theory contained in
\cite{Murre} is now easy to prove.

\begin{Proposition}\label{descent-main.thm}
Let $f:X\rightarrow Y$ be a morphism of $S$-schemes,
and set $X_{\bullet}^+:=\Nerve(f)$, and
$X_{\bullet}:=X^+_{\bullet}|_{(\Delta)}$.
Let $\Cal M\in\Mod(X_\bullet)$.
Then we have the following.
\begin{description}
\item[\bf 0] The counit of adjunction
\[
\varepsilon:(R_{(\Delta)}\Cal M)_{(\Delta)}\rightarrow \Cal M
\]
is an isomorphism.
\item[\bf1] If $f$ is quasi-compact separated and $\Cal M\in
\LQco(X_{\bullet})$, then
$R_{(\Delta)}\Cal M\in\LQco(X^+_{\bullet})$.
\item[\bf2] If $f$ is faithfully flat quasi-compact separated
and $\Cal M\in \Qco(X_{\bullet})$, then we have $R_{(\Delta)}\Cal M\in
\Qco(X^+_{\bullet})$.
\item[\bf 3] If $f$ is faithfully flat quasi-compact separated,
$\Cal N\in\EM(X^+_{\bullet})$, and $\Cal N|_{(\Delta)}\in\Qco(X_\bullet)$,
then the unit of adjunction
\[
u:\Cal N\rightarrow R_{(\Delta)}(\Cal N_{(\Delta)})
\]
is an isomorphism.
In particular, $\Cal N$ is quasi-coherent.
\item[\bf 4] If $f$ is faithfully flat quasi-compact separated, then
the restriction functor
\[
(?)_{(\Delta)}:\Qco(X^+_\bullet)\rightarrow
\Qco(X_\bullet|_{(\Delta)})
\]
is an equivalence, with
$R_{(\Delta)}$ its quasi-inverse.
\end{description}
\end{Proposition}

\proof The assertion {\bf 0} follows from Lemma~\ref{induction-restrict.thm}.

We prove {\bf 1}.
By {\bf 0}, it suffices to prove that 
\[
(R_{(\Delta)}\Cal M)_{-1}=
\Ker( e(0)_*\beta_{\delta_0(1)}- e(0)_*\beta_{\delta_1(1)})
\]
is quasi-coherent.
This is obvious by \cite[(9.2.2)]{EGA-I}.

Now we assume that $f$ is faithfully flat quasi-compact separated, 
to prove the assertions {\bf 2}, {\bf 3}, and {\bf 4}.

We prove {\bf 2}.
As we already know that $R_{(\Delta)}\Cal M$ is locally quasi-coherent,
it suffices to show that it is equivariant.
As $(d_0^+):(X_\bullet|_{(\Delta)})'\rightarrow X_\bullet$ is faithfully
flat, it suffices to show that $(d_0^+)^*R_{(\Delta)}\Cal M$ is
equivariant, by Lemma~\ref{inv-im-qco.thm}.
Now the assertion is obvious by Lemma~\ref{inv-im-right-ind.thm} and
Corollary~\ref{equiv-descent.thm}, as $\Cal M'$ is quasi-coherent.

We prove {\bf 3}.
Note that the composite map
\begin{equation}\label{ukyarapa.eq}
(d_0^+)^*\Cal N\specialarrow{(d_0^+)^* u}
(d_0^+)^*R_{(\Delta)}(\Cal N_{(\Delta)})
\rightarrow
R_{(\Delta)}(d_0)^*\Cal N_{(\Delta)}
\cong
R_{(\Delta)}((d_0^+)^*\Cal N)_{(\Delta)}
\end{equation}
is nothing but the unit of adjunction $u((d_0^+)^*\Cal N)$.
As $(\alpha^+):(d_0^+)^*\Cal N\rightarrow (\Cal N_{(\Delta)})'$ is
an isomorphism since $\Cal N$ is equivariant, we have that
$u((d_0^+)^*\Cal N)$ is an isomorphism by Lemma~\ref{unit-isom.thm}.
As the second arrow in (\ref{ukyarapa.eq}) is an isomorphism
by Lemma~\ref{inv-im-right-ind.thm}, we have that
$(d_0^+)^*u$ is an isomorphism.
As $(d_0^+)$ is faithfully flat, we have that
$u:\Cal N\rightarrow R_{(\Delta)}(\Cal N_{(\Delta)})$ is an isomorphism,
as desired.
The last assertion is obvious by {\bf 2}, and {\bf 3} is proved.

The assertion {\bf 4} is a consequence of {\bf 0}, 
{\bf 2} and {\bf 3}.
\qed

\begin{Corollary}
Let $f:X\rightarrow Y$ be a faithfully flat quasi-compact separated
morphism of $S$-schemes, and $\Cal M\in\Mod(Y)$.
Then $\Cal M$ is quasi-coherent if and only if $f^*\Cal M$ is.
\end{Corollary}

\proof Only if part is trivial.

If $f^*\Cal M$ is quasi-coherent, then $\Cal N\in L_{-1}\Cal M$ satisfies
the assumption of {\bf 2} of the proposition, as can be seen easily.
So $\Cal M\cong (\Cal N)_{-1}$ is quasi-coherent.
\qed

\begin{Corollary}\label{ascent-descent.thm}
Let the notation be as in the proposition,
and assume that 
$f$ is faithfully flat quasi-compact separated.
The composite functor
\[
\Bbb A:= (?)_{(\Delta)}\circ L_{-1}:\Qco(Y)\rightarrow \Qco(X_\bullet)
\]
is an equivalence with
\[
\Bbb D:=(?)_{-1}\circ R_{(\Delta)}
\]
its quasi-inverse.
\end{Corollary}

\proof Follows immediately by the proposition and
Lemma~\ref{initial-equivariant.thm}, {\bf 2}, since
$[-1]$ is the initial object of $(\Delta^+)$.
\qed

We call $\Bbb A$ in the corollary the ascent functor,
and $\Bbb D$ the descent functor.

\begin{Corollary}\label{lqco-qco.thm}
Let the notation be as in the proposition.
Then the composite functor
\[
\Bbb A\circ \Bbb D:\LQco(X_\bullet)\rightarrow\Qco(X_\bullet)
\]
is the right adjoint functor 
of the inclusion $\Qco(X_\bullet)\hookrightarrow\LQco(X_\bullet)$.
\end{Corollary}

\proof Note that $\Bbb D:\LQco(X_\bullet)\rightarrow\Qco(Y)$ is
a well-defined functor, and hence $\Bbb A\circ\Bbb D$ is a functor
from $\LQco(X_\bullet)$ to $\Qco(X_\bullet)$.

For $\Cal M\in\Qco(X_\bullet)$ and $\Cal N\in\LQco(X_\bullet)$,
we have
\lrsplit
\Hom_{\Qco(X_\bullet)}(\Cal M,\Bbb A\Bbb D\Cal N)
\cong
\Hom_{\Qco(Y)}(\Bbb D\Cal M,\Bbb D\Cal N)
\cong
\Hom_{\LQco(X_\bullet^+)}(R_{(\Delta)}\Cal M,R_{(\Delta)}\Cal N)\\
\cong
\Hom_{\LQco(X_\bullet)}((R_{(\Delta)}\Cal M)_{(\Delta)},\Cal N)
\cong
\Hom_{\LQco(X_\bullet)}(\Cal M,\Cal N)
\endlrsplit
by the proposition, Corollary~\ref{ascent-descent.thm}, and
Lemma~\ref{initial-equivariant.thm}, {\bf 1}.
\qed

\begin{Corollary}\label{qco-right-adjoint.thm}
Let $X_\bullet$ be a simplicial groupoid of $S$-schemes,
and assume that $d_0(1)$ and $d_1(1)$ are separated and quasi-compact.
Then
\[
(d_0)_*^{\Qco}\circ \Bbb A:\Qco(X_0)\rightarrow\Qco(X_\bullet)
\]
is a right adjoint of
$(?)_0:\Qco(X_\bullet)\rightarrow\Qco(X_0)$,
where $\Bbb A:\Qco(X_0)\rightarrow \Qco(X_\bullet'|_{(\Delta)})$ is the
ascent functor defined in Corollary~\ref{ascent-descent.thm}.
\end{Corollary}

\proof Note that $(d_0)_*^{\Qco}$ is well-defined,
because $(d_0)$ is quasi-compact separated of fiber type
and $X_\bullet$ has flat $(\Delta)\mon$-arrows by Lemma~\ref{groupoid.thm}
and the assumption.
It is obvious that $\Bbb D\circ (d_0)^*_{\Qco}$ is the left adjoint
of $(d_0)_*^{\Qco}\circ \Bbb A$ by Corollary~\ref{ascent-descent.thm}.
On the other hand, we have $(?)_0\cong \Bbb D\circ (d_0)^*_{\Qco}$ by
Corollary~\ref{equiv-descent.thm}.
Hence, $(d_0)_*^{\Qco}\circ \Bbb A$ is a right adjoint of $(?)_0$,
as desired.
\qed

\section{Local noetherian property}\label{noetherian.sec}

An abelian category $\Cal A$ is called {\em locally noetherian}
if it is a $\Cal U$-category, 
satisfies the (AB5) condition, and has a small set of 
noetherian generators \cite{Gabriel}.
For a locally noetherian category $\Cal A$, we denote the full subcategory
of $\Cal A$ consisting of its noetherian objects by $\Cal A_f$.

\begin{Lemma}\label{locally-noetherian.thm}
Let $\Cal A$ be an abelian $\Cal U$-category which satisfies the
{\rm(AB3)} condition, and $\Cal B$ a locally noetherian category.
Let $F:\Cal A\rightarrow \Cal B$ be a faithful exact functor,
and $G$ its right adjoint.
If $G$ preserves filtered inductive limits,
then the following hold.
\begin{description}
\item[\bf 1] $\Cal A$ is locally noetherian
\item[\bf 2] $a\in\Cal A$ is a noetherian object if and only if 
$Fa$ is a noetherian object.
\end{description}
\end{Lemma}

\proof The `if' part of {\bf 2} is obvious, as $F$ is faithful and  exact.
Note that $\Cal A$ satisfies the (AB5) condition, as $F$ is 
faithful exact and colimit preserving,
and $\Cal B$ satisfies the (AB5) condition.

Note also that, for $a\in \Cal A$, the set of subobjects of $a$ is small,
because the set of subobjects of $Fa$ is small \cite{Tohoku} 
and $F$ is faithful exact.

Let $S$ be a small set of noetherian generators of $\Cal B$.
As any noetherian object is a quotient of a finite sum of objects in $S$,
we may assume that any noetherian object in $\Cal B$ is isomorphic to
an element of $S$, replacing $S$ by some larger small set, if necessary.
For each $s\in S$, the set of subobjects of $Gs$ is small by the last
paragraph.
Hence, there is a small subset $T$ of $\ob(\Cal A)$ such that,
any element $t\in T$ admits a monomorphism $t\rightarrow Gs$ for some
$s\in S$, $Ft$ is noetherian, and that 
if $a\in \Cal A$ admits a monomorphism
$a\rightarrow Gs$ for some $s\in S$ and $Fa$ noetherian then
$a\cong t$ for some $t\in T$.

We claim that any $a\in \Cal A$ 
is a filtered inductive limit $\indlim a_\lambda$
of subobjects $a_\lambda$
of $a$, with each $a_\lambda$ is isomorphic to some element in $T$.

If the claim is true, then {\bf 1} is obvious, as $T$ is a small set of
noetherian generators of $\Cal A$, and $\Cal A$ satisfies the (AB5) condition,
as we have already seen.

The `only if' part of {\bf 2} is also true if the claim is true,
since if $a\in\Cal A$ is noetherian, then it 
is a quotient of a finite sum of elements of $T$, and hence
$Fa$ is noetherian.

It suffices to prove the claim.
As $\Cal B$ is locally noetherian, we have $Fa=\indlim b_\lambda$,
where $(b_\lambda)$ is the filtered inductive system of noetherian subobjects
of $Fa$.

Let $u:\Id\rightarrow GF$ be the unit of adjunction,
and $\varepsilon:FG\rightarrow\Id$ be the counit of adjunction.
It is well-known that we have $(\varepsilon F)\circ(Fu)=\id_F$.
As $Fu$ is a split monomorphism, $u$ is also a monomorphism.
We define $a_\lambda:=u(a)^{-1}(G b_\lambda)$.
As $G$ preserves filtered inductive limits and $\Cal A$ satisfies the
(AB5) condition,
we have
\[
\indlim a_\lambda=u(a)^{-1}(G\indlim b_\lambda)=u(a)^{-1}(GFa)=a.
\]
Note that $a_\lambda\rightarrow G b_\lambda$ is a monomorphism,
with $b_\lambda$ being noetherian.

It remains to show that $Fa_\lambda$ is noetherian.
Let $i_\lambda:a_\lambda\hookrightarrow a$ be the inclusion map,
and $j_\lambda:b_\lambda\rightarrow Fa$ the inclusion.
Then the diagram
\[
\begin{array}{ccccc}
Fa & \specialarrow{Fu(a)} & FGFa & \specialarrow{\varepsilon F(a)}& Fa\\
\hbox to 0pt{\hss\small$Fi_\lambda$}\uparrow &  &
\uparrow \hbox to 0pt{\small$FGj_\lambda$\hss}& &
\uparrow \hbox to 0pt{\small$j_\lambda$\hss}\\
Fa_\lambda & \specialarrow{\hphantom{Fu(a)}} 
& FGb_\lambda & \specialarrow{\,\varepsilon (b_\lambda)\,}
& b_\lambda
\end{array}
\]
is commutative.
As the composite of the first row is the identity map and $Fi_\lambda$
is a monomorphism,
we have that the composite of the second row
$Fa_\lambda\rightarrow b_\lambda$ is a monomorphism.
As $b_\lambda$ is noetherian, we have that $Fa_\lambda$ is also noetherian,
as desired.
\qed

\paragraph
Let $S$ a scheme, and $X_\bullet\in \Cal P((\Delta),\Sch/S)$.

\begin{Lemma}\label{fe.thm}
The restriction functor $(?)_0:\EM(X_\bullet)\rightarrow\Mod(X_0)$ is
faithful exact.
\end{Lemma}

\proof This is obvious, because for any $[n]\in(\Delta)$, 
there is a morphism $[0]\rightarrow [n]$.

\begin{Lemma} Let $X_\bullet$ be a simplicial groupoid of $S$-schemes,
and assume that $d_0(1)$ and $d_1(1)$ are quasi-compact and separated.
Assume moreover that $\Qco(X_0)$ is locally noetherian.
Then we have
\begin{description}
\item[\bf 1] $\Qco(X_\bullet)$ is locally noetherian.
\item[\bf 2] $\Cal M\in\Qco(X_\bullet)$ is a noetherian object
if and only if $\Cal M_0$ is a noetherian object.
\end{description}
\end{Lemma}

\proof Let $F:=(?)_0:\Qco(X_\bullet)\rightarrow \Qco(X_0)$ be the
restriction.
By Lemma~\ref{fe.thm}, $F$ is faithful exact.
Let $G:=(d_0)_*^{\Qco}\circ \Bbb A$ be the right adjoint of $F$,
see Corollary~\ref{qco-right-adjoint.thm}.
As $\Bbb A$ is an equivalence and $(d_0)_*^{\Qco}$ preserves filtered
inductive limits by Lemma~\ref{direct-image-flim.thm},
$G$ preserves filtered inductive limits.
As $\Qco(X_\bullet)$ satisfies (AB3) by Lemma~\ref{indlim-qco.thm},
the assertion is obvious by Lemma~\ref{locally-noetherian.thm}.
\qed

The following is well-known, see \cite[pp.~126--127]{Hartshorne}.

\begin{Corollary}\label{qco-ln-single.thm}
Let $Y$ be a noetherian separated scheme.
Then $\Qco(Y)$ is locally noetherian, and $\Cal M\in\Qco(Y)$ is
a noetherian object if and only if it is coherent.
\end{Corollary}

\proof This is obvious if $Y=\Spec A$ is affine.
Now consider the general case.
Let $(U_i)_{1\leq i\leq r}$ be an affine open covering of $Y$, and
set $X:=\coprod_i U_i$.
Let $p:X\rightarrow Y$ be the canonical map, and set $X_\bullet:=
\Nerve(f)|_{(\Delta)}$.
Note that $p$ is faithfully flat, quasi-compact and separated.
By assumption and the lemma, we have that $\Qco(X_\bullet)$
is locally noetherian, and $\Cal M\in\Qco(X_\bullet)$ is a noetherian
object if and only if $\Cal M_0$ is noetherian, i.e., coherent.
As $\Bbb A:\Qco(Y)\rightarrow\Qco(X_\bullet)$ is an equivalence,
we have that $\Qco(Y)$ is locally noetherian, and
$\Cal M\in\Qco(Y)$ is a noetherian object if and only if
$(\Bbb A\Cal M)_0=p^*\Cal M$ is coherent if and only if $\Cal M$ is coherent.
\qed

Now the following is obvious.

\begin{Corollary}\label{ln-groupoid.thm}
Let $X_\bullet$ be a simplicial groupoid of $S$-schemes, 
with $d_0(1)$ and $d_1(1)$ quasi-compact and separated.
If $X_0$ is separated noetherian, then
$\Qco(X_\bullet)$ is locally noetherian, and $\Cal M\in\Qco(X_\bullet)$
is a noetherian object if and only if $\Cal M_0$ is coherent.
\end{Corollary}

\begin{Lemma} Let $I$ be a finite category, $S$ a scheme,
and $X_\bullet\in\Cal P(I,\Sch/S)$.
If $X_\bullet$ is separated noetherian, then $\Mod(X_\bullet)$ and
$\LQco(X_\bullet)$ are locally noetherian.

\end{Lemma}

\proof Let $J$ be the discrete subcategory of $I$ such that $\ob(J)=\ob(I)$.
Obviously, the restriction $(?)_J$ is faithful and exact.
For $i\in\ob(I)$, there is an isomorphism of functors
$(?)_iR_J\cong \prod_{j\in\ob(J)}\prod_{\phi\in I(i,j)} (X_\phi)_*(?)_j$.
The product is a finite product, as $I$ is finite.
As each $X_\phi$ is quasi-compact, $(X_\phi)_*(?)_j$ preserves filtered
inductive limits.
Hence $R_J$ preserves filtered inductive limits.
Note also that $R_J$ preserves local quasi-coherence.

Hence we may assume that $I$ is a discrete finite category, which
case is trivial by \cite[Theorem~II.7.8]{Hartshorne2} and
Corollary~\ref{qco-ln-single.thm}.
\qed

\section{Groupoid of schemes}

\paragraph
Let $\Cal C$ be a category with finite limits.
A $\Cal C$-groupoid $X_*$ is a functor from $\Cal C\op$ to the
category of $\Cal U$-small groupoids (i.e., $\Cal U$-small category
all of whose morphisms are isomorphisms) such that the set
valued functors
$X_0:=\ob\circ X_*$ and $X_1:=\Mor\circ X_*$ are representable.

Let $X_*$ be a $\Cal C$-groupoid.
Let us denote the source (resp.\ target) $X_1\rightarrow X_0$ by
$d_1$ (resp.\ $d_0$).
Then $X_2:=X_1\ttimes X_1$ represents the functor of
pairs $(f,g)$ of morphisms of $X_*$ such that $f\circ g$ is defined.
Let $d_0':X_2\rightarrow X_1$ (resp.\ $d_2':X_2\rightarrow X_1$) be the
first (resp.\ second) projection , and $d_1':X_2\rightarrow X_1$ be the
composition.

By Yoneda's lemma, $d_0$, $d_1$, $d_0'$, $d_1'$, and $d_2'$ are
morphisms of $\Cal C$.
$d_1:X_1(T)\rightarrow X_0(T)$ is surjective for any $T\in\ob(\Cal C)$.
Note that the squares
\begin{equation}\label{three.eq}
\begin{array}{ccccccccccc}
  X_2 & \specialarrow{d_0'} & X_1 &~~~~~&  
             X_2 & \specialarrow{d_1'} & X_1 &~~~~~& 
             X_2 & \specialarrow{d_1'} & X_1\\
  \sdarrow{d_2'} & & \sdarrow{d_1}&    &
             \sdarrow{d_2'} &        &\sdarrow{d_1} & &
             \sdarrow{d_0'} &  & \sdarrow{d_0}\\
  X_1 & \specialarrow{d_0} & X_0  &     & 
             X_1  & \specialarrow{d_1}  & X_0 &     & 
             X_1 & \specialarrow{d_0} & X_0
\end{array}
\end{equation}
are fiber squares.
In particular,
\begin{equation}\label{standard-mumford.eq}
X_*:=X_2
\begin{array}{c}
\mathop{\verylongrightarrow}\limits^{d_0'}\\
\mathop{\verylongrightarrow}\limits^{d_1'}\\
\mathop{\verylongrightarrow}\limits^{d_2'}
\end{array}
X_1
\begin{array}{c}
\mathop{\verylongrightarrow}\limits^{d_0}\\
\mathop{\verylongrightarrow}\limits^{d_1}
\end{array}
X_0
\end{equation}
forms an object of $\Cal P(\Delta_M,\Cal C)$.
Finally, by the associativity,
\begin{equation}\label{composition-associativity.eq}
\circ(\circ\times 1)=\circ(1\times\circ),
\end{equation}
where $\circ: X_1 \ttimes X_1\rightarrow X_1$ 
denotes the composition,
or $\circ$ is the composite
\[
X_1\ttimes X_1\cong X_2\specialarrow{d_1'}X_1.
\]

Conversely, a diagram $X_*\in \Cal P(\Delta_M,\Cal C)$ as in 
(\ref{standard-mumford.eq}) such that the squares in 
(\ref{three.eq}) are fiber squares, $d_1(T):X_1(T)\rightarrow X_0(T)$ 
are surjective for all $T\in\ob(\Cal C)$, and the associativity
(\ref{composition-associativity.eq}) holds gives a $\Cal C$-groupoid
\cite{SGA3-V}.
In the sequel, we mainly consider that a $\Cal C$-groupoid is an
object of $\Cal P(\Delta_M,\Cal C)$.

Let $S$ be a scheme.
We say that $X_\bullet\in \Cal P(\Delta_M,\Sch/S)$ is an $S$-groupoid,
if $X_\bullet$ is a $(\Sch/S)$-groupoid with flat arrows.

\paragraph Let $X_\bullet$ be a $(\Sch/S)$-groupoid,
and set $X_n:=X_1\ttimes X_1\ttimes\cdots
\ttimes X_1$ ($X_1$ appears $n$ times) for $n\geq 2$.
For $n\geq 2$, 
$d_i:X_n\rightarrow X_{n-1}$ is defined by
$d_0(x_{n-1},\ldots,x_1,x_0)=(x_{n-1},\ldots,x_1)$, 
$d_n(x_{n-1},\ldots,x_1,x_0)=(x_{n-2},\ldots,x_0)$, and
$d_i(x_{n-1},\ldots,x_1,x_0)=(x_{n-2},\ldots,x_i\circ x_{i-1},\ldots,x_0)$
for $0<i<n$.
$s_i:X_n\rightarrow X_{n+1}$ is defined by
$s_i(x_{n-1},\ldots,x_1,x_0)=(x_{n-1},\ldots,x_i,\id,x_{i-1},\ldots,x_0)$.

It is easy to see that this gives  a simplicial $S$-scheme 
$\Sigma(X_\bullet)$ such that $\Sigma(X_\bullet)|_{\Delta_M}=X_\bullet$.

For any simplicial $S$-scheme
$Z_\bullet$ and $\psi_\bullet:Z_\bullet|_{\Delta_M}\rightarrow X_\bullet$,
there exists some
unique $\varphi_\bullet:Z_\bullet\rightarrow \Sigma(X_\bullet)$
such that $\varphi|_{\Delta_M}:Z_\bullet|_{\Delta_M}\rightarrow 
\Sigma(X_\bullet)|_{\Delta_M}=X_\bullet$ equals $\psi$.
Indeed, $\varphi$ is given by
\[
\varphi_n(z)=(\psi_1(Q_{n-1}(z)),\ldots,\psi_1(Q_{0}(z))),
\]
where $q_i:[1]\rightarrow[n]$ is the injective monotone map such that
$\Image q_i=\{i,i+1\}$ for $0\leq i<n$, and $Q_i:Z_n\rightarrow Z_1$
is the associated morphism.
This shows that $\Sigma(X_\bullet)\cong \cosk_{\Delta_M}^{(\Delta)}X_\bullet$,
and the counit map $(\cosk_{\Delta_M}^{(\Delta)}X_\bullet)|_{\Delta_M}
\rightarrow X_\bullet$ is an isomorphism.

Note that under the identification $\Sigma(X_\bullet)_{n+1}\cong 
X_n \mathbin{{}_{r_0}\!\times_{d_0}} X_1$, the morphism $d_0:
\Sigma(X_\bullet)_{n+1}\rightarrow
\Sigma(X_\bullet)_{n}$ is nothing but the first projection.
So $(d_0):\Sigma(X_\bullet)'\rightarrow\Sigma(X_\bullet)$ is of fiber type.
If moreover $d_0(1)$ is flat, then $(d_0)$ is faithfully flat.

We construct an isomorphism 
$h_\bullet:
\Sigma(X_\bullet)'\rightarrow \Nerve(d_1(1))$.
Define $h_{-1}=\id$ and $h_0=\id$.
Define $h_1$ to be the composite
\[
X_1\ttimes X_1\specialarrow{(d_0\boxtimes d_2)^{-1}}X_2
\specialarrow{d_1\boxtimes d_2} X_1\dtimes X_1.
\]
Now define $h_n$ to be the composite
\lrsplit
X_1\ttimes X_1\ttimes \cdots\ttimes X_1\ttimes X_1
\specialarrow{1\times h_1}
X_1\ttimes X_1\ttimes \cdots\ttimes X_1\dtimes X_1\\
\specialarrow{\via h_1}
\cdots
\specialarrow{h_1\times 1}
X_1\dtimes X_1\dtimes\cdots\dtimes X_1\dtimes X_1.
\endlrsplit
It is straightforward to check that this gives a well-defined isomorphism
$h_\bullet:
\Sigma(X_\bullet)'\rightarrow \Nerve(d_1(1))$.
In conclusion, we have

\begin{Lemma}\label{groupoid-coskelton.thm}
If $X_\bullet$ is an $S$-groupoid, then 
$\cosk_{\Delta_M}^{(\Delta)}X_\bullet
$ is a simplicial $S$-groupoid, and the
counit $\varepsilon:(\cosk_{\Delta_M}^{(\Delta)}X_\bullet
)|_{\Delta_M}\rightarrow
X_\bullet$ is an isomorphism.
\end{Lemma}

Conversely, the following hold.

\begin{Lemma}
If 
$Y_\bullet$ is a simplicial $S$-groupoid, then $Y_\bullet|_{\Delta_M}$ is
an $S$-groupoid, and the unit map $u:Y_\bullet\rightarrow
\cosk_{\Delta_M}^{(\Delta)}(Y_\bullet|_{\Delta_M})$ is an
isomorphism.
\end{Lemma}

\proof It is obvious that $Y_\bullet|_{\Delta_M}$ has flat arrows.
So it suffices to show that $Y_\bullet|_{\Delta_M}$ is a 
$(\Sch/S)$-groupoid.
Since $Y_\bullet'$ is isomorphic to $\Nerve d_1$, the square
\[
\begin{array}{ccc}
  Y_2 &\specialarrow{d_2} & Y_1\\
\sdarrow{d_1} &  & \sdarrow{d_1}\\
  Y_1 & \specialarrow{d_1} & Y_0
\end{array}
\]
is a fiber square.
Since $(d_0):Y_\bullet'|_{(\Delta)}\rightarrow Y_\bullet$ is a
morphism of fiber type, the squares
\[
\begin{array}{ccccccc}
  Y_2 &\specialarrow{d_1} & Y_1 &
~~~~~~~~~~~~~ & Y_2 & \specialarrow{d_2} & Y_1\\
\sdarrow{d_0} &  & \sdarrow{d_0} &  & \sdarrow{d_0} & & \sdarrow{d_0}\\
  Y_1 & \specialarrow{d_0} & Y_0 &  & Y_1 & \specialarrow{d_1} & Y_0
\end{array}
\]
are fiber squares.

As $d_1s=\id$, $d_1(T):X_1(T)\rightarrow X_0(T)$ is surjective for
any $S$-scheme $T$.

Let us denote the composite
\[
X_1\ttimes X_1\cong X_2\specialarrow{d_1}X_1
\]
by $\circ$.
It remains to show the associativity.

As the three squares in the diagram
\[
\begin{array}{ccccc}
X_3 & \specialarrow{d_3} & X_2 & \specialarrow{d_2} & X_1\\
\sdarrow{d_0} & & \sdarrow{d_0} & & \sdarrow{d_0}\\
X_2 & \specialarrow{d_2} & X_1 & \specialarrow{d_1} & X_0\\
\sdarrow{d_0} & & \sdarrow{d_0}\\
X_1 & \specialarrow{d_1} & X_0
\end{array}
\]
are all fiber squares, the canonical map $Q:=Q_2\boxtimes Q_1\boxtimes Q_0:
X_3\rightarrow X_1\ttimes X_1\ttimes X_1$ is
an isomorphism.
So it suffices to show that the maps
\[
X_3\specialarrow{Q}
X_1\ttimes X_1\ttimes X_1
\specialarrow{\circ\times 1}
X_1\ttimes X_1
\specialarrow{\circ}
X_1
\]
and
\[
X_3\specialarrow{Q}
X_1\ttimes X_1\ttimes X_1
\specialarrow{1\times\circ}
X_1\ttimes X_1
\specialarrow{\circ}
X_1
\]
agree.
But it is not so difficult to show that the first map is $d_1d_2$, while
the second one is $d_1d_1$.
So $Y_\bullet|_{\Delta_M}$ is an $S$-groupoid.

Set $Z_\bullet:=\cosk_{\Delta_M}^{(\Delta)}(Y_\bullet|_{\Delta_M})$, and
we are to show that the unit $u_\bullet: Y_\bullet\rightarrow Z_\bullet$ 
is an isomorphism.
Since $Y_\bullet|_{\Delta_M}$ is an $S$-groupoid,
$\varepsilon_\bullet:Z_\bullet|_{\Delta_M}\rightarrow Y_\bullet|_{\Delta_M}$
is an isomorphism.
It follows that $u_\bullet|_{\Delta_M}:Y_\bullet|_{\Delta_M}
\rightarrow Z_\bullet|_{\Delta_M}$ is also an isomorphism.
Hence $\Nerve(d_1(1)(u_\bullet)):\Nerve(d_1(1)(Y_\bullet))
\rightarrow\Nerve(d_1(1)(Z_\bullet))$ is also an isomorphism.
As both $Y_\bullet$ and $Z_\bullet$ are simplicial $S$-groupoids by
Lemma~\ref{groupoid-coskelton.thm}, $u_\bullet':Y_\bullet'\rightarrow
Z_\bullet'$ is an isomorphism.
So $u_n:Y_n\rightarrow Z_n$ are all isomorphisms, and we are done.
\qed

\begin{Lemma}\label{ln-groupoid2.thm}
Let $S$ be a scheme, and $X_\bullet$ an $S$-groupoid,
with $d_0(1)$ and $d_1(1)$ quasi-compact and separated.
If $X_0$ is noetherian, then
$\Qco(X_\bullet)$ is locally noetherian, and $\Cal M\in\Qco(X_\bullet)$
is a noetherian object if and only if $\Cal M_0$ is coherent.
\end{Lemma}

\proof This is immediate by Corollary~\ref{ln-groupoid.thm} and
Lemma~\ref{simplicial-M.thm}.
\qed

\paragraph 
Let $f:X\rightarrow Y$ be a faithfully flat quasi-compact separated
$S$-morphism.
Set $X_\bullet^+:=(\Nerve(f))|_{\Delta_M^+}$
and $X_\bullet:=(X_\bullet^+)|_{\Delta_M}$.
We define the descent functor 
\[
\Bbb D: \LQco(X_\bullet)\rightarrow \Qco(Y)
\]
to be the composite $(?)_{[-1]}R_{\Delta_M}$.
The left adjoint $(?)_{\Delta_M}L_{[-1]}$ is denoted by $\Bbb A$, 
and called the ascent functor.

\begin{Lemma} Let the notation be as above.
Then $\Bbb D:\Qco(X_\bullet)\rightarrow \Qco(Y)$ is an equivalence, 
with $\Bbb A$ its quasi-inverse.
The composite
\[
\Bbb A\circ \Bbb D: \LQco(X_\bullet)\rightarrow \Qco(X_\bullet)
\]
is the right adjoint of the inclusion
$\Qco(X_\bullet)\hookrightarrow \LQco(X_\bullet)$.
\end{Lemma}

\proof Follows easily from
Corollary~\ref{ascent-descent.thm} and
Corollary~\ref{lqco-qco.thm}.
\qed

\begin{Lemma} Let $X_\bullet$ be an $S$-groupoid,
and assume that $d_0(1)$ and $d_1(1)$ are 
separated and quasi-compact.
Set $Y_\bullet^+:=((\cosk_{\Delta_M}^{(\Delta)}X_\bullet)')|_{\Delta_M^+}$,
and $Y_\bullet:=(Y_\bullet^+)|_{\Delta_M}$.
Let $(d_0):Y_\bullet\rightarrow X_\bullet$ be the canonical map
\[
Y_\bullet=((\cosk_{\Delta_M}^{(\Delta)}X_\bullet)')|_{\Delta_M}
\specialarrow{(d_0)|_{\Delta_M}}
(\cosk_{\Delta_M}^{(\Delta)}X_\bullet)|_{\Delta_M}
\cong X_\bullet.
\]
Then $(d_0)$ is quasi-compact separated faithfully flat of fiber type,
and
\[
(d_0)_*^{\Qco}\circ \Bbb A:\Qco(X_0)\rightarrow\Qco(X_\bullet)
\]
is a right adjoint of
$(?)_0:\Qco(X_\bullet)\rightarrow\Qco(X_0)$,
where $\Bbb A:\Qco(X_0)\rightarrow \Qco(Y_\bullet)$ is the
ascent functor.
\end{Lemma}

\proof Follows easily from Corollary~\ref{qco-right-adjoint.thm}.
\qed

\section{Comparison of local Ext sheaves}

\paragraph Let $S$ be a scheme, and $X_\bullet$ an $S$-groupoid.
Assume that $d_0(1)$ and $d_1(1)$ are affine faithfully flat of finite type,
and $X_0$ is noetherian.
Let $i:D(\Qco(X_\bullet))\rightarrow D(\Mod(X_\bullet))$ be the inclusion.
Note that $\Qco(X_\bullet)$ is locally noetherian 
(Lemma~\ref{ln-groupoid2.thm}).
So any complex in $C(\Qco(X_\bullet))$ admits a $K$-injective resolution
with each term an injective object.
The following is a generalization of \cite[Theorem~II.1.1.12]{Hashimoto}.

\begin{Lemma} Let $\Bbb F\in K_{\Coh(X_\bullet)}^-(\Mod(X))$ and $\Bbb G
\in K^+(\Qco(X_\bullet))$.
If $\Bbb G$ is a bounded below 
complex consisting of injective objects of $\Qco(X_\bullet)$, then
$\Bbb G$ is $\uHom_{\Cal O_{X_\bullet}}^\bullet(\Bbb F,?)$-acyclic.
\end{Lemma}

\proof By the way-out lemma, we may assume that $\Bbb F$ is a single
coherent sheaf, and 
$\Bbb G$ is a single injective object of $\Qco(X_\bullet)$.

To prove this case, it suffices to show that
\[
\uExt^i_{\Cal O_{X_\bullet}}(\Bbb F,\Bbb G)=0
\]
for $i>0$.

Let $\{U_1,\ldots,U_n\}$ be a finite affine open covering of $X_0$,
and set $Y:=\coprod_{i=1}^n U_i$.
Let $\varphi:Y\rightarrow X_0$ be the associated morphism.
Set $f': Y\times_{X_0}\,{}_{d_1(1)}X_1\rightarrow Y$ to be the first
projection, and set $Y_\bullet^+:=\Nerve(f')|_{\Delta_M^+}$
and $Y_\bullet:=\Nerve(f')|_{\Delta_M}$.
Let $\psi_\bullet:Y_\bullet\rightarrow X_\bullet'$ be the canonical map,
where $X_\bullet':=\Nerve(d_1(1))|_{\Delta_M}$.
As $(d_0)_*\Bbb A\varphi_*\cong (d_0)_*(\psi_\bullet)_*\Bbb A$ has
a faithful left adjoint $\varphi^*(?)_0$, there exists some injective object
$\Bbb I$ of $Y$ such that $\Bbb G$ is a direct summand
of $(d_0)_*(\psi_\bullet)_*\Bbb A\Bbb I$.
Set $q_\bullet:=(d_0)\circ \psi_\bullet$.
Note that $q_\bullet$ is affine faithfully flat of finite type and 
fiber type.
We may assume that $\Bbb G=(q_\bullet)_*\Bbb A\Bbb I$.
As
\[
R\uHom^\bullet_{\Cal O_{X_\bullet}}(\Bbb F,R(q_\bullet)_*\Bbb A\Bbb I)
\cong R(q_\bullet)_*
R\uHom^\bullet_{\Cal O_{Y_\bullet}}(q_\bullet^*\Bbb F,
\Bbb A\Bbb I),
\]
$q_\bullet$ is affine, and 
$R\uHom^\bullet_{\Cal O_{Y_\bullet}}(q_\bullet^*\Bbb F,
\Bbb A\Bbb I)$ has quasi-coherent cohomology groups,
we may assume that $X_\bullet=\Nerve(f)|_{\Delta_M}$ 
for some faithfully flat morphism $f:X\rightarrow Y$ of finite type
between affine noetherian $S$-schemes.
For each $l$, we have that
\begin{multline*}
(?)_l R\uHom^\bullet_{\Cal O_{X_\bullet}}(\Bbb F,\Bbb G)
\cong
R\uHom^\bullet_{\Cal O_{X_l}}((?)_l\Bbb A\Bbb D\Bbb F,(?)_l\Bbb A\Bbb D
\Bbb G)
\cong\\
R\uHom^\bullet_{\Cal O_{X_l}}(e(l)^*\Bbb D\Bbb F,e(l)^*\Bbb D\Bbb G)
\cong
e(l)^*R\uHom^\bullet_{\Cal O_Y}(\Bbb D\Bbb F,\Bbb D\Bbb G).
\end{multline*}
As $\Bbb D:\Qco(X_\bullet)\rightarrow \Qco(Y)$ is an equivalence, 
$\Bbb D\Bbb G$ is an injective object of $\Qco(Y)$.
Hence $\uExt^i_{\Cal O_Y}(\Bbb D\Bbb F,\Bbb D\Bbb G)=0$ for $i>0$
by Lemma~\ref{I-B-N.thm},
as desired.
\qed

\section{Group schemes flat of finite type}

\paragraph Let $S$ be a noetherian separated scheme.

\paragraph 
Let $\Cal F$ (resp.\ $\Cal F_M$)
denote the subcategory of $\Cal P((\Delta),\Sch/S)$
(resp.\ subcategory of $\Cal P(\Delta_M,\Sch/S)$)
consisting of objects separated of finite type over $S$ with flat arrows
and morphisms of finite type and of fiber type.

Let $G$ be a flat $S$-group scheme separated of finite type.
Note that $G$ is faithfully flat over $S$.
Set $\Cal A_G$ to be the category of
$S$-schemes separated of finite type with left $G$-actions.

For $X\in \Cal A_G$, we associate a simplicial scheme 
$B_G(X)$ by $B_G(X)_n=G^{n}\times X$.
For $n\geq 1$,
$d_i(n):G^{n}\times X\rightarrow G^{n-1}\times X$ is 
the projection $p\times 1_{G^{n-1}\times X}$ if
$i=n$, where $p:G\rightarrow S$ is the structure morphism.
While $d_i(n)=1_{G^{n-1}}\times a$ if $i=0$, where $a:G\times X\rightarrow X$ 
is the action.
If $0<i<n$, then $d_i(n)= 1_{G^{n-i-1}}\times
\mu\times 1_{G^{i-1}\times X}$,
where $\mu:G\times G\rightarrow G$ is the product.
For $n\geq 0$, 
$s_i(n):G^n\times X\rightarrow G^{n+1}\times X$
is given by
\[
s_i(n)(g_n,\ldots,g_1,x)=(g_n,\ldots,g_{i+1},e,g_i,\ldots,g_1,x),
\]
where $e:S\rightarrow G$ is the unit element.

Note that $(B_G(X)')|_{(\Delta)}$ is canonically isomorphic to 
$B_G(G\times X)$, where $G\times X$ is viewed as a principal $G$-action.
Note also that there is an isomorphism from $B_G(X)'$ 
to $\Nerve(p_2:G\times X\rightarrow X)$ given by
\begin{multline*}
B_G(X)'_n=G^{n+1}\times X\rightarrow (G\times X)\times_X\cdots\times_X
(G\times X)=\Nerve(p_2)_n
\\
(g_n,\ldots,g_0,x)\mapsto((g_n\cdots g_0,x),(g_{n-1}\cdots g_0,x),\ldots,
(g_0,x)).
\end{multline*}
Hence we have

\begin{Lemma}\label{group-groupoid.thm}
Let $G$ and $X$ be as above.
Then $B_G(X)$ is a simplicial $S$-groupoid with 
$d_0(1)$ and $d_1(1)$ 
faithfully flat separated of finite type.
\end{Lemma}

We denote the restriction $B_G(X)|_{\Delta_M}$ by $B_G^M(X)$.
Obviously, $B_G^M(X)$ is an $S$-groupoid with
$d_0(1)$ and $d_1(1)$ 
faithfully flat separated of finite type.

For a morphism $f:X\rightarrow Y$ in $\Cal A_G$, we define 
$B_G(f):B_G(X)\rightarrow B_G(Y)$ by $(B_G(f))_n=1_{G^n}\times f$.
It is easy to check that $B_G$ is a functor from $\Cal A_G$ to $\Cal F$.
Thus $B_G^M$ is a functor from $\Cal A_G$ to $\Cal F_M$.

By definition, a $(G,\Cal O_X)$-module (or a $G$-linearlized 
$\Cal O_X$-module \cite{GIT} or $G$-equivariant $\Cal O_X$-module) 
is an equivariant $\Cal O_{B_G^M(X)}$-module,
see \cite{BL} and \cite{Hashimoto}.
Note that the category of equivariant $\Cal O_{B_G(X)}$-modules and
the category of equivariant $\Cal O_{B_G^M(X)}$-modules are equivalent
(Lemma~\ref{simplicial-M.thm}).
Thus a quasi-coherent (resp.\ coherent) $(G,\Cal O_X)$-module is nothing
but a quasi-coherent (resp.\ coherent) $\Cal O_{B_G^M(X)}$-module.
We denote the category of $(G,\Cal O_X)$-modules by $\EM(G,X)$.
The category of quasi-coherent (resp.\ coherent) $(G,\Cal O_X)$-modules
is denoted by $\Qco(G,X)$ (resp.\ $\Coh(G,X)$).

The discussion on derived categories of categories of sheaves over
diagrams of schemes are interpreted to the derived categories of 
the categories of $G$-equivariant modules.

By Lemma~\ref{group-groupoid.thm} and
Lemma~\ref{ln-groupoid2.thm}, we have

\begin{Lemma} Let $X\in\Cal A_G$.
Then $\Qco(G,X)$ is a locally noetherian abelian category.
\end{Lemma}

Let $\Cal M$ be a $(G,\Cal O_X)$-module.
If there is no danger of confusion, we may write $\Cal M_0$ 
instead of $\Cal M$.
For example, $\Cal O_X$ means $\Cal O_{B_G^M(X)}$,
since $\Cal O_{B_G^M(X)}$ is equivariant and $(\Cal O_{B_G(X)})_0=\Cal O_X$.
This abuse of notation is what we always do when $S=X=\Spec k$
and $G$ is an affine algebraic group over $k$.
A $G$-module and its underlying vector space are denoted by the same symbol.
Similarly, an object of $D(\Mod(B_G^M(X)))$ and its restriction to
$B_G^M(X)_{0}=X$ are sometimes denoted by the same symbol.
Moreover, for a morphism $f$ in $\Cal A_G$, we denote for example 
$R(B_G^M(f))_*$ by $Rf_*$, and $B_G^M(f)^!$ by $f^!$.

Thus, as a corollary to Theorem~\ref{egd.thm}, we have

\begin{Theorem}[$G$-Grothendieck's duality]\label{ggd.thm}
Let $S$ be a noetherian separated scheme, and 
$G$ a flat $S$-group scheme separated of finite type.
Let $X$ and $Y$ be $S$-schemes separated of finite type with $G$-actions,
and $f:X\rightarrow Y$ a proper $G$-morphism.
Then the composite
\begin{multline*}
\Theta(f):
R f_* R\uHom_{\Mod(B_G^M(X))}^\bullet(\Bbb F,f^\times \Bbb G)
\specialarrow{H}
R\uHom_{\Mod(B_G^M(Y))}^\bullet(R f_*\Bbb F,R f_*f^\times \Bbb G)\\
\specialarrow{\varepsilon}
R\uHom_{\Mod(B_G^M(Y))}^\bullet(R f_*\Bbb F,\Bbb G)
\end{multline*}
is an isomorphism in $D^+_{\Qco(G,Y)}(B_G^M(Y))$
for $\Bbb F\in D_{\Coh(G,X)}^-(B_G^M(X))$ and
$\Bbb G\in  D^+_{\Qco(G,Y)}(B_G^M(Y))$.
\end{Theorem}

\section{Compatibility with derived $G$-invariance}

\paragraph\label{G-invariance.par}
Let $S$ be a separated noetherian scheme, and $G$ a flat separated
$S$-group scheme of finite type.
Let $X$ be an $S$-scheme separated of finite type with a trivial $G$-action.
That is, $a:G\times X\rightarrow X$ agrees with the second projection $p_2$.
In other words, $d_0(1)=d_1(1)$ in $B_G^M(X)$.

For an object $\Cal M$ of $\Mod(B_G^M(X))$, we define the {\em $G$-invariance}
of $\Cal M$ to be the kernel of the natural map
\[
\beta_{d_0(1)}-\beta_{d_1(1)}:
\Cal M_{0}\rightarrow d_0(1)_*\Cal M_{1}=
d_1(1)_*\Cal M_{1},
\]
and we denote it by $\Cal M^G$.

\paragraph Let $X$ be as in (\ref{G-invariance.par}).
Define $\tilde B_G^M(X)$ to be the augmented diagram
\[
G\times_S G\times_S X
\begin{array}{c}
\mathop{\verylongrightarrow}\limits^{1_G\times a}\\
\mathop{\verylongrightarrow}\limits^{\mu\times 1_X}\\
\mathop{\verylongrightarrow}\limits^{p_{23}}
\end{array}
G\times_S X
\begin{array}{c}
\mathop{\verylongrightarrow}\limits^{a}\\
\mathop{\verylongrightarrow}\limits^{p_2}
\end{array}
X
\,\specialarrow{\id}\,
X.
\]
Note that $\tilde B_G^M(X)$ is an object of $\Cal P(\tilde \Delta_M,S)$.

\begin{Lemma}\label{invariance-right-induction.thm}
The functor $(?)^G:\Mod(B_G^M(X))\rightarrow \Mod(X)$
agrees with $(?)_{-1}R_{\Delta_M}$.
\end{Lemma}

\proof Follows easily from (\ref{right-induction.par}).
\qed

\paragraph We say that an object $\Cal M$ of $B_G^M(X)$ is $G$-trivial if
$\Cal M$ is equivariant, and the canonical inclusion $\Cal M^G\rightarrow
\Cal M_{0}$ is an isomorphism.
Note that $(?)_{\Delta_M} L_{-1}$ is the exact left adjoint of 
$(?)^G$.
Note also that $\Cal M$ is $G$-trivial if and only if 
the counit of adjunction $\varepsilon: (?)_{\Delta_M} L_{-1}\Cal M^G
\rightarrow \Cal M$ is an isomorphism if and only if $\Cal M\cong
\Cal N_{\Delta_M}$ for some $\Cal N\in\EM_{\tilde B_G^M(X)}$.

Let $\triv(G,X)$ denote the full subcategory of $\Mod(B_G(X))$ consisting
of $G$-trivial objects.
Note that $(?)^G: \triv(G,X)\rightarrow \Mod(X)$ is an equivalence, 
whose quasi-inverse is $(?)_{\Delta_M}L_{-1}$.

Note that if $\Cal M$ is locally quasi-coherent, then $\Cal M^G$ is
quasi-coherent.
Thus, we get a derived functor $R(?)^G:D^+_{\LQco(B_G^M(X))}(\Mod(B_G^M(X)))
\rightarrow D^+_{\Qco(X)}(\Mod(X))$.

\begin{Proposition}\label{twisted-G-cohomology.thm}
Let $X$ and $Y$ be $S$-schemes 
separated of finite type with trivial $G$-actions,
and $f:X\rightarrow Y$ an $S$-morphism, which is automatically a 
$G$-morphism.
Then there is a canonical isomorphism
\[
f^! R(?)^G\cong R(?)^G f^!
\]
between functors from $D_{\LQco}^+(B_G^M(Y))$ to $D_{\Qco}^+(X)$.
\end{Proposition}

\proof As $(?)_{-1}$ is exact, we have $R(?)^G\cong (?)_{-1}RR_{\Delta_M}$ 
by Lemma~\ref{invariance-right-induction.thm}.
Thus, we have a composite isomorphism
\[
f^! R(?)^G\cong f^! (?)_{-1}RR_{\Delta_M}
\specialarrow{\bar\xi^{-1}}
(?)_{-1}(\tilde B^G_M(f))^! RR_{\Delta_M}
\specialarrow{\bar\xi}
(?)_{-1}RR_{\Delta_M}f^!
\cong
R(?)^G f^!
\]
by Proposition~\ref{twisted-restrict.thm} and 
Lemma~\ref{bar-xi-rind.thm}.
\qed

\section{Equivariant dualizing complexes and canonical modules}

\paragraph Let $I$ be a finite ordered category, $S$ a scheme, 
and $X_\bullet\in\Cal P(I,\Sch/S)$.

\begin{Lemma} Assume that $X_\bullet$ has flat arrows.
Then an object $\Bbb I\in K(\Mod(X_\bullet))$ is $K$-injective
if and only if $\Bbb I_i$ is $K$-injective for any $i\in\ob(I)$.
\end{Lemma}

\proof The \lq only if' part is obvious, since $(?)_i$ has an exact left
adjoint $L_i$.

We prove the converse by induction on the number of objects of $I$.
We may assume that $I$ has at least two objects.
Let $i$ be a minimal element of $\ob(I)$, and set $J:=\ob(I)\setminus
\{i\}$.
Let $\Bbb F$ be an exact complex in $K(\Mod(X_\bullet))$.
Consider a distinguished triangle
\[
L_i\Bbb F_i\specialarrow{\varepsilon}\Bbb F\rightarrow \Bbb C
\specialarrow{(1)}\Sigma L_i\Bbb F_i.
\]
As
\[
\Hom_{\Mod(X_\bullet)}^\bullet(L_i\Bbb F_i,\Bbb I)
\cong
\Hom_{\Mod(X_i)}^\bullet(\Bbb F_i,\Bbb I_i)
\]
is exact by assumption, 
it suffices to show that $\Hom_{\Mod(X_\bullet)}^\bullet(\Bbb C,\Bbb I)$
is exact.
Note that $\Bbb C$ is exact.

Consider the distinguished triangle
\[
L_J\Bbb C_J\specialarrow{\varepsilon}\Bbb C\specialarrow{u} R_i\Bbb C_i
\specialarrow{(1)}\Sigma L_J\Bbb C_J.
\]
As
\[
\Hom_{\Mod(X_\bullet)}^\bullet(L_J\Bbb C_J,\Bbb I)
\cong
\Hom_{\Mod(X_\bullet|_J)}^\bullet(\Bbb C_J,\Bbb I_J),
\]
is exact by induction assumption, it suffices to show that
$\Hom_{\Mod(X_\bullet)}(R_i\Bbb C_i,\Bbb I)$ is exact.
To verify this, it suffices to show that $\Bbb C_i$ is null-homotopic.
But this is trivial, since
\[
(?)_i L_i\Bbb F_i\specialarrow{(?)_i\varepsilon}\Bbb F_i\rightarrow \Bbb C_i
\specialarrow{(1)}(?)_i\Sigma L_i\Bbb F_i
\]
is a distinguished triangle, and $(?)_i\varepsilon$ is an isomorphism.
\qed

\begin{Corollary}\label{finite-injective1.thm}
Let $I$, $S$, and $X_\bullet$ be as in the lemma.
Let $J$ be a full subcategory of $I$ such that $\ob(J)$ is a filter
of $\ob(I)$.
If $\Bbb J$ is a $K$-injective complex of $K(\Mod(X_\bullet|_{J}))$,
then so is $L_{J}\Bbb J$.
\end{Corollary}

\begin{Corollary}\label{fid.thm}
Let $I$, $S$, and $X_\bullet$ be as in the lemma.
Let $\Bbb F$ be an object of $K(\Mod(X_\bullet))$.
Then $\Bbb F$ has finite injective dimension if and only if 
$\Bbb F_i$ has finite injective dimension for any $i\in\ob(I)$.
\end{Corollary}

\proof We prove the corollary by induction on the number of objects of $I$.
We may assume that $I$ has at least two objects.
Let $i$ be a minimal element of $\ob(I)$, and $J$ the full subcategory
whose object set is $\ob(I)\setminus\{i\}$.
Then
\[
L_J\Bbb F_J\specialarrow{\varepsilon}\Bbb F\specialarrow{u}
R_i\Bbb F_i
\specialarrow{(1)}\Sigma L_J\Bbb F_J
\]
is a distinguished triangle.
As $\Bbb F_J$ has finite injective dimension by induction assumption,
$L_J\Bbb F_J$ has finite injective dimension by 
Corollary~\ref{finite-injective1.thm}.
Since $R_i$ is exact and has an exact left adjoint, 
$R_i\Bbb F_i$ has finite injective dimension.
Hence $\Bbb F$ also has finite injective dimension.
\qed

\paragraph\label{dualizing-setting.par}
Let the notation be as in Theorem~\ref{finite-data.thm}.
Let $X_\bullet$ be an object of $\Cal F$ (i.e., an $I\op$-diagram of
separated $S$-schemes of finite type with flat arrows).
We say that $\Bbb F\in D(\Mod(X_\bullet))$ is a {\em dualizing complex} of 
$X_\bullet$ if $\Bbb F\in D^b_{\Coh(X_\bullet)}(\Mod(X_\bullet))$, 
$\Bbb F$ has
finite injective dimension, and the canonical map
\[
\Cal O_{X_\bullet}\specialarrow{\trace}R\uHom^\bullet_{\Cal O_{X_\bullet}}
(\Bbb F,\Bbb F)
\]
is an isomorphism.

\begin{Lemma}\label{dualizing-restrict.thm}
Let the notation be as in {\rm(\ref{dualizing-setting.par})}.
An object $\Bbb F\in D(\Mod(X_\bullet))$ is a dualizing complex of
$X_\bullet$ if and only if $\Bbb F$ has equivariant cohomology groups
and $\Bbb F_i\in D(\Mod(X_i))$ is a dualizing complex of $X_i$ for 
any $i\in\ob(I)$.
\end{Lemma}

\proof This is obvious by Lemma~\ref{H_J-isom.thm} and
Corollary~\ref{fid.thm}.
\qed

\begin{Corollary}\label{Gor-dualizing.thm}
Let the notation be as in {\rm (\ref{dualizing-setting.par})}.
If $X_\bullet$ is Gorenstein with finite Krull dimension, then 
$\Cal O_{X_\bullet}$ is a dualizing complex of $X_\bullet$.
\end{Corollary}

\proof This is clear by the lemma and \cite[(V.10)]{Hartshorne2}.
\qed

\begin{Lemma} Let the notation be as in {\rm (\ref{dualizing-setting.par})}.
If $X_\bullet$ has a dualizing complex $\Bbb F$, then $X_\bullet$ has
finite Krull dimension, and $X_\bullet$ has Gorenstein arrows.
\end{Lemma}

\proof As $\Bbb F_i$ is a dualizing complex of $X_i$ for each $i\in\ob(I)$,
$X_\bullet$ has finite Krull dimension by \cite[Corollary~V.7.2]{Hartshorne2}.

Let $\phi:i\rightarrow j$ be a morphism of $I$.
As $X_\phi$ is flat, $\alpha_\phi:X_\phi^*\Bbb F_i\rightarrow \Bbb F_j$ 
is an isomorphism of $D(\Mod(X_j))$.
As $X_\phi^*\Bbb F_i$ is a dualizing complex of $X_j$, $X_\phi$ is
Gorenstein by \cite[(5.1)]{AF}.
\qed

\begin{Proposition} Let the notation be as above, and $\Bbb I$ a
dualizing complex of $X_\bullet$.
Let $\Bbb F\in D_{\Coh(X_\bullet)}(\Mod(X_\bullet))$.
Then we have 
\begin{description}
\item[1] $R\uHom^\bullet_{\Cal O_{X_\bullet}}(\Bbb F,\Bbb I)\in
D_{\Coh(X_\bullet)}(\Mod(X_\bullet))$.
\item[2] The canonical map
\[
\Bbb F\rightarrow R\uHom^\bullet_{\Cal O_{X_\bullet}}(
  R\uHom^\bullet_{\Cal O_{X_\bullet}}(\Bbb F,\Bbb I),\Bbb I)
\]
is an isomorphism for $\Bbb F\in D_{\Coh(X_\bullet)}(\Mod(X_\bullet))$.
\end{description}
\end{Proposition}

\proof {\bf 1} As $\Bbb I$ has finite finite injective dimension,
$R\uHom^\bullet_{\Cal O_{X_\bullet}}(?,\Bbb I)$ is way-out in
both directions.
Hence by \cite[Proposition~I.7.3]{Hartshorne2}, we may assume that
$\Bbb F$ is bounded.
This case is trivial by Lemma~\ref{H_J-isom.thm} and 
Lemma~\ref{ext-equivariant.thm}.

{\bf 2} Using Lemma~\ref{H_J-isom.thm} twice, we may assume that
$X_\bullet$ is a single scheme.
This case is \cite[Proposition~V.2.1]{Hartshorne2}.
\qed

\begin{Lemma} Let the notation be as in {\rm (\ref{dualizing-setting.par})}.
Let $f_\bullet:X_\bullet\rightarrow Y_\bullet$ be a morphism in 
$\Cal F$, and let $\Bbb I$ be a dualizing complex of $Y_\bullet$.
Then $f_\bullet^!(\Bbb I)$ is a dualizing complex of $X_\bullet$.
\end{Lemma}

\proof By Corollary~\ref{preservation.thm}, $f_\bullet^!(\Bbb I)$ has
coherent cohomology groups.
By Lemma~\ref{fid.thm}, $f_\bullet^!(\Bbb I)$ has finite injective dimension.

By Lemma~\ref{dualizing-restrict.thm} and
Proposition~\ref{twisted-restrict.thm},
we may assume that $f:X\rightarrow Y$ is a morphism of single schemes,
and it suffices to show that 
\[
\Cal O_{X}\specialarrow{\trace}R\uHom^\bullet_{\Cal O_{X}}
(f^!\Bbb I,f^!\Bbb I)
\]
is an isomorphism.
As the question is local both on $Y$ and $X$,
We may assume that both $Y$ and $X$ are affine, and $f$ is either an
affine $n$-space or a closed immersion.
This case is done in \cite[Chapter~V]{Hartshorne2}.
\qed

\begin{Lemma} Let the notation be as in {\rm (\ref{dualizing-setting.par})},
and $\Bbb I$ and $\Bbb J$ dualizing complexes on $X_\bullet$.
If $X_\bullet$ is connected and $X_i$ is non-empty for some $i\in\ob(I)$, 
then there exist a unique invertible sheaf
$\Cal L$ and a unique integer $n$ such that
\[
\Bbb J\cong \Bbb I\otimes^{\bullet,L}_{\Cal O_{X_\bullet}}\Cal L[n].
\]
Such $\Cal L$ and $n$ are determined by
\[
\Cal L[n]\cong R\uHom^\bullet_{\Cal O_{X_\bullet}}(\Bbb I,\Bbb J).
\]
\end{Lemma}

\proof Use \cite[Theorem~V.3.1]{Hartshorne2}.

\begin{Def}
Let the notation be as in (\ref{dualizing-setting.par}),
and $\Bbb I$ a fixed dualizing complex of $X_\bullet$.
For any object $f_\bullet:Y_\bullet\rightarrow X_\bullet$ of
$\Cal F/X_\bullet$, we define {\em the} dualizing complex of $Y_\bullet$
(or better, of $f_\bullet$) to be $f_\bullet^!\Bbb I$.
If $Y_\bullet$ is connected and $Y_i$ is non-empty for some $i\in\ob(I)$,
then we define the {\em canonical sheaf} $\omega_{Y_\bullet}$
of of $Y_\bullet$ (or better, $f_\bullet$) to be $H^s(f_\bullet^!\Bbb I)$,
where $s$ is the smallest $i$ such that $H^i(f_\bullet^!\Bbb I)\neq 0$.
If $Y_\bullet$ is disconnected, then we define $\omega_{Y_\bullet}$
componentwise.
\end{Def}

Let $R$ be a noetherian ring, 
$S:=\Spec R$, and $G$ a flat $R$-group scheme separated of finite type.

\begin{Lemma}\label{c-i.thm}
$G\rightarrow S$ is a \(flat\) local complete intersection morphism.
That is, \(it is flat and\) all fibers are locally complete intersection.
\end{Lemma}

\proof 
We may assume that $S=\Spec k$, with $k$ a field.
Then by \cite[Theorem~1]{Avramov}, we may assume that $k$ is algebraically
closed.

First assume that the characteristic is $p>0$.
Then there is some $r\gg0$ such that the image of the Frobenius map
$F^r: G\rightarrow G^{(r)}$ is reduced (or equivalently, $k$-smooth)
and agrees with $G^{(r)}_{\red}$.
Note that the induced morphism $G\rightarrow G^{(r)}_{\red}$ is flat,
since the flat locus is a $G$-stable open subset of $G$ 
\cite[Lemma~2.1.10]{Hashimoto}, and the morphism is flat at the generic point.

As the group $G_{\red}$ acts on $G$ transitively, it suffices to show that
$G$ is locally a complete intersection at the unit element $e$.
So by \cite[Theorem~2]{Avramov}, it suffices to show that the 
$r$th Frobenius kernel $G_r$ is a complete intersection.
As $G_r$ is finite connected, this is well known \cite[(14.4)]{Waterhouse}.

Now consider the case that $G$ is of characteristic zero.
We are to prove that $G$ is $k$-smooth.
Take a finitely generated $\Bbb Z$-subalgebra $R$ of $k$ 
such that $G$ is defined.
We may take $R$ so that $G_R$ is $R$-flat.
Set $H:=(G_R)_{\red}$.
We may take $R$ so that $H$ is also $R$-flat.
Then $H$ is a closed subgroup scheme over $R$, since $H\times_R H$ is
reduced.
As a reduced group scheme over a field of characteristic zero is smooth,
we may localize $R$ if necessary, and we may assume that $H$ is $R$-smooth.

Let $\Cal J$ be the defining ideal sheaf of $H$ in $G_R$.
There exists some $s\geq 0$ such that $\Cal J^{s+1}=0$.
Note that $\Cal G:=\bigoplus_{i=0}^{s}\Cal J^i/\Cal J^{i+1}$
is a coherent $(H,\Cal O_H)$-module.
Applying Corollary~\ref{ascent-descent.thm} to the case that
$Y=\Spec R$  and $X_\bullet=B_H^M(H)$, the coherent $(H,\Cal O_H)$-module
$\Cal G$
is of the form $f^*(V)$, where $V$ is a finite $R$-module,
and $f:H\rightarrow \Spec R$ is the structure map.
Replacing $R$ if necessary, we may assume that $V\cong R^u$.
Now we want to prove that $u=1$ so that $H=G_R$, which implies $G$ is
$k$-smooth.

There exists some prime number 
$p>u$ and a maximal ideal $\frak m$ of $R$ such that
$R/\frak m$ is a finite field of characteristic $p$.
Let $\kappa$ be the algebraic closure of $R/\frak m$, and consider the base
change $(\bar ?):=?\otimes_R \kappa$.
Note that $\bar {\Cal G}
=\bigoplus_i\bar{\Cal J}^i/\bar {\Cal J}^{i+1}$, since $R$ is $\Bbb Z$-flat
and $H$ and $G_R$ are $R$-flat.
Let $\Cal I^{[p^r]}$ denote 
the defining ideal of the $r$th Frobenius kernel of $\bar G_R$.
By \cite[(14.4)]{Waterhouse} again, $\dim_k (\Cal O_{\bar G_R}/\Cal I^{[p^r]})
_e$ is a power of $p$, say $p^{v(r)}$.
Similarly, the $k$-dimension of 
the coordinate ring $(\Cal O_{\bar G_R}/(\bar {\Cal J}+\Cal I^{[p^r]}))_e$ 
of the $r$th Frobenius kernel of $\bar H$ is a power of $p$, say 
$p^{w(r)}$.
Note that
\[
p^{w(r)}\leq p^{v(r)}\leq \dim_k
(\Cal O_{\bar G_R}/\Cal I^{[p^r]}\otimes_{\Cal O_{\bar G_R}}\bar {\Cal G})_e
= p^{w(r)}u < p^{w(r)+1}.
\]
Hence $w(r)\leq v(r)<w(r)+1$, and we have $\bar 
{\Cal J}_e\subset \Cal I^{[p^r]}_e$ 
for any $r$.
By Krull's intersection theorem, $\bar {\Cal J}_e\subset
\bigcap_r (\Cal I_e)^{p^r}=0$.
This shows that $\bar G_R$ is reduced at $e$, which shows that
$\bar G_R$ is $\kappa$-smooth everwhere.
So the nilpotent ideal $\bar {\Cal J}$ must be zero, and this shows $u=1$.
\qed

Now assume that $R$ is a Gorenstein local ring of dimension $d$.

\begin{Lemma}  $B_G(S)$ is Gorenstein of finite Krull dimension.
In particular, $\Cal O_S[d]$ is an equivariant dualizing complex of $X$
\(i.e., $\Cal O_{B_G(S)}$ is a dualizing complex of $B_G(S)$\).
\end{Lemma}

\proof As $S=\Spec R$ is Gorenstein by assumption and
$G$ is Gorenstein over $S$ by Lemma~\ref{c-i.thm}, the assertions are
trivial.
\qed

\begin{Def} We choose and fix the equivariant dualizing complex $\Cal O_S[d]$.
For an object $X\in\Cal A_G$, 
the corresponding equivariant dualizing complex is called the equivariant 
dualizing complex of $X$.
The corresponding equivariant canonical sheaf is called the equivariant
canonical sheaf of $X$, and is denoted by $\omega_X$.
\end{Def}

\begin{Lemma}\label{G-trivial-dualizing.thm}
Let $X\in\Cal A_G$, and assume that $G$ acts on $X$ trivially.
Then the dualizing complex $\Bbb I_X:=f^!(\Cal O_S[d])$ has $G$-trivial
cohomology groups, where $f:X\rightarrow S$ is the structure map.
In particular, $\omega_X$ is $G$-trivial.
\end{Lemma}

\proof By Proposition~\ref{twisted-restrict.thm}, 
\[
f^!(\Cal O_S[d])\cong f^!((\Cal O_{\tilde B_G^M(X)})_{B_G^M(X)})
\cong
(?)_{B_G^M(X)}(\tilde B_G^M(f)^!(\Cal O_{\tilde B_G^M(X)})).
\]
By Corollary~\ref{preservation.thm}, 
$\tilde B_G^M(f)^!(\Cal O_{\tilde B_G^M(X)})$ has coherent cohomology 
groups.
Hence, $f^!(\Cal O_S[d])$ has $G$-trivial cohomology groups.
\qed

\section{A generalization of Watanabe's theorem}

\begin{Lemma}\label{connected-inv.thm}
Let $R$ be a noetherian commutative ring, and 
$G$ a finite group which acts on $R$.
Set $A=R^G$, and assume that $\Spec A$ is connected.
Then $G$ permutes the connected components of $R$ transitively.
\end{Lemma}

\proof Since $\Spec R$ is a noetherian space, 
$\Spec R$ has only finitely many connected components, say $X_1,\ldots,X_n$.
Then $R=R_1\times\cdots\times R_n$, and each $R_i$ is of the form $Re_i$,
where $e_i$ is a primitive idempotent.
Note that $E:=\{e_1,\ldots,e_n\}$ is the set of primitive idempotents of $R$,
and $G$ acts on $E$.
Let $E_1$ be an orbit of this action.
Then $e=\sum_{e_i\in E_1}e_i$ is in $A$.
As $A$ does not have any nontrivial idempotent, $e=1$.
This shows that $G$ acts on $E$ transitively, and we are done.
\qed

\begin{Lemma}\label{transitive.thm}
Let $R$ be a noetherian commutative ring, and 
$G$ a finite group which acts on $R$.
Set $A=R^G$, and assume that the inclusion 
$A\hookrightarrow R$ is finite.
If $\frak p\in\Spec A$, then $G$ acts transitively on the set of
primes of $R$ lying over $\frak p$.
\end{Lemma}

\proof Note that $A$ is noetherian by Eakin-Nagata theorem 
\cite[Theorem~3.7]{CRT}.
Let $A'$ be the $\frak pA_{\frak p}$-adic completion of $A_{\frak p}$,
and set $R':=A'\otimes_A R$.
As $A'$ is $A$-flat, $A'=(R')^G$.
It suffices to prove that $G$ acts transitively on the maximal ideals of
$R'$.
But $R'$ is the direct product $\prod_i R'_i$ of complete local rings $R'_i$.
Consider the corresponding primitive idempotents.
Since $A'$ is a local ring, $G$ permutes these idempotents
transitively by Lemma~\ref{connected-inv.thm}.
It is obvious that this action induces a transitive action on the maximal
ideals of $R'$.
\qed

\paragraph
Let $k$ be a field, and $G$ a finite $k$-group scheme.
Let $S=\Spec R$ be an affine $k$-scheme of finite type
with a left $G$-action.
It gives a $k$-algebra automorphism action of $G$ on $R$.
Let $A:=R^G$ be the ring of invariants.

\begin{Proposition}
Assume that $G$ is linearly reductive \(i.e., any $G$-module is semisimple\).
Then the following hold.
\begin{description}
\item[1] If $R$ satisfies Serre's $(S_r)$ condition,
then the $A$-module $R$ satisfies $(S_r)$, and 
$A$ satisfies $(S_r)$.
\item[2] If $R$ is Cohen-Macaulay, then 
$R$ is a maximal Cohen-Macaulay $A$-module,
and $A$ is also a Cohen-Macaulay ring.
\item[3] If $R$ is normal, then so is $A$.
\item[4] If $R$ is Cohen-Macaulay, then $\omega_R^G\cong \omega_A$ as
$A$-modules.
\item[5]
Assume that $R$ is Gorenstein and $\omega_R\cong R$ as $(G,R)$-modules.
Then $A=R^G$ is Gorenstein and $\omega_A\cong A$.
\end{description}
\end{Proposition}

\proof Note that the associated morphism $\pi:S=\Spec R\rightarrow \Spec A$ 
is finite surjective.

To prove the proposition, we may assume that $\Spec A$ is connected. 

Set $\bar G:=G\otimes_k \bar k$, and $\bar R=R\otimes_k \bar k$,
where $\bar k$ is the algebraic closure of $k$.
Let $G_0$ be the identity component (or the Frobenius kernel for sufficiently
high Frobenius maps, if the characteristic is nonzero) of $\bar G$, 
which is a normal subgroup scheme of $\bar G$.
Note that $\Spec \bar R\rightarrow \Spec \bar R^{G_0}$ is finite and is 
a homeomorphism, 
since $G_0$ is trivial if the characteristic is zero, and
$\bar R^{G_0}$ contains some sufficiently high Frobenius power of $\bar R$,
if the characteristic is positive.
On the other hand, the finite group
$\bar G(\bar k)=(\bar G/G_0)(\bar k)$ acts on $\bar R^{G_0}$, and
the ring of invariants under this action is $A\otimes_k\bar k$.
By Lemma~\ref{transitive.thm},
for any prime ideal $\frak p$ of $A\otimes_k\bar k$, 
$\bar G(\bar k)$ acts transitively on the set of prime ideals
of $\bar R$ (or $\bar R^{G_0}$) lying over $\frak p$.
It follows that for any prime ideal $\frak p$ of $A$ and a
prime ideal $\frak P$ of $R$ lying over $\frak p$, we have 
$\height \frak p=\height \frak P$.

Let $M$ be the sum of all non-trivial simple $G$-submodules of $R$.
As $G$ is linearly reductive, $R$ is the direct sum of $M$ and $A$ as
a $G$-module.
It is easy to see that $R=M\oplus A$ is a direct sum decomposition as
a $(G,A)$-module.

{\bf 1} Since $A$ is a direct summand of $R$ as an $A$-module,
it suffices to prove that the $A$-module $R$ satisfies the 
$(S_r)$-condition.
Let $\frak p\in\Spec A$ and assume that $\depth_{A_{\frak p}}R_{\frak p}<r$.
Note that
$\depth_{A_{\frak p}} R_{\frak p}=\inf_{\frak P}\depth R_{\frak P}$,
where $\frak P$ runs through the prime ideals lying over $\frak p$.
So there exists some $\frak P$ such that
$\depth R_{\frak P}\leq \depth_{A_{\frak p}}R_{\frak p}<r$.
As $R$ satisfies Serre's $(S_r)$-condition, we have that $R_{\frak P}$ 
is Cohen-Macaulay.
So
\[
\height \frak p=\height \frak P=\depth R_{\frak P}
\leq \depth_{A_{\frak p}}R_{\frak p}
\leq\depth A_{\frak p}\leq \height\frak p,
\]
and all $\leq$ must be $=$.
In particular, $R_{\frak p}$ is a maximal Cohen-Macaulay 
$A_{\frak p}$-module.
This shows that the $A$-module $R$ satisfies Serre's $(S_r)$-condition.

{\bf 2} is obvious by {\bf 1}.

{\bf 3} Consider the case that $k$ is of characteristic zero.
If so, then we may assume that $k$ is algebraically closed,
$G$ is a finite group, and $\Spec A$ is connected.
The action of $G$ on $R$ is extended to an action on the total
quotient ring $Q(R)$.
We have $A=R\cap Q(R)^G$.
On the other hand, any nonzero element of
$A$ is a nonzerodivisor of $R$.
In fact, if $0\neq a\in A$ is a zerodivisor of $R$, then there exists some
primitive idempotent $e_0$ of $R$ such that $e_0a=0$.
Let $e_0,e_1,\ldots,e_r$ be the primitive idempotents of $R$.
As $G$ acts transitively on these, there exists some 
$g_i\in G$ such that $e_i=g_ie_0$ for each $i$.
Then $a=1a=(e_0+\cdots+e_r)a=(g_0+\cdots+g_r)(e_0a)=0$, which is
a contradiction.
So $Q(A)\subset Q(R)^G$ (in fact $=$ holds), and $A$ is integrally
closed, as can be seen easily.

So assume that $k$ is of characteristic $p>0$.
We need to show that $A$ satisfies the $(R_1)$-condition.
Let $\frak p$ be a height one prime ideal of $A$.
As $R_{\frak p}$ is a one-dimensional (semilocal) regular ring,
the direct summand subring $A_{\frak p}$ is strongly $F$-regular,
see \cite{HH}.
Namely, $A_{\frak p}$ is a DVR.
It follows that $A$ satisfies the $(R_1)$-condition.

{\bf 4} 
Note that $\pi: S=\Spec R\rightarrow \Spec A$ is a finite 
$G$-morphism.
Set $d=\dim R=\dim A$.
As $A$ is Cohen-Macaulay and $\Spec A$ is connected, $A$ is
equidimensional of dimension $d$.
So $\height \frak m=d$ for all maximal ideals of $A$.
The same is true of $R$, and hence $R$ is also equidimensional.
So $\omega_R[d]$ and $\omega_A[d]$ are
the equivariant dualizing complexes of $R$ and $A$, respectively.
In particular, we have $\pi^!\omega_A\cong\omega_R$.
By Lemma~\ref{G-trivial-dualizing.thm}, $\omega_A$ is $G$-trivial.
By Theorem~\ref{ggd.thm}, we have
\begin{multline*}
\omega_R^G \cong  R(?)^G R\pi_* R\uHom_{\Cal O_{\Spec R}}
(\Cal O_{\Spec R},\pi^!\omega_A)\\
\cong
R(?)^G R\uHom_{\Cal O_{\Spec A}}(R\pi_*\Cal O_{\Spec R},\omega_A).
\end{multline*}
As $\pi$ is affine, $R\pi_*\Cal O_{\Spec R}=R$.
As $R$ is a maximal Cohen-Macaulay $A$-module and $\omega_A$ is a 
finitely generated $A$-module which is of finite injective dimension,
we have that $\Ext^i_A(R,\omega_A)=0$ ($i>0$).
As $G$ is linearly reductive, we have that
\[
R(?)^G R\uHom_{\Cal O_{\Spec A}}(R\pi_*\Cal O_{\Spec R},\omega_A)
\cong
\Hom_A(R,\omega_A)^G=\Hom_{G,A}(R,\omega_A).
\]
As $G$ is linearly reductive, there is a canonical direct
sum decomposition $R\cong R^G\oplus U_R$, where $U_R$ is the sum of 
all non-trivial simple $G$-submodules of $R$.
As $\omega_A$ is $G$-trivial, $\Hom_G(U_R,\omega_A)=0$.
In particular, $\Hom_{G,A}(U_R,\omega_A)=0$.

On the other hand, we have that
\[
\Hom_A(R^G,\omega_A)^G=\Hom_A(A,\omega_A)^G=\omega_A^G=\omega_A.
\]
Hence $\omega_R^G\cong \omega_A$.

{\bf 5} follows from {\bf 2} and {\bf 4} immediately.
\qed

\begin{Corollary} Let $k$ be a field, $G$ a linearly reductive finite
$k$-group scheme, and $V$ a $G$-module.
Assume that the representation $G\rightarrow \GL(V)$ factors through
$\SL(V)$.
Then the ring of invariants $A:=(\Sym V)^G$ is Gorenstein, and 
$\omega_A\cong A$.
\end{Corollary}

\proof Set $R:=\Sym V$.
As $R$ is $k$-smooth, we have that $\omega_R\cong \ext^n\Omega_{R/k}
\cong R\otimes\ext^n V$, where $n=\dim_k V$.
By assumption, $\ext^n V\cong k$, and we have that $\omega_R\cong R$,
as $(G,R)$-modules.
By the proposition, $A$ is Gorenstein and $\omega_A\cong A$.
\qed

\end{document}